\documentclass[11pt, a4paper,reqno]{amsart}
\pdfoutput=1
\usepackage{tikz,genyoungtabtikz,enumitem,rotating,latexsym,bm,stmaryrd,caption}
\usetikzlibrary{positioning,intersections,decorations}
 
\usetikzlibrary{shapes}

 \usepackage{amsmath,amsthm,amsfonts,amssymb,mathrsfs,pb-diagram}
\usepackage[%pagebackref=true, 
bookmarks=true,colorlinks=true,linktoc=page,citecolor=darkgreen,linkcolor=blue,urlcolor=mediumblue]{hyperref}
\usepackage{caption}
\usepackage{lipsum,wasysym}
\usepackage{mathtools}
\usepackage[a4paper,margin=0.9in]{geometry}
\usepackage{cleveref}
 \usepackage{amsmath}
\mathchardef\mhyphen="2D
\usepackage{color}
\usepackage{xcolor}
\usepackage{ifthen}
\usepackage{sidecap}   
\definecolor{mediumblue}{rgb}{0.0, 0.0, 0.8}
 \colorlet{darkgreen}{green!50!black}
\newcommand{\Res}{{\rm Res}}

\newcommand{\Ind}{{\rm Ind}}
\newcommand{\hstar}{\mathfrak{h}^*}
\newcommand\mptn[2]{\mathscr{P}^{#1}_{#2}}
\renewcommand{\geq}{\geqslant}
\renewcommand{\leq}{\leqslant}
\renewcommand{\trianglerighteq}{\trianglerighteqslant}
\renewcommand{\trianglelefteq}{\trianglelefteqslant}

\tikzset{wei/.style= 
{red,double=red,double
distance=0.5pt}}

\newcommand{\fS}{\mathfrak{S}}

\newcommand{\fc}{\mathfrak{c}}

\newcommand{\fh}{\mathfrak{h}}

\newcommand{\hd}{\operatorname{hd}}
\newcommand{\cA}{\mathcal{A}}
\newcommand{\cB}{\mathcal{B}}

\newcommand{\cF}{\mathcal{F}}

\newcommand{\cO}{\mathcal{O}}
\newcommand{\cP}{\mathcal{P}}

\newcommand{\Z}{\mathbb{Z}}

\newcommand{\N}{\mathbb{N}}
\newcommand{\C}{\mathbb{C}}

\DeclareMathOperator{\reg}{reg}
\DeclareMathOperator{\image}{im}

\tikzset{wei2/.style={red,double=red,double
distance=0.5pt}}

\allowdisplaybreaks
\numberwithin{equation}{section}
\usepackage{scalefnt}

\newtheorem{thm}{Theorem}[section]
\newtheorem{cor}[thm]{Corollary}

\newtheorem{lem}[thm]{Lemma}
\newtheorem{prop}[thm]{Proposition}

\newtheorem*{prop*}{Proposition}
\newtheorem*{thmA*}{Theorem A}
\newtheorem*{thmB*}{Theorem B}
\newtheorem*{thmC*}{Theorem C}\newtheorem*{thm*}{Theorem D}
\newtheorem*{cor*}{Corollary}

\newtheorem*{conj*}{Conjecture A}

\newtheorem*{mainA}{Theorem A} 
\newtheorem*{main}{Theorem B} 

\newtheorem*{conj1*}{Conjecture B}
\newtheorem*{Acknowledgements*}{Acknowledgements}

\theoremstyle{rmk}

\theoremstyle{defn}
\newtheorem{rmk}[thm]{Remark}
\newtheorem{defn}[thm]{Definition}
\newtheorem{eg}[thm]{Example}

\newcommand{\great}{>}
\newcommand{\less}{<}

\newcommand{\triv}{\mathrm{triv}}

\newcommand{\rad}{\mathrm{rad}}
\newcommand{\res}{\mathrm{res}}
\newcommand{\ik}{{k}}
\newcommand{\Std}{{\rm Std}}
\newcommand{\SStd}{{\rm SStd}}

\newcommand{\Shape}{\operatorname{Shape}} 
\newcommand{\Path}{{\rm Path}}
\newcommand{\la}{\lambda}
\newcommand{\I}{i}
\newcommand{\J}{j}

\newcommand{\M}{m}

\newcommand{\SSTS}{\mathsf{S}}

\newcommand{\SSTT}{\mathsf{T}}  %semistandard index
\newcommand{\SSTU}{\mathsf{U}}  %semistandard index 
\newcommand{\sts}{\mathsf{s}}  %standard index
\newcommand{\stt}{\mathsf{t}}  %standard index
\newcommand{\stu}{\mathsf{u}}  %standard index
\newcommand{\ZZ}{{\mathbb Z}}
\newcommand{\NN}{{\mathbb N}}
\newcommand{\g}{\ell}
 
\newcommand{\CC}{{\mathbb{C}}}
\newcommand{\RR}{{\mathbb R}}

\DeclareMathOperator{\Hom}{Hom}

\let\<=\langle
\let\>=\rangle

\tikzset{
ultra thin/.style= {line width=0.05pt},
very thin/.style=  {line width=0.2pt},
thin/.style=       {line width=0.1pt},
semithick/.style=  {line width=0.6pt},
thick/.style=      {line width=0.8pt},
very thick/.style= {line width=1.2pt},
ultra thick/.style={line width=1.6pt}
}

\crefname{defn}{Definition}{Definitions}
\crefname{thm}{Theorem}{Theorems}
\crefname{prop}{Proposition}{Propositions}
\crefname{lem}{Lemma}{Lemmas}
\crefname{cor}{Corollary}{Corollaries}
\crefname{conj}{Conjecture}{Conjectures}
\crefname{section}{Section}{Sections}
\crefname{subsection}{Subsection}{Subsections}
\crefname{eg}{Example}{Examples}
\crefname{figure}{Figure}{Figures}
\crefname{rem}{Remark}{Remarks}
\crefname{rmk}{Remark}{Remarks}
\crefname{equation}{equation}{equation}

\Crefname{defn}{Definition}{Definitions}
\Crefname{thm}{Theorem}{Theorems}
\Crefname{prop}{Proposition}{Propositions}
\Crefname{lem}{Lemma}{Lemmas}
\Crefname{cor}{Corollary}{Corollaries}
\Crefname{conj}{Conjecture}{Conjectures}
\Crefname{section}{Section}{Sections}
\Crefname{subsection}{Subsection}{Subsections}
\Crefname{eg}{Example}{Examples}
\Crefname{figure}{Figure}{Figures}
\Crefname{rem}{Remark}{Remarks}
\Crefname{rmk}{Remark}{Remarks}

  \newcommand{\Mull}{{\rm M}}
    
\usepackage[hang,flushmargin]{footmisc}

\def\Item{\item\abovedisplayskip=0pt\abovedisplayshortskip=5pt~\vspace*{-\baselineskip}} 

\begin{document}
  
% \title[On simple modules for quiver Hecke and  Cherednik algebras]{Characteristic-free bases  and BGG resolutions  \\  of     unitary simple modules  for  \\ quiver Hecke   and  Cherednik algebras}

 \title[On simple modules for quiver Hecke and  Cherednik algebras]{Characteristic-free bases  and BGG resolutions  \\  of     unitary simple modules  for  \\ quiver Hecke   and  Cherednik algebras}

% \title[On simple modules for quiver Hecke and  Cherednik algebras]{BGG resolutions of  unitary representations of \\ Cherednik and quiver Hecke algebras}

\author{C. Bowman} 
\address{\scalefont{0.9}School of Mathematics, Statistics and Actuarial Science, University of Kent, 
CT2 7NF, UK}
\email{\scalefont{0.9}C.D.Bowman@kent.ac.uk}

\author {E. Norton}
\address{\scalefont{0.9} Hausdorff Center for Mathematics, Universit\"{a}t Bonn, 53115 Bonn, Germany}
\email{\scalefont{0.9}enorton@mpim-bonn.mpg.de}

\author {J. Simental}
 \address{\scalefont{0.9} Department of Mathematics, UC Davis. One Shields Ave, Davis, CA 95616, USA.}
\email{\scalefont{0.9}jsimental@ucdavis.edu}
  
  \!\!\!\!\!\!\!\!\!\!\!\!\!\! 
\begin{abstract} 
We  provide a homological construction 
 of unitary simple modules of 
 Cherednik and Hecke algebras of type $A$ via BGG resolutions, solving a conjecture of Berkesch--Griffeth--Sam. 
 We vastly generalize the conjecture and its solution to cyclotomic Cherednik and Hecke algebras over arbitrary ground fields, compute characteristic-free bases for this family of simple modules,  and calculate the Betti numbers and   Castelnuovo--Mumford regularity  of certain   symmetric linear subspace arrangements. 

 \end{abstract}
 \maketitle
 
\!\!\!\!\!\!\!\!\!\!\!\!\!\! 
 \section*{Introduction}

\!\noindent In \cite{bgg}, Bernstein--Gelfand--Gelfand  utilise    resolutions of simple modules by Verma modules to prove certain beautiful properties 
  of  finite-dimensional  Lie algebras.  
  Such resolutions (now known as BGG resolutions) 
  have had spectacular  applications in the study of the  Laplacian on Euclidean space   \cite{MR2180410}, 
% the  
complex representation theory and homology of     Kac--Moody   algebras   \cite{MR0414645},  
 statistical mechanics and conformal  field theories \cite{GJSV13,MS94,MW03}, and they provide 
  graded free resolutions (in the sense of commutative algebra) for determinantal varieties \cite{MR520233,MR2037715}.  
   Remarkably, such resolutions have never been used in the study of 
  symmetric and general  linear groups in positive characteristic --- or indeed anywhere in modular representation theory!

%  Unitary representations hold a distinguished position in complex  Lie theory.
%  Given a Lie theoretic object, % over $\mathbb{C}$, 
  One of the most important %open 
  problems in %representation
  Lie  theory  is to classify and construct %the 
  unitary simple representations. %; this problem was originally motivated by harmonic analysis {\color{red}cite ?}. % of Lie theoretic objects.  
%     For Lie groups, classifying and constructing such representations 
% is one of most important open problems in representation theory 
% 
% with applications in harmonic analysis 
For Lie groups, this ongoing project draws on techniques from %topology \cite{MR2409701}, 
 Dirac cohomology \cite{MR3714563}, Kazhdan--Lusztig theory \cite{MR3263031},  and the Langlands Program \cite{MR1770725}, and 
has provided profound insights into relativistic quantum mechanics \cite{MR1503456}. 
% The Cherednik algebra of a complex reflection group $W$ is a Lie theoretic  object whose representation theory is built on structural parallels with classical Lie theory and a direct relationship with the representation theory of the Hecke algebra of $W$ at a root of unity \cite{GGOR}. 
 The Cherednik algebra of a complex reflection group, $W$, is an important Lie theoretic   object which possesses %all the 
 hallmarks from the classical theory% for semi-simple Lie algebras
: a triangular decomposition and a category $\mathcal{O}$ with a highest weight theory \cite{GGOR}, analogues of translation functors \cite{Losev2017}, induction and restriction functors \cite{MR2511190} with associated Harish--Chandra series \cite{MR3782434}, and Kazhdan--Lusztig theory 
%in the case of cyclotomic groups $W$ \cite{RSVV}.
   \cite{RSVV} (for  $W=G(\ell,1,n)$).
%The huge importance of  unitary representations
%   of  both  real reductive groups \cite{MR0399356,MR0439993,MR0439994} 
%     and of Cherednik algebras \cite{MR3810569}
%     in harmonic analysis .  
Both the   unitary representations
   of     real reductive groups \cite{MR0399356,MR0439993,MR0439994} 
     and  those of Cherednik algebras \cite{MR3810569}
     are of huge importance in algebraic harmonic analysis.  
     
       For   Cherednik algebras of symmetric groups, $H_{1/e}(\mathfrak{S}_n)$, %the study of unitary simple representations is motivated by harmonic analysis \cite{MR3810569}. The
  the simple unitary representations $L(\lambda)$ of $H_{1/e}(\mathfrak{S}_n)$ were classified by Etingof, Griffeth and Stoica   \cite{MR2534594} by a combinatorial condition on the partition $\lambda$ of $n$ labeling the ``highest weight" of $L(\lambda)$. In the spirit of classical results in Lie theory, Berkesch, Griffeth, and Sam subsequently conjectured that any unitary simple $L(\lambda)$ admits a BGG resolution \cite[Conjecture 4.5]{conjecture}.
 
 The primary purpose of this paper is to prove Berkesch--Griffeth--Sam's conjecture %\cite[Conjecture 4.5]{conjecture}
   and thus {\em homologically construct} the unitary simple $ H_{1/e}(\mathfrak{S}_n)$-modules:

%   
%The original motivation of this paper was to provide  graded free resolutions (in the sense of commutative algebra) for the algebraic varieties [Jos\'e -- insert  some sexiness about these varieties here]
%   $$
%X_{e, k, n} := \mathfrak{S}_{n}\{(z_{1}, \dots, z_{n}) \in \C^{n} : z_{ie+1} = \cdots=z_{(i+1)e} \text{ for } 0\leq i < k \} 
%   $$  
%and to provide an explicit homological  construction of all unitary simple $H_{1/e}(\mathfrak{S}_n)$-modules via BGG resolutions (and hence   verify Berkesch--Griffeth--Sam's 
%conjecture) and to  compute Dirac cohomology of Cherednik algebras. 
% 

       \begin{mainA}
Associated to any simple unitary $  H_{1/e}(\mathfrak{S}_n)$-module, $L(\lambda)$,   we have a 
  complex 
  $   {\bf C}_\bullet (\lambda)=   \bigoplus_{
\begin{subarray}c 
 \la\trianglerighteq  \nu
\end{subarray}
}\Delta(\nu)\langle \ell(\nu)\rangle
  $
with differential given by an alternating  sum over all ``one-column homomorphisms".  
This complex 
 is exact except in degree zero, where    $H_0({\bf C}_\bullet(\la))=L(\la).$    
   \end{mainA} 
%We hence provide an explicit homological  construction of all unitary simple $H_{1/e}(\mathfrak{S}_n)$-modules via BGG resolutions, prove Berkesch--Griffeth--Sam's 
%conjecture and 
% compute Dirac cohomology of Cherednik algebras.  
%Our representation theoretic resolutions also carry a geometric interpretation.
%Surprisingly,  our  proof  only depends on  elementary combinatorics and Frobenius reciprocity.   
In contrast to classical papers on BGG resolutions and unitary representations, which usually employ ideas from algebraic geometry, our methods are completely algebraic and moreover, yield several geometric results. Namely, each standard module $\Delta(\nu)$ is a free $\C[x_1,\dots,x_n]$-module, and as a consequence we obtain $\mathfrak{S}_{n}$-equivariant, graded free resolutions (in the sense of commutative algebra) for the $e$-equals variety
  $$X_{e, 1, n} := \mathfrak{S}_{n}\{(z_{1}, \dots, z_{n}) \in \C^{n} : z_{1} = \dots = z_{e}\},$$
  and for the following algebraic varieties when $n=ke$:
   $$
%   	X_{e, 1, n} := \mathfrak{S}_{n}\{(z_{1}, \dots, z_{n}) \in \C^{n} : z_{1} = \dots = z_{e}\}.
%\qquad 
X_{e, k, n} := \mathfrak{S}_{n}\{(z_{1}, \dots, z_{n}) \in \C^{n} : z_{ie+1} = \cdots=z_{(i+1)e} \text{ for } 0\leq i < k \}
   .$$   
We hence provide formulae for the  graded Betti numbers and calculate the Castelnuovo--Mumford regularity of these varieties -- a notoriously difficult problem in general  \cite{MR1942401, MR3299723}. Moreover, we also provide formulae for these invariants in the cyclotomic case, where the equalities in the equations defining the above varieties become equalities \emph{up to multiplication by an $\ell$th root of unity.}
Finally, we remark that the Cherednik algebra approach to geometric resolutions  was inspired by the 
 Lie theoretic construction of Lascoux's  resolutions of determinantal  varieties 
(via parabolic BGG resolutions of unitary modules)
   \cite{MR2037715,conjecture}; it would be interesting to find a   purely geometric proof of the resolutions of our varieties by analogy with \cite{MR520233}.

%\cite{MR0414645,GJSV13,MS94,MW03,MR520233,MR2037715}.  

%The CM regularity of $X_{e, 1, n}$ was computed in \cite{conjecture} for $e \geq (n+1)/2$ but the computation crucially depends on the fact that, in this case, $X_{e,1,n}$ is Cohen-Macaulay. For $e < (n+1)/2$, the variety $X_{e,1,n}$ is not Cohen-Macaulay and a different approach is needed.

%However,  our interest soon evolved beyond the usual confines of the word ``unitary".  
%Within  the context Hecke algebra 
%Applying the KZ functor to the  resolutions of Theorem A, one immediately obtains resolutions of the corresponding simple modules for the Hecke algebra over $\mathbb{C}$.  Furthermore, we obtain explicit bases 
%and representing matrices. %     and completely describe  the branching rule for these simples!     

A key ingredient to our proofs is to work in the 2-categorical  setting of {\em diagrammatic Cherednik algebras} of \cite{Webster}.  The diagrammatic calculus is easier for calculation and benefits from a graded structure.  
% We construct our resolutions using the (Morita equivalent) diagrammatic Cherednik algebra introduced by Webster \cite{Webster}. 
 The diagrammatic approach allows us to generalize the original conjecture to higher levels and arbitrary ground fields; we prove this more general version.
 We recast the combinatorial condition in type $A$ for $L(\lambda)$ to be unitary \cite{MR2534594} as, the partition $\lambda$ lies in the \textit{fundamental alcove} of the dominant chamber in an affine type $A$ alcove geometry.  In our BGG resolution,  $\Delta(\nu)$ appears in  homological degree $d$ if and only if 
 $\nu$ is obtained from  $\lambda$ by  reflecting across $d$ walls (increasing the distance from the fundamental alcove by 1 at each step).
    This alcove model  vastly generalizes to  the set of all  $\ell$-partitions whose components {\em each have at most $h$ columns}, $\mptn \ell n (h)$.  For any multipartition lying in the fundamental alcove we then construct a BGG resolution of the corresponding simple $H_c(G(\ell,1,n))$-module. We remark that  Griffeth has obtained a combinatorial description of the $\ell$-partitions that label unitary irreducible modules for $H_{c}(G(\ell, 1, n))$, \cite{griffeth}, and it would be interesting to compare this condition to the one arising from the alcove model.

  Working with  quiver Hecke algebras furthermore allows us to obtain our results over fields, $\Bbbk$, of arbitrary characteristic. 
 The search for an effective description of 
 the dimensions of simple representations of symmetric groups over arbitrary fields is a centre of gravity for much research in modular Lie theory % and is currently undergoing a revolution %breakthroughs by Williamson, Riche, Elias, Lusztig and others
   \cite{w16,MR3766576,MR3755781}.   
 We  construct explicit bases and representing matrices of unitary simple modules (as modules for the quiver Hecke algebra over $\mathbb{C}$) at the same time as we establish the properties of their resolutions.  
 We show that  our  bases, representing matrices, and  resolutions for %$e$-
 unitary simples %of quiver Hecke algebras
  remain  stable under reduction modulo $p$ --- in other words,   the many beautiful properties of unitary modules extend beyond the confines of   characteristic zero (a necessary condition for the  definition of unitary modules via bilinear forms to make sense) to arbitrary fields.
% 
% 
% 
%  
% %
%%
%%
%Undeterred,  a dreamer could perhaps ask  
% an even deeper and more structural question: {\em is it possible to construct explicit bases 
%    of modular 
% simples  
%  of   symmetric groups?}   
%   
%  
This makes these  ``$p$-unitary simples" the most well-understood family of    simple
  modules for  symmetric groups in positive characteristic.
   In particular,  
all   results   on  non-unitary  
 simples  
 hold  on the level of 
   {\em dimensions or characters}   (and  
   proceed by calculating decomposition matrices) and so 
   our results are on a higher structural level than these. 
     Finally, in   \cref{mull,easypeasy} we obtain a simple closed form for  the Mullineux involution, $\Mull$, on   unitary simples 
      and  {\em explicitly construct this  isomorphism}  --- to our knowledge, this is the first time such an isomorphism has been explicitly constructed (outside of the semisimple case).  %In \cref{mulling}, we define an involution, $\Mull$, on the  set of $p$-restricted  tableaux.  

   This pivots the impact of our result from Cherednik algebras and geometry of subspace arrangements, to modular representations of the symmetric group.     As our main result is the first of its kind  
 for  symmetric groups %(or indeed anywhere in modular Lie theory) 
   we 
 state it now in this simplified form.  
 For the far more general statement concerning all cyclotomic 
  quiver    Hecke %and Cherednik
   algebras, %of type $G(\ell,1,n)$,
     see \cref{maintheomre,action}.  

      \begin{main}Let $\Bbbk$ be a   field of characteristic $p>0$.  
 For $D^\Bbbk_n(\la)$    a  ``$p$-unitary simple" we have an associated 
  $\Bbbk \mathfrak{S}_n$-complex 
   $C_\bullet(\la)  =   \bigoplus_{
\begin{subarray}c 
 \la\trianglerighteq  \nu
\end{subarray}
}S _n(\nu)\langle \ell(\nu)\rangle
$ 
with differential given by an alternating  sum over all ``one-column homomorphisms".  
This complex 
 is exact except in degree zero, where %it is isomorphic to the simple module 
 $$H_0(C_\bullet(\la))=D^\Bbbk_n(\la).$$   
Moreover,  the simple $\Bbbk \mathfrak{S}_n$-module $D^\Bbbk_n(\la)$ is free as a $\mathbb{Z}$-module with basis  
$ 
\{c_{\sts} \mid \sts \in \Std_p(\la)  \} 
$ where $\Std_p(\la)\subseteq \Std(\la) $ is the  set of $p$-restricted  tableaux.  
The action on this basis is given in \cref{action}.    
We have that $D^\Bbbk_n(\la) \otimes {\rm sgn}  \cong D^\Bbbk_n(\la_\Mull)$ under the map $:c_\sts \mapsto c_{\sts_\Mull}.$
 \end{main}

We thus provide the first instances of BGG resolutions anywhere in modular representation theory and in particular the first homological construction of a family of simple modules for symmetric groups.  
For the symmetric groups and their Hecke algebras, our   bases lift a combinatorial result  of Kleshchev \cite{MR1383482}  to   a structural level
%too repetitive now that we added it on p.1: 
% and  our  construction of  BGG resolutions 
  % settles  Berkesch--
    %  Griffeth--Sam's conjecture over $\mathbb{C}$  \cite{conjecture}  (and generalises this work to arbitrary fields!)  
    and our resolutions provide  a structural lift of a character-theoretic result of Ruff   
     \cite{MR2266877}.   
For   Hecke algebras of type $B$, the simplest examples of our resolutions have   appeared in work of mathematical physicists concerning Virasoro and blob algebras  \cite{GJSV13,MS94,MW03}
  but our bases (and representing matrices) are entirely new.    
%We vastly  generalise these results and conjectures   to arbitrary fields and higher levels: in particular, to simple modules labelled by the set of  $\ell$-partitions whose components {\em each have at most $h$ columns}, $\mptn \ell n (h)$. 
 We remark that our results for the  Hecke algebras depend only on the quantum parameter $e\in \mathbb{N}$ and are entirely independent of the characteristic of the underlying field (for $\ell=1$ we set $e=p$ in 
  \cref{maintheomre,action} to obtain the above result for symmetric groups).  
%  Furthermore 
% This is   the first characteristic-free description of a family of simple modules of general quiver Hecke  algebras.   

%As $n\rightarrow \infty$, the composition series of  our  Specht modules     quickly become incalculable. 
%Let $p>0$, our theorem has no restriction on the $p$-weight of a $p$-unitary partition; 

The partitions and multipartitions we consider (namely those lying in the fundamental alcove)  
have    no restriction on their $e$-weight;    calculating the composition series of the corresponding  Specht modules for symmetric groups and Hecke algebras  in positive characteristic    is 
 far beyond  
  the current realms of conjecture 
 (which at present have been stretched as far as $w(\la)<p^2$ for $h=3$ by Lusztig--Williamson  \cite{MR3766576}).  
Over $\mathbb{C}$ calculating the  composition series of these Specht modules 
 is theoretically possible using Kazhdan--Lusztig theory --- however it quickly becomes computationally impossible --- we provide examples of 
  series of   Specht modules (of rank $n$ as $n\rightarrow \infty$) for which   the length of   the composition series   tends to infinity.  
  Thus our two descriptions (homological and via bases) of unitary simple modules provide the {\em only} contexts in which we can hope to understand  these simple modules.  See \cref{section4} for more details.

\section{Combinatorics of reflection groups and their deformations }\label{sec:3}

We first recall the basic combinatorics controlling Hecke and Cherednik algebras and their diagrammatic analogues.  
In the spirit of higher representation theory, we begin by introducing the necessary quiver-theoretic notation.  
 Fix  $e\in \{2,3,4,\dots\}$ $\cup\{\infty\}$. If $e=\infty$ then we set $I=\ZZ$, while if $e<\infty$ then we set $I=\ZZ/e\ZZ$. 
 We  let $\Gamma_e$ be the
  quiver with vertex set $I $ and edges $i\longrightarrow i+1$, for
  $i\in I$.   %Hence, we are considering either the linear quiver~$\ZZ$ ($e=\infty$)
%  or a cyclic quiver ($e<\infty$):
   To the quiver $\Gamma_e$ we attach the symmetric Cartan matrix with  entries 
    $(a_{i,j})_{i,j\in I}$   defined by $a_{ij}=2\delta_{ij}-\delta_{i(j+1)}-\delta_{i(j-1)}$.
    Following \cite[Chapter~1]{Kac}, let $\widehat{\mathfrak{sl}}_e$ be the
  Kac-Moody algebra of~$\Gamma_e$  with simple roots
  $\{\alpha_i \mid i\in I\}$, fundamental weights $\{\Lambda_i\mid i\in I\}$,
  positive weight lattice $P^+=\bigoplus_{i\in I}\ZZ_{\geq 0}\Lambda_i$ and
  positive root lattice $Q^+=\bigoplus_{i\in I}\ZZ_{\geq 0}\alpha_i$. Let
  $(\cdot,\cdot)$ be the usual invariant form associated with this data,
  normalised so that $(\alpha_i,\alpha_j)=a_{ij}$ and
  $(\Lambda_i,\alpha_j)=\delta_{ij}$, for $i,j\in I$.   
  Fix a sequence $\kappa=(\kappa_1,\dots,\kappa_\ell)\in I^\ell$, the
  $e$-{\sf multicharge}, and define
  $\Lambda=\Lambda(\kappa) =
      \Lambda_{{\kappa}_1}+\dots+\Lambda_{{\kappa}_\ell}      $. Then $\Lambda\in P^+$ is
  dominant weight of {\sf level}~$\ell$.  
%   Throughout this paper we shall 
Antecedents  to this paper \cite{Webster,CoxBowman,bcs15,bs15}
 feature a  parameter $\theta\in\ZZ^\ell$; we have   fixed 
      $\theta=(1,2,\dots ,\ell) $  and dropped this from our notation.  
% Upon fixing a dominance order in \cref{}, we fixed a particular weighting 
% (in the parlance of 
%\cite{Webster,CoxBowman,bcs15,bs15})  
%and so we have dropped this from our notation 

We define a {\sf partition}, $\lambda$,  of $n$ to be a   finite weakly decreasing sequence  of non-negative integers $ (\lambda_1,\lambda_2, \ldots)$ whose sum, $|\lambda| = \lambda_1+\lambda_2 + \dots$, equals $n$.   
An    {\sf $\ell $-partition}  $\lambda=(\lambda^{(1)},\dots,\lambda^{(\ell)})$ of $n$ is an $\ell $-tuple of     partitions  such that $|\lambda^{(1)}|+\dots+ |\lambda^{(\ell)}|=n$. 
We will denote the set of   $\ell$-partitions of $n$ by $\mptn {\ell}n$.
Given  $\lambda=(\lambda^{(1)},\lambda^{(2)},\ldots ,\lambda^{(\ell)}) \in \mptn {\ell}n$, the {\sf Young diagram} of $\lambda$    is   the set of nodes, 
\[
\{(r,c,m) \mid  1\leq  c\leq \lambda^{(m)}_r\}.
\]
%We do not distinguish between the $\ell$-partition and its Young diagram.  
  We refer to a node $(r,c,m)$ as being in the $r$th row and $c$th column of the $m$th component of $\lambda$.  
 Given a node, $(r,c,m)$,  
we define the {\sf residue} of this node  to be  ${\rm res}(r,c,m)  = \kappa_m+  c - r\pmod e$.  
We refer to a node of residue  $i\in I$ as an $i$-node.

Given $\lambda\in \mptn \ell n$,  the associated  {\sf   Russian array} 
is defined as follows.
For each $1\leq m\leq \ell$, we place a point on the real line at $m$ and consider the region bounded by  half-lines at angles $3\pi/4 $ and $\pi/4$.
		We tile the resulting quadrant with a lattice of squares, each with diagonal of length $2 \ell $.  
We place a box  $(1,1,m)\in \lambda$   at the point  $m$ on the real line, with rows going northwest from this node, and columns going northeast. 
We do not distinguish between $\la$ 
and its  Russian array.

There are many different orderings on $\mptn \ell n$, each gives rise to a different (diagrammatic) Cherednik algebra and a different lens through which to study the quiver Hecke algebra \cite{manycell}.  
 We shall restrict our attention to the most natural of these orderings, which we now define.

\begin{defn} \label{domdef}
Let  $(r,c,m), (r',c',m')    $ be two boxes.  
We write $(r,c,m) \rhd   (r',c',m')$ if 
either 
\begin{itemize}
\item[$(i)$] $  \ell(r-c)+m <  \ell(r'-c') + {m'} $ 
or
\item[$(ii)$] $  \ell(r-c) +m =  \ell(r'-c')+m'$
and $r+c < r'+c'$.  
\end{itemize}
Given  $\la,\mu\in \mptn \ell n$, we say that  
$\la$  {\sf dominates} $\mu$ (and write $\mu \trianglelefteq   \la$) 
if for every $i$-box  $(r,c,m) \in \mu$, there exist  at least as many $i$-boxes $(r',c',m') \in \la$ which  dominate $(r,c,m)$ than  there do 
$i$-boxes $(r'',c'',m'') \in \mu$ which  dominate $(r,c,m)$.  
\end{defn}

%\begin{rmk}
%In \cite{CoxBowman}, it is demonstrated that the above ordering is the most natural from the viewpoint of modular representation theory of quiver Hecke algebras.   
%\end{rmk}
\begin{rmk}
This dominance ordering    is a coarsening of the usual $c$-function ordering on the Fock spaces   of Foda--Leclerc--Okado--Thibon--Welsh \cite{Webster,manycell}.  
This  is the only ordering for which we have a 
closed form for a labelling of the simple modules for the quiver Hecke algebra  \cite{MR1671762} (in other words,   a labelling of the   component of the $\widehat{\mathfrak{sl}}_e$   crystal containing the empty $\ell$-partition).  
% (in other words, the component labelling  simple modules of the Hecke algebra).   %KZ-component of the   crystal.  
\end{rmk}

\begin{defn}  \label{stndrd}
Given $\lambda\in \mptn \ell n$, we define a {\sf tableau} of shape $\lambda$ to be a filling of the boxes of 
the  Russian array of  $\la $ with the numbers 
$\{1,\dots , n\}$.  We define a {\sf  standard tableau} to be a tableau  in which    the entries increase along the rows and columns of each component.  
We let $\Std (\lambda)$ denote the set of all standard tableaux of shape $\lambda\in\mptn\ell n$. 
 \end{defn}

  \!\!\!\!\!\!\!\!\!\!\!\!\!\! 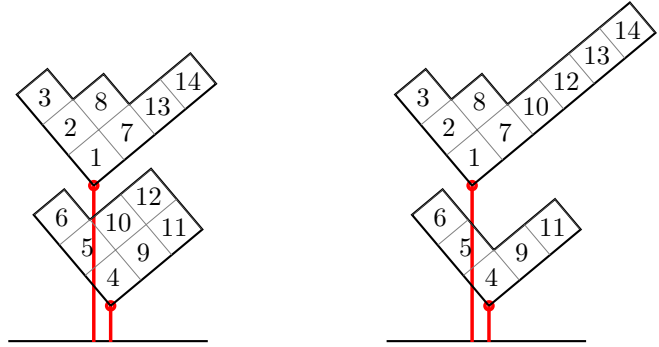
\begin{figure}[ht]\captionsetup{width=0.9\textwidth}\scalefont{0.9}
\[
 \begin{tikzpicture}[scale=0.75]
\begin{scope}
\draw[thick](2,-1.95)--(-1.5,-1.95); 
{     \path (0,0.8) coordinate (origin);     \draw[wei2] (origin)   circle (2pt);
 
\clip(-2.2,-1.95)--(-2.2,4.3)--(3.5,4.3)--(3.5,-1.95)--(-2.2,-1.95); 
\draw[wei2] (origin)   --(0,-2); 
\draw[thick] (origin)  
--++(130:3*0.7)
--++(40:1*0.7)
--++(-50:1*0.7)	
--++(40:1*0.7) --++(-50:1*0.7)
--++(40:2*0.7) --++(-50:1*0.7) 
--++(-50-90:4*0.7);
 \clip (origin)  
--++(130:3*0.7)
--++(40:1*0.7)
--++(-50:1*0.7)	
--++(40:1*0.7) --++(-50:1*0.7)
--++(40:2*0.7) --++(-50:1*0.7) 
--++(-50-90:4*0.7);
\path  (origin)--++(40:0.35)--++(130:0.35)  coordinate (111); 
   \node at  (111)  {  $1$};
\path  (origin)--++(40:0.35)--++(130:0.35+0.7)  coordinate (121); 
   \node at  (121)  {  $2$};
     \path  (origin)--++(40:0.35)--++(130:0.35+0.7*2)  coordinate (131); 
   \node at  (131)  {  $3$};
     \path  (origin)--++(40:0.35+0.7)--++(130:0.35 )  coordinate (211); 
   \node at  (211)  {  $7$};
            \path  (origin)--++(40:0.35+0.7)--++(130:0.35+0.7 )  coordinate (221); 
   \node at  (221)  {  $8$};
  %%%%%%%%%%%%
           \path  (origin)--++(40:0.35+0.7*2)--++(130:0.35 )  coordinate (211); 
   \node at  (211)  {  $13$};
              \path  (origin)--++(40:0.35+0.7*3)--++(130:0.35 )  coordinate (211); 
   \node at  (211)  {  $14$};
\path (40:1cm) coordinate (A1);
\path (40:2cm) coordinate (A2);
\path (40:3cm) coordinate (A3);
\path (40:4cm) coordinate (A4);
\path (130:1cm) coordinate (B1);
\path (130:2cm) coordinate (B2);
\path (130:3cm) coordinate (B3);
\path (130:4cm) coordinate (B4);
\path (A1) ++(130:3cm) coordinate (C1);
\path (A2) ++(130:2cm) coordinate (C2);
\path (A3) ++(130:1cm) coordinate (C3);
\foreach \i in {1,...,19}
{
\path (origin)++(40:0.7*\i cm)  coordinate (a\i);
\path (origin)++(130:0.7*\i cm)  coordinate (b\i);
\path (a\i)++(130:4cm) coordinate (ca\i);
\path (b\i)++(40:4cm) coordinate (cb\i);
\draw[thin,gray] (a\i) -- (ca\i)  (b\i) -- (cb\i); } 
}
\end{scope}
\begin{scope}
{   
\path (0,-1.3)++(40:0.3*0.7)++(-50:0.3*0.7) coordinate (origin);  
\draw[wei2]  (origin)   circle (2pt);
\clip(-2.2,-1.95)--(-2.2,3)--(3,3)--(3,-1.95)--(-2.2,-1.95); 
\draw[wei2] (origin)   --++(-90:4cm); 
\draw[thick] (origin) 
--++(130:3*0.7)
--++(40:1*0.7)
--++(-50:1*0.7)	
--++(40:2*0.7)   --++(-50:2*0.7) 
--++(-50-90:3*0.7);
 \clip (origin) 
--++(130:3*0.7)
--++(40:1*0.7)
--++(-50:1*0.7)	
--++(40:2*0.7)   --++(-50:2*0.7) 
--++(-50-90:3*0.7);
  \path  (origin)--++(40:0.35)--++(130:0.35)  coordinate (111); 
   \node at  (111)  {  $4$};
\path  (origin)--++(40:0.35)--++(130:0.35+0.7)  coordinate (121); 
   \node at  (121)  {  $5$};
     \path  (origin)--++(40:0.35)--++(130:0.35+0.7*2)  coordinate (131); 
   \node at  (131)  {  $6$};
     \path  (origin)--++(40:0.35+0.7)--++(130:0.35 )  coordinate (211); 
   \node at  (211)  {  $9$};
            \path  (origin)--++(40:0.35+0.7)--++(130:0.35+0.7 )  coordinate (221); 
   \node at  (221)  {  $10$};
           \path  (origin)--++(40:0.35+2*0.7)--++(130:0.35  )  coordinate (224); 
   \node at  (224)  {  $11$};       
   \path  (origin)--++(40:0.35+2*0.7)--++(130:0.35+0.7 )  coordinate (223); 
   \node at  (223)  {  $12$};
\path (40:1cm) coordinate (A1);
\path (40:2cm) coordinate (A2);
\path (40:3cm) coordinate (A3);
\path (40:4cm) coordinate (A4);
\path (130:1cm) coordinate (B1);
\path (130:2cm) coordinate (B2);
\path (130:3cm) coordinate (B3);
\path (130:4cm) coordinate (B4);
\path (A1) ++(130:3cm) coordinate (C1);
\path (A2) ++(130:2cm) coordinate (C2);
\path (A3) ++(130:1cm) coordinate (C3);
\foreach \i in {1,...,19}
{
\path (origin)++(40:0.7*\i cm)  coordinate (a\i);
\path (origin)++(130:0.7*\i cm)  coordinate (b\i);
\path (a\i)++(130:4cm) coordinate (ca\i);
\path (b\i)++(40:4cm) coordinate (cb\i);
\draw[thin,gray] (a\i) -- (ca\i)  (b\i) -- (cb\i); } 
}  \end{scope}
\end{tikzpicture}
\qquad 
\begin{tikzpicture}[scale=0.75]
\begin{scope}
\draw[thick](2,-1.95)--(-1.5,-1.95); 
{     \path (0,0.8) coordinate (origin);     \draw[wei2] (origin)   circle (2pt);
 
\clip(-2.2,-1.95)--(-2.2,4.3)--(3.5,4.3)--(3.5,-1.95)--(-2.2,-1.95); 
\draw[wei2] (origin)   --(0,-2); 
\draw[thick] (origin)  
--++(130:3*0.7)
--++(40:1*0.7)
--++(-50:1*0.7)	
--++(40:1*0.7) --++(-50:1*0.7)
--++(40:4*0.7) --++(-50:1*0.7) 
--++(-50-90:6*0.7);
 \clip (origin)  
--++(130:3*0.7)
--++(40:1*0.7)
--++(-50:1*0.7)	
--++(40:1*0.7) --++(-50:1*0.7)
--++(40:4*0.7) --++(-50:1*0.7) 
--++(-50-90:6*0.7);
\path  (origin)--++(40:0.35)--++(130:0.35)  coordinate (111); 
   \node at  (111)  {  $1$};
\path  (origin)--++(40:0.35)--++(130:0.35+0.7)  coordinate (121); 
   \node at  (121)  {  $2$};
     \path  (origin)--++(40:0.35)--++(130:0.35+0.7*2)  coordinate (131); 
   \node at  (131)  {  $3$};
     \path  (origin)--++(40:0.35+0.7)--++(130:0.35 )  coordinate (211); 
   \node at  (211)  {  $7$};
            \path  (origin)--++(40:0.35+0.7)--++(130:0.35+0.7 )  coordinate (221); 
   \node at  (221)  {  $8$};
  %%%%%%%%%%%%
           \path  (origin)--++(40:0.35+0.7*2)--++(130:0.35 )  coordinate (211); 
   \node at  (211)  {  $10$};
              \path  (origin)--++(40:0.35+0.7*3)--++(130:0.35 )  coordinate (211); 
   \node at  (211)  {  $12$};
                 \path  (origin)--++(40:0.35+0.7*4)--++(130:0.35 )  coordinate (211); 
   \node at  (211)  {  $13$};
                 \path  (origin)--++(40:0.35+0.7*5)--++(130:0.35 )  coordinate (211); 
   \node at  (211)  {  $14$};
\path (40:1cm) coordinate (A1);
\path (40:2cm) coordinate (A2);
\path (40:3cm) coordinate (A3);
\path (40:4cm) coordinate (A4);
\path (130:1cm) coordinate (B1);
\path (130:2cm) coordinate (B2);
\path (130:3cm) coordinate (B3);
\path (130:4cm) coordinate (B4);
\path (A1) ++(130:3cm) coordinate (C1);
\path (A2) ++(130:2cm) coordinate (C2);
\path (A3) ++(130:1cm) coordinate (C3);
\foreach \i in {1,...,19}
{
\path (origin)++(40:0.7*\i cm)  coordinate (a\i);
\path (origin)++(130:0.7*\i cm)  coordinate (b\i);
\path (a\i)++(130:4cm) coordinate (ca\i);
\path (b\i)++(40:4cm) coordinate (cb\i);
\draw[thin,gray] (a\i) -- (ca\i)  (b\i) -- (cb\i); } 
}
\end{scope}
\begin{scope}
{   
\path (0,-1.3)++(40:0.3*0.7)++(-50:0.3*0.7) coordinate (origin);  
\draw[wei2]  (origin)   circle (2pt);
\clip(-2.2,-1.95)--(-2.2,3)--(3,3)--(3,-1.95)--(-2.2,-1.95); 
\draw[wei2] (origin)   --++(-90:4cm); 
\draw[thick] (origin) 
--++(130:3*0.7)
--++(40:1*0.7)
--++(-50:2*0.7)	
--++(40:2*0.7)   --++(-50:1*0.7) 
--++(-50-90:3*0.7);
 \clip (origin) 
--++(130:3*0.7)
--++(40:1*0.7)
--++(-50:2*0.7)	
--++(40:2*0.7)   --++(-50:1*0.7) 
--++(-50-90:3*0.7);
  \path  (origin)--++(40:0.35)--++(130:0.35)  coordinate (111); 
   \node at  (111)  {  $4$};
\path  (origin)--++(40:0.35)--++(130:0.35+0.7)  coordinate (121); 
   \node at  (121)  {  $5$};
     \path  (origin)--++(40:0.35)--++(130:0.35+0.7*2)  coordinate (131); 
   \node at  (131)  {  $6$};
     \path  (origin)--++(40:0.35+0.7)--++(130:0.35 )  coordinate (211); 
   \node at  (211)  {  $9$};
            \path  (origin)--++(40:0.35+0.7)--++(130:0.35+0.7 )  coordinate (221); 
   \node at  (221)  {  $3$};
           \path  (origin)--++(40:0.35+2*0.7)--++(130:0.35  )  coordinate (224); 
   \node at  (224)  {  $11$};       
   \path  (origin)--++(40:0.35+2*0.7)--++(130:0.35+0.7 )  coordinate (223); 
   \node at  (223)  {  $2$};
\path (40:1cm) coordinate (A1);
\path (40:2cm) coordinate (A2);
\path (40:3cm) coordinate (A3);
\path (40:4cm) coordinate (A4);
\path (130:1cm) coordinate (B1);
\path (130:2cm) coordinate (B2);
\path (130:3cm) coordinate (B3);
\path (130:4cm) coordinate (B4);
\path (A1) ++(130:3cm) coordinate (C1);
\path (A2) ++(130:2cm) coordinate (C2);
\path (A3) ++(130:1cm) coordinate (C3);
\foreach \i in {1,...,19}
{
\path (origin)++(40:0.7*\i cm)  coordinate (a\i);
\path (origin)++(130:0.7*\i cm)  coordinate (b\i);
\path (a\i)++(130:4cm) coordinate (ca\i);
\path (b\i)++(40:4cm) coordinate (cb\i);
\draw[thin,gray] (a\i) -- (ca\i)  (b\i) -- (cb\i); } 
}  \end{scope}
\end{tikzpicture}
\]
 \caption{The tableau $\stt^\lambda$ for two bi-partitions $\la,\mu\in  \mptn 2 {14}(3)$ such that $\la \rhd \mu$.  Note that the latter is obtained from the former by removing a strip of nodes from one column and adding them in another column {\em further to the right}. 	 	}  
\label{cyclic tableau}
\end{figure}

\begin{defn}\label{nicetab}
Given $\la \in \mptn \ell n$ we let $\stt^\la \in \Std (\la)$ be the tableau   obtained by placing the entry $n$ in the least dominant removable box $(r,c,m) \in \lambda $  and then placing the entry $n-1$ in the least dominant removable box of $\la  \setminus \{(r,c,m) \}$ and continuing in this fashion.  
 \end{defn}

Given $\la \in \mptn \ell n$, we let   ${\rm Rem} (\la)$ (respectively ${\rm Add}  (\la)$) 
denote the set of all    removable  respectively addable) boxes  of the Young diagram of $\la$ so that the resulting diagram is the Young diagram of a $\ell$-partition.  Given $i\in I$, we let  ${\rm Rem}_i  (\la)\subseteq {\rm Rem}  (\la)$ (respectively  ${\rm Add}_i  (\la)\subseteq {\rm Add}  (\la)$) denote the subset of boxes of residue $i\in I$.

     \begin{defn}\label{admissible}  For $h\in \mathbb{Z}_{>0}$  we say that $\kappa\in I^\ell $ is  
 $h$-{\sf admissible}   if  
$(\Lambda,\alpha_i+\alpha_{i+1}+\dots+\alpha_{i+h-1})\le 1$  for all 
 $ i\in I$.  
 \end{defn}

  Given $h\in \NN$, we  let $\mptn \ell n (h)$ denote the subset of $\mptn \ell n$ consisting of those $\ell$-partitions which have at most $h$ columns in each component, that is
  $$
\mptn \ell n (h)=\{
 \lambda= (\lambda^{(1)},\lambda^{(2)},\dots, \lambda^{(\ell)})
 \mid\lambda^{(m)}_1 \leq h \text{ for }1\leq m \leq \ell\}.
  $$
Given $\la \in \mptn \ell n$, we define its residue sequence, ${\rm res}(\la)$ to be the sequence obtained by recording the  residues of the boxes of $\la$ according to the dominance ordering on boxes.

%\begin{eg}
%Let $\kappa=(0,3)$ and $e=7$.  
%The residue sequence of $((3,2,1^2), (3,2^2)) \in \mptn 2 {14}$  is given by 
%$(2,5,1,4,0,0,3,3,6,2,2,5,1,4)$.
%\end{eg}

\section{Diagrammatic    algebras}
  We now recall the diagrammatic Cherednik and Hecke algebras of \cite{Webster} and some of the
  necessary results concerning their  representation theory   from \cite{CoxBowman,bcs15,bs15}.  
  We first tilt  the  Russian array of $\lambda\in \mptn\ell n$   ever-so-slightly in the clockwise direction so that the top vertex of the box $(r,c,m)\in \lambda$  
   has $x$-coordinate $\mathbf{I}_{(r,c,m)}= m + \ell(r-c) + (r+c)\epsilon$   (using standard small-angle identities to approximate the coordinate to order $\epsilon^2$)
   providing  
$\epsilon\ll \tfrac{1}{n}$. 
 %coordinate and the residue of the box $(r,c,m)$.  
Given $\la\in\mptn \ell n$, we let ${\bf I}_\la$ denote   the disjoint union over the $\mathbf{I}_{(r,c,m)}$ for $(r,c,m)\in \la$.

\begin{defn}
We define a   {\sf diagram} {of type} $G(\ell,1,n)$ to  be a  {\sf frame}
$\mathbb{R}\times [0,1]$ with  distinguished  solid points on the northern and southern boundaries
given by   $\mathbf{I}_\mu$ and $\mathbf{I}_\lambda$ for  
$\lambda,\mu  \in \mptn \ell n$ and a 
collection of {\sf solid strands} each starting at a northern point, ${\bf I}_{(r,c,m)}$ for $(r,c,m)\in \mu$,  
 and ending at a southern
  point, ${\bf I}_{(r',c',m')}$ for $(r',c',m')\in \la$. %,  for $(r,c,m)$ and $(r',c',m')$ two boxes of the same residue, $i\in I$  say, which    we refer to this as a {\sf solid $i$-strand}.  
  Each strand carries a residue, $i \in I$ say, and is referred to as a {\sf solid $i$-strand}.  
 We   require that each solid strand has a 
mapping diffeomorphically to $[0,1]$ via the projection to the $y$-axis.  
Each solid strand is allowed to carry any number of dots.  We draw
\begin{itemize}[leftmargin=*]
\item  a dashed line $\ell$ units to the left of each solid $i$-strand, which we call a {\sf ghost $i$-strand} or $i$-{\sf ghost};
\item vertical red lines at $ m   \in 
\ZZ$ each of which carries a 
residue $\kappa_m$ for $1\leq m\leq \ell$ which we call a {\sf red $\kappa_m$-strand}. 
\end{itemize}
Finally, we require that there are no triple points or tangencies involving any combination of strands, ghosts or red lines and no dots lie on crossings. \end{defn}

\!\!\!
\begin{figure}[ht!]
 \[
\begin{tikzpicture}[baseline, thick,yscale=1,xscale=4]\scalefont{0.8}

  \draw[wei]  (-1.1, -2)  to[out=90,in=-90] (-1.1, 2) ;
  \draw[red] (-1.1, -2.05) [below] node{$\ \kappa_2$};

  \draw[wei]  (-1.2, -2) to[out=90,in=-90] (-1.2, 2);
   \draw[red] (-1.2, -2.05) [below] node{$\ \kappa_1$};
 
 \draw(-1.5,-2) rectangle (1.4,2);

  \draw[gray, densely dotted] (-0.4+-.5,-2) to[out=90,in=-90]   (-0.4+-1,-1)  to[out=90,in=-90](-0.4+1,.3) to[out=90,in=-90]  
    (-0.4+0,2);
  \draw[gray, densely dotted] (-0.4+.5,-2) to[out=90,in=-90]  (-0.4+.5,0) to[out=90,in=-90] (-0.4+1,2);
  \draw[gray, densely dotted]  (-0.4+1,-2) to[out=90,in=-90]  (-0.4+-1,1) to[out=90,in=-90] (-0.4+.5,2);
  \draw[gray, densely dotted] (-0.4+0,-2) to[out=90,in=-90]  
  (-0.4+0.1,0) to[out=90,in=-90]
  (-0.4+-.5,2);
  \draw[gray, densely dotted]  (-0.4+-1, -2) to[out=95,in=-90]  
  (-0.4+-.5,-1.1) to[out=90,in=-90] (-0.4+-1,0) to[out=90,in=-90] (-0.4+-.5,1)
  to[out=90,in=-90]   (-0.4+-1,2);

  \draw (-.5,-2) to[out=90,in=-90] node[below,at start]{$\ i_2$} (-1,-1)  to[out=90,in=-90](1,.2) to[out=90,in=-90]    (0,2);
  \draw (.5,-2) to[out=90,in=-90] node[below,at start]{$\ i_4$} (.5,0) to[out=90,in=-90] node[midway,circle,fill=black,inner sep=2pt]{} (1,2);
  \draw  (1,-2) to[out=90,in=-90]  node[below,at start]{$\ i_5$} (-1,1) to[out=90,in=-90]
   (.5,2);
  \draw (0,-2) to[out=90,in=-90] node[below,at start]{$\ i_3$}node[midway,below,circle,fill=black,inner sep=2pt]{}
  (0.05,0) to[out=90,in=-90]
  (-.5,2);
  \draw  (-1, -2) to[out=90,in=-90] node[below,at start]{$\ i_1$}
  (-.5,-1) to[out=90,in=-90] (-1,0) to[out=90,in=-90] (-.5,1)
  to[out=90,in=-90]   (-1,2);

\end{tikzpicture} 
\]

         \!\!\!\!\!\! 
\caption{A  diagram  for  $\ell=2$  with northern and southern loading
  ${\bf I}_\omega$ where $\omega= (\varnothing,(1^5))$. }
\end{figure}

\begin{defn}[Definition 4.1 \cite{Webster}]\label{defintino2}
We let   ${\bf A}(n,\kappa)$ denote  the  $\Bbbk$-algebra spanned by   all  diagrams modulo the following local relations (here a local relation means one that 
can be applied on a small region of the diagram).  The product $d_1 d_2$ of two diagrams $d_1,d_2 \in A(n, \kappa)$ is given by putting $d_1$ on top of $d_2$.
This product is defined to be $0$ unless the southern border of $d_1$ matches the northern border of $d_2$, in which case we obtain a new diagram in the obvious fashion.  

\begin{enumerate}[label=(2.\arabic*)] 
\item\label{rel1}  Any diagram may be deformed isotopically; that is,
 by a continuous deformation
 of the diagram which %at no point introduces or removes any crossings of 
%strands.    
 %We consider these diagrams equivalent if they are related by an isotopy that 
 avoids  tangencies, double points and dots on crossings. 
\item\label{rel2} 
For $i\neq j$ we have that dots pass through crossings. 
\[   \scalefont{0.8}\begin{tikzpicture}[scale=0.5,baseline]
  \draw[very thick](-4,0) +(-1,-1) -- +(1,1) node[below,at start]
  {$i$}; \draw[very thick](-4,0) +(1,-1) -- +(-1,1) node[below,at
  start] {$j$}; \fill (-4.5,.5) circle (5pt);
    \node at (-2,0){$=$}; \draw[very thick](0,0) +(-1,-1) -- +(1,1)
  node[below,at start] {$i$}; \draw[very thick](0,0) +(1,-1) --
  +(-1,1) node[below,at start] {$j$}; \fill (.5,-.5) circle (5pt);
  \node at (4,0){ };
\end{tikzpicture}\]
 \item\label{rel3}  For two like-labelled strands we get an error term.
\[
\scalefont{0.8}\begin{tikzpicture}[scale=.5,baseline]
  \draw[very thick](-4,0) +(-1,-1) -- +(1,1) node[below,at start]
  {$i$}; \draw[very thick](-4,0) +(1,-1) -- +(-1,1) node[below,at
  start] {$i$}; \fill (-4.5,.5) circle (5pt);
   \node at (-2,0){$=$}; \draw[very thick](0,0) +(-1,-1) -- +(1,1)
  node[below,at start] {$i$}; \draw[very thick](0,0) +(1,-1) --
  +(-1,1) node[below,at start] {$i$}; \fill (.5,-.5) circle (5pt);
  \node at (2,0){$+$}; \draw[very thick](4,0) +(-1,-1) -- +(-1,1)
  node[below,at start] {$i$}; \draw[very thick](4,0) +(0,-1) --
  +(0,1) node[below,at start] {$i$};
\end{tikzpicture}  \quad \quad \quad
\scalefont{0.8}\begin{tikzpicture}[scale=.5,baseline]
  \draw[very thick](-4,0) +(-1,-1) -- +(1,1) node[below,at start]
  {$i$}; \draw[very thick](-4,0) +(1,-1) -- +(-1,1) node[below,at
  start] {$i$}; \fill (-4.5,-.5) circle (5pt);
       \node at (-2,0){$=$}; \draw[very thick](0,0) +(-1,-1) -- +(1,1)
  node[below,at start] {$i$}; \draw[very thick](0,0) +(1,-1) --
  +(-1,1) node[below,at start] {$i$}; \fill (.5,.5) circle (5pt);
  \node at (2,0){$+$}; \draw[very thick](4,0) +(-1,-1) -- +(-1,1)
  node[below,at start] {$i$}; \draw[very thick](4,0) +(0,-1) --
  +(0,1) node[below,at start] {$i$};
\end{tikzpicture}\]
\item\label{rel4} For double-crossings of solid strands with $i\neq j$, we have the following.
\[
\scalefont{0.8}\begin{tikzpicture}[very thick,scale=0.5,baseline]
\draw (-2.8,-1) .. controls (-1.2,0) ..  +(0,2)
node[below,at start]{$i$};
\draw (-1.2,-1) .. controls (-2.8,0) ..  +(0,2) node[below,at start]{$i$};
\node at (-.5,0) {$=$};
\node at (0.4,0) {$0$};
\end{tikzpicture}
\hspace{.7cm}
\scalefont{0.8}\begin{tikzpicture}[very thick,scale=0.5,baseline]
\draw (-2.8,-1) .. controls (-1.2,0) ..  +(0,2)
node[below,at start]{$i$};
\draw (-1.2,-1) .. controls (-2.8,0) ..  +(0,2)
node[below,at start]{$j$};
\node at (-.5,0) {$=$}; 

\draw (1.8,-1) -- +(0,2) node[below,at start]{$j$};
\draw (1,-1) -- +(0,2) node[below,at start]{$i$}; 
\end{tikzpicture}
\]
\end{enumerate}
\begin{enumerate}[resume, label=(2.\arabic*)]  
\item\label{rel5} If $j\neq i-1$,  then we can pass ghosts through solid strands.
\[\begin{tikzpicture}[very thick,xscale=1,yscale=0.5,baseline]
\draw (1,-1) to[in=-90,out=90]  node[below, at start]{$i$} (1.5,0) to[in=-90,out=90] (1,1)
;
\draw[densely dashed] (1.5,-1) to[in=-90,out=90] (1,0) to[in=-90,out=90] (1.5,1);
\draw (2.5,-1) to[in=-90,out=90]  node[below, at start]{$j$} (2,0) to[in=-90,out=90] (2.5,1);
\node at (3,0) {$=$};
\draw (3.7,-1) -- (3.7,1) node[below, at start]{$i$}
;
\draw[densely dashed] (4.2,-1) to (4.2,1);
\draw (5.2,-1) -- (5.2,1) node[below, at start]{$j$};
\end{tikzpicture} \quad\quad \quad \quad 
%%%
%%%
\scalefont{0.8}\begin{tikzpicture}[very thick,xscale=1,yscale=0.5,baseline]
\draw[densely dashed] (1,-1) to[in=-90,out=90] (1.5,0) to[in=-90,out=90] (1,1)
;
\draw (1.5,-1) to[in=-90,out=90] node[below, at start]{$i$} (1,0) to[in=-90,out=90] (1.5,1);
\draw (2,-1) to[in=-90,out=90]  node[below, at start]{$\tiny j$} (2.5,0) to[in=-90,out=90] (2,1);
\node at (3,0) {$=$};
\draw (4.2,-1) -- (4.2,1) node[below, at start]{$i$}
;
\draw[densely dashed] (3.7,-1) to (3.7,1);
\draw (5.2,-1) -- (5.2,1) node[below, at start]{$j$};
\end{tikzpicture}
\]
\end{enumerate} 
\begin{enumerate}[resume, label=(2.\arabic*)]  
\item\label{rel6} On the other hand, in the case where $j= i-1$, we have the following.
\[\scalefont{0.8}\begin{tikzpicture}[very thick,xscale=1,yscale=0.5,baseline]
\draw (1,-1) to[in=-90,out=90]  node[below, at start]{$i$} (1.5,0) to[in=-90,out=90] (1,1)
;
\draw[densely dashed] (1.5,-1) to[in=-90,out=90] (1,0) to[in=-90,out=90] (1.5,1);
\draw (2.5,-1) to[in=-90,out=90]  node[below, at start]{$\tiny i\!-\!1$} (2,0) to[in=-90,out=90] (2.5,1);
\node at (3,0) {$=$};
\draw (3.7,-1) -- (3.7,1) node[below, at start]{$i$}
;
\draw[densely dashed] (4.2,-1) to (4.2,1);
\draw (5.2,-1) -- (5.2,1) node[below, at start]{$\tiny i\!-\!1$} node[midway,fill,inner sep=2.5pt,circle]{};
\node at (5.75,0) {$-$};

\draw (6.2,-1) -- (6.2,1) node[below, at start]{$i$} node[midway,fill,inner sep=2.5pt,circle]{};
\draw[densely dashed] (6.7,-1)-- (6.7,1);
\draw (7.7,-1) -- (7.7,1) node[below, at start]{$\tiny i\!-\!1$};
\end{tikzpicture}\]

\item\label{rel7} We also have the relation below, obtained by symmetry.  \[
\scalefont{0.8}\begin{tikzpicture}[very thick,xscale=1,yscale=0.5,baseline]
\draw[densely dashed] (1,-1) to[in=-90,out=90] (1.5,0) to[in=-90,out=90] (1,1)
;
\draw (1.5,-1) to[in=-90,out=90] node[below, at start]{$i$} (1,0) to[in=-90,out=90] (1.5,1);
\draw (2,-1) to[in=-90,out=90]  node[below, at start]{$\tiny i\!-\!1$} (2.5,0) to[in=-90,out=90] (2,1);
\node at (3,0) {$=$};
\draw[densely dashed] (3.7,-1) -- (3.7,1);
\draw (4.2,-1) -- (4.2,1) node[below, at start]{$i$};
\draw (5.2,-1) -- (5.2,1) node[below, at start]{$\tiny i\!-\!1$} node[midway,fill,inner sep=2.5pt,circle]{};
\node at (5.75,0) {$-$};

\draw[densely dashed] (6.2,-1) -- (6.2,1);
\draw (6.7,-1)-- (6.7,1) node[midway,fill,inner sep=2.5pt,circle]{} node[below, at start]{$i$};
\draw (7.7,-1) -- (7.7,1) node[below, at start]{$\tiny i\!-\!1$};
\end{tikzpicture}\]
\end{enumerate}
\begin{enumerate}[resume, label=(2.\arabic*)]
\item\label{rel8} Strands can move through crossings of solid strands freely.
\[
\scalefont{0.8}\begin{tikzpicture}[very thick,scale=0.5 ,baseline]
\draw (-2,-1) -- +(-2,2) node[below,at start]{$k$};
\draw (-4,-1) -- +(2,2) node[below,at start]{$i$};
\draw (-3,-1) .. controls (-4,0) ..  +(0,2)
node[below,at start]{$j$};
\node at (-1,0) {$=$};
\draw (2,-1) -- +(-2,2)
node[below,at start]{$k$};
\draw (0,-1) -- +(2,2)
node[below,at start]{$i$};
\draw (1,-1) .. controls (2,0) ..  +(0,2)
node[below,at start]{$j$};
\end{tikzpicture}
\]
\end{enumerate}
 for any $i,j,k\in I$.    Similarly, this holds for triple points involving ghosts, except for the following relations when $j=i-1$.
\begin{enumerate}[resume, label=(2.\arabic*)]  \Item\label{rel9}
\[
\scalefont{0.8}\begin{tikzpicture}[very thick,xscale=1,yscale=0.5,baseline]
\draw[densely dashed] (-2.6,-1) -- +(-.8,2);
\draw[densely dashed] (-3.4,-1) -- +(.8,2); 
\draw (-1.1,-1) -- +(-.8,2)
node[below,at start]{$\tiny j$};
\draw (-1.9,-1) -- +(.8,2)
node[below,at start]{$\tiny j$}; 
\draw (-3,-1) .. controls (-3.5,0) ..  +(0,2)
node[below,at start]{$i$};

\node at (-.75,0) {$=$};

\draw[densely dashed] (.4,-1) -- +(-.8,2);
\draw[densely dashed] (-.4,-1) -- +(.8,2);
\draw (1.9,-1) -- +(-.8,2)
node[below,at start]{$\tiny j$};
\draw (1.1,-1) -- +(.8,2)
node[below,at start]{$\tiny j$};
\draw (0,-1) .. controls (.5,0) ..  +(0,2)
node[below,at start]{$i$};

\node at (2.25,0) {$-$};

\draw (4.9,-1) -- +(0,2)
node[below,at start]{$\tiny j$};
\draw (4.1,-1) -- +(0,2)
node[below,at start]{$\tiny j$};
\draw[densely dashed] (3.4,-1) -- +(0,2);
\draw[densely dashed] (2.6,-1) -- +(0,2);
\draw (3,-1) -- +(0,2) node[below,at start]{$i$};
\end{tikzpicture}
\]
\Item\label{rel10}
\[
\scalefont{0.8}\begin{tikzpicture}[very thick,xscale=1,yscale=0.5,baseline]
\draw[densely dashed] (-3,-1) .. controls (-3.5,0) ..  +(0,2);  
\draw (-2.6,-1) -- +(-.8,2)
node[below,at start]{$i$};
\draw (-3.4,-1) -- +(.8,2)
node[below,at start]{$i$};
\draw (-1.5,-1) .. controls (-2,0) ..  +(0,2)
node[below,at start]{$\tiny j$};

\node at (-.75,0) {$=$};

\draw (0.4,-1) -- +(-.8,2)
node[below,at start]{$i$};
\draw (-.4,-1) -- +(.8,2)
node[below,at start]{$i$};
\draw[densely dashed] (0,-1) .. controls (.5,0) ..  +(0,2);
\draw (1.5,-1) .. controls (2,0) ..  +(0,2)
node[below,at start]{$\tiny j$};

\node at (2.25,0) {$+$};

\draw (3.4,-1) -- +(0,2)
node[below,at start]{$i$};
\draw (2.6,-1) -- +(0,2)
node[below,at start]{$i$}; 
\draw[densely dashed] (3,-1) -- +(0,2);
\draw (4.5,-1) -- +(0,2)
node[below,at start]{$\tiny j$};
\end{tikzpicture}
\]
\end{enumerate}
 The ghost strands may pass through red strands freely.
For $i\neq j$, the solid $i$-strands may pass through red $j$-strands freely. If the red and solid strands have the same label, a dot is added to the solid strand when straightening.  Diagrammatically, these relations are given by 
\begin{enumerate}[resume, label=(2.\arabic*)] \Item\label{rel11}
\[
\scalefont{0.8}\begin{tikzpicture}[very thick,baseline,scale=0.5]
\draw[wei] (-1.2,-1)  .. controls (-2.8,0) ..  +(0,2)
node[below,at start]{$i$};
\draw (-2.8,-1)  .. controls (-1.2,0) ..  +(0,2)
node[below,at start]{$i$};
 
\node at (-.3,0) {$=$};

\draw[wei] (2.8,-1) -- +(0,2)
node[below,at start]{$i$};
\draw (1.2,-1) -- +(0,2)
node[below,at start]{$i$};
\fill (1.2,0) circle (5pt);

\draw[wei] (6.8,-1)  .. controls (5.2,0) ..  +(0,2)
node[below,at start]{$j$};
\draw (5.2,-1)  .. controls (6.8,0) ..  +(0,2)
node[below,at start]{$i$};

\node at (7.7,0) {$=$};

\draw (9.2,-1) -- +(0,2)
node[below,at start]{$i$};
\draw[wei] (10.8,-1) -- +(0,2)
node[below,at start]{$j$};
\end{tikzpicture}
\]
\end{enumerate}
for $i\neq j$ and their mirror images. All solid crossings and dots can pass through  red strands, with a
correction term.
\begin{enumerate}[resume, label=(2.\arabic*)]
\Item\label{rel12}
\[
\scalefont{0.8}\begin{tikzpicture}[very thick,baseline,scale=0.5]
\draw[wei] (-3,-1) .. controls (-4,0) ..  +(0,2)
node[at start,below]{$k$};
\draw (-2,-1)  -- +(-2,2)
node[at start,below]{$i$};
\draw (-4,-1) -- +(2,2)
node [at start,below]{$j$};
 
\node at (-1,0) {$=$};
\draw[wei] (1,-1) .. controls (2,0) .. +(0,2)
node[at start,below]{$k$};

\draw (2,-1) -- +(-2,2)
node[at start,below]{$i$};
\draw (0,-1) -- +(2,2)
node [at start,below]{$j$};
 
\node at (2.8,0) {$+ $};
\draw[wei] (6.5,-1) -- +(0,2)
node[at start,below]{$k$};

\draw (7.5,-1) -- +(0,2)
node[at start,below]{$i$};
\draw (5.5,-1) -- +(0,2)
node [at start,below]{$j$};
 \node at (3.8,-.2){$\delta_{i,j,k}$};
\end{tikzpicture}
\]
\Item\label{rel13}
\[
\scalefont{0.8}\begin{tikzpicture}[scale=0.5,very thick,baseline=2cm]
\draw[wei] (-2,2) -- +(-2,2);
\draw (-3,2) .. controls (-4,3) ..  +(0,2);
\draw (-4,2) -- +(2,2);

\node at (-1,3) {$=$};

\draw[wei] (2,2) -- +(-2,2);
\draw (1,2) .. controls (2,3) ..  +(0,2);
\draw (0,2) -- +(2,2);

\draw (6,2) -- +(-2,2);
\draw (5,2) .. controls (6,3) ..  +(0,2);
\draw[wei] (4,2) -- +(2,2);

\node at (7,3) {$=$};

%\draw (10,2) -- +(-2,2);
%\draw (9,2) .. controls (10,3) ..  +(0,2);
\draw[wei] (8,2) -- +(2,2);

\draw (9,2) .. controls (8,3) ..  +(0,2);
\draw (10,2) -- +(-2,2);

\end{tikzpicture}
\]

\Item\label{rel14}
\[
\scalefont{0.8}\begin{tikzpicture}[very thick,baseline,scale=0.5]
\draw[wei](-3,0) +(1,-1) -- +(-1,1);\draw(-3,0) +(-1,-1) -- +(1,1);
\draw[wei](1,0) +(1,-1) -- +(-1,1);

\fill (-3.5,-.5) circle (5pt);
\node at (-1,0) {$=$};
\draw(1,0) +(-1,-1) -- +(1,1);
 \fill (1.5,.5) circle (5pt);

\draw[wei](1,0) +(3,-1) -- +(5,1);
\draw(1,0) +(5,-1) -- +(3,1);
\fill (5.5,-.5) circle (5pt);
\node at (7,0) {$=$};
\draw[wei](5,0) +(3,-1) -- +(5,1);
\draw(5,0) +(5,-1) -- +(3,1);
\fill (8.5,.5) circle (5pt);
\end{tikzpicture}
\]
\end{enumerate}
Finally, we have the following non-local idempotent relation.
\begin{enumerate}[resume, label=(2.\arabic*)]
\item\label{rel15}
Any idempotent in which a solid strand is $\ell n$ units to the left of the leftmost red-strand is referred to as unsteady and set to be equal to zero.
\end{enumerate}
 \end{defn}

\begin{rmk}\label{agrade}
The algebra ${\bf A}(n, \kappa)$  admits a $\ZZ$-grading \cite[Defintion 4.2]{Webster}\label{grsubsec}.  
We do not recall this explicitly here, but in \cref{sec3} we  shall encode this grading using the path combinatorics of \cite{CoxBowman}.  
We let $t$ be an indeterminate over $\ZZ_{\geq 0}$.
If $M=\oplus_{k \in \ZZ} M_k$  is a graded $\Bbbk$-module, we  write   
$\dim_t(M)= \sum_{k \in \ZZ}(\dim_\Bbbk(M_k))t^k$.  
\end{rmk}

    \begin{defn} \label{semistandard:defn} Let $\lambda,\mu \in \mptn {\ell}n$.  A $\lambda$-tableau of weight $\mu$  is a
  bijective map $\SSTT :  \lambda  \to  {\bf I}_\mu$. %from the nodes of the Young diagram  of $\la$ to those of $\mu$. 
   We let $\mathcal{T}(\lambda,\mu)$ denote the set of all tableaux of shape $\lambda$ and weight $\mu$.  
%
% which respects residues.  
%  In other words,  if $
%  \SSTT(r,c,m)= \mathbf{i}_{(r',c',m')} \in {\mathbf i}_\mu$ for   $(r,c,m) \in  \lambda  $ 
%  and $(r',c',m')\in \mu$, then  $\kappa_m+c-r = \kappa_{m'}+c'-r' \pmod e$.    
%    \end{defn}
%
%
% \begin{defn}
   We say that a tableau $\SSTT$ is {\sf semistandard} if it  satisfies the following additional properties   \begin{itemize}
\item[$(i)$]     $\SSTT(1,1,m)> m$,
\item[$(ii)$]    $\SSTT(r,c,m)> \SSTT(r-1,c,m)  +\g$,
\item[$(iii)$] 
$\SSTT(r,c,m)> \SSTT(r,c-1,m) -\g$.   
\end{itemize}
We  denote the set of all  semistandard tableaux of shape $\lambda$
and weight $\mu$ by $\SStd(\lambda,\mu)$.   Given $\SSTT \in 
\SStd(\lambda,\mu)$, we   write $\Shape(\SSTT)=\lambda$.  
 We let $\SStd^+_n(\lambda,\mu) \subseteq \SStd_n(\lambda,\mu) $ denote the subset of 
 tableaux which respect residues.  
  In other words,  if $
  \SSTT(r,c,m)= {(r',c',m')} $ for   $(r,c,m) \in  \lambda  $ 
  and $(r',c',m')\in \mu$, then  $\kappa_m+c-r = \kappa_{m'}+c'-r' \pmod e$.    

    \end{defn}

 Associated to $\SSTT\in\mathcal{T}(\la,\mu)$ we have a  diagram $C_\SSTT$   with 
   solid points on the  northern and southern  borders given by ${\bf I}_\mu$ and ${\bf I}_\lambda$ respectively;  
 the $n$ solid strands each connect  a northern and southern   point and  
  they trace out the bijection determined by $\SSTT$  using the minimal number of
  crossings   (this can be chosen arbitrarily);  
 the strand terminating at southern point   ${(r,c,m)} $ for $(r,c,m) \in \lambda$ carries residue equal to  $\res(r,c,m)\in I$.  
 We  let ${C}_{\SSTS  \SSTT}=C  _\SSTS C^\ast_\SSTT$ where $C^\ast_\SSTT$
is the diagram obtained from $C_\SSTT $ by flipping it through the
horizontal axis.    
 
 Given $\lambda\in \mptn \ell n$ and  $\underline{i} \in I^n$, we have an associated idempotent 
 ${\sf 1}_\la^{ \underline{i}}$ given by the diagram with northern and southern points 
 $\bf{I}_\la$, no crossing strands,  and northern (or equivalently southern) residue sequence of the diagram given by  $\underline{i} \in I^n$.    
 If the residue sequence is equal to that of the partition, ${ \res(\la)}$,  then we let ${\sf 1}_\la:=
  {\sf 1}_\la^{ \res(\la)}$.  
%  \begin{defn}
 We define the {\sf diagrammatic Cherednik algebra}, $A(n,\kappa)$ to be the  algebra
 $$
 A(n,\kappa):=
{\sf E}^+
 {\bf A}(n,\kappa)
{\sf E}^+ \quad \text{where} \quad {\sf E}^+=\textstyle \sum_{\la \in \mptn \ell n} {\sf 1}_\la . $$
%  \end{defn}

% 
% \begin{defn}
% We let $\SStd^+_n(\lambda,\mu) \subseteq \SStd_n(\lambda,\mu) $ denote the subset of 
% tableaux which respect residues.  
%  In other words,  if $
%  \SSTT(r,c,m)= {(r',c',m')} $ for   $(r,c,m) \in  \lambda  $ 
%  and $(r',c',m')\in \mu$, then  $\kappa_m+c-r = \kappa_{m'}+c'-r' \pmod e$.    
% \end{defn}
%The following theorem  might seem  overstated; however, 
 
\begin{thm} \label{cellularitybreedscontempt-2}
The $R$-algebra   $ {\bf A} (n,\kappa)$ %a%nd  $ A (n,\kappa)$  are 
is a quasi-hereditary graded cellular algebra 
 with  basis 
\[
 \{C_{\SSTS  \SSTT} \mid \SSTS \in \SStd (\lambda,\mu), \SSTT\in \SStd  (\lambda,\nu), 
\lambda,\mu, \nu \in \mptn {\ell}n\}
\]
and the subalgebra $A(n,\kappa)$ is a quasi-hereditary graded cellular algebra 
 with  basis 
 \[
 \{C_{\SSTS  \SSTT} \mid \SSTS \in \SStd^+(\lambda,\mu), \SSTT\in \SStd^+ (\lambda,\nu), 
\lambda,\mu, \nu \in \mptn {\ell}n\}.  
\]
%with respect to cell-data $((\mptn \ell n,\rhd), \ast, \SStd, C) $  and $((\mptn \ell n,\rhd), \ast, \SStd^+, C) $ respectively.  
For both algebras, the involution is given by $\ast$ and the ordering on $\mptn \ell n(h)$ is that of \cref{domdef}.  
   We denote the corresponding left cell-modules for  $ {\bf A} (n,\kappa)$ and  $ A (n,\kappa)$  by  $   {\bm\Delta}(\lambda) $   and    $   {\Delta}(\lambda) $ 
%   $$
% {\bm\Delta}(\lambda) = \{ C_\SSTS \mid \SSTS \in \SStd  (\lambda, \mu ), \mu \in \mptn \ell n \}  \quad
%   {\Delta}(\lambda) = \{ C_\SSTS \mid \SSTS \in \SStd^+ (\lambda, \mu ), \mu \in \mptn \ell n \}
%   $$
respectively.  The modules    $   {\bm\Delta}(\lambda) $   and    $   {\Delta}(\lambda) $  have    simple heads, denoted  by  ${ \bf L}(\la) $ and ${  L}(\la) $ respectively. % denote the simple head of this module. %quotient by the radical of the cellular bi-linear form.  
 \end{thm}

 The standard  tableaux $\Std (\lambda)$ form the predictable  subset of semistandard tableaux of weight $\omega=(1^n)$ as follows.  Let  $\lambda \in \mptn \ell n$.
If $\sts \in \Std (\lambda)$ is such that  $\sts(r_k,c_k,m_k)=k$ for $1\leq k \leq n$, then 
we let $\SSTS \in\mathcal{T}(\lambda,\omega)$ denote the tableau 
  $\SSTS: \la \to \omega$   determined by $\SSTS(r_k,c_k,m_k)={\bf I}_{(k,1,\ell)}$.
We have a 
 bijective map 
$ 
\varphi : \Std (\lambda) \rightarrow  \mathcal{T}(\lambda,\omega) .
$ 
given by $\varphi(\sts)=\SSTS$.

\begin{defn} We define the    {\sf Schur  idempotent}, ${\sf E}_\omega$, and 
  quiver Hecke algebra, $R_n(\kappa)$, as follows 
 $$
 {\sf E}_\omega= \sum_{\underline{i} \in I^n}{\sf 1}^{\underline{i}}_{\omega} 
\quad \text{ and } \quad % $$
% We define the quiver Hecke algebra, $R_n(\kappa)$,  to be the algebra
%$$
R_n(\kappa):={\sf E}_\omega {\bf A}(n,\kappa) {\sf E}_\omega.  
$$
\end{defn}
%When proving results on homomorphisms, the algebra $ A (n,\kappa)$ is smaller than  $ {\bf A} (n,\kappa)$ and much easier for computation.  We shall then (trivially) inflate these results to $ {\bf A} (n,\kappa)$ and apply the Schur functor to obtain the corresponding result for $R_n(\kappa)$.  

%\begin{defn}  
%Given $\lambda\in \mptn \ell n$, we define a {\sf tableau} of shape $\lambda$ to be a filling of the boxes of 
%the  Russian array of  $\la $ with the numbers 
%$\{1,\dots , n\}$.  We define a {\sf  standard tableau} to be a tableau  in which    the entries increase along the rows and columns of each component.  
%We let $\Std (\lambda)$ denote the set of all standard tableaux of shape $\lambda\in\mptn\ell n$. 
%Given 
%$\stt\in \Std (\lambda)$, we set $\Shape(\stt)=\la$.  
%Given $1\leq k \leq n$, we let $\stt{\downarrow}_{\{1,\dots ,k\}}$ be the subtableau of $\stt$ whose  entries belong to the set
%$\{1,\dots,k\}$.  \end{defn}
%

\begin{thm}[\cite{manycell}]
 The algebra  $R_n(\kappa)$   admits a  graded cellular  structure with respect to the poset $(\mptn \ell n,\rhd) $, the   basis
   $$ 
   \{ c_{\sts\stt}:=C_{ \varphi(\sts),\varphi(\stt)}  \mid \lambda \in \mptn \ell n, \sts,\stt \in \Std (\lambda)\},
   $$  
 and the involution $\ast$.     We denote the   left cell-module by 
   $    S_n(\lambda) = \{ c_\sts \mid \sts \in \Std (\lambda)\}
   $ 
   % and let ${  D}_n(\la)$ denote the simple head of $  S_n(\lambda)$.  % the quotient by the radical of the cellular bi-linear form.  
%The module $     S_n(\lambda) $ is cyclically generated by $c_{\stt^\lambda}$ and has simple head, ${  D}_n(\la)$.  
\end{thm} 
 
 When proving results on homomorphisms, the algebra $ A (n,\kappa)$ is smaller than  $ {\bf A} (n,\kappa)$ and much easier for computation.  We shall then (trivially) inflate these results to $ {\bf A} (n,\kappa)$ and apply the Schur functor to obtain the corresponding result for $R_n(\kappa)$.

  We have an embedding   $R_{n}(\kappa)\hookrightarrow R_{n+1}(\kappa)$  and $R_{n+1}(\kappa)$ is free as a $R_{n}(\kappa)$-module.  
%If $M$ is a $R_{n}(\kappa)$-module, 
We let
\[
{\rm res}^{n+1}_{n}: R_{n+1}(\kappa) \text{-}{\rm mod}
\to 	R_{n}(\kappa)\text{-}{\rm mod}
\qquad
{\rm ind}^{n+1}_{n } :R_n(\kappa) \text{-}{\rm mod}
\to 	R_{n+1}(\kappa)\text{-}{\rm mod} \] 
denote the obvious restriction and induction functors.   We   post-compose these functors  with the projection onto a block   in the standard fashion.  This amounts to multiplying by an idempotent  
$$
{\sf E}_r  = \sum_{\underline i = (i_1, \dots, i_{n-1},r )} {\sf 1}^{\underline i}_\omega
$$
for $r\in I$.  We hence   decompose these  restriction/induction functors into  
 \[
%{\rm res}^{n}_{n-1} = \bigoplus_{r \in I} r\text{-}{\rm res}^{n}_{n-1}
%\quad
%{\rm ind}^{n}_{n-1} = \bigoplus_{r \in I} r\text{-}{\rm ind}^{n}_{n-1}
%\qquad 
 r\text{-}{\rm res}^{n+1}_{n} = {\sf E}_r \circ  {\rm res}^{n+1}_{n}\quad
 r\text{-}{\rm ind}^{n}_{n-1} = {\sf E}_r \circ  {\rm ind}^{n}_{n-1}
\]
Finally we recall  a simple case of  \cite[Theorem 6.1]{manycell}.  
% Let  $\lambda \in \mptn \ell n$ and let $A_1\rhd  A_2 \rhd \dots \rhd A_z$ denote the $r$-removable boxes of $\lambda$, totally ordered according to the dominance ordering.  
%  The module 
%$r\text{-}{\rm res}^{n}_{n-1}( S_n(\la))$  has an $R_{n-1}(\kappa)$-module filtration 
%\begin{align}\label{restriction}
%0= S^{z+1,\lambda} \subset 
%S^{z,\lambda}
%\subset \dots 
%\subset S^{1,\lambda}= r\text{-}{\rm res}^{n}_{n-1} ( S_n(\la))
%\end{align}such that 
%\begin{align}\label{restriction2}
%S (\la-A_r)\langle  	\rangle  
%\cong
%S^{r,\lambda} / S^{r+1,\lambda} 
%\end{align} 
%for $1\leq r\leq z$.  In particular, 
If  $\la$ has precisely one removable $r$-box, $\square \in {\rm Rem}_r(\la)$, then we set 
${\sf E}_r(\la) = \la -\square$ and we have that $ r\text{-}{\rm res}^{n}_{n-1} (S_n(\la))= S_{n-1}( \la-\square)$.

% \subsection{Morita equivalences}\label{MOO}
 Finally, we briefly recall from \cite[Theorem 4.9]{CoxBowman} that for the above three algebras, there are graded Morita equivalences relating the subcategories of 
 ${\bf  A}(n,\kappa)\text{-}{\rm mod}$,  
  ${A}(n,\kappa)\text{-}{\rm mod}$, 
  and 
   $R_n(\kappa)\text{-}{\rm mod}$  
whose simple constituents are labelled by  $\mptn \ell n(h)$.  
   In particular
\begin{equation}\label{MOO}
\Hom_{A(n,\kappa)}
({  {\Delta}}(\mu), 
 {\Delta}(\la))
\cong 
\Hom_{{\bm A}(n,\kappa)}
({\bf \Delta}(\mu), {\bm \Delta}(\la))
\cong 
\Hom_{R_n(\kappa)}
(S_n(\mu), S_n(\la))
\end{equation}
   for $\la,\mu \in \mptn \ell n(h)$.   
  This   allows  us to focus on  the algebra 
  $A (n,\kappa)$ where we have the benefit of a highest weight theory which is intimately related to the underlying alcove geometry.  
Both isomorphisms are simply given by idempotent truncation from ${\bf    A}(n,\kappa)$.

 \section{Alcove geometries and bases of diagrammatic algebras}  \label{sec3}
 \renewcommand{\sts}{\mathsf{S}}   
\renewcommand{\stt}{\mathsf{T}}

In \cref{3.1}, we recall the alcove geometry controlling the subcategories  of representations for 
  quiver Hecke and Cherednik algebras of interest in this paper.
%   For  $h$-admissible $\kappa\in I^\ell$,   we     define the set of fundamental $\ell$-partitions in the context of this geometry.
  In \cref{3.2} we cast the semistandard tableaux  for diagrammatic Cherednik algebras in this geometry; 
  this path-combinatorial framework  will be essential for our proofs. 
    In \cref{3.3} we cast the  standard tableaux combinatorics of the quiver Hecke algebra in this geometry ---
     this allows us to define the 
    $e$-restricted tableaux which we will  prove provide bases of simple   $R_n(\kappa)$-modules in \cref{section4}.

 \subsection{The alcove geometry}\label{3.1}

Fix integers $h,\ell \in \ZZ_{>0}$ and $e \geq h\ell$.
For each $ 1\leq \I \leq h$ and $ 0\leq \M < \ell$ we let
$\varepsilon_{h\M+\I}$ denote a
formal symbol, and define an   $\ell h$-dimensional real vector space, 
\[
{\mathbb E}_{h,\ell }
=\bigoplus_{
\begin{subarray}c 
1\leq \I \leq  h   \\ 
0 \leq \M <  \ell  
 \end{subarray}
} \mathbb{R}\varepsilon_{h\M+\I}
\]
%to be the associated $\ell h$-dimensional real vector space.
We have an inner product $\langle \, , \, \rangle$ given by extending
linearly the relations 
\[
\langle \varepsilon_{h\M+\I} , \varepsilon_{ht+\J} \rangle= 
\delta_{\I,\J}\delta_{t,\M}
\]
for all $1\leq \I, \J \leq h$ and $0 \leq  \M , t< \ell $, where
$\delta_{i,j}$ is the Kronecker delta.  We let $\Phi$ and 
$\Phi_0$  denote the root systems of type $A_{\ell h -1}$ and  $A_{h -1} \times A_{h -1}  \times \dots A_{h -1} $ respectively which  consist  of the
roots 
$$\{\varepsilon_{h\M+ \I}-\varepsilon_{ht+ \J}\mid 1\leq \I,\J
\leq h \text{ and } 0\leq \M,t < \ell\text{ with } (i,m)\neq (j,t)\} \text{ and }$$
$$\{\varepsilon_{h\M+ \I}-\varepsilon_{h\M+ \J}\mid 1\leq \I,\J \leq h
\text{ with } \I\neq \J \text{ and } 0\leq \M < \ell\}$$
 respectively.   
    We identify $\lambda \in \mptn \ell n(h)$ with a point  in $\mathbb{E}_{h\ell} $ via the transpose   map 
\begin{equation}\label{lhfskadj} ( \lambda^{(1)} ,\ldots, \lambda^{(\ell )} )\mapsto \sum_{\begin{subarray}c
1\leq i \leq \ell   	\\
1\leq j \leq h
\end{subarray}}( \lambda^{(\M)})^T_\I \varepsilon_{h(m-1)+i},
\end{equation}
(where the $T$ denotes the transpose partition).  
   Given    $r\in \ZZ$ and  $\alpha \in \Phi$ we let  
  $ s_{\alpha,re }$ denote the reflection which acts on ${\mathbb E}_{h,\ell}$ by 
  \[s_{\alpha, re}x =x-(\langle x, \alpha \rangle -re) \alpha\]
  \noindent and we let $W^e $ and $W^e_0$  denote the groups 
  generated by the reflections 
  $$\mathcal{S} = \{s_{\alpha, re} \mid   \alpha \in \Phi, r\in\mathbb{Z}\} \qquad \mathcal{S}_0=\{s_{\alpha, 0} \mid   \alpha \in \Phi_0\} $$ respectively. 
%  If $e=\infty$ then we let $W^e$ be the finite reflection group
%  generated by the reflections 
%  $$\mathcal{S} = \{s_{\alpha, 0} \mid   \alpha \in \Phi\}$$
%  and let $W^e_0$  denote the   parabolic subgroup generated by 
%  $$\mathcal{S}_0=\{s_{\alpha, 0} \mid   \alpha \in \Phi_0\}.$$
   For $e\in\mathbb{Z}_{>0}$  we assume that 
 $\kappa \in I^\ell $ is $h$-admissible.    
We shall consider a shifted action of the Weyl group $W ^e$  on ${\mathbb E}_{hl}$ %^{\circledcirc}$
by the element
 $$\rho:= (\rho_{1}, \rho_{2}, \ldots, \rho_{\ell}) \in I^{h\ell}, \; \rho_i := (e-\kappa_i, e-\kappa_i-1, \dots,e-\kappa_i-h+1)\in I^{h},$$   %(e^{h\ell}) 
 %and $\rho = (\rho_1,\rho_2, \ldots,\rho_\ell) \in I^{h }.$  
%$$\rho=(e,e-1,e-2, \dots , 1)$$
\noindent that is, given an an element $w\in W  ^e$,  we set 
$
w\cdot_{\rho}x=w(x+\rho)-\rho.
$
%and define the ``dot'' action of $s_{\alpha,re}$ on ${\mathbb
%  E}_{h,\ell}$ by 
%$$s_{\alpha,re}\cdot x=s_{\alpha,re}(x+\rho)-\rho.$$
%(RELIES ON DEFINITION OF $\rho$)
We let ${\mathbb E} ({\alpha, re})$ denote the affine hyperplane
consisting of the points  
$${\mathbb E} ({\alpha, re}) = 
\{ x\in{\mathbb E}_{h,\ell} \mid  s_{\alpha,re} \cdot x = x\} .$$
Note that our assumption that $\kappa\in I^\ell$ is $h$-admissible implies that the origin does not lie on any hyperplane.   Given a
hyperplane ${\mathbb E} ( \alpha,re)$ we 
remove the hyperplane from ${\mathbb E}_{h,\ell}$ to obtain two
distinct subsets ${\mathbb 
  E}^{\great}(\alpha,re)$ and ${\mathbb E}^{\less}(\alpha,re)$
where the origin $\circledcirc \in {\mathbb E}^{\less }(\alpha,re)$.  
%We define
%${\mathbb E}^{\greatoreq }(\alpha,re) = {\mathbb E} (\alpha,re)\cup
%{\mathbb E}^{\great }(\alpha,re)$ and similarly ${\mathbb
%  E}^{\lessoreq }(\alpha,re) = {\mathbb E} (\alpha,re)\cup {\mathbb
%  E}^{\less}(\alpha,re)$. 
  The dominant Weyl chamber, denoted
${\mathbb E}_{h,\ell }^\circledcirc $, is set to be 
$${\mathbb E}_{h,\ell }^{\circledcirc}=\bigcap_{ \begin{subarray}c
 \alpha \in \Phi_0
 \end{subarray}
 } {\mathbb E}^{\less} (\alpha,0).$$

	\newcommand{\Po}{\operatorname{Po}}
\begin{defn}
Let $\lambda \in {\mathbb E }_{h,\ell} $.     There are only finitely many hyperplanes   
lying   between  the   point  $\lambda  \in {\mathbb E} _{h,\ell}$
and the point  $\nu \in {\mathbb E} _{h,\ell}$.  
For  a root $\varepsilon_i - \varepsilon_j  \in \Phi$, we let $\ell_{\varepsilon_i-\varepsilon_j}  (\la,\nu )$ denote the total number of these hyperplanes which are perpendicular to $\varepsilon_i - \varepsilon_j\in \Phi$ (including any hyperplanes upon which $\lambda$ or $\nu$  lies).
  We let $\ell (\la,\nu )= \sum_{1\leq i < j\leq h\ell}\ell_{\varepsilon_i - \varepsilon_j}(\la,\nu)$.  
  We let $\ell(\lambda):=\ell(\la,\circledcirc)$ for $\circledcirc$ the origin and refer to $\ell(\la)$ simply as the {\sf length} of the point $\la \in {\mathbb E} _{h,\ell}$.    
\end{defn}

%  Let $\lambda,\nu \in {\mathbb E }_{h,\ell} $ and $1\leq i < j\leq h\ell$.     There are only finitely many hyperplanes   
%parallel to ${\mathbb E}^{\less} (\varepsilon_i-\varepsilon_j,0)  $ lying   between  these   point   $\lambda  \in {\mathbb E} _{h,\ell}$
%%and the point  $\nu\in {\mathbb E} _{h,\ell}$.  
%We let $\ell_{\varepsilon_i-\varepsilon_j}  (\la )$ denote the total number of these hyperplanes (including any hyperplanes upon which $\lambda$ or $\nu$ 
%lies).  We let $\ell(\lambda):=\ell(\la,\astrosun)$ for $\astrosun$ the origin and refer to $\ell(\la)$ simply as the {\sf length} of the point $\la \in {\mathbb E} _{h,\ell}$.    

\begin{defn} Given $n\in \mathbb{N}$, we define the {\sf fundamental alcove} %  or  {\sf unitary alcove}
  to be 
 $$\mathcal{F}^\ell_n(h)=\{\lambda \in \mptn \ell n \mid \ell(\lambda)=0\}. $$
%There is at most one  $\lambda\in \mathcal{F}^\ell_n(h)$  in any given block and so for $\la \in \mathcal{F}^\ell_n(h)$, 
%we set
%$$\Po_e(\lambda)= \widetilde{\mathfrak{S}}_{h\ell} \cdot \la =\{\mu \in \mptn \ell n(h)\mid  \lambda \trianglerighteq \mu\}.$$
  \end{defn}

  The name ``fundamental alcove" is clearly inspired by classical Lie theory.  
  We stress that    Specht modules from the fundamental alcove can become arbitrarily complicated and that understanding 
the decomposition numbers   $d_{\lambda\mu}=[S_n(\lambda):D_n(\mu)]$ for $\lambda \in \mathcal{F}^\ell_n(h)$  and $\mu \lhd \lambda$
is, in general, a hopelessly difficult task   (see \cref{aremark,aremark2}).  
 This might surprise classical Lie theorists, but this is because we are working in the Ringel dual setting, see \cite{CoxBowman} for more details.

%\subsection{An abaci theoretic description for the Cherednik algebra of the symmetric group} 

% \section{Paths in the geometry}
%    We now consider paths in our space ${\mathbb E}_{h,\ell}$ and define a
%degree function on such paths in terms of the hyperplanes in our
%geometry.
%

\begin{defn}  
 Given   a map  $s: \{1,\dots
  , n\}\to \{1,\dots , \ell h\}$ we define points $\sts(k)\in
  {\mathbb E}_{h,\ell}$ by
$
\sts(k)=\sum_{1\leq i \leq k}\varepsilon_{s(i)} 
$ 
for $1\leq k \leq n$. We define the associated path of length $n$ in
 ${\mathbb E}_{h,\ell}$ by $$\sts=(\sts(0),\sts(1),\sts(2), \ldots, \sts(n)),$$  where we fix all
paths to begin at the origin, so that $\sts(0)=\circledcirc \in
{\mathbb E}_{h,\ell}$.  We let $\sts_{\leq \ik}$ denote the subpath
of $\sts$ of length $\ik$ corresponding to the restriction of the
map $s$ to the domain $\{1,\dots , k\} \subseteq \{1,\dots , n\} $.
We let $\Shape(\SSTS)$ denote the point in $  {\mathbb E}_{h,\ell}$ at which $\SSTS$ terminates.  
   \end{defn}

\begin{rmk}\label{rmkrmk}
Let $\sts $ be a path which passes through a hyperplane ${\mathbb
  E}_{\alpha,re}$ at point $\sts(\ik)$ (note that $\ik$ is not
necessarily unique).  Let $\stt$ be the path obtained from
$\sts$ by applying the reflection $s_{\alpha,re}$ to all the steps in
$\sts $ after the point $\sts(\ik)$.  In other words,
$\stt(i)=\sts(i)$ for all $1\leq i \leq \ik$ and
$\stt(i)=s_{\alpha,re} \cdot \sts(i)$ for $\ik \leq i \leq n$.  We
refer to the path $\stt$ as the reflection of $\sts$ in ${\mathbb
  E}_{\alpha,re}$ at point $\sts(\ik)$ and denote this by
$s_{\alpha,re}^{\ik}\cdot \sts$.  We write $\sts \sim \stt$ if the
path $\stt$ can be obtained from $\sts$ by a series of such
reflections.
\end{rmk}

\begin{defn}
Let $\stt$ denote a fixed path from $\circledcirc$ to $\nu\in \mptn \ell n(h)$.  
  We let $\Path_n(\lambda,\stt )$
  denote the set of  paths from the origin to $\lambda$ obtainable by applying repeated reflections to $\stt $, in other
  words
 $$\Path_n(\lambda,\stt ) = \{\sts \mid \sts(n)=\lambda, \sts \sim
  \stt \}.$$ 
We let 
 $\Path_n^+(\lambda,\stt)
 \subseteq \Path_n(\lambda,\stt) 
 $ denote the  set of  paths  which at no point 
 leave the dominant Weyl chamber, in other words 
 $$\Path_n^+(\lambda,\stt)
 =
  \{\sts \in \Path_n(\lambda,\stt)
  \mid \sts(k) \in {\mathbb E}_{h,\ell}^\circledcirc
  \text{ for all }1\leq k \leq n
  \}.$$
 \end{defn}

\begin{defn}\label{Soergeldegreee}
Given a path $\sts=(\sts(0),\sts(1),\sts(2), \ldots, \sts(n))$, we define $d(\sts(k),\sts(k-1))$ as follows. 
For $\alpha\in\Phi$ we set $d_{\alpha}(\sts(k),\sts(k-1))$ to be
\begin{itemize}
\item $+1$ if $\sts(k-1) \in 
   {\mathbb E}(\alpha,re)$ and 
   $\sts(k) \in 
   {\mathbb E}^{\less}(\alpha,re)$;
   
\item $-1$ if $\sts(k-1) \in 
   {\mathbb E}^{\great}(\alpha,re)$ and 
   $\sts(k) \in 
   {\mathbb E}(\alpha,re)$;
\item $0$ otherwise.  
   \end{itemize}
We set 
$\deg(\sts(0))=0$ and define  
 \[
d (\sts(k-1),\sts(k))= \sum_{\alpha \in \Phi}d_\alpha(\sts(k-1),\sts(k))\quad \text{ and } \quad \deg(\sts ) = \sum_{1\leq k \leq n} d (\sts(k),\sts(k-1)).
 \]
\end{defn}

    \begin{figure}[ht] $$
 \scalefont{0.8}  \begin{tikzpicture}[scale=1.2]%[scale=1.4]       
     \clip (-3.2,-0.3) rectangle ++(7.4,7.7);
      \path (120:2.4cm)++(60:2.4cm) coordinate (BB);                    
       \path (0,0) coordinate (origin);
         \path   (origin) ++ (90:6.2 cm)     coordinate (BOOM);   
                  \path   (BOOM) ++(0:2.8 cm)     coordinate (BOOM1);                
  \draw[->](BOOM1) to  ++ (0:0.5 cm);
  \draw[->](BOOM1) to  ++ (120:0.5 cm);
   \draw[->](BOOM1) to  ++ (240:0.5 cm);
 \path(BOOM1)    ++ (0:0.7 cm) coordinate (BOOM2);
 \draw (BOOM2) node {$\varepsilon_2$};
  \path(BOOM1)    ++ (120:0.7 cm) coordinate (BOOM3);
 \draw (BOOM3) node {$\varepsilon_1$};
  \path(BOOM1)    ++ (240:0.7 cm) coordinate (BOOM4);
 \draw (BOOM4) node {$\varepsilon_3$};

   \path[dotted] (origin)++(60:6.4cm) coordinate (aa);
      \path[dotted] (origin)++(60:4.8cm) coordinate (ab);
      
            \path[dotted] (origin)++(60:4.8cm)++(120:1.6) coordinate (ba);
   \path[dotted] (origin)++(120:6.4cm) coordinate (bb);
   \path[dotted] (origin)++(120:4.8cm) coordinate (bb1);
   \path[dotted] (bb1)++(60:3.2cm) coordinate (bb2);
      \path[dotted] (bb2)++(180:1.6cm) coordinate (bb25);
            \path[dotted] (ab)++(0:1.6cm) coordinate (bb325);
      \path[dotted] (bb25)++(-60:1.6cm) coordinate (bb252);

      \path[dotted] (bb2)++(-60:1.6cm) coordinate (bb3);
      \path[dotted] (origin)++(120:3.2cm) coordinate (cc);
             \path[dotted] (cc)++(-120:1.6cm) coordinate (cc2);
                          \path[dotted] (cc2)++(-60:1.6cm) coordinate (cc88);
                        \path[dotted] (cc2)++(+0:1.6cm) coordinate (cc3);
            \path[dotted] (origin)++(120:1.6cm)++(60:1.6cm) coordinate (dd);
                        \path[dotted] (origin)++(60:3.2cm)++(0:1.6cm) coordinate (de);

            \path[dotted] (dd)++(-120:1.6cm)  coordinate (dd1);
            \path[dotted] (dd1)++(0:1.6cm)  coordinate (dd2);
                        \path[dotted] (dd)++(0:1.6cm)  coordinate (dd3);

   \draw[thick] (bb325)-- (ab)--(ba)--(bb3)--(bb2)--(bb25)--(bb252)--(bb1)--(cc)--(cc2)--(cc88)--(cc3)--(dd1)--(dd2)--(dd3)
   --(de)--(bb325);
   \clip (bb325)-- (ab)--(ba)--(bb3)--(bb2)--(bb25)--(bb252)--(bb1)--(cc)--(cc2)--(cc88)--(cc3)--(dd1)--(dd2)--(dd3)
   --(de)--(bb325);

      \path[dotted] (ab)++(0:0.4  cm) coordinate (XV);
            \path[dotted] (dd)++(0:0.4+2  cm) coordinate (YV);
\draw[thin,gray!30](XV)--(YV); 

  \path[dotted] (ab)++(0:0.8  cm) coordinate (XV);
            \path[dotted] (dd)++(0:0.8+2.4  cm) coordinate (YV);
\draw[thin,gray!30](XV)--(YV);

  \path[dotted] (ab)++(0:1.2  cm) coordinate (XV);
            \path[dotted] (dd)++(0:2.8  cm) coordinate (YV);
\draw[thin,gray!30](XV)--(YV);

  \path[dotted] (ab)++(0:1.6  cm) coordinate (XV);
            \path[dotted] (dd)++(0:3.2  cm) coordinate (YV);
\draw[thick](XV)--(YV); 
%
% \path[dotted] (ab)++(0:0.8+.4  cm) coordinate (XV);
%            \path[dotted] (dd)++(0:0.8+2.4  cm) coordinate (YV);
%\draw(XV)--(YV); 

         \path[dotted] (ab)++(0:0.4  cm) coordinate (XV);
            \path[dotted] (de)++(0:0.4  cm) coordinate (YV);
\draw[thin,gray!30]
(XV)--(YV); 

      \path[dotted] (ab)++(0:0.8  cm) coordinate (XV);
            \path[dotted] (de)++(0:0.8  cm) coordinate (YV);
\draw[thin,gray!30](XV)--(YV); 

    \path[dotted] (ab)++(0:1.2  cm) coordinate (XV);
            \path[dotted] (de)++(0:1.2  cm) coordinate (YV);
\draw[thin,gray!30](XV)--(YV); 
 
\draw[thick] (cc2)--(cc3);
%   \draw[very thick,red] (ab)--(ba)--(bb3)--(bb2)--(bb25)--(bb252)--(bb1)--(cc)--(cc2)--(cc3)--(dd)--(de)--(ab);

       \foreach \i in {-12,-11,...,60}
  {
    \path  (origin)++(60:0.2*\i cm)  coordinate (a\i);
    \path (origin)++(120:0.2*\i cm)  coordinate (b\i);
    }
       \foreach \i in {0,2,...,60}{    
    \path[dotted] (a\i)++(120:7cm) coordinate (ca\i);
    \path[dotted] (b\i)++(60:4.8cm) coordinate (cb\i);
   }
      \foreach \i in {-8,...,48}
     {
    \path[dotted] (origin)++(60:0.2*\i cm)  coordinate (a\i);
    \path[dotted] (origin)++(120:0.2*\i cm)  coordinate (b\i);
%    \draw[thin,gray!30]  (a\i) -- (b\i) ;   
     }
     \foreach \i in {0,8,16,24,32,40}
{  \draw[thick,black]       (a\i) -- (b\i) ;   
  \draw[thick,black]       (a\i) -- (ca\i);
    \draw[thick,black]       (b\i) -- (cb\i) ; } ;
        \foreach \i 
        in {0,2,...,40}
{  \draw[thin,gray!30]       (a\i) -- (b\i) ;   
  \draw[thin,gray!30]       (a\i) -- (ca\i);
    \draw[thin,gray!30]       (b\i) -- (cb\i) ; } ;
 {
      \path[dotted] (b13)++(60:0.2*1 cm)++(0:1*0.2cm)  coordinate (2dr);
         \path[dotted] (b15)++(60:0.2*3 cm)++(0:1*0.2cm)   coordinate (1dr);
  \path[dotted] (b17)++(60:0.2*2 cm)++(120:1*0.2cm)   coordinate (2r);
         \path[dotted] (b22)++(60:0.2*7 cm)++(-120:1*0.2cm)   coordinate (3r);
   \fill (2dr)  circle (1.6pt);
   \fill (1dr)  circle (1.6pt);
      \fill (2r)  circle (1.6pt);
         \fill (3r)  circle (1.6pt);
   }
    {
      \path[dotted] (a9)++(120:0.2*21 cm)++(0:1*0.2cm)   coordinate (1u);
      \path[dotted] (1u)++(120:0.2*4 cm) ++(60:0.2*1 cm)++(120:1*0.2cm)++(-180:1*0.2cm)    coordinate (1uuuuu);

            \path[dotted]  (1uuuuu)          --++(90:0.24)coordinate(1uuuuuR)  ;
      
            \path[dotted]  (1uuuuu)   ++(-60:0.4*1 cm) ++(0:0.4*1 cm) coordinate (harrow);
      
         \path[dotted] (a10)++(120:0.2*10 cm)  coordinate (zero);
  \path[dotted] (a13)++(120:0.2*10 cm)++(-120:1*0.2cm)   coordinate (2ur);
         \path[dotted] (a14)++(120:0.2*11 cm)++(120:1*0.2cm)   coordinate (2ul);
         \path[dotted] (a18)++(120:0.2*9 cm)++(120:1*0.2cm)   coordinate (2uy);
   \fill (2ur)  circle (1.6pt);
    \fill (2ul)  circle (1.6pt);
        \fill (harrow)  circle (1.6pt);
       \fill (1u)  circle (1.6pt);
       \fill (1uuuuu)  circle (1.6pt) ;
     \draw(1uuuuuR)  node[above]
     {$\bm{  1^8\mid \varnothing \mid  \varnothing }$}; 
          \fill (2uy)  circle (1.6pt);
         \draw (zero)  node  {$\astrosun$};     
%         \path(zero)--(120:6*0.4)--++(-120:2*0.4) coordinate (602);
                  \path(zero)--++(120:6*0.4)--++(-120:2*0.4)  --++(-90:0.24) coordinate (602);
                  \path(zero)--++(-120:6*0.4)--++(120:2*0.4)  --++(-90:0.24) 
--++(-120:1*0.4)  --++(-60:1*0.4)                  coordinate (206);         
                  
%                  \draw (602) node [below] {$	\bm{  1^6\mid \varnothing \mid  1^2 }	$};
                  \draw (206) node   {$	\bm{   \varnothing\mid \varnothing \mid  1^8 }	$};
         
         \path(zero) ++(-60:0.4*4 cm)  coordinate (new6);
                  \path(new6) ++(-60:0.4*1 cm) ++(180:0.4*2 cm)  coordinate (new8);
         \path(zero) ++(0:0.4*2 cm)  coordinate (new1);
                  \path(new1) ++(120:0.4*1 cm) ++(0:0.4*2 cm)  coordinate (new2);
                  \path(new1) ++(0:0.4*6 cm)   coordinate (new3);                     
 \path(new2) ++(120:0.4*2 cm) ++(0:0.4*1 cm)  coordinate (new4);
      \path(new3) ++(-120:0.4*2 cm) ++(120:0.4*1 cm)  coordinate (new5);    
                  
            \fill (new1)  circle (1.6pt);            \fill (new2)  circle (1.6pt);       
                 \fill (new3)  circle (1.6pt);  
             \fill (new4)  circle (1.6pt);        \fill (new5)  circle (1.6pt);     
                  \fill (new8)  circle (1.6pt);                    \fill (new6)  circle (1.6pt);  
            
\draw[very thick, red](zero)--++(120:2*0.4)--++(-120:4*0.4)--++(120:2*0.4)
coordinate(hiyer!) circle (2pt)circle (1.2pt)circle (.4pt);

\draw[very thick, blue](zero)--++(0:1*0.4)--++(120:6*0.4)--++(-120:1*0.4)
coordinate(hiyer2!) circle (2pt)circle (1.2pt)circle (.4pt);

\draw[very thick, violet](zero)--++(-120:5*0.4)--++(0:1*0.4)--++(120:2*0.4)
coordinate(hiyer3!)circle (2pt)circle (1.2pt)circle (.4pt);

    \path(new3)  --++(90:0.24)coordinate (1uuuuuR);    
    \path(hiyer!)  --++(90:0.24)coordinate (hiyer!!);    
  \draw(hiyer!!)  node[above]
     {$\bm{   1^4\mid   \varnothing \mid 1^4}$}; 

  \draw(1uuuuuR)  node[above]
     {$\bm{  \varnothing \mid 1^8\mid   \varnothing }$}; 
                  
     }
               \foreach \i in {0,2,...,24}
     {
    \path[dotted] (origin)++(60:0.2*\i cm)  coordinate (a\i);
    \path[dotted] (a\i)++(-60:0.2*\i cm)  coordinate (b\i);
    \draw[thin,gray!30]  (a\i) -- (b\i) ;   
     \path[dotted] (a\i)++(-0:1.6 cm)  coordinate (c\i);
    \draw[thin,gray!30]  (a\i) -- (c\i) ;         
      \path   (origin)++(60:0.2*16 cm)  coordinate (a1);
          \path  (a1)++(-60:0.2*8 cm)  coordinate (b1);
                    \path  (b1)++(-180:0.2*8 cm)  coordinate (c1);
    \draw  (a1) -- (b1) --(c1) ;   
      \path   (origin)++(60:0.2*24 cm)  coordinate (a2);
          \path  (a2)++(-60:0.2*8 cm)  coordinate (b2);
                    \path  (b2)++(-180:0.2*8 cm)  coordinate (c2);
    \draw  (a2) -- (b2) --(c2) ;   
           \path  (origin)++(60:12 cm)  coordinate (GHY);
    \draw  (origin)--(GHY) ;   
}
                  \foreach \i in {-4,...,8,10,12,14,13,16,18,20,22}
     {
    \path[dotted] (origin)++(120:0.4*\i cm)  coordinate (a\i);
    \path[dotted] (a\i)++(-120:0.4*\i cm)  coordinate (b\i);
    \draw[thin,gray!30]  (a\i) -- (b\i) ;   
     \path[dotted] (a\i)++(180:1.6 cm)  coordinate (c\i);
    \draw[thin,gray!30]  (a\i) -- (c\i) ;         
      \path   (origin)++(120:0.4*16 cm)  coordinate (a1);
          \path  (a1)++(-120:0.4*8 cm)  coordinate (b1);
                    \path  (b1)++(0:0.4*8 cm)  coordinate (c1);
    \draw  (a1) -- (b1) --(c1) ;   
      \path   (origin)++(120:0.4*24 cm)  coordinate (a2);
          \path  (a2)++(-120:0.4*8 cm)  coordinate (b2);
                    \path  (b2)++(0:0.4*8 cm)  coordinate (c2);
    \draw  (a2) -- (b2) --(c2) ;   
           \path  (origin)++(120:12 cm)  coordinate (GHY);
    \draw  (origin)--(GHY) ;   
}
 \foreach \i in {0,2,4,6,...,16}
{
         \path[dotted] (origin)++(-120:0.2*\i cm)  coordinate (x\i);
     \path[dotted] (x\i)++(120:6 cm)  coordinate (y\i);
     \draw[thin,gray!30]  (x\i) -- (y\i) ;         
     }
 \path (origin)++(120:1.6 cm)  coordinate (AA);
  \path (origin)++(120:3.2 cm)  coordinate (BB);
    \path (origin)++(120:4.8 cm)  coordinate (CC);
      \path (origin)++(-120:1.6)++(120:4.8 cm)  coordinate (DD);
        \draw   (origin)++(120:1.6 cm)  -- (0:2 cm)  ;
 \path[dotted] (origin)++(120:13*0.2)++(180:0.2 cm)++(-120:1*0.2cm)   coordinate (XON);  
  \path[dotted] (XON)++(-120:4*0.4)++(0:0.4*2 cm) coordinate (XON2);  
 \fill (XON) circle (1.6pt);              \fill (XON2) circle (1.6pt);    
          
    \end{tikzpicture} 
$$

\!\!\!\!\!\!   \caption{The black points label the $3$-partitions of a  block of $R_8({0,1,2})$ with $e=4$.       The origin is labelled as $\astrosun$.  
There are three separate paths drawn on the diagram
belonging to  $  \Path^+((1^4 \mid \varnothing\mid1^4), 
\SSTT^
{(1^8\mid \varnothing \mid \varnothing )}  
   )$,
   $  \Path^+((1^6 \mid 1 \mid 1), 
\SSTT^
{(  \varnothing \mid  1^8 \mid \varnothing )}  
   )$,
and    
  $  \Path^+((1^2 \mid 1 \mid 1^5), 
\SSTT^
{(     \varnothing \mid   \varnothing\mid   1^8  )}  
   )$.  These are coloured red, blue, and violet respectively.     
We have labelled some of the points in the diagram for reference.     }
\label{2diag}
\end{figure}

\subsection{Semistandard tableaux as paths}   \label{3.2}
Let $e>h\ell$.  
%We now  encode the bases of cell-modules    in terms of paths.  
We now provide path-theoretic bases for the diagrammatic Cherednik algebra.  
Let  $  \mu \in \mptn \ell n (h)$. % and  $\SSTT\in\SStd^+_n(\lambda,\mu)$.
We define the {\sf component word} of $\mu$  to be the series of $\ell$-partitions
\[
\varnothing=  \mu^{(0)}  \xrightarrow{+X_1}
\mu^{(1)}  \xrightarrow{+X_2}
\mu^{(2)}  \xrightarrow{+X_3} 
\cdots 
\xrightarrow{+X_{n-1}} 
\mu^{(n-1)}  \xrightarrow{+X_n}
\mu^{(n)}  = \mu
\]
where   $X_{k}=(r_k,c_k,m_k)$ is the least  dominant removable node of the partition $\mu^{(k)}\in\mptn \ell {k}(h)$. 
Using the component word of $\mu$, we define a distinguished path $\stt^\mu$ from the origin to $\mu$   as follows 
$$
\stt^\mu=(+\varepsilon_{X_1},+\varepsilon_{X_2}, \dots, +\varepsilon_{X_n}).   
$$
%%$$
%\quad \text{ and let }\quad
%%$$
For $\lambda\in \mptn \ell n (h)$, we let  $$  
\SSTS =(+\varepsilon_{Y_1},+\varepsilon_{Y_2},\dots , +\varepsilon_{Y_n})\in \Path (\lambda,\SSTT^\mu).  $$
From $\SSTS$, we  obtain a tableau    $ \overline{\SSTS}\in \mathcal{T}(\la,\mu)$ by setting 
$ 
 {\overline{\SSTS}}({X}_k) = {\bf I}_{{Y}_k}.  
$ 
We   freely identify paths and tableaux in this manner (and so we drop the overline). 
Under this identification, we obtain a   bijection
 $
\SStd^+(\la,\mu) \leftrightarrow \Path^+(\la,\SSTT^\mu)
 $
and hence we can rewrite the basis of \cref{cellularitybreedscontempt-2} in terms of paths (see \cite[Theorem 5.21]{CoxBowman}) as follows.  For $\lambda \in \mptn \ell n (h)$ we have that 
\begin{equation}\label{anequation}
\Delta(\lambda)=\{C_\SSTS \mid \SSTS \in \Path^+(\lambda,\SSTT^\mu), \mu \in \mptn \ell n (h)\}.  
\end{equation}   

%Finally, we following paths will be essential in constructing homomorphisms in the later sections.  

\begin{defn}
Let  $\la,\mu \in \mptn \ell n(h)$ and suppose that $\la \rhd \mu$.  Then we let 
$ 
\SSTT^\mu_\lambda \in \Path(\lambda,\SSTT^\mu)
$ 
denote the unique path satisfying 
$$\deg(\SSTT^\mu_\la) = \ell (\mu)-\ell(\la).$$
\end{defn}

The above definition is well-defined by \cite[Proposition 7.4]{CoxBowman} and these paths will be very useful later on.  Examples of this path/tableau for three distinct pairs $(\la,\mu)$ are given in \cref{2diag}.  

\begin{rmk}
%In \cite{CoxBowman}, we make the additional assumption that $e\neq h\ell$ 
%so that the degree function on tableaux and paths can be matched up uniformly.  
% This is unnecessary for our purposes (as we only care about the  degree of a homomorphism, not that of a basis element). % and so we allow $e=h\ell$. 
%   In more detail, suppose  $e=h\ell$ and $\mu\in \mptn \ell n$ is in the image of the map $\det_h$ defined in \cite[Section 8]{CoxBowman}. If  $\SSTS \in \Path(\lambda,\SSTT^\mu)$, then $\deg(C_{\SSTT^\mu})=0$ whereas $\deg(\SSTT^\mu)\neq 0$.  However $\deg(\SSTS)-\deg(\SSTT^\mu)=\deg(C_\SSTS)-\deg(C_{\SSTT^\mu})$.  
%
%In \cite{CoxBowman}, we make the additional assumption that $e\neq h\ell$ 
%so that the degree function on tableaux and paths can be matched up uniformly.  
% This is unnecessary for our purposes (as we only care about the  degree of a homomorphism, not that of a basis element). % and so we allow $e=h\ell$. 
%   In more detail, suppose  $e=h\ell$ and $\mu\in \mptn \ell n$ is in the image of the map $\det_h$ defined in \cite[Section 8]{CoxBowman}.
If $e= h\ell$, 
 then all the results of this paper go through unchanged modulo   minor edits to the proofs.  
 Annoyingly, the definition of the degree and reflections of paths require some   tinkering   
  (akin to the  case $e=\infty$ case covered in detail in \cite[Section 6.4]{CoxBowman}).  %This is essentially because if $e=h\ell$ and $\mu=((n),(n),\dots,(n))$,  then 
% $(a)$ the unique path in   $\Path^+(\mu,\SSTT^\mu)$ has degree 1
%  whereas the unique tableau in 
%  $\SStd^+(\mu,\SSTT^\mu)$ has degree 0
% and  $(b)$  $\Path^+(((n,1),(n),(n)\dots,(n),(n-1)),\SSTT^\mu)\neq 
%  \emptyset$ 
%   and 
%  $\SStd^+(((n,1),(n),(n)\dots,(n),(n-1)),\SSTT^\mu)=
%  \emptyset$. 
  In what follows, we only discuss the case $e>h\ell$  explicitly. 
 For Cherednik algebras of symmetric groups, we provide an explicit and independent proof of our  main result in quantum characteristic $e=h$ in     \cref{Cahracteristicchange}.  
%For $e=h\ell$ the proofs require only trivial modifications to the language. All the statements of results are given in terms of  
%  semistandard tableaux and so work equally well for $e=  h\ell$.  
%  In     \cref{Cahracteristicchange}, we provide an explicit and separate  
%   proof of our main result for  $\ell=1$   (our main interest) in characteristic $e=h$.  
\end{rmk}

\renewcommand{\sts}{\mathsf{s}}   
\renewcommand{\stt}{\mathsf{t}}

\subsection{Standard   tableaux as paths}\label{3.3}
%It is easier to identify standard tableaux (rather than semistandard tableaux)
% with paths in $\mathbb{E}_{h,\ell}^+$.  
 Given $\lambda \in \mptn \ell n (h)$, a  tableau $\stt\in \Std(\lambda)$
  is  easily identified with the series of partitions 
  $\stt(k)$ for $0\leq k \leq n$, which in turn determine a path  in $\mathbb{E}_{h,\ell}^+$ 
 via the map in \cref{lhfskadj}. This provides path-theoretic bases of Specht modules $S_n(\la)$ for  $\lambda \in \mptn \ell n (h)$.  
We now restrict our attention to  $\lambda \in \mathcal{F}^\ell_n(h)\subseteq \mptn \ell n (h)$
 and  define the subset of $e$-restricted  standard $\lambda$-tableaux which will index the basis 
of the simple module $D_n(\la)$ for $\lambda \in \mathcal{F}^\ell_n(h)$.   

%   We now  wish to define the subset of $e$-restricted  standard tableaux which will index the bases 
%of simple modules for quiver Hecke algebras.   

\begin{defn}
Given $\lambda \in \mathcal{F}^\ell_n(h)$, we say that $\sts \in\Std(\lambda)$ is $e$-restricted if
 $
\sts(k) \in \mathcal{F}^\ell_k(h)
 $
  for all $1\leq k \leq n$.  We let $\Std_e(\lambda)$ denote the set of all $e$-restricted tableaux of shape $\lambda$.  
\end{defn}

Given $\la \in \mathcal{F}^\ell_n(h)$, we say that a  node $\square \in {\rm Rem} (\la)$ is {\sf good} if 
$\la - \square \in \mathcal{F}^\ell_{n-1}(h)$ (we remark that this is easily seen to coincide with the classical definition of a good node).  
We let $ \mathcal{F}_h(\la)$ denote the set of all {\sf good} removable nodes of $\lambda$.  
The following result is obvious, but will be essential for the  proof of our main theorem.  
\begin{prop}\label{tableau}
Given $\lambda \in \mathcal{F}^\ell_n(h)$, we have that $\langle c_\sts , c_\stt\rangle =\delta_{\sts,\stt}$ for $\sts,\stt \in \Std_e(\la)$.  Furthermore, 
$$
\Bbbk\{c_\sts %{\sf1}^{\res(\sts)}_{\omega}
  \mid \sts \in \Std_e(\lambda)\} \subseteq D_n(\lambda)
\quad\text{ and }\quad
\Std_e(\lambda) \leftrightarrow \bigsqcup_{\square \in  \mathcal{F}_h(\la)}
\Std_e(\lambda-\square) .
$$
%and $c_\sts = {\sf1}^{\res(\sts)}_{\omega}$.
\end{prop}
\begin{proof}
For $\la\in \mathcal{F}^\ell_{n}(h)$, we have that $\sts\in \Std_e(\lambda)$ if and only if
 $\sts{\downarrow}_{n-1}\in \Std_r(\nu)$ for some $\nu\in \mathcal{F}^\ell_{n-1}(h)$;  the bijection follows.  
 To see that $\{c_\sts \mid \sts \in \Std_e(\la)\}\subseteq D^\Bbbk_n(\la)$ and that $\langle c_\sts,c_\stt\rangle = \delta_{\sts,\stt}$ for $\sts,\stt\in\Std_e(\la)$, % and that $c_\sts={\sf1}^{\res(\sts)}_{\omega}$ ,
  it is enough to show that 
 \begin{equation}\label{eqnhg}
 {\sf 1}^{\res(\sts)}_{\omega} S_n(\nu) \neq 0 \text{ implies }\nu
 % \not \vartriangleleft 
\vartriangleright   \la \text{ or }\nu=\la
\text{ and } {\sf 1}^{\res(\sts)}_{\omega} S_n(\la)=c_\sts
%\text{ and }
% {\sf 1}^{\res(\sts)}_{\omega} S_n(\la)= c_\sts
 \end{equation} 
 for $\sts \in\Std_e(\la)$.  
 To see this, assume that  $c_\sts$ for $\sts \in \Std_e(\la)$ belongs  to some simple composition factor $L(\nu)$ of $S_n(\la)$ for $\nu \neq \la$;  
in which case  $  \lambda\rhd \nu$ and $${\sf 1}^{\res(\sts)}_{\omega}L(\nu) \subseteq {\sf 1}^{\res(\sts)}_{\omega}S_n(\nu)\neq 0$$    
which gives us our required contradiction.  
% This    implies that $1_\omega^{\res(\sts)}\in  S_n(\la)$.  
% Furthermore, if  
%\begin{equation}\label{eqnhg2} {\sf 1}^{\res(\sts)}_{\omega} S_n(\la)= c_\sts
%\end{equation}
%for $\sts \in\Std_e(\la)$, then we can also deduce that $c_\sts = {\sf1}^{\res(\sts)}_{\omega}$.  
  Now we turn to the  proof of  \cref{eqnhg}.  If $\nu \trianglelefteq \la$,  then $\nu \in \mptn \ell n(h)$.   
 Given $\stt \in \Std(\nu)$ with  $\stt (r_k,c_k,m_k) = k $, 
   we identify $\stt$ with the path 
   $$(+\varepsilon_{hm_1+c_1},+\varepsilon_{hm_2+c_2}, + \dots, 
   +\varepsilon_{hm_n+c_n}).$$
   Given $\stt \in \Std(\nu)$, we have that 
    $$\Path^+(\lambda,\stt )=\{ \stu \in \Std(\lambda) \mid \res(\stu) =\res( \stt) \}. $$
 Given any $\sts \in \Std_e(\lambda)$, we have that 
 $\sts(k) \in \mathcal{F}^\ell_k(h)$ for all $1\leq k \leq n$ and hence $\sts(k) \not \in \mathbb{E}(\alpha,me)$ for any $\alpha\in \Phi$, $m\in \ZZ$.  Hence   
\begin{equation}\label{slhfdjksalfkhjasfhlkhjsadflkjsahdflaksjdfh}
  \sts \not \in \bigcup_{
  \begin{subarray}c
  \nu\in \mptn \ell n(h)%\setminus\{\la\}
  \\
  \stt\in \Std(\nu), 
 \stt \neq \sts
  \end{subarray}
  } \Path^+(\lambda,\stt )
\end{equation}   and the result follows.  
 \end{proof}

\begin{eg}\label{eg1}
Let $h=1$ and $\ell=3$ and $\kappa =(0,1,2)\in (\ZZ/4\ZZ)^3$ as in \cref{2diag}.  
The unique $\la\in \mathcal{F}^3_n(h)   $  is given by 
$\la=((1^3),(1^3), (1^2))\}$.   
The tableau $\stt^\la$  is the unique element of $\Std_e(\la)$   and hence $D_8(\la)$ is 1-dimensional.   
\end{eg}

 \section{One column homomorphisms}

\renewcommand{\j}{i}
\renewcommand{\i}{j}

In  \cref{onecoller} we 
 construct the maps which will provide the backbone  of our BGG complexes.  
 We then consider how these homomorphisms compose (in terms of   ``diamonds") and it is 
 these in-depth diagrammatic calculations that provide the technical crux of the paper:   
 In  
  \cref{dimon} we classify the diamonds in terms of pairs of reflections in the alcove geometry;  
  in \cref{dimonpath} we localise to consider the $\mu$ weight-spaces of  cell-modules $\Delta(\lambda)$ for $\mu, \lambda$ two points in a given diamond; and finally 
  in \cref{dimonpath2} we use these results to prove that, within a diamond,  composition of  homomorphisms  commutes   (up to scalar multiplication by $\pm1$) or is zero  (for degenerate diamonds).
% also prove the     background results 
%which will be  necessary for proving that the complex is, indeed, a complex (in other words, we consider the composition of these one column homomorphisms).   

\subsection{One column homomorphisms}\label{onecoller}Let $e>h\ell$.   
Given $1\leq i < j\leq h\ell$ and $\alpha,\beta \in \mptn \ell n (h)$,
we suppose that $\ell(\alpha)=\ell(\beta)-1$ and that $\beta \rhd \alpha$.  
Then there exists a unique hyperplane
$  \mathbb{E} (\varepsilon_i-\varepsilon_j,\mu_{ij}e)$ 
for $1\leq i,j \leq h\ell$ and $\mu_{ij}\in \ZZ$ such that  $s_{i-j, \mu_{ij}e}(\alpha)=\beta.$ 
 By definition, this amounts to removing a series of nodes from the $j$th column of $\alpha$ 
and adding them in the $i$th column of $\alpha$ to obtain $\beta \in \mptn \ell n$ or vice versa.  
By not assuming that $i<j$, we can  use the notation 
$$s_{i-j, \mu_{ij}e}(\alpha)=\beta$$
to always mean that $\beta$ is obtained  by removing a series of nodes from the $j$th column of $\alpha$ 
and adding them in the $i$th column of $\beta$.  
 There are two distinct cases to consider.  The most familiar case (to many Lie theorists) is that in which 
 $\ell_{\varepsilon_{i} - \varepsilon_j}(\alpha,\beta)=1.$   In other words   
$\mu_{ij}\in \mathbb{Z}$ is the {\em unique} value such that 
$$\alpha\in \mathbb{E}^>(\varepsilon_i-\varepsilon_j,\mu_{ij}e)
\qquad 
\beta\in \mathbb{E}^<(\varepsilon_i-\varepsilon_j,\mu_{ij}e).$$
We refer to such pairs $(\alpha,\beta)$ as {\sf maximal} pairs.    These    pairs  include those  related by a reflection ``through  a common   alcove wall".  
The other case (which should be familiar to those who study blob and Virasoro algebras)  is that in which   $$\ell_{\varepsilon_i-\varepsilon_j}(\alpha,\beta)=2\ell_{\varepsilon_i-\varepsilon_j}(\alpha)-1$$ and so $\mathbb{E} (\varepsilon_i-\varepsilon_j,\mu_{ij}e)$ is just one of  many  hyperplanes  lying between $\alpha$ and $\beta$; these pairs $(\alpha,\beta)$  correspond to pairs which 
 are as far away as possible in the alcove geometry.  We refer to such pairs $(\alpha,\beta)$ as  {\sf minimal pairs}   (and they only exist for $\ell>1$).  
We wish to   distinguish between such pairs.  Therefore, for a minimal (respectively maximal)
pair we set $m_{ij}:=\mu_{ij}$ (respectively $M_{ij}:=\mu_{ij}$). 
We have that $m_{ij}\in \{0,1\}$ for any pair $\alpha,\beta \in\mptn\ell n(h)$.

\begin{eg}\label{big reflection}
Let $h=1$ and $\ell=3$ and $\kappa =(0,1,2)$ as in \cref{2diag}.  The pair  
$ 
((1^8 \mid \varnothing\mid \varnothing),
(1^2\mid \varnothing\mid 1^6))
 $ 
is a minimal pair. There are three hyperplanes parallel to $\mathbb{E}_{\varepsilon_1-\varepsilon_3,0}$ separating these two points.   \end{eg}

\begin{rmk}\label{maxmin}
Note that, near the origin, it is possible that a reflection is both maximal and minimal.  
For example, consider the pair $(1^6\mid \varnothing \mid 1^2)$  and 
$(1^4\mid \varnothing \mid 1^4)$ pictured in \cref{2diag}.  
\end{rmk}

%\begin{eg}\label{familiar}
%Let $h=1$ and $\ell=3$ and $\kappa =(0,1,2)$ as in \cref{2diag}.  
%Let 
%$$
%\alpha=(\varnothing \mid\varnothing \mid1^8)
%\quad
%\beta=(1^2\mid 1 \mid1^5) 
%\quad
%\gamma=(\varnothing \mid1^3\mid 1^5)
%\quad
%\delta=(1^2\mid \varnothing \mid1^6).
%$$
%The pairs 
%$$(\alpha,\gamma)
%\quad
%(\alpha,\delta)\quad
%( \gamma,\beta)\quad
%( \delta,\beta)$$
% are all examples of maximal pairs and should be familiar to classical Lie theorists.  
%\end{eg}

\begin{thm} \label{onecol} Let $e>h\ell$ and suppose that $\kappa \in I^\ell$ is $h$-admissible.   
Let $\alpha,\beta \in \mptn \ell n(h)$ be such that   $\ell(\beta)= \ell( \alpha)-1$.
Let $1\leq i,j \leq h\ell$ and $\mu_{ij}\in \NN$ be such that 
$s_{i-j, \mu_{ij}e}(\alpha)=\beta$.      
We have that   
$$\alpha\setminus \alpha\cap \beta
=\{X_1, X_2,\dots, X_k\} 
\quad \text{and}\quad 
  \beta\setminus \alpha\cap \beta
=\{Y_1, Y_2,\dots, Y_k\}  $$
with $X_a \rhd X_{a+1}$  (respectively $Y_a \rhd Y_{a+1}$) for $1\leq a < k$   is  a sequence   of nodes   belonging to the $j$th column of  $\alpha$ (respectively $i$th column of $\beta$).
There is a unique   $\SSTT^\alpha_\beta\in \SStd^+(  \beta , \alpha) $, as follows
% $ 
%| \SStd^+(  \beta , \alpha)  | =1 
% $ 
%and this unique tableau, $\SSTT^\alpha_\beta$, can be described as follows:
$$
\SSTT^\alpha_\beta (\square ) =
\begin{cases}
\square  & \text{if } \square  \in \alpha\cap \beta \\
 Y_k  & \text{if } \square  = X_k.
\end{cases} 
 $$
%belonging to  $ 
%| \SStd^+(  \beta , \alpha)  | =1 
% $.   
We have 
 $   \Hom_{A(n,\kappa)}({\Delta}(\alpha),{\Delta}(\beta) )=\Bbbk\{ \varphi^\alpha_\beta\}$  
 where $  \varphi^\alpha_\beta$  is determined by
  $ 
 \varphi^\alpha_\beta ( C_{\SSTT^\alpha}) = C_{\SSTT^\alpha_\beta}.
 $ 
  We define 
$\stt^\alpha_{\alpha\cap\beta}\in \Std(\alpha)$ and 
$\stt^\beta_{\alpha\cap\beta}\in \Std(\beta)$   to be the unique standard  tableaux of given shape determined by 
 $$\stt^\alpha_{\alpha\cap\beta}
 (r,c,m)
 =
 \stt^{\alpha\cap\beta}(r,c,m)=\stt^\beta_{\alpha\cap\beta}
 (r,c,m)$$ for
  $(r,c,m) \in \alpha\cap \beta $.  We have that 
   $ 
 \phi^\alpha_\beta ( c_{\stt^\alpha_{\alpha\cap\beta}}) = c_{\stt^\alpha_{\alpha\cap\beta}} 
 $  
determines the corresponding unique homomorphism in $\Hom_{R_n(\kappa)}(S_n(\alpha),S_n(\beta))$.  
The homomorphisms $\varphi^\alpha_\beta$   and $\phi^\alpha_\beta$ are both of degree $t^1$.  
  \end{thm}
\begin{proof}
For the statement  for $A(n,\kappa)$ see {\cite[Corollary 10.12]{CoxBowman}}.  
 By the definition of 
 $\stt^\alpha_{\alpha\cap\beta}$, 
  we have that $\Path(\lambda, \stt^\alpha_{\alpha\cap\beta})=\emptyset$  unless $\lambda \trianglerighteq \alpha$.  Therefore
  $ 
  e(\stt^\alpha_{\alpha\cap\beta})\Delta(\lambda) =0
  $ 
 unless $\lambda \trianglerighteq \alpha$.  Therefore $c_{\stt^\alpha_{\alpha\cap\beta}}\in L(\alpha)$ and thus  it is enough to define a homomorphism, $\phi^\alpha_\beta$ say,
 by where it sends $ c_{\stt^\alpha_{\alpha\cap\beta}}$.  
 Now, we have that 
$$
C_{\varphi(\stt^\alpha_{\alpha\cap \beta})}
C_{\SSTT^\alpha_\beta} =
C_{\varphi(\stt^\beta_{\alpha\cap \beta})}    \in {\bf  A}(n,\kappa)
$$
and so the result follows by applying the Schur idempotent.  
%By   \cref{MOO}, we have  that $\dim_\Bbbk(\Hom_{R_n(\kappa)}(S_n(\alpha),S_n(\beta)))=1$.
%  By the definition of 
% $\stt^\alpha$, 
%  we have that 
%  $ 
%  e(\stt^\alpha)\Delta(\lambda) =0
%  $ 
% unless $\lambda \trianglerighteq \alpha$.  Therefore $c_{\stt^\alpha_\beta}\in L(\alpha)$ and thus  it is enough to define a homomorphism, $\phi^\alpha_\beta$ say,
% by where it sends $ c_{\stt^\alpha_\beta}$.  
%have that 
%  $ e(\stt^\alpha)\Delta(\beta) = c_{ \stt^{\alpha}_{\beta}}$ and so any homomorphism (and we already know one exists)
%   must send $c_{\stt^\alpha}$ to $ c_{ \stt^{\alpha}_{\beta}}$.   
%   Alternatively, one can directly argue for $R_n(\kappa)$ by mimicking the proof of {\cite[Corollary 10.12]{CoxBowman}} (as $R_n(\kappa)$ is quasi-hereditary).    
 \end{proof}

% 
%\begin{eg}\label{big reflection}
%Let $h=1$ and $\ell=3$ and $\kappa =(0,1,2)$ as in \cref{2diag}.  The pair  
%$$
%((1^8 \mid \varnothing\mid \varnothing),
%(1^2\mid \varnothing\mid 1^6))
% $$
%is a minimal pair. There are three hyperplanes parallel to $\mathbb{E}_{\varepsilon_1-\varepsilon_3,0}$ separating these two points.   \end{eg}
%
%\begin{rmk}\label{maxmin}
%Note that, near the origin, it is possible that a reflection is both maximal and minimal.  
%For example, consider the pair $(1^6\mid \varnothing \mid 1^2)$  and 
%$(1^4\mid \varnothing \mid 1^4)$ pictured in \cref{2diag}.  
%\end{rmk} 
% 

\subsection{Diamonds formed by pairs of one-column morphisms}\label{dimon} We wish to consider all possible ways of composing a pair of such one-column homomorphisms.    %Simply by inspection, we have the following.   
%
%
%\begin{Diamonds}
Let $\alpha,\beta \in \mptn \ell n(h)$ 
be  such that $\beta \rhd \alpha$ and $\ell(\alpha)=\ell(\beta)+2$.  
There are six such  cases to consider which we now list.  The first five cases of homomorphisms should be familiar to all Lie theorists.  
We first consider the cases in which $\alpha,\beta$  differ in precisely three columns.  
In other words, $\alpha,\beta$ belong to a plane 
$\mathbb{R}\{\varepsilon_j-\varepsilon_i,\varepsilon_k-\varepsilon_\j
 \}$ 
for some $1\leq i,j,k \leq  h\ell$ and (without loss of generality) we can assume  that 
\begin{equation}\label{assump}
\langle \alpha, \varepsilon_\j \rangle  >  \langle \alpha, \varepsilon_\i \rangle >  \langle \alpha ,\varepsilon_k \rangle. 
\end{equation}

\begin{figure}[ht!]$$ \begin{tikzpicture}[scale=1.5]
{\scalefont{0.8}\path(0,0) --++(120:1cm)--++(180:2.2) coordinate (origin);
\draw[->](origin)--++(120:0.5cm);
\path(origin) --++(120:0.7cm) node {$+\varepsilon_\j$};
\draw[->](origin)--++(0:0.5cm);

\path(origin) --++(0:0.7cm) node {$+\varepsilon_\i$};
\draw[->](origin)--++(-120:0.5cm);
\path(origin) --++(-120:0.7cm) node {$+\varepsilon_k$};

}\fill[gray!30](0,0)--++(120:2cm)--++(-120:1cm)--++(-60:1cm)--++(0:1cm);
 \draw(0,0)--++(120:2cm)--++(00:1cm)--++(-60:1cm)
 --++(-120:1cm)--++(180:1cm)--++(120:1cm)--++(60:1cm)
 --++(0:1cm)--++(-120:2cm)--++(120:1cm)--++(0:2cm);
 \path(0,0)--++(120:1cm)--++(90:0.443cm) node {$\alpha$};

 \path(0,0)--++(120:1cm)--++(30:0.443cm) node {$\gamma$};
 
  \path(0,0)--++(120:1cm)--++(150:0.443cm) node {$\delta$};

 \path(0,0)--++(120:1cm)--++(-30:0.443cm) node {$\beta$};

\draw[very thick, white]
(0,0)--++(60:1cm)--++(120:1cm)--++(180:1cm)
--++(-120:1cm)--++(-60:1cm)--++(0:1cm);

 \draw[very thick, white]
(0,0)--++(60:1cm)--++(120:1cm)--++(180:1cm)
--++(-120:1cm)--++(-60:1cm)--++(0:1cm);
\end{tikzpicture}
 \quad 
 \begin{tikzpicture}[scale=1.5]
 \fill[gray!30](0,0)--++(60:1cm)--++(120:1cm)--++(-120:2cm)--++(0:1cm);
 \draw(0,0)--++(120:2cm)--++(00:1cm)--++(-60:1cm)
 --++(-120:1cm)--++(180:1cm)--++(120:1cm)--++(60:1cm)
 --++(0:1cm)--++(-120:2cm)--++(120:1cm)--++(0:2cm);
 \path(0,0)--++(120:1cm)--++(90:0.443cm) node {$\alpha$};
 \path(0,0)--++(120:1cm)--++(-150:0.443cm) node {$\beta$};
% \path(0,0)--++(90:0.443cm)--++(180:1cm) node {$\beta$};
 \path(0,0)--++(120:1cm)--++(30:0.443cm) node {$\delta$ };
 
  \path(0,0)--++(120:1cm)--++(150:0.443cm) node {$\gamma $ };

 \draw[very thick, white]
(0,0)--++(60:1cm)--++(120:1cm)--++(180:1cm)
--++(-120:1cm)--++(-60:1cm)--++(0:1cm);
\end{tikzpicture}
  \quad 
 \begin{tikzpicture}[scale=1.5]
 \draw(0,0)--++(120:2cm)--++(00:1cm)--++(-60:1cm)
 --++(-120:1cm)--++(180:1cm)--++(120:1cm)--++(60:1cm)
 --++(0:1cm)--++(-120:2cm)--++(120:1cm)--++(0:2cm);
 \path(0,0)--++(120:1cm)--++(150:0.443cm) node {$\alpha$};

 \path(0,0)--++(120:1cm)--++(-90:0.443cm) node {$\beta$};
 \path(0,0)--++(120:1cm)--++(-30:0.443cm) node {$\delta$ };
 
  \path(0,0)--++(120:1cm)--++(-150:0.443cm) node {$\gamma$ };

 \draw[very thick, white]
(0,0)--++(60:1cm)--++(120:1cm)--++(180:1cm)
--++(-120:1cm)--++(-60:1cm)--++(0:1cm);
\end{tikzpicture}
 \quad 
 \begin{tikzpicture}[scale=1.5]
 \draw(0,0)--++(120:2cm)--++(00:1cm)--++(-60:1cm)
 --++(-120:1cm)--++(180:1cm)--++(120:1cm)--++(60:1cm)
 --++(0:1cm)--++(-120:2cm)--++(120:1cm)--++(0:2cm);
 \path(0,0)--++(120:1cm)--++(30:0.443cm) node {$\alpha$};

 \path(0,0)--++(120:1cm)--++(-90:0.443cm) node {$\beta$};
% \path(0,0)--++(120:1cm)--++(-30:0.443cm) node {$\delta$ };
% 
%  \path(0,0)--++(120:1cm)--++(-150:0.443cm) node {$\gamma$ };
 \path(0,0)--++(120:1cm)--++(-30:0.443cm) node {$\gamma$ };
 
  \path(0,0)--++(120:1cm)--++(-150:0.443cm) node {$\delta$ }; \draw[very thick, white]
(0,0)--++(60:1cm)--++(120:1cm)--++(180:1cm)
--++(-120:1cm)--++(-60:1cm)--++(0:1cm);
\end{tikzpicture}
 $$
 \caption{The first four cases of diamond pairs $(\alpha,\beta)$.}
 \label{fourcases}
 \end{figure}
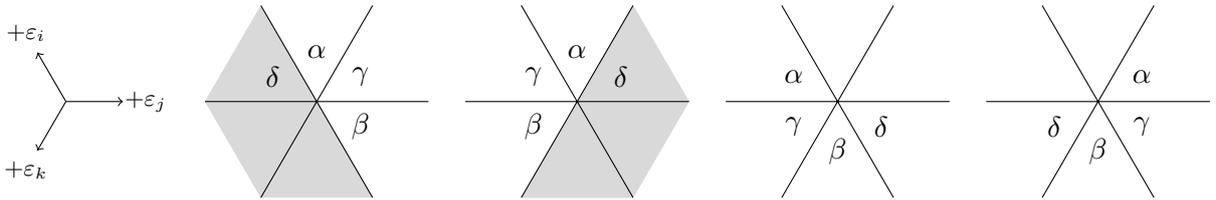

\begin{itemize}
\item[$(1)$]  
We have 
$\beta:=s_{k-\j,\mu_{k \j}e}s_{\i-\j,\mu_{\i \j}e} (\alpha) $
and 
 $\gamma:= s_{\i-\j,\mu_{\i \j}e}(\alpha)$.  There are two subcases 
\begin{itemize}
\item[$(a)$]    $\delta:= s_{k-\i,\mu_{k \i}e}(\alpha)\in \mptn \ell n$;
\item[$(b)$]      $\delta:=s_{k-\i,\mu_{k\i}e}(\alpha)\not \in  \mptn \ell n$;
 
\end{itemize}

\item[$(2)$]     
we have 
$\beta:=s_{k-i,\mu_{ki}e}s_{k-j,\mu_{kj}e}(\alpha)  
$ 
and $\gamma:= s_{k-j,\mu_{kj}e}(\alpha)$.    
There are two subcases 
\begin{itemize}
\item[$(a)$]    $\delta:= s_{j-i,\mu_{\i\j}e}(\alpha)\in \mptn \ell n$;
\item[$(b)$]   
 $\delta:= s_{j-i,\mu_{\i\j}e}(\alpha)\not \in  \mptn \ell n$;
 
\end{itemize}

\item[$(3)$]  
 
   $\delta:=s_{\i-\j,\mu_{\i\j}e}(\alpha)$ and  $\gamma:=s_{k-\j,\mu_{k\j}e}(\alpha)$ 
 $\beta:=s_{k-j,\mu_{kj}e} (\delta)  = s_{\i-i,\mu_{\i i}e}           (\gamma) $, 
  all belong to $\mptn \ell n$;

\item[$(4)$]  
 $\delta:=s_{k-\i,\mu_{k\i}e}(\alpha)$, and 
$\gamma:=s_{k-\j,\mu_{k\j}e}(\alpha)$ 
$\beta:=s_{ j-i,\mu_{ ji}e}(\delta)=   	s_{k-j,\mu_{kj}e}	(\gamma)	$, 
all belong to $\mptn \ell n$.

\end{itemize}
We  now assume that  $\alpha$ and $\beta$ differ in four columns (so that we cannot picture them belonging to a plane). Without loss of generality, we assume that 
$$\langle \alpha, \varepsilon_\j \rangle  >  \langle\alpha, \varepsilon_\i \rangle \qquad 
		  \langle \alpha ,\varepsilon_k \rangle		>   \langle \alpha ,\varepsilon_l \rangle  .$$
  This is the  case  in which 
\begin{itemize}

\item[$(5)$]   %$\alpha$,   
$\gamma:=s_{\i-\j,\mu_{\i\j}e}(\lambda)$, 
				$\delta:=s_{l-k,\mu_{lk}e}(\alpha)$, and  
 		 $\beta:=s_{l-k,\mu_{lk}e}(\gamma)= s_{\i-\j,\mu_{\i\j}e}(\delta)$  
 all belong to $\mptn \ell n$.

\end{itemize}
Finally, we have   one additional case to consider in which 
$\alpha$ and $\beta$ differ {\em only} in 2 columns.  
In other words $\alpha$ and $\beta$ belong to a line $\RR\{\varepsilon_i-\varepsilon_j\}$.  %\begin{itemize}
% \item[$(6a)$]   
% $\gamma:=s_{j-i,M_{ji}}(\alpha)$, 
%  $\beta:=s_{j-i,M_{ji}}  (\gamma)
%=
%  s_{i-j,m_{ij}}(\delta)$      all belong to $\mptn \ell n$ and 
%  \begin{itemize}
%\item[(i)] $\delta:=s_{j-i,m_{ji}}(\alpha)$  does not  belong to $\mptn \ell n$;   
%\item[(ii)]  for $\ell >1$ we have that $\delta:=s_{j-i,m_{ji}}(\alpha)$  does   belong to $\mptn \ell n$;   
%\end{itemize}
\begin{itemize}
 \item[$(6a)$]   
%  $\beta:=s_{j-i,M_{ji}}  (\gamma)
%=
%  s_{i-j,m_{ij}}(\delta)$  
We have $\beta=  s_{i-j,(1-m_{ji})e}   s_{j-i,m_{ji}e}(\alpha)$  
     and $\gamma:=s_{j-i,M_{ji}e}(\alpha)$ belong to $\mptn \ell n$    and 
  \begin{itemize}
\item[$(i)$] 
$\delta:=s_{j-i,m_{ji}e}(\alpha)$  does not  belong to $\mptn \ell n$;   
\item[$(ii)$]   For $\ell >1$ we have that $\delta:=s_{j-i,m_{ji}e}(\alpha)$  does   belong to $\mptn \ell n$;   
\end{itemize}

 \item[$(6b)$]   For $\ell >1$ we have that 
 $\beta:=s_{i-j,m_{ji}e}    s_{j-i,M_{ji}e}  (\alpha)$,  and 
    $\gamma:=s_{j-i,M_{ji}e}(\alpha)$ and $\delta:=s_{j-i,m_{ji}e}(\alpha)$ 
     belong to $\mptn \ell n$ and cannot be written in the form specified in case $(6a)$.

 \end{itemize}
 In the ``does belong to $\mptn \ell n$" cases, we get a diamond in the complex and so we refer to these 4-tuples as {\sf diamonds} and $(\alpha,\beta)$ as a {\sf diamond pair}.  
 In the ``does not  belong to $\mptn \ell n$" case, we get a single strand in the complex,
 and so we  refer to these 3-tuples as {\sf degenerate-diamonds} or {\sf strands}.  
% \end{Diamonds} 
  The first four cases  can be pictured by projecting into the plane $\mathbb{R}\{\varepsilon_ j-\varepsilon_ i,\varepsilon_k-\varepsilon_\j
 \}$   as depicted in \cref{fourcases}.

 In the first two cases, the lightly coloured-in region denotes the ``missing" region of $(1b)$ and $(2b)$.  
For diamonds formed entirely of maximal pairs, 
 the pictures in \cref{fourcases} consists only of six $e$-alcoves and their walls; thus the hyperplanes
  pictured are the {\em only} hyperplanes between $\alpha$ and $\beta$.
  See \cref{2diag} and \cref{familiar} for such an example.    
For (degenerate) diamonds involving one or two minimal pairs,   there can be  many other  hyperplanes between $\alpha$ and $\beta$ which are not pictured. See  \cref{2diag} and \cref{big reflection2} for such an example.    
The fifth case arises from a pair of orthogonal reflections and cannot be pictured in 2-dimensional space, however it is also the easiest case and so we do not lose much by being unable to picture it.  
The subcases of $(6)$ for which $\ell>1$  are easily pictured and 
should be familiar to those who work with Virasoro and blob algebras.  See  \cref{2diag} and \cref{big reflection} for such an example (many further examples can be found in \cite{bcs15}).    

\begin{eg} 
Let $h=1$ and $\ell=3$ and $\kappa =(0,1,2)$ as in \cref{2diag}. 
 The diamond consisting of  
\begin{equation}\label{big reflection}
\alpha=(1^8 \mid \varnothing\mid \varnothing)\quad
\beta=(1^2\mid \varnothing\mid 1^6)
\quad \gamma=(1^6 \mid \varnothing\mid 1^2)
\qquad
\delta=(1^2 \mid \varnothing\mid 1^6).  
\end{equation}
  is as in case $(6a)$.    
  The diamond consisting of the $3$-partitions, 
\begin{equation}\label{big reflection2}
\alpha=(\varnothing\mid  1^8 \mid  \varnothing)\quad
\beta=(1^6\mid 1\mid 1 )
\quad \gamma=(1^7 \mid \varnothing\mid 1  )
\qquad
\delta=( \varnothing\mid 1^7 \mid  1).  
\end{equation}
  is as in case $(4)$  
 and is a mixture of minimal and maximal pairs.  
 Let 
\begin{equation}\label{familiar}
\alpha=(\varnothing \mid\varnothing \mid1^8)
\quad
\beta=(1^2\mid 1 \mid1^5) 
\quad
\gamma=(\varnothing \mid1^3\mid 1^5)
\quad
\delta=(1^2\mid \varnothing \mid1^6).
\end{equation}
The diamond $(\alpha,\beta,\gamma,\delta)$ 
  is as in case $(4)$   and consists solely  of maximal pairs.   \end{eg}

\begin{rmk}
We   added a clause so that  cases $(6a)$ and $(6b)$ are mutually exclusive.  Without that clause, these cases would have a non-trivial intersection for points  near the   origin (see \cref{maxmin}).  
We have added this clause as these two subcases are genuinely different, 
see \cref{tableaux} below.  
\end{rmk}
 
\begin{defn}\label{dimon1}  Let $(\alpha,\beta)$ be a  non-degenerate   diamond pair. 
   We define the $(\alpha,\beta)$-{\sf vertex} to be  
      $$\xi =( \alpha \cap \beta \cap \gamma\cap \delta) \in \mptn \ell n (h).$$ 
In case $(6)$, 
%,  we let $W_\xi$ denote the copy of $\mathfrak{S}_2\times \mathfrak{S}_2$ generated by the reflections 
%through the hyperplanes 
%$  \mathbb{E}({j-i},m_{ji}e)$  and $  \mathbb{E}(j-i,M_{ji}e)$.  
 there exists $0< y <e$ and $x \geq 0$  such that 
$  \langle \alpha - \xi , \varepsilon_i\rangle =  xe+y. $ 
In cases $ (1)$, $(2)$  and $ (4)$,  we let $W_\xi$ denote the copy of $\mathfrak{S}_3$ generated by the reflections 
through the hyperplanes 
$  \mathbb{E}({j-i},\mu_{ji}e)$, $ \mathbb{E}({k-j},\mu_{kj}e)$, and  $ \mathbb{E}({k-i},\mu_{ki}e)$.  
 Given   $s\in W_\xi$ we let 
$$x=\langle \alpha - \xi , \varepsilon_i\rangle = \langle s(\alpha) - \xi   ,\varepsilon_{s(i)}\rangle 
\qquad y=\langle \alpha- \xi , \varepsilon_j\rangle = 
 \langle s(\alpha)- \xi , \varepsilon_{s(j)} \rangle
.$$  
 We let  
$\{X^{s(\alpha)}_1,X^{s(\alpha)}_2,\dots, X^{s(\alpha)}_x\}$ denote the final $x$ nodes of the 
  $s(i)$th column of $s(\alpha)$ and   let  
 $\{Y^{s(\alpha)}_1,Y^{s(\alpha)}_2,$ $\dots, Y^{s(\alpha)}_y\}$ 
 denote the final $y$ nodes of the  $s(j)$th column   of  $s(\alpha)$.

\end{defn}
 
\begin{rmk}In cases $(1)$, $(2)$  and $(4)$   we have
 that  $0 < y<e$ and  $\res(X_k)=\res(Y_k)$ for $1\leq k \leq \min\{x, y\}$.  
 \end{rmk}

 \subsection{Paths in diamonds}\label{dimonpath}
 We shall now consider reflections of the corresponding paths in the hyperplanes 
described in our 6 cases above.  We remark that each of these paths passes through 
each hyperplane {\em  at  most once}.  Therefore, we simplify our notation
 of \cref{rmkrmk} by dropping the superscript on the reflection.  %For brevity, we also drop the mention of $e$ in the subscript.  
   We now consider the (dominant) paths in $\Path(\beta,\alpha)$.  
 In   case   $(1a)$  there at   two paths 
  $$ 
\SSTS^\alpha_\beta  := s_{j-k,\mu_{kj}e} s_{k-\j,\mu_{k\j}e}(\SSTT^\alpha)
 \text{ and }   
\SSTT^\alpha_\beta :=s_{\i-\j,\mu_{\i\j}e} s_{k-\i,\mu_{k\i}e}(\SSTT^\alpha)  
$$ of degrees 0 and 2 respectively, 
  which are both dominant.  
  Generic examples of such paths (drawn from the point at which they meet the hyperplane $ \mathbb{E}({k-j},\mu_{kj}e)$ onwards) are pictured below
 $$ 
  \begin{tikzpicture}[scale=1.5]
{\scalefont{0.8}\path(0,0) --++(120:1cm)--++(180:2.2) coordinate (origin);
\draw[->](origin)--++(120:0.5cm);
\path(origin) --++(120:0.7cm) node {$+\varepsilon_\j$};
\draw[->](origin)--++(0:0.5cm);

\path(origin) --++(0:0.7cm) node {$+\varepsilon_\i$};
\draw[->](origin)--++(-120:0.5cm);
\path(origin) --++(-120:0.7cm) node {$+\varepsilon_k$};

}
  \draw(0,0)--++(120:2cm)--++(00:1cm)--++(-60:1cm)
 --++(-120:1cm)--++(180:1cm)--++(120:1cm)--++(60:1cm)
 --++(0:1cm)--++(-120:2cm)--++(120:1cm)--++(0:2cm);
 \path(0,0)--++(180:0.2) coordinate (neworigin);
 
   \path(neworigin)--++(0:0.1)--++(120:0.1)  --++(0:0.1)coordinate(neworigin);
 
  \draw[red,thick](neworigin) 
--++(120:0.1)
  --++(0:0.1)--++(120:0.1)
  --++(0:0.1)--++(120:0.1)--++(0:0.1)--++(120:0.1)
  --++(0:0.1) 
 --++(120:0.5)
  --++(-120:0.4)
   --++(0:0.3) circle (0.8pt) circle (0.4pt) ; 
 
 \draw[very thick, white]
(0,0)--++(60:1cm)--++(120:1cm)--++(180:1cm)
--++(-120:1cm)--++(-60:1cm)--++(0:1cm);
\end{tikzpicture}
 \quad 
  \begin{tikzpicture}[scale=1.5]

 \draw(0,0)--++(120:2cm)--++(00:1cm)--++(-60:1cm)
 --++(-120:1cm)--++(180:1cm)--++(120:1cm)--++(60:1cm)
 --++(0:1cm)--++(-120:2cm)--++(120:1cm)--++(0:2cm);
 \path(0,0)--++(180:0.2) coordinate (neworigin);
 
  \path(neworigin)--++(0:0.1)--++(120:0.1)  --++(0:0.1)coordinate(neworigin);
 
  \draw[red,thick](neworigin) 
--++(120:0.1)
  --++(-120:0.1)
  --++(120:0.1)
  --++(-120:0.1)--++(120:0.1)--++(-120:0.1)--++(120:0.1)
  --++(-120:0.1) 
 --++(120:0.5)
  --++(0:0.4)
   --++(0:0.3)circle (0.8pt) circle (0.4pt) ;  
 
 \draw[very thick, white]
(0,0)--++(60:1cm)--++(120:1cm)--++(180:1cm)
--++(-120:1cm)--++(-60:1cm)--++(0:1cm);
\end{tikzpicture}
$$
and are of degree 0 and 2 respectively.  
In   case  $(2a)$ there is a  unique path  
  $$ 
 \SSTT^\alpha_\beta:=  s_{k-\j,\mu_{k\j}e}s_{k-\i,\mu_{k\i}e}(\SSTT^\alpha)
$$ 
 which is  of degree  2 and dominant.  A generic example of such a path  is pictured below. %  in the $(a)$ case but not in the $(b)$ case.   
 $$ 
  \begin{tikzpicture}[scale=1.5]
{\scalefont{0.8}\path(0,0) --++(120:1cm)--++(180:2.2) coordinate (origin);
\draw[->](origin)--++(120:0.5cm);
\path(origin) --++(120:0.7cm) node {$+\varepsilon_i$};
\draw[->](origin)--++(0:0.5cm);

\path(origin) --++(0:0.7cm) node {$+\varepsilon_j$};
\draw[->](origin)--++(-120:0.5cm);
\path(origin) --++(-120:0.7cm) node {$+\varepsilon_k$};

} 

 \draw(0,0)--++(120:2cm)--++(00:1cm)--++(-60:1cm)
 --++(-120:1cm)--++(180:1cm)--++(120:1cm)--++(60:1cm)
 --++(0:1cm)--++(-120:2cm)--++(120:1cm)--++(0:2cm);
 \path(0,0)--++(180:0.2) coordinate (neworigin);
 
  \path(neworigin)--++(0:0.1)--++(120:0.1)  --++(0:0.1)coordinate(neworigin);
 
  \draw[red,thick](neworigin) 
--++(120:0.1)
  --++(-120:0.1)
  --++(120:0.1)
  --++(-120:0.1)--++(120:0.1)--++(-120:0.1)--++(120:0.1)
  --++(-120:0.1) 
 --++(120:0.5)
  --++(120:0.4)
   --++(-120:0.3) circle (0.8pt) circle (0.4pt) ; 
 
 \draw[very thick, white]
(0,0)--++(60:1cm)--++(120:1cm)--++(180:1cm)
--++(-120:1cm)--++(-60:1cm)--++(0:1cm);
\end{tikzpicture}
$$
    In each of cases $(1b)$ and $(2b)$ there is a  single path 
   $$ 
s_{j-k,\mu_{jk}e} s_{k-\j,\mu_{k\j}e}(\SSTT^\alpha)\text{ and }    s_{\i-\j,\mu_{\i\j}e}s_{k- i,\mu_{ki}e}(\SSTT^\alpha)$$
of degree  2, {\em neither} of which is dominant.  These are pictured below
$$ 
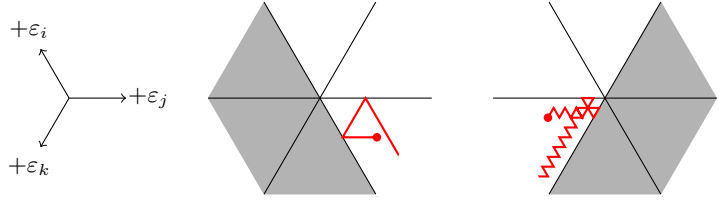
\begin{tikzpicture}[scale=1.5]
{\scalefont{0.8}\path(0,0) --++(120:1cm)--++(180:2.2) coordinate (origin);
\draw[->](origin)--++(120:0.5cm);
\path(origin) --++(120:0.7cm) node {$+\varepsilon_i$};
\draw[->](origin)--++(0:0.5cm);

\path(origin) --++(0:0.7cm) node {$+\varepsilon_j$};
\draw[->](origin)--++(-120:0.5cm);
\path(origin) --++(-120:0.7cm) node {$+\varepsilon_k$};

}\fill[gray!60](0,0)--++(120:2cm)--++(-120:1cm)--++(-60:1cm)--++(0:1cm);  \draw(0,0)--++(120:2cm)--++(00:1cm)--++(-60:1cm)
 --++(-120:1cm)--++(180:1cm)--++(120:1cm)--++(60:1cm)
 --++(0:1cm)--++(-120:2cm)--++(120:1cm)--++(0:2cm);
 \path(0,0)--(60:0.4cm) coordinate (neworigin);
\draw[red,thick](neworigin)
 --++(120:0.6)
   --++(-120:0.4)
    --++(0:0.3) circle (0.8pt) circle (0.4pt) ;
   
 \path(0,0)--++(0:2.5cm) coordinate(newneworigin);
 
  \fill[gray!60](newneworigin)--++(60:1cm)--++(120:1cm)--++(-120:2cm)--++(0:1cm);

 \draw(newneworigin)--++(120:2cm)--++(00:1cm)--++(-60:1cm)
 --++(-120:1cm)--++(180:1cm)--++(120:1cm)--++(60:1cm)
 --++(0:1cm)--++(-120:2cm)--++(120:1cm)--++(0:2cm);
 
  \path(newneworigin)--++(180:1cm)--++(120:0.1cm) coordinate(newnewneworigin);

\draw[red,thick](newnewneworigin)--++
(120:0.1cm)--++(0:0.1cm)
--++(120:0.1cm)--++(0:0.1cm)
--++(120:0.1cm)--++(0:0.1cm)
--++(120:0.1cm)--++(0:0.1cm)
--++(120:0.1cm)--++(0:0.1cm)
--++(120:0.1cm)--++(0:0.1cm)
--++(120:0.1cm)--++(0:0.1cm)
--++(120:0.1cm)--++(0:0.1cm)
--++(120:0.1cm)--++(0:0.1cm)
--++(-120:0.1cm)--++(0:0.1cm)
%  --++(-60:0.1)
%   --++(180:0.3); 
--++(-120:0.1cm)--++(120:0.1cm) 
--++(-120:0.1cm)--++(120:0.1cm) 
--++(-120:0.1cm)--++(120:0.1cm) --++(-120:0.1cm)--++(120:0.1cm) --++(-120:0.1cm)
circle (0.8pt) circle (0.4pt) ; %--++(-120:0.2cm) ;

\draw[very thick, white]
(newneworigin)--++(60:1cm)--++(120:1cm)--++(180:1cm)
--++(-120:1cm)--++(-60:1cm)--++(0:1cm);

 \draw[very thick, white]
(0,0)--++(60:1cm)--++(120:1cm)--++(180:1cm)
--++(-120:1cm)--++(-60:1cm)--++(0:1cm);
\end{tikzpicture}
 $$
 In each of cases  (3) and (4) there is a unique path 
 $$
 \SSTT^\alpha_\beta:=s_{k-\i,\mu_{k\i}e}s_{\i-\j,\mu_{\i\j}e}(\SSTT^\alpha)
 \qquad 
 \SSTT^\alpha_\beta:=s_{\i-\j,\mu_{\j\i}e}  s_{k-\i,\mu_{k\i}e	}(\SSTT^\alpha)
$$
respectively, of degree  2.  
Generic examples of such paths are pictured below 
$$ 
  
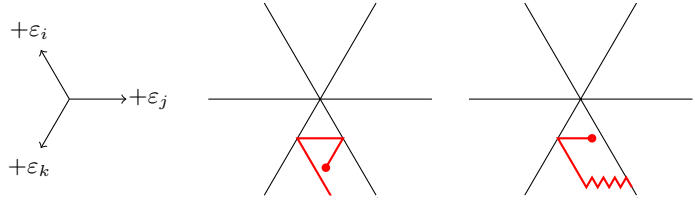
\begin{tikzpicture}[scale=1.5]
{\scalefont{0.8}\path(0,0) --++(120:1cm)--++(180:2.2) coordinate (origin);
\draw[->](origin)--++(120:0.5cm);
\path(origin) --++(120:0.7cm) node {$+\varepsilon_i$};
\draw[->](origin)--++(0:0.5cm);

\path(origin) --++(0:0.7cm) node {$+\varepsilon_j$};
\draw[->](origin)--++(-120:0.5cm);
\path(origin) --++(-120:0.7cm) node {$+\varepsilon_k$};

}

 \draw(0,0)--++(120:2cm)--++(00:1cm)--++(-60:1cm)
 --++(-120:1cm)--++(180:1cm)--++(120:1cm)--++(60:1cm)
 --++(0:1cm)--++(-120:2cm)--++(120:1cm)--++(0:2cm);
 \path(0,0)--++(180:0.4) coordinate (neworigin);
 
  \draw[red,thick](neworigin)
  --++(120:0.6)  --++(0:0.4)--++(-120:0.3)circle (0.8pt) circle (0.4pt) ;
 
 \draw[very thick, white]
(0,0)--++(60:1cm)--++(120:1cm)--++(180:1cm)
--++(-120:1cm)--++(-60:1cm)--++(0:1cm);
\end{tikzpicture}
 \quad 
  \begin{tikzpicture}[scale=1.5]

 \draw(0,0)--++(120:2cm)--++(00:1cm)--++(-60:1cm)
 --++(-120:1cm)--++(180:1cm)--++(120:1cm)--++(60:1cm)
 --++(0:1cm)--++(-120:2cm)--++(120:1cm)--++(0:2cm);
 \path(0,0)--++(180:0.2) coordinate (neworigin);
 
   \path(neworigin)--++(0:0.1)--++(120:0.1)  --++(0:0.1)coordinate(neworigin);
 
  \draw[red,thick](neworigin) 
--++(120:0.1)
  --++(-120:0.1)
  --++(120:0.1)
  --++(-120:0.1)--++(120:0.1)--++(-120:0.1)--++(120:0.1)
  --++(-120:0.1) 
 --++(120:0.5)
  --++(0:0.3)circle (0.8pt) circle (0.4pt) ; 
 
 \draw[very thick, white]
(0,0)--++(60:1cm)--++(120:1cm)--++(180:1cm)
--++(-120:1cm)--++(-60:1cm)--++(0:1cm);
\end{tikzpicture}
$$
%In case $(4)$, there is also a unique path of degree 0, 
%$$\SSTS^\alpha_\beta :=
%s_{j-k,\mu_{jk}e}
%s_{i-k,\mu_{ik}e}\SSTT^\alpha.  
%$$
  In case $(5)$, the reflections are orthogonal and there is a unique (dominant) path 
 and 
 if we assume (without loss of generality) that $\langle \alpha , \varepsilon_i\rangle >
 \langle \alpha , \varepsilon_k\rangle $, then this path is given by 
 $$
 \SSTT^\alpha_\beta:= s_{l-k,\mu_{lk}e} s_{j-i,{\mu_{ji}e}}(\SSTT^\alpha) 
 $$  and is of degree  2.  
     In case $(6a)$   
we have   $(x-1)$ distinct dominant paths of degree 0 given as follows, 
 $$
 \SSTS_\chi = 
 \begin{cases}
  s_{i-j, (M_{ij}-\chi-1)e   }
  s_{i-j, (M_{ij}-\chi)e  }
( \SSTT^\alpha) &\text{for $m_{ij}=1$}
 \\
  s_{i-j, (M_{ij}+\chi+1)e   }
  s_{i-j, (M_{ij}+\chi )e }
 (\SSTT^\alpha) &\text{for $m_{ij}=0$} 
 \end{cases}
 $$
 for $1\leq \chi < x$ (for $x$ as in \cref{dimon1}); we also  have a  unique path of degree 2 given by  
 $$
 \SSTT^\alpha_\beta = 
 s_{j-i,(1-m_{ji})e }
 s_{i-j,m_{ij}e }(\SSTT^\alpha) \in \Path(\beta,\SSTT^\alpha) 
$$ 
which is dominant if and only if we are in case $(6a)(ii)$.      In case $(6b)$  we have a unique (dominant)   path 
% $$ \SSTT^\alpha_\beta = 
% s_{j-i,M_{ij}-m_{ji} }
% s_{i-j,M_{ij} }(\SSTT^\alpha)$$
$$ \SSTT^\alpha_\beta = 
 s_{i-j,(2m_{ji}-M_{ji})e	 }
 s_{j-i,m_{ji}e }(\SSTT^\alpha)$$
 of 
 degree 2. 
% We have that $ \SSTT^\alpha_\beta\in  \Path^+(\beta,\SSTT^\alpha) $  if and only if we are in case $(6a)(ii)$ (and in particular $\ell>1$).  
% In case $ (6b)$  we have no tableaux of degree zero.  
Using \cref{anequation}, we now summarise the above as follows. 
\begin{prop}\label{tableaux}
Let  $\alpha,\beta \in \mptn \ell n(h)$ 
be a pair such that  $\beta \rhd \alpha$ and $\ell(\alpha)=\ell(\beta)+2$.  We have that 
$$\dim_t({\sf1}_{\alpha}\Delta(\beta))=
\begin{cases}
0  		&\text{in cases }(1b) \text{ and }   (2b) 		\\
t^2 + 1		&\text{in case  }(1a) 
\\

t^2  		&\text{in cases  }(2a),   (3),   (4),    		(5)  \text{ and } (6b) %	(6a)(ii)
\\
  	x- 1  		&\text{in case }		(6a)(i)		\\ 
 t^2 + x -1 	 &\text{in case  } (6a)(ii) 
\end{cases} $$
where $x\in \ZZ_{>0}$ is defined in \cref{dimon1} and $t$ is the indeterminate over $\ZZ_{\geq 0}$ from \cref{agrade}. 

\end{prop} 

%\begin{eg}
%In \cref{2diag}, we have that 
%\begin{itemize}
%\item the tableau $\SSTT^{(1^8\mid \varnothing\mid \varnothing)}_{(1^4\mid \varnothing\mid 1^4)}
%=   s_{1-3,M_{13}}^{(6)}\circ s_{3-1,M_{31}}^{(2)}  (\SSTT^{(1^8\mid \varnothing\mid \varnothing)})
%$ is as in case $(6a)$.
%\item the tableau $\SSTT^{(1^8\mid \varnothing\mid \varnothing)}_{(1^4\mid \varnothing\mid 1^4)}
%= s_{3-1,m_{31}}^{(7)} \circ s_{1-2,M_{12}}^{(1)}(\SSTT^{(1^8\mid \varnothing\mid \varnothing)})
%$ is as in case $(4)$.
\begin{eg}
In \cref{2diag}, we have that   the tableaux 
$$\SSTT^{(1^8\mid \varnothing\mid \varnothing)}_{(1^4\mid \varnothing\mid 1^4)}
=   s_{1-3,M_{13}}^{(6)}\circ s_{3-1,M_{31}}^{(2)}  (\SSTT^{(1^8\mid \varnothing\mid \varnothing)})
\quad \SSTT^{(1^8\mid \varnothing\mid \varnothing)}_{(1^4\mid \varnothing\mid 1^4)}
= s_{3-1,m_{31}}^{(7)} \circ s_{1-2,M_{12}}^{(1)}(\SSTT^{(1^8\mid \varnothing\mid \varnothing)})$$
$$\text{ and }\quad 
\SSTT^{(1^8\mid \varnothing\mid \varnothing)}_{(1^4\mid \varnothing\mid 1^4)}
= s_{1-2,m_{12}}^{(6)} \circ s_{3-2,m_{32}}^{(5)}(\SSTT^{(1^8\mid \varnothing\mid \varnothing)})$$
are as in cases $(6a)$,   $(4)$, and   $(4)$ respectively and are all of degree $t^2$.   
 \end{eg}

\subsection{Compositions of one-column homomorphisms in diamonds} \label{dimonpath2}
We now consider the composition of the one-column homomorphisms in terms of the path basis constructed in  \cref{tableaux}.   Let $\SSTT \in \mathcal{T}(\la,\mu)$ and 
$\SSTT(X)={\bf I}_Y\in \ZZ[\epsilon]$ for   $X\in \lambda$, $Y\in \mu$;
 we  abuse notation by writing either $\SSTT(X)= Y$ or $\SSTT(X)={\bf I}_Y$.    
 From \cref{tableaux}, we deduce the immediate corollary.  
\begin{cor}\label{cooooor}

The composition of two one-column homomorphisms is   zero  in the degenerate cases, in other words cases (1b), (2b) and (6a)(i). 
\end{cor}
  
  \begin{proof}
  Cases  (1b) and (2b) are clear.  Case (6a)(i) follows because 
   the  composition of two homomorphism of degree $t^1$ must be a vector 
   of degree $t^2$  and no such vector exists (by \cref{tableaux}).  
  \end{proof}

\begin{prop}\label{propwithaproof}
Let  $\alpha,\beta \in \mptn \ell n(h)$ 
be a pair such that  $\beta \rhd \alpha$ and $\ell(\alpha)=\ell(\beta)+2$. 
Assume that we are not in one of   degenerate cases of \cref{cooooor}.  
We have that 
%$$  
%C_{\SSTT_\beta^\alpha }= C_{\SSTT ^\alpha_{\delta} }
%C_{\SSTT  ^{\delta}_{\beta}}   
%$$
%in  all remaining cases.  Furthermore, we have  that 
%$$ C_{\SSTT_\beta^\alpha }
%=   C_{\SSTT ^\alpha_{\gamma} }  C_{\SSTT  ^{\gamma}_{\beta}} $$
\begin{itemize} 
\item $C_{\SSTT_\beta^\alpha }= C_{\SSTT ^\alpha_{\delta} }
C_{\SSTT  ^{\delta}_{\beta}}   
$ 
in  all remaining cases;
\item $ C_{\SSTT_\beta^\alpha }
=   C_{\SSTT ^\alpha_{\gamma} }  C_{\SSTT  ^{\gamma}_{\beta}} $ 
in all remaining cases except  $(1a)$, $(6)$, and $(4)$.  
\end{itemize}
\end{prop}

\begin{proof}
% We first consider  cases $(1)$ to  $(4)$.  In   cases $(1a)$, $(2)$,   $(3)$ and $(4)$ we claim that
%$$
% C_{\SSTT^\alpha_\beta}
%=
%  C_{\SSTT_{\beta}^{\delta}}C_{   \SSTT^{\alpha}_{\delta}} 
%$$
%and in cases $(2)$ and  $(3)$ we claim that 
%$$
% C_{\SSTT^\alpha_\beta}
%=
%  C_{\SSTT_{\beta}^{\gamma}}C_{   \SSTT^{\alpha}_{\gamma}} 
%$$
%and moreover, that both of these  can be seen without even applying any relations.  
For $\gamma \in \mptn \ell n$ (similarly for  $\delta\in \mptn \ell n$) it is  clearly enough to show that  
\begin{equation}\label{prop1}
  \SSTT^{\alpha}_{\gamma} \SSTT_{\beta}^{\gamma}    =\SSTT^\alpha_\beta
\in \mathcal{T}(\beta,\alpha) 
\end{equation}
on the level of  bijective maps  $:\beta \to \alpha$,  and furthermore that 
if 
\begin{equation}\label{prop2} (r,c,m) \rhd
(r',c',m')\text{ 
and 
}
 \SSTT_{\beta}^{\gamma}(r,c,m)
 \lhd 
 \SSTT_{\beta}^{\gamma}(r',c',m')
%$$ then this implies that $$
\text{  implies }
{   \SSTT^{\alpha}_{\gamma}} \SSTT_{\beta}^{\gamma}(r,c,m)
\lhd 
{   \SSTT^{\alpha}_{\gamma}} \SSTT_{\beta}^{\gamma}(r',c',m')
\end{equation}    
for any two nodes  $(r,c,m), (r',c',m') \in \beta\setminus  \xi$ of the same or adjacent residue. 
This is simply by the definition of the bases elements corresponding to these tableaux
(and the fact that double-crossings between strands of non-adjacent distinct residues can be removed by relation \ref{rel5}).  
The cases listed in the above proposition are precisely those for which \cref{prop1,prop2} are both true 
(in other words,   \ref{prop1} and \ref{prop2}  both hold in all cases except  in cases $(1a)$,  $(4)$, and $(6)$  for the product  $
\SSTT_\gamma^\alpha
\SSTT^\gamma_\beta$ --- which will be discussed separately).

  We shall consider case  $(2)$, as the other cases are identical.  
  It is clear that 
 $\SSTS(r,c,m) = (r,c,m)$ if $(r,c,m) \in \xi$ for $\SSTT ^\lambda_\mu$ for $\lambda, \mu \in\{\alpha,\beta,\gamma,\delta\}$.    
Thus it remains to consider the restriction of these bijections 
to   $:\beta\setminus \xi \to \alpha \setminus  \xi$ (via both $\gamma \setminus \xi $ and $\delta\setminus \xi $).   
 We have that $$
\SSTT^\gamma_\beta(X^\beta_{p})= Y^\gamma_{p}\quad 
\SSTT^\gamma_\beta(X^\beta_{y+q})=X^\gamma_{y+q}\quad 
\SSTT^\gamma_\beta(Y^\beta_{p})=X^\gamma_{p} 
$$
$$
\SSTT_ \gamma^\alpha(Y^\gamma_{p})= Y^\alpha_{p}\quad 
\SSTT _\gamma^\alpha(X^\gamma_{y+q})=X^\alpha_{y+q}\quad 
\SSTT_\gamma^\alpha(X^\gamma_{p})=X^\alpha_{p} 
$$
 and 
 $$
\SSTT^\delta_\beta(X^\beta_{p})= X^\delta_{p}\quad 
\SSTT^\delta_\beta(X^\beta_{y+q})=X^\delta_{y+q}\quad 
\SSTT^\delta_\beta(Y^\beta_{p})=Y^\delta_{p} 
$$
$$
\SSTT_\delta^\alpha(X^\delta_{p})= Y^\alpha_{p}\quad 
\SSTT_\delta^\alpha(X^\delta_{y+q})=X^\alpha_{y+q}\quad 
\SSTT_\delta^\alpha(Y^\beta_{p})=X^\alpha_{p} 
$$
and 
$$
\SSTT_\beta^\alpha(X^\beta_{p})= Y^\alpha_{p}\quad 
\SSTT_\beta^\alpha(X^\beta_{y+q})=X^\alpha_{y+q}\quad 
\SSTT_\beta^\alpha(Y^\beta_{p})=X^\alpha_{p} 
$$
for $1\leq p \leq y$ and $1\leq q \leq x-y$.  Therefore \cref{prop1} holds. 
To see that \cref{prop2} holds, one requires the following observation 
%We have that $$
%\SSTT^\gamma_\beta(X^\beta_{p})= Y^\gamma_{p}\quad 
%\SSTT^\gamma_\beta(X^\beta_{y+q})=X^\gamma_{y+q}\quad 
%\SSTT^\gamma_\beta(Y^\beta_{r})=X^\gamma_{r} 
%$$
%$$
%\SSTT_ \gamma^\alpha(Y^\gamma_{p})= Y^\alpha_{p}\quad 
%\SSTT _\gamma^\alpha(X^\gamma_{y+q})=X^\alpha_{y+q}\quad 
%\SSTT_\gamma^\alpha(X^\gamma_{r})=X^\alpha_{r} 
%$$
% and 
% $$
%\SSTT^\delta_\beta(X^\beta_{p})= X^\delta_{p}\quad 
%\SSTT^\delta_\beta(X^\beta_{y+q})=X^\delta_{y+q}\quad 
%\SSTT^\delta_\beta(Y^\beta_{r})=Y^\delta_{r} 
%$$
%$$
%\SSTT_\delta^\alpha(X^\delta_{p})= Y^\alpha_{p}\quad 
%\SSTT_\delta^\alpha(X^\delta_{y+q})=X^\alpha_{y+q}\quad 
%\SSTT_\delta^\alpha(Y^\beta_{r})=X^\alpha_{r} 
%$$
%and 
%$$
%\SSTT_\beta^\alpha(X^\beta_{p})= Y^\alpha_{p}\quad 
%\SSTT_\beta^\alpha(X^\beta_{y+q})=X^\alpha_{y+q}\quad 
%\SSTT_\beta^\alpha(Y^\beta_{r})=X^\alpha_{r} 
%$$
%for $1\leq p \leq y$ and $1\leq q \leq x-y$ and $1\leq r \leq y$.  Therefore \cref{prop1} holds. 
%To see that \cref{prop2} holds, one requires the following observation 
 $$
X^\beta_j \rhd Y^\beta_j \qquad
Y^\gamma_j \rhd X^\gamma_j \qquad
X^\delta_j \rhd Y^\delta_j \qquad
Y^\alpha_j \rhd X^\alpha_j 
$$
for all   $1\leq j \leq y$; one can apply this observation to each of the above tableaux in turn.   Thus \cref{prop2}   holds, as required.  
\end{proof}

  It remains to consider the $\gamma$ subcases of $(1)$, $(4)$, and $(6)$ not considered above.  
%  Namely, we must show that 
%  $$
% C_{\SSTT^\alpha_\beta}
%=
% (-1)^{y+1} C_{\SSTT_{\beta}^{\gamma}}C_{   \SSTT^{\alpha}_{\gamma}} 
%\qquad 
% C_{\SSTT^\alpha_\beta}
%=
%(-1)^{x-y}  C_{\SSTT_{\beta}^{\gamma}}C_{   \SSTT^{\alpha}_{\gamma}} 
%$$ in cases  $(1a)$ and $(6b)$   in the former case and $(4)$ in the latter case.  
 In all these cases, we shall see that \cref{prop1,prop2} fail.  Thus, we must apply some relations in order to 
rewrite  each  product-diagram in the required form.   
We let $A^{\rhd \beta}(n,\kappa)= \sum_{\lambda\rhd \beta}A (n,\kappa){\sf 1}_\la A (n,\kappa)$.  
 Given $(r,c,m)\in \alpha$, we let  ${\sf y}(r,c,m){\sf 1}_\alpha$ denote the diagram ${\sf 1}_\alpha$ with a dot added on the 
vertical solid strand with $x$-coordinate given by ${\bf I}_{(r,c,m)}$.  
Following \cite{MR2551762}, we set  $y_k={\sf y}(k,1,\ell)$.

%We shall consider diagrams $C_\SSTS$ in which any strand labelled by a node  from $\xi=\calphca\cap\beta\cap\gamma\cap \delta$ is vertical;  therefore  we shall focus on the strands relating nodes from 
%  $\alpha \setminus \xi$,
%  $\beta \setminus \xi$, $\gamma \setminus \xi$, and $\delta\setminus \xi$; furthermore, we  use the notation 
% from \cref{dimon} to label these nodes.  
%As in the proof of \cref{propwithaproof}, we shall continue to 
%use the notation of  \cref{dimon} to focus on strands  connecting nodes from 
%% Each strand connects a node on the northern or southern 
%% edge of the diagram
%Each diagram  we consider has a northern and southern loading given by a pair of distinct elements from the set $\{\alpha,\beta,\gamma,\delta\}$; 

\begin{prop}\label{ausefulprop}
Let $(\alpha,\beta)$ be a diamond pair as in case  $(1a)$.  
Then 
\begin{equation}\label{fashdljfdlshak}
   C_{   \SSTT^{\alpha}_{\gamma}}  C_{\SSTT_{\beta}^{\gamma}}=
 - {\sf y}({X^\alpha_y}){\sf 1}_\alpha
C_{\SSTS^{\alpha}_\beta}
=(-1)^{y+1}C_{\SSTT^{\alpha}_\beta}+ A^{\rhd \beta}(n,\kappa).
\end{equation}\end{prop}
\begin{proof}
 In case  $(1a)$, we have    
 $ 
 {   \SSTT^{\alpha}_{\gamma}}\circ  {\SSTT_{\beta}^{\gamma}} =\SSTS^\alpha_\beta \in \mathcal{T}(\beta,\alpha) 
  $ as bijective maps.  
 However, the corresponding product of   diagrams  has   a single double-crossing of non-zero degree; 
 this is between   the strand from $X^\beta_{y+1}$ on the southern edge  to $X^\alpha_{y+1}$ on the northern edge
 and the strand 
from  $Y^\beta_{y}$ on the southern edge  to $X^\alpha_{y}$ on the northern edge.  In particular, $\res(X_{y+1}^\beta)= \res(Y_{y}^\beta)-1$, and  
 $$
 X_{y+1}^\beta\lhd Y_y^\beta
\qquad
 {\SSTT_{\beta}^{\gamma}}( X_{y+1}^\beta)=X_{y+1}^\gamma\rhd 
 Y^\gamma_y= {\SSTT_{\beta}^{\gamma}}(Y_y^\beta)
 \qquad
  {\SSTT^{\alpha}_{\gamma}}(Y_{y}^\gamma)=X_{y}^\alpha
  \rhd 
 X_{y+1}^\alpha= {\SSTT^{\alpha}_{\gamma}}(X_{y+1}^\gamma).
 $$
%and $\res(X_{y+1}^\beta)= \res(Y_{y}^\beta)-1$.  
For  $1\leq p ,p' \leq y$ the strand 
from $X^\beta_{y+p}$ on the southern edge to $X^\alpha_{y+p}$ 
 on the northern edge double-crosses  with  the strand 
from  $Y^\beta_{p'}$ on the southern edge to $X^\alpha_{ p'}$ on the northern edge;  
 since $y<e$, we can remove  all of the   double-crossings  for $(p,q) \neq (1, 1)$ using relation \ref{rel5}.  
%There are a total of $y$ double-crossings in the diagram (involving each 
%of $y$ the strands 
%from $X^\beta_{y+p}$ to $X^\alpha_{y+p}$ 
%and 
%from $X^\beta_{y+p}$ to $X^\alpha_{y+p}$ 
%) 
We now resolve  the  final double crossing (for $(p,q)=(1,1)$)  using relation \ref{rel6} and hence obtain 
$$
  C_{   \SSTT^{\alpha}_{\gamma}}  C_{\SSTT_{\beta}^{\gamma}}={\sf y}(X_{y+1}^\alpha	)	{\sf 1}_\alpha C_{\SSTS^\alpha_\beta}
-
{\sf y}(X_{y}^\alpha	)	{\sf 1}_\alpha C_{\SSTS^\alpha_\beta}.
$$
Concerning the former diagram: we pull the dot down the strand and encounter no like-crossings on the way; hence this term is equal to zero.  
%The latter diagram is equal to  $ 
% (-1)^{y+1} 
%C_{\SSTT^\alpha_\beta}
%$ as stated in the first equality in the statement of the proposition.  
It remains to prove the second equality in \cref{fashdljfdlshak}.  
We let  $\SSTU_y  \in \mathcal{T}(\alpha,\beta)$ denote  the map 
$$
\SSTU_y(r,c,m) = 
\begin{cases}
{X^\alpha_y}			&\text{for } (r,c,m) = X^\beta_y\\
{Y^\alpha_y}			&\text{for } (r,c,m) = Y^\beta_y\\
\SSTS^\alpha_\beta( r,c,m) &\text{otherwise}.
\end{cases}.
$$
We claim that 
$${\sf y}({X^\alpha_y}){\sf 1}_\alpha
C_{\SSTS^{\alpha}_\beta}
=-C_{\SSTU_y} + A^{\rhd \beta}(n,\kappa).$$
To see this,   pull the dot at the top of the diagram ${\sf y}({X^\alpha_y}){\sf 1}_\alpha
C_{\SSTS^{\alpha}_\beta}$ 
down   the strand on which it lies (from   $  {X^\alpha_y} $ on northern edge to $ {Y^\beta_y} $ on the southern edge) towards the bottom of the diagram. 
By \cref{defintino2}, we can do this freely until we encounter a like-crossing of the form in 
relation \ref{rel3}.
Such a crossing involves %some vertical strand along with 
  the aforementioned 
strand  (between points   $ {X^\alpha_y} $  and ${Y^\beta_y}$ on the northern and southern edges)
  and some vertical strand of the same residue.  Such a vertical strand  either $(i)$ corresponds to a step   of the form $+\varepsilon _m$ for $m\not \in \{i,j,k\}$ 
 or $(ii)$  is   the vertical  strand 
from $ {X^\beta_y} $ on the southern edge to $ {Y^\alpha_y} $ on the northern edge.
    In the former case, the resulting error term belongs to  $A^{\rhd \beta}(n,\kappa)$.  
%In the latter case, we apply relation \ref{rel3} to move the dot past the crossing at the expense of an  error term, which   
%is equal to $-C_{\SSTU_y}$.  
In the latter case, we apply relation \ref{rel3} to move the dot past the crossing at the expense of acquiring an  error term, which   
is equal to $-C_{\SSTU_y}$.  
Finally (in the diagram which has a dot) we continue pulling the  dot reaches the bottom of the diagram,  the resulting diagram again belongs to $A^{\rhd \beta}(n,\kappa)$.   
Thus the only non-zero term acquired in this process
is $-C_{\SSTU_y}$ and the claim holds.   
 If $y=1$, then $\SSTU_y =\SSTT^\alpha_\beta$ and we are done.  Suppose that $y>1$.
Consider
\begin{itemize}
\item[$(i)$] the solid strand 
from ${Y^\beta_{y-1}}$ on the southern edge to ${X^\alpha_{y-1}}$ on the northern edge
\item[$(ii)$] the solid strand 
from ${X^\beta_{y-1}}$ on the southern edge to ${Y^\alpha_{y-1}}$ on the northern edge
\item[$(iii)$] the ghost  strand from 
 ${Y^\beta_{y}}$ on the southern edge to ${Y^\alpha_{y}}$ on the northern edge.  
 \end{itemize}
 These three strands together form a triple-crossing as on the right-hand side of  relation \ref{rel10}. Applying  relation \ref{rel10}, we can undo the crossing (at the expense of multiplication by  minus one and an error term with the same number of crossings).
 %The error term is the diagram on the left-hand side of relation \ref{rel10};
Consider the error term: We are free to pull the ghost   strand
 (of the strand connecting ${Y^\beta_{y}}$ and  ${Y^\alpha_{y}}$) 
  to the left 
 to obtain a diagram which belongs to $ A^{\rhd \beta}(n,\kappa)$.  
That leaves one remaining non-zero diagram   which differs from $-C_{\SSTU_y}$ in that we have undone the aforementioned  triple-crossing; to summarise
$${\sf y}({X^\alpha_y}){\sf 1}_\alpha
C_{\SSTS^{\alpha}_\beta}
= C_{\SSTU_{y-1}} + A^{\rhd \beta}(n,\kappa)
%$$
%where $\SSTU_{y-1}  \in \mathcal{T}(\alpha,\beta)$ is the map 
%$$
\quad\text{with}\quad  \SSTU_y(r,c,m) = 
\begin{cases}
{X^\alpha_{y-1}}			&\text{for } (r,c,m) =  X^\beta_{y-1}\\
{Y^\alpha_{y-1}}			&\text{for } (r,c,m) =   Y^\beta_{y-1}\\
\SSTU_{y}( r,c,m) &\text{otherwise}.
\end{cases} 
$$
Repeat this argument until all $y$ crossings have been resolved, the results follows.   
 \end{proof}

     \begin{prop}\label{ausefulprop3}
Let $(\alpha,\beta)$ be a diamond pair as in case  $(4)$.  
Then 
$$
  C_{   \SSTT^{\alpha}_{\gamma}}  C_{\SSTT_{\beta}^{\gamma}} =
(-1)^{x} C_{\SSTT^{\alpha}_\beta}+ A^{\rhd \beta}(n,\kappa).
$$
\end{prop}

\begin{proof}
We have that 
$$
\SSTT^{\alpha}_\gamma\circ \SSTT_{\beta}^\gamma
(r,c,m)
=
\begin{cases}
X_p^\alpha = \SSTT^{\alpha}_ {\beta}(Y_p^\beta)       & \text{if } (r,c,m) = X_p^\beta \text{ for  }1\leq p \leq x		\\
Y_p^\alpha = \SSTT^{\alpha}_ {\beta}(X_p^\beta)       & \text{if } (r,c,m) = Y_p^\beta	\text{ for  }1\leq p \leq x		\\
Y_q^\alpha 
= \SSTT^{\alpha}_ {\beta}(Y_q^\beta)       & \text{if } (r,c,m) = Y_q^\beta	\text{ for  }x+1\leq q \leq y		\\
\SSTT^{\alpha}_ {\beta}(r,c,m)      & \text{for  }(r,c,m) \in \xi
\end{cases}
$$
Consider
\begin{itemize}
\item[$(i)$] the solid strand 
from ${Y^\beta_{x}}$ on the southern edge to ${X^\alpha_{x}}$ on the northern edge;
\item[$(ii)$] the solid strand 
from ${X^\beta_{x}}$ on the southern edge to ${Y^\alpha_{x}}$ on the northern edge;
\item[$(iii)$] the ghost strand of the strand from 
 ${Y^\beta_{x+1}}$  on the southern edge to  ${Y^\alpha_{x+1}}$ on the northern edge.  
 \end{itemize}
These strands together form a crossing as on the righthand-side of relation \ref{rel10}.
Undoing this crossing we obtain an error term (corresponding to the diagram on the lefthand-side of  relation \ref{rel10}) which belongs to $A^{\rhd \beta}(n,\kappa)$ and another (non-zero) term.  One can then repeat the above argument with the 
latter diagram (except replacing the subscript $`x' $ with $`x-1'$).  Continuing in this fashion, we obtain the required result.  
\end{proof}

     \begin{prop}\label{ausefulprop3}
Let $(\alpha,\beta)$ be a diamond pair as in case  $(6b)$.  
Then 
$$
   C_{   \SSTT^{\alpha}_{\gamma}} C_{\SSTT_{\beta}^{\gamma}} =
  (-1)^{y}C_{\SSTT^{\alpha}_\beta}+ A^{\rhd \beta}(n,\kappa).
$$
\end{prop}

\begin{proof}
 We let $\{X^\alpha_1,X^\alpha_2,\dots, X^\alpha_{ex},Y^\alpha_{1},\dots, Y^\alpha_{y}\}$ 
   denote the final $xe+y$  nodes of the $i$th column of $\alpha$.  
% We let $\{X^\alpha_1,X^\alpha_2,\dots, X^\alpha_{ex +y}\}$ 
%   denote the final $xe+y$  nodes of the $i$th column of $\alpha$.  
 We let 
    $\{X^\beta_1,X^\beta_2,\dots, X^\beta_{ex}\}$ denote the final 
    $xe$ nodes of the $j$th column of $\beta$ and 
 $   \{Y_1^\beta,\dots, Y_y^\beta\}$ 
 denote the final $y$ nodes of the $i$th column of $\beta$.  
We have that 
$$
\SSTT^{\alpha}_\gamma\circ \SSTT_{\beta}^\gamma
(r,c,m)
=
\begin{cases}
\SSTT^{\alpha}_\beta(Y^\beta_p)=X^\alpha_p   	 	&\text{for }(r,c,m) = X^\beta_p \text{ and }1\leq p \leq y			\\
\SSTT^{\alpha}_\beta(X^\beta_p)=Y^\alpha_{p}   	 	&\text{for }(r,c,m) = Y^\beta_p \text{ and }1\leq p \leq y			\\
\SSTT^{\alpha}_\beta(X^\beta_q)= X^\alpha_{q}
&\text{for }(r,c,m) = X^\beta_q \text{ and }y< q \leq ex			\\
\SSTT^{\alpha}_\beta(r,c,m)		&\text{for }(r,c,m) \in \xi.			\\
\end{cases}
$$
Consider
\begin{itemize}
\item[$(i)$] the solid strand 
from ${Y^\beta_{y}}$ on the southern edge to ${X^\alpha_{y}}$ on the northern edge;
\item[$(ii)$] the solid strand 
from ${X^\beta_{y}}$ on the southern edge to ${Y^\alpha_{y}}$ on the northern edge;
\item[$(iii)$] the ghost   of the strand from 
 ${X^\beta_{y+1}}$ on the southern edge to ${X^\alpha_{y+1}}$  on the northern edge.    
 \end{itemize}
These strands together form a crossing as on the righthand-side of relation \ref{rel10}.
Undoing this crossing we obtain an error term (corresponding to the diagram on the lefthand-side of  relation \ref{rel10}) which belongs to $A^{\rhd \beta}(n,\kappa)$ and another (non-zero) term.  One can then repeat the above argument with the 
latter diagram (except replacing the subscript $`y' $ with $`y-1'$).  Continuing in this fashion, we obtain the required result.   \end{proof}

  \begin{prop}\label{ausefulprop4}
Let $(\alpha,\beta)$ be a diamond pair as in case  $(6a)(ii)$.  
Then 
$$
  C_{   \SSTT^{\alpha}_{\gamma}}  C_{\SSTT_{\beta}^{\gamma}} =
(-1)^{e(x+1)+y}
 C_{\SSTT^{\alpha}_\beta}+ A^{\rhd \beta}(n,\kappa).
$$
\end{prop}
\begin{proof}
We first fix some notation.   We denote the final $e$ nodes at the end of the 
  $j$th column of $\beta$    
by   $X^\beta_{1}, \dots, X^\beta_{e}$.
 We denote the final $e(x-1)+y$ nodes at the end of the 
  $i$th column of $\beta$    
by  $X^\beta_{e+1}, \dots, X^\beta_{ex}, Y^\beta_{1}, \dots, Y^\beta_{y} $.  
 We let $X^\alpha_{1},X^\alpha_{2},\dots, X^\alpha_{ex },
 Y^\alpha_{1}, \dots, Y^\alpha_{y}$ denote the final $ex+y$ nodes at the end of the $i$th column of $\alpha$.   
Given $\sigma \in \mathfrak{S}_x$ we define  $\SSTU_\sigma \in \mathcal{T}(\beta,\alpha)$ as follows,
\begin{equation}\label{fhdsallahjskf}
\SSTU_\sigma 
(X^\beta_{ep-q})
=
\SSTU_\sigma  
(X^\beta_{e\sigma(p)-q})   \qquad
 \qquad
\SSTU _\sigma
(Y^\beta_{t}) = 
 { Y^\alpha_{t} }
\end{equation}
  for  $1\leq  p,  \leq  x$, $0\leq q < e$,   $0\leq t < y$ and     such that 
$\SSTU_\sigma (r,c,m) =  {(r,c,m) } $  for $(r,c,m) \in \xi$.  
%We refer to such tableaux as {\sf brick tableaux} as they permute the $e$-bricks of nodes from $\alpha \setminus (\alpha\cap\beta)$ and $\beta\setminus (\alpha\cap\beta)$.   
% permute the $e$-bricks of nodes from $\alpha \setminus (\alpha\cap\beta)$ and $\beta\setminus (\alpha\cap\beta)$.   
We have that  $\SSTT^\alpha_\beta=\SSTU_{\rm id}$ for ${\rm id}\in \mathfrak{S}_x$ and 
 $\SSTS_{\chi}=\SSTU_{\sigma}$ for $\sigma=s_1s_2\dots s_\chi$ for $1
 \leq \chi< x$ (and so any element of $\SStd(\beta,\alpha)$ can be written in the form of  \cref{fhdsallahjskf}).  
 
 \newcommand{\yr}{r}
 
We now state a claim that will provide the crux of the proof.   
 Set  $\sigma=s_{1}s_{2}\dots s_{\chi}  $ for $1\leq \chi \leq x$.  
%We refer to the strand  originating at the point ${X^\beta_e}$ on the southern edge   as
% the {\sf principal strand} of the diagram $C_{\SSTU_\sigma}$. 
Given  $1\leq \yr < \chi $, we refer to the 
  strand  in $C_{\SSTU_\sigma}  $ 
 from $X^\beta_e$ on the southern edge 
 to
 $X^\alpha_{\chi e}$ on the northern  as the {\sf principal strand}.  
 Let 
 $C^\yr_{\SSTU_\sigma}  $ denote the diagram obtained from  $C_{\SSTU_\sigma}  $ by 
  placing a dot on the principal strand 
  at any point in the interval 
 $({\bf I}_{X^\alpha_{e\yr}  }, {\bf I}_{X^\alpha_{e\yr+1}  })\times [0,1]$. 
 %; in other words between the $r$th and $(r+1)$th bricks.
For $\sigma \neq s_1$, we claim that 
\begin{align}
\label{1}
& C^\yr_{\SSTU_\sigma}  = C^{\yr-1}_{\SSTU_\sigma}  + 
(-1)^{e+1}C^{\yr-1}_{\SSTU_{\sigma'}}  \\
 \label{2}
&C^\yr_{\SSTU_{s_1}}  = C^{\yr-1}_{\SSTU_{s_1}}  + 
(-1)^{e}C_{\SSTU_{{s_1}}}  
 \end{align}
modulo $ A^{\rhd \beta}(n,\kappa)$ where $\sigma'=s_{1}s_{2}\dots {s}_{\yr-1}  {s}_{\yr+1}\dots s_\chi %p
 	   $.  
Diagrammatically, we can think of our claim as simply a beefed-up version of 
relation \ref{rel3} in which we consider   crossings involving collections  of strands (each of  size $e>1$). 
  We let $i=\res(X^\beta_{e})$.% We let $r=\res(X^\beta_{y})$.

We now prove the claim. First apply relation \ref{rel3} to pull the dot through the crossing $i$-strands and hence obtain 
 $ C^{\yr-1}_{\SSTU_\sigma}$  plus another term
 with a minus sign. 
   For this latter diagram, the ghost of the principal $i$-strand   
can be pulled to the left through the crossing solid $(i+1)$-strands as in relation 
 \ref{rel10}.  We  hence obtain two diagrams: one with the same number of crossings, and one 
 in which the crossing of $(i+1)$-strands has been undone.  
 The former is zero modulo the stated ideal.  The latter diagram now has a crossing of two solid $(i+2)$-strands and a ghost $(i+1)$-strand as in  relation 
 \ref{rel10}.  
 Repeating as necessary, this process terminates with a diagram 
(occurring with coefficient  $(-1)^{e}$) which traces out the bijection of 
 ${\SSTU_{\sigma'}}
$ but with many  double-crossings.  
\begin{itemize}[leftmargin=*]
\item If $\sigma= s_1$, then all of these double-crossings are of degree zero;
\item If $\sigma \neq s_1$, then precisely  one of these double-crossings has non-zero degree: that between the solid  strand  from $X_{e}^\beta$
 on the southern edge   to $X^\alpha_{re}$ on the northern edge 
 and the ghost of the strand 
from   $ X^\beta_{  \yr e-e+1} $ on the southern edge  to 
 $  X^\alpha_{\chi e-e+1}  $ on the northern edge.  
\end{itemize}
In the latter case, we resolve this double-crossing as in relation \ref{rel10} and obtain two diagrams: one is of the required form and the other belongs to the stated ideal.  In either case,  the claim holds.   
Having proven our claim, we are now ready to   prove  the result.   We have that 
 $$
 \SSTT^\alpha_\gamma
 \circ 
 {\SSTT^\gamma_\beta}
 (r,c,m)
 =
 \begin{cases}
 \SSTU_{s_1\dots s_{x-1}} (X^\beta_{t})=X^\alpha_{xe-e+t}   &  \text{for $(r,c,m)=Y^\beta_{t} $ and }1\leq t \leq y \\
  
   \SSTU_{s_1\dots s_{x-1}} (Y^\beta_{t})= Y^\alpha_{t}   &  \text{for $(r,c,m)=X^\beta_{t} $ and }1\leq t \leq y \\

\SSTU_{s_1\dots s_{x-1}}  (r,c,m)  &\text{otherwise.}
 
  \end{cases}
 $$
Therefore, using $y$ applications of 
\ref{rel10} we obtain 
a diagram which traces out the same bijection as $\SSTU_{s_1\dots s_{x-1}}$ (modulo error terms).  
However the resulting diagram contains a single degree 2 double-crossing 
between the solid   strand from $X^\beta_e$ to $X^\alpha_{ex}$ (on the southern and northern edges, respectively) with the ghost of the strand from 
$Y^\beta_1$ to $Y^\alpha_{1}$ (on the southern and northern edges, respectively). 
Resolving this crossing using relation \ref{rel6}, we obtain 
that 
$$C_{\SSTT^\alpha_\gamma}C_{\SSTT^\gamma_\beta}
=(-1)^{y+1}C_{\SSTU_{s_1\dots s_{x-1}}} +A^{\rhd \beta}(n,\kappa).$$   
%  modulo $ A^{\rhd \beta}(n,\kappa)$.  
We now successively apply  \cref{1}  a total of $x-1$  times, 
  followed by  a single application of \cref{2}.  The error terms all belong to $A^{\rhd\beta}(n,
 \kappa)$ and the result follows.
\end{proof}

 \begin{thm}\label{key ingredient}
Let  $\alpha,\beta \in \mptn \ell n(h)$ 
be a non-degenerate diamond pair. % such that  $\beta \rhd \gamma,\delta \rhd \alpha$ belong to  $\mptn \ell n(h)$. 
 We have that  
$$    \varphi^\alpha_{\delta} \circ  \varphi^{\delta}_{\beta}
=
 \varphi_\beta^\alpha 
=  \varepsilon_{\alpha,\beta,\gamma,\delta} 
  \varphi^\alpha_{\gamma}  \circ  \varphi^{\gamma}_{\beta} 
\quad \text{where}\quad\varepsilon_{\alpha,\beta,\gamma,\delta} 
=
\begin{cases}

  (-1) ^ {y+1} 		&\text{in cases $(1a) $  }\\
          (-1) ^ {x} 		&\text{in cases $(4) $  }\\
    (-1) ^ {y} 		&\text{in cases $(6b) $  }\\
(-1)^{e(x+1)+y}	&\text{in case  $(6a)(ii)$}  						
  \\
  1 		&\text{otherwise.} 

  \end{cases}$$  Moreover the map  $\varphi^\alpha_\beta$ is determined by
   $\varphi^\alpha_\beta(C_{\SSTT^\alpha})= C_{\SSTT^\alpha_\beta}$  and 
$\dim_\Bbbk( \Hom_{{\bf    A} (n,\kappa)}({\bm\Delta}(\alpha),{\bm\Delta}(\beta)))=t^2.
$
 \end{thm}
\begin{proof}
 This theorem is mostly a restatement of the earlier results of this section
 (proved in the $A(n,\kappa)$ setting) using \cref{MOO}.  
To verify that the homomorphism space is 1-dimensional, it remains to check that $C_\SSTS \in L(\beta)$ for 
$\SSTS \in \SStd^+(\beta,\alpha)$ for each $\SSTS$ such that $\deg(\SSTS)=0$.  
We will not need  the dimension result in what follows and so we leave this as an exercise for the reader.  
\end{proof}

\section{A characteristic-free BGG-complex in the quiver Hecke algebra \\ and characteristic-free bases of simple modules }
\label{section4}
 
We are now ready to prove the main result from the introduction over $\Bbbk$  an arbitrary field.    
Given $\alpha\in \mathcal{F}^\ell_n(h)$, we define  an associated $R_n(\kappa)$-complex  and  show that this complex forms a BGG resolution of $D_n(\alpha)$.     We simultaneously construct bases and representing matrices for $D_n(\alpha)$ and completely determines its restriction along the tower of cyclotomic quiver Hecke algebras.

\begin{prop}Let $e>h\ell$.  For $\lambda \in \mathcal{F}^\ell_n(h)$, %a unitary $\ell$-partition, 
 we    set $${\bf C}_\bullet(\lambda) := %{\bf C}^\ell_n(\lambda) =
  \bigoplus_{
\begin{subarray}c 
%\mu \in \mptn \ell n(h)
%\\
\lambda\trianglerighteq  \mu
%, \ell(\mu)=\ell
\end{subarray}
}{ \bm\Delta} _n(\mu)[\ell(\mu)]
$$
where for $\ell\geq 0$ we define  the differential is the homomorphism of graded degree $t^1$   given as follows 
$$
 {\bm \delta}_{\ell}=\sum_{
 \begin{subarray}c
 \la\trianglerighteq  \mu \rhd \nu% \\
% \ell(\mu)=\ell(\nu)-1
 \end{subarray}
 } \varepsilon(\mu,\nu)\varphi^{\mu}_{\nu}%{\bf C}^{\ell+1}_n(\lambda) \to {\bf C}^\ell_n(\lambda)
$$
%to be an alternating  sum over all one-column homomorphisms 
for    $\varepsilon(\mu,\nu)\in\Bbbk\setminus\{0\}$  assigned arbitrarily, providing that each diamond in the complex satisfies 
\begin{equation}\label{crucial}\varepsilon_{\alpha,\beta,\gamma,\delta}=
-\varepsilon(\alpha,\gamma)
\varepsilon(\alpha,\delta)
\varepsilon( \gamma,\beta)
\varepsilon(  \delta,\beta).   
\end{equation}
In which case, we have that 
$ 
{\rm Im}({\bm\delta}_{\ell+1}) \subseteq {\rm ker}({\bm\delta}_{\ell})%\qquad 
%{\rm Im}(\delta_{\ell+1}) \subseteq {\rm ker}(\delta_{\ell})
$,  
in other words 
  $ {\bf C}_\bullet(\lambda)$  
is a complex. 

\end{prop}
\begin{proof}
This is a standard argument using \cref{key ingredient} and the fact that if $\ell(\alpha)=\ell(\beta)+2$,
then there exists at most two $\ell$-partitions, $\gamma$ and $\delta$ say, such that 
$\alpha \rhd \gamma,\delta \rhd \beta$.  
\end{proof} 

 We now apply the Schur functor  to the above to obtain  a complex of modules in the quiver Hecke algebra as follows, 
$$
C_\bullet(\lambda) :=
 {\sf E}_\omega {\bf C}_\bullet(\lambda)
%=
%C^\ell_n(\lambda)
= 
 \bigoplus_{
\begin{subarray}c 
 \lambda \trianglerighteq  \mu
 \end{subarray}
}S _n(\mu)[\ell(\mu)] 
 \qquad \text{ with }\qquad 
{\sf E}_\omega {\bm \delta}_{\ell }=\delta_{\ell}. 
$$
%where $\delta_\ell$ is a linear combination of the $\phi$
%By \cite[Theorem 4.9]{CoxBowman}, this simply amounts to pushing the complex through a Morita equivalence and so this is mostly vacuous.   
%However, the quiver Hecke algebra comes equipped with restriction functors, which will be instrumental in the proof of the following result.

\begin{thm}\label{maintheomre}
Let $e>h\ell$, let $\kappa \in I^\ell$ be $h$-admissible,  let  $\Bbbk$ be  a field, and $\lambda \in \mathcal{F}^\ell_n(h)$.  
The $R_n(\kappa)$-complex $C_\bullet(\lambda)$ is exact except in degree zero, where $$H_0(C_\bullet(\lambda))=D_n(\lambda).$$  We have  $ D_n(\lambda)=\Bbbk
\{c_{\sts} %={\sf 1}_{\omega}^{\res(\sts)}
  \mid \sts \in \Std_e(\lambda)\}  
${ and }$\rad(S_n(\lambda))=\Bbbk
\{c_{\sts} %={\sf 1}_{\omega}^{\res(\sts)}
  \mid \sts  \in \Std (\lambda)\setminus \Std_e(\lambda)\}. 
 $     
Furthermore,  $$\res^{n}_{n-1}  (D_n(\lambda))= \bigoplus_{ \square \in \mathcal{F}_h(\lambda) }   D_n(\lambda-\square).$$
   \end{thm}

\begin{proof}
We assume, by induction, that if $\lambda\in \mathcal{F}^\ell_  {n-1}(h)$, then the complex 
$C_\bullet(\lambda)$ %(for which $\varepsilon(\mu,\nu) \in \Bbbk\setminus\{0\}$ for all $\mu,\nu \trianglerighteq \lambda$ only need satisfy \cref{crucial}) 
 forms a BGG resolution and that $\{c_{\sts} \mid \sts \in \Std_e(\lambda)\}$ forms a basis of the simple module $D_n(\lambda)$.   
% We prove by induction on n that any complex C in which the coordinates \phi^\mu_\nu are all non-zero is exact
We now assume that $ \lambda\in \mathcal{F}^\ell_  {n}(h)$ and consider   the complex  $C_\bullet(\lambda)$.  We have that 
 $$\res^n_{n-1}(C_\bullet(\lambda))=\bigoplus_{r\in I} {\sf E}_r(C_\bullet(\lambda)).$$
We now consider one residue at a time. 
 As $\la$ belongs to an  alcove, %$e$-regular block,  
 we have that $\la$ (and any  $\mu \lhd \la$) has either 0 or 1 removable $r$-boxes for each $r\in I$. We let  ${\sf E}_r(\la)$ denote the unique $\ell$-composition (respectively  $\ell$-partition) 
  which differs from $\la$ by removing an $r$-node.  
For each residue, there are two possible cases.
\begin{itemize}[leftmargin=*]
\item  We have that ${\sf E}_r(\lambda) $ lies on an alcove wall or ${\sf E}_r(\lambda)\not\in \mptn \ell n$.   
By   restriction, we have that 
${\rm Im}({\sf E}_r\delta_{\ell+1})\subseteq \ker({\sf E}_r\delta_{\ell}) $ and so 
${\sf E}_r(C_\bullet(\la))$ forms  a complex. 
We have that ${\sf E}_r(\lambda)$ is fixed by reflection through some hyperplane 
and the $\ell$-compositions of $n$ which dominate $\lambda\in \mathcal{F}^\ell_n(h)$  come in pairs $(\mu^+,\mu^-)$ with $\mu^- \rhd \mu^+$ and $\ell(\mu^+)=\ell(\mu^-)+1$ and furthermore such that  
$${\sf E}_r(\mu^+)={\sf E}_r(\mu^-)=\mu \in \widetilde{\mathfrak{S}}_{h\ell} \cdot ({\sf E}_r(\lambda)) .$$
%In order that we can discuss things solely in terms of $\mu$ (and drop the superscripts where possible) we formally set $\ell(\mu)=\ell(\mu^+)-1/2=\ell(\mu^+)+1/2$.  
We have that 
%$$
%\ell(\mu^+)=\ell(\mu-)+1 \qquad \mu^- \rhd \mu^+
%$$
%and
$$
{\sf E}_r(S_n(\mu^+))
={\sf E}_r(S_n(\mu^-))
=
\begin{cases}
0		&\text{if either $\mu^+\not \in \mptn \ell n$ or $\mu^-  \not\in \mptn \ell n$ } \\
S_{n-1}(\mu)		&\text{otherwise}
\end{cases}
$$
Thus our complex ${\sf E}_r(C_\bullet(\la))$ decomposes into two chains of identical modules as follows,
\begin{equation}\label{hiyer}
{\sf E}_r(C_\bullet (\la)) =  \bigoplus_{\begin{subarray}c\mu \lhd \lambda-\square %\\ \mu \in \mptn \ell {n-1}
 \end{subarray}} S_n(\mu)[\ell(\mu)-1])   
 \bigoplus_{\begin{subarray}c\mu \lhd \lambda-\square %\\ \mu \in \mptn \ell {n-1}
  \end{subarray}} S_n(\mu)[\ell(\mu)].
\end{equation}
%We now consider the homomorphisms ${\sf E}_r\delta_\ell$ for $\ell\geq 0$.  
Given $\mu \lhd \lambda-\square$, the restriction of  
$
\phi^{\mu^+}_{\mu^-}\in \Hom_{R_n}(S_n(\mu^+),S_n(\mu^-))
$
is %non-zero (as restriction does not kill any submodule) and so (by highest weight theory) is 
 equal to  
\begin{equation}\label{shfjahljslhjkfads}  {\sf 1}_\mu \in {\rm End}_{R_n(\kappa)}(S_n(\mu))\end{equation}   
by the construction  of $\phi^{\mu^+}_{\mu^-}$ in \cref{onecol} 
 and \cite[Theorem 6.1]{manycell}.  
By restriction, we have  
$$
{\rm Im}({\sf E}_r(\delta_{\ell+1}))
\subseteq 
{\rm ker}({\sf E}_r(\delta_{\ell}))
$$
and by \cref{shfjahljslhjkfads}, we have that 
 $
{\sf E}_r\delta_{\ell+1}=\sum_{\ell(\mu)=\ell+1} {{\sf 1}_\mu}+ \dots
 $
 and so the complex is exact.  
We conclude that $H  ({\sf E}_r (C_\bullet (\lambda))=0$.

\item  We have that ${\sf E}_r(\lambda)\in \mathcal{F}^\ell_{n-1}(h) $. % lies in the fundamental/unitary alcove %(equivalently, $\square \in {\rm Rem}_\mathcal{F}(\lambda)$).
  Then 
 ${\sf E}_r(C_\bullet(\lambda))$ is given by 
 $$
{\sf E}_r\bigg(\bigoplus_{\lambda  \trianglerighteq  \mu}(S_n(\mu))[ \ell(\mu)] \bigg)
 $$
 with  homomorphisms 
 ${\sf E}_r \delta_{\ell}:{\sf E}_r C^{\ell}_n\to {\sf E}_r C^{\ell-1}_n$.   
We have that 
  $$
{\sf E}_r  S_n(\mu)[ \ell(\mu )]  
= S_{n-1}(\mu-\square) [ \ell(\mu-\square)] 
 $$if $ {\rm Rem}_i(\mu)\neq \emptyset$ and is  zero  otherwise.  In the non-zero case, this is simply because $ \mu-\square$ belongs to the same alcove  as $\mu$
 (and therefore  the   lengths coincide) for   $\mu \trianglelefteq \lambda$.  
 Now,   for a  pair $\mu,\mu'$ with $ \square\in {\rm Rem}_i(\mu)$
  and $ \square'\in {\rm Rem}_i(\mu')$, 
  we have that 
${\sf E}_r \phi^\mu_{\mu'}=   \phi^{\mu-\square}_{\mu'-\square'}$  
(again  by the construction  of $\phi^{\mu}_{\mu'}$ in \cref{onecol} 
 and \cite[Theorem 6.1]{manycell}).  
Thus ${\sf E}_r(C_\bullet(\lambda))= 
C_\bullet(\lambda-\square)$  
  and the righthand-side is exact except 
$H_0(C_\bullet(\lambda-\square))=D_n(\lambda-\square)$ by our inductive assumption.  
 Thus ${\sf E}_r (H_0(C_\bullet(\lambda)))=H_0(C_\bullet(\lambda-\square))= D_n(\lambda- \square)$
 and ${\sf E}_r(H_j(C_\bullet(\lambda)))=0$ for all  $j>0$.
      \end{itemize}
Putting all of the above together, we have shown that 
\begin{equation}\label{44}
\res^{n}_{n-1} (H_j (C_\bullet(\lambda)))=
\begin{cases}
\bigoplus_{ \square \in \mathcal{F}_h(\lambda) }   D_n(\lambda-\square)   &\text{if }j=0 \\
0			& \text{otherwise}.
\end{cases}
\end{equation}
Now, since ${\rm Head}(S_n(\lambda))=D_n(\lambda)\not \subset  {\rm Im}(\delta_1)$, we are able to conclude  that 
\begin{equation}\label{11} \res^{n}_{n-1}  (D_n(\lambda)) \subseteq \bigoplus_{ \square \in \mathcal{F}_h(\lambda) }   D_n(\lambda-\square).
\end{equation}
Conversely, we have that 
\begin{equation}\label{22}
|\Std_e(\lambda)| = \sum_{\square \in \mathcal{F}_h(\lambda) }|\Std_e(\lambda-\square)|
\end{equation}
by \cref{tableau}.  By induction,  the righthand-side of \cref{22} is equal to the dimension of the righthand-side of \cref{11}.  The lefthand-side of \cref{22} is a lower bound for the dimension of the lefthand-side of 
\cref{11}.  Putting these two things together, we deduce that 
\begin{equation}\label{33} \res^{n}_{n-1}  (D_n(\lambda))= \bigoplus_{ \square \in \mathcal{F}_h(\lambda) }   D_n(\lambda-\square) 
\end{equation}
and   furthermore, the  set $\{c_{\sts} \mid \sts \in \Std_e(\lambda)\} $ does indeed form a basis of $D_n(\lambda)$; to obtain the basis of the radical, recall that $e_{\stt}L(\mu)=0$ for $\la\rhd \mu $ and $\stt\in \Std_e(\lambda)$.  
Putting \cref{44} and \cref{33} together, we have that
$$
\res^{n}_{n-1} (H_j (C_\bullet(\lambda)))=
\begin{cases}
 \res^{n}_{n-1} D_n(\lambda )   &\text{if }j=0 \\
0			& \text{otherwise} 
\end{cases}
\quad\qquad  H_j (C_\bullet(\lambda) )=
\begin{cases}
D_n(\lambda) &\text{if }j=0 \\
0			& \text{otherwise} 
\end{cases}
$$
where the  second equality follows because   $\res^{n}_{n-1} D_n(\mu)\neq0$ for any $\lambda\trianglerighteq\mu$  (even though ${\sf E}_r (D_n(\mu))=0$ is possible for a given $r\in I$, as seen above).  %Finally, we remark that   $\res^{n}_{n-1} D_n(\mu)\neq0$ for any $\lambda\trianglerighteq\mu$  (even though ${\sf E}_r (D_n(\mu))=0$ is possible for a given $r\in I$, as seen above) and therefore $$
% H_j (C_\bullet(\lambda) )=
%\begin{cases}
%D_n(\lambda) &\text{if }j=0 \\
%0			& \text{otherwise} 
%\end{cases}
%$$ as required.  
  \end{proof}

 Note that the restriction rule  was used as the starting point in  \cite{MR1383482}, where Kleshchev   obtains results concerning the dimensions of simple modules.  
 Weirdly, our proof deduces that the homology of the complex is equal to $D^\Bbbk_n(\la)$, that the basis $D^\Bbbk_n(\la)$ is of the stated form, and the restriction of the simple module is of the stated form all at once!  % We invite the reader to compare the theorem below with \cite[Theorem 2.13$(vi)$ and Definition 4.1]{MR2266877} and \cite[Definition 2.3]{MR2534594}.  

\begin{thm}\label{action}
For $\lambda \in \mathcal{F}^\ell_n(h)$ the action of $R_n(\kappa)$ on  $D_n(\lambda)=\Bbbk\{c_\sts \mid \sts \in \Std_e(\lambda)\}$ is as follows:  
$$
y_k(c_\sts)=0\qquad {\sf 1}_\omega^{\underline{i}} (c_\sts)= \delta_{\underline{i},{\rm res}(\sts)}\qquad
\psi_r(c_\sts)=
\begin{cases}
c_{\sts_{k\leftrightarrow k+1}} &\text{if }|\res(\sts(r)) - \res(\sts(r+1))|>1 			\\
0   &\text{otherwise }
\end{cases}$$
where ${\sts_{k\leftrightarrow k+1}}$ is the tableau obtained from $\sts$ by swapping 
the entries $k$ and $k+1$.  %for $\sts \in \Std_e(\lambda)$.  
In particular, the  subalgebra
   $\langle y_k, {\sf 1}_\omega^{\underline{i}} \mid 1 \leq k\leq r, \underline{i}\in I^n\rangle \leq R_n(\kappa)$ acts semisimply on % any unitary simple module 
    $D_n(\lambda)$. The weight-spaces of $D_n(\lambda)$ are all 1-dimensional and  $D_n(\lambda)$ is concentrated in degree zero only.   
Finally, the cellular bilinear form is given by 
 $
\langle c_\sts , c_\stt \rangle =\delta_{\sts,\stt}
$ 
for $\sts,\stt \in \Std_e(\la)$. % (and therefore is positive definite if  $\Bbbk=\mathbb{R}$).  
\end{thm}

\begin{proof}
The  statements not  relating to the action  follow  from \cref{tableau} and \cref{slhfdjksalfkhjasfhlkhjsadflkjsahdflaksjdfh}  and \cref{maintheomre}.  
The action of the idempotents  is obvious.  The other zero-relations all follow because the product has non-zero degree (and the module $D_n(\lambda)$ is concentrated in degree 0).  
Finally, assume $|\res(\sts(r)) - \res(\sts(r+1))|>1$.  
The strands terminating at ${(r,1,\ell)} $ and ${(r+1,1,\ell)} $ on the northern edge either do or do not cross.  In the former case, we can resolve the double crossing in $\psi_r c_\sts$ without cost by our assumption on the residues and the result follows.  The latter case is trivial.  
Finally, notice that ${\sts_{k\leftrightarrow k+1}}\in \Std_e(\lambda)$  under the assumption that 
$|\res(\sts(r)) - \res(\sts(r+1))|>1 			$.  
 \end{proof}

%\begin{rmk}\label{aremark}
%Let $p>0$.  
%  Computing the composition series of    $S_n(\lambda)$ for $\lambda \in  \mptn 1 n(h)$ 
%  for arbitrary primes seems to be an impossible task    \cite{w13}.  
%Roughly speaking, if we assume that $p\gg h$  is suitably large
%% so  that the  conjectures of  Lusztig and Andersen  hold true   (as proven in \cite{w16}), 
%  then we can use Kazhdan--Lusztig theory to 
%  calculate $\dim_\Bbbk(D_n(\lambda))$
%  provided {\bf all partitions} $\mu\in  \mptn 1 n(h)$ such that $\mu \lhd \lambda$ belong to the first $p^2$-alcove \cite{w16}.  This is  equivalent to the $p$-weight of $\lambda$ (defined in \cref{abacus}) being less than $p$.  
% For $h=3$ this combinatorics has been conjecturally  extended (in terms of billiards in an alcove geometry) to the first $p^3$-alcove    \cite{MR3766576}.   
%We stress that there is no restriction on the $p$-weight of   $\lambda\in \mathcal{F}^1_n(h)$.  Therefore understanding the  composition series of unitary Specht modules  is   well beyond   the current state of the art.  Thus our two descriptions of the simple   modules $D_n(\lambda)$ for $\lambda\in \mathcal{F}^1_n(h)$ provide  the only contexts in which these modules can currently be hoped to be understood.  
%
%  \end{rmk}

\begin{rmk}\label{aremark}
Let $p>0$.  
  Combinatorially computing the composition series of    $S_n(\lambda)$ for $\lambda \in  \mptn 1 n(h)$ 
  for arbitrary primes seems to be an impossible task    \cite{w13}.  
 If we assume that $p\gg h$  is suitably large
   then we can use Kazhdan--Lusztig theory to 
 combinatorially calculate $\dim_\Bbbk(D_n(\lambda))$,
  this requires (as a minimum) that {\bf all partitions} $\mu\in  \mptn 1 n(h)$ such that $\mu \lhd \lambda$ belong to the first $p^2$-alcove \cite{w16}.  
  This is  equivalent to the requirement that the $p$-weight of $\lambda$ (defined in \cref{abacus}) is less than $p$.  
 For $h=3$ this combinatorics has been conjecturally  extended (in terms of billiards in an alcove geometry) to the first $p^3$-alcove    \cite{MR3766576}.   
We stress that there is no restriction on the $p$-weight of   $\lambda\in \mathcal{F}^1_n(h)$.  Therefore understanding the  composition series of unitary Specht modules  is   well beyond   the current state of the art.  Thus our two descriptions of the simple   modules $D_n(\lambda)$ for $\lambda\in \mathcal{F}^1_n(h)$ provide  the only contexts in which these modules can currently be hoped to be understood.  

  \end{rmk}

%\begin{rmk}\label{aremark}
%Let    $\ell=1$ and $\Bbbk$ be a field of positive characteristic.  
%%We stress that there is no restriction on the 
% Computing the composition series of   Specht modules $S_n(\lambda)$ for $\lambda \in  \mptn \ell n(h)$ 
%  for arbitrary primes seems to be an impossible task, as  shown by Williamson in    \cite{w13}.  
%  Even if we assume that $p\gg h$  is suitably large
%   that the  conjectures of  Lusztig and Andersen  hold true for $w(\lambda)<p^2$ 
%  (as proven in \cite{w16}), 
%  then we only have a combinatorial understanding of $D_n(\lambda)$ if all partitions $\mu \lhd \lambda$ belong to the first $p^2$-alcove.  
%  
%  understanding the  composition series of unitary Specht modules  is still  well-beyond even the theoretical limits of our current understanding
% --- see for example recent work of Lusztig--Williamson for the case $h=3$ and $w(\lambda)<p^3$ \cite{MR3766576}.   
%  \end{rmk}

 \begin{eg}[{\cite[Proposition 7.6]{CoxBowman}}]\label{aremark2}
 Let $\ell\geq 2$, $e=\ell+1$, and  $\kappa=(0,1,2,\dots , \ell-1)\in (\ZZ/e\ZZ)^\ell$, and $\Bbbk$ be arbitrary.  
 We have that 
$ 
 \lambda:= ((n),(n), \dots,(n))\in \mathcal{F}^\ell_{n\ell}(h)
$ 
and that 
\begin{equation}
\label{sgkhjfsghjlhsgdf}
[S_{n\ell} { ((n),(n), \dots,(n)) } : D_{n\ell}(\nu) ]=t^{\ell(\nu)}+ \dots
 \end{equation}
modulo terms of lower order degree.  Therefore every simple module $D(\nu)$ for $\la \rhd \nu \in \mptn \ell {n\ell}(1)$   appears with multiplicity at least 1.  
% For $\ell>2$  some of these multiplicities will   be  strictly greater than 1.  
Therefore as $n\rightarrow \infty$, the number of composition factors  of  $S_{n\ell}((n),(n), \dots,(n))$ tends to infinity and so is impossible to compute.  
In contrast, the module $D_{n\ell}(\nu)$ is  1-dimensional and   easily seen to be spanned by $c_{\stt^\lambda}$ for $\stt^\la$ as in \cref{nicetab}.  
\end{eg}

\section{Symmetric group combinatorics: $e$-abaci} \label{abacus}

  We now discuss how the combinatorial description of resolutions simplifies for (diagrammatic)  Cherednik algebras of symmetric groups.  
  In this case, we choose to emphasize the abacus presentation of partitions. We first recall this classical combinatorial approach, then flesh out the notion of homological degree introduced in \cite{conjecture} that is key to \cite[Conjecture 4.5]{conjecture}, and finally identify all this as the level $1$ case of the alcove geometry already studied in the previous sections.
  \subsection{The abacus of a partition}\label{subsect:abacus}
  Let $\lambda\in \mptn 1 n(h)$. Then $\lambda$ can be encoded by an abacus with at least $h$ beads, where each bead stands for a column of $\lambda$.
  This is simply a sequence of spaces and beads which records the shape of the border of $\lambda$, since knowing the border of $\lambda$ is the same as knowing $\lambda$. We form the $\Z$-abacus $\cA^h_\Z(\lambda)$ with $h$ beads by walking along the border from the top right corner to the bottom left corner of the Young diagram of $\lambda$, writing a space every time we walk down and a bead every time we walk left.

	\begin{eg}  
	The $\mathbb{Z}$-abaci of  $ (3^4,1)  , (3^3,2,1^2) \in \mptn 1 {13}(3)$ with 3 beads are as follows 
	
		\begin{center}
			\begin{tikzpicture}[scale=1.3,every text node part/.style={align=center}]
			
\foreach\x in{0,0.5,1,...,4.5}{
\draw[thick] (\x,-0.1+0.5)--(\x,0.1+0.5);
%\draw(\x,-0.2) coordinate {$\x$};
 }
 \draw[thick](0,0.5)--(4.5,0.5);
   \draw[thick,densely dashed](5,0.5)--(4.5,0.5);

        \node[circle,shading=ball,minimum width=3pt] (ball) at (2,+0.5) {};
                \node[circle,shading=ball,minimum width=3pt] (ball) at (2.5,+0.5) {};
                        \node[circle,shading=ball,minimum width=3pt] (ball) at (3.5,+0.5) {};

			\end{tikzpicture}\qquad\qquad
\begin{tikzpicture}[scale=1.3,every text node part/.style={align=center}]
			
\foreach\x in{0,0.5,1,...,4.5}{\draw[thick] (\x,-0.1+0.5)--(\x,0.1+0.5);
 }
 \draw[thick](0,0.5)--(4.5,0.5);
   \draw[thick,densely dashed](5,0.5)--(4.5,0.5);

        \node[circle,shading=ball,minimum width=3pt] (ball) at (1.5,+0.5) {};
                \node[circle,shading=ball,minimum width=3pt] (ball) at (2.5,+0.5) {};
                        \node[circle,shading=ball,minimum width=3pt] (ball) at (4,+0.5) {};

			\end{tikzpicture}			 \end{center}	\end{eg}

	Fix $e\geq 2$. We obtain an $e$-abacus $\cA_e(\lambda)$ by looping the $\ZZ$-abacus around $e$ horizontal runners. This can be described as follows: subdivide $\cA^h_\ZZ(\lambda)$ into segments of length $e$ starting from the leftmost position, then rotate each segment counterclockwise by ninety degrees so that it is vertical. The partition is now written on $e$ horizontal runners. Thus the runners of our $e$-abacus resemble a musical staff, and $\cA_e(\lambda)$ resembles sheet music. 
% We label the decorated runners of  $\mathcal{A}_e(\lambda)$ by placing the integer  $1$
% on the $(w+1)$th decorated strand (reading upwards from the southernmost runner and cycling round once we reach the $h$th runner) and then 
%continuing in this fashion until all runners are decorated.
	Like a staff, the runners of $\cA_e(\lambda)$ are bounded to the left. We let them extend infinitely to the right, because we want to think of being able to move beads in that direction by adding boxes or $e$-strips at the bottom of the Young diagram of a partition.  In French a musical score is called a \textit{partition}, so we may say that our abaci are written in the French style.

\begin{eg}\label{Debussy abacus} A $5$-abacus.\\
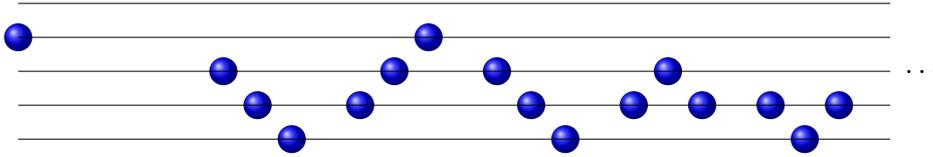
	\begin{center}
			\begin{tikzpicture}[scale=.45,every text node part/.style={align=center}]
			\node at (0,1)[circle,shading=ball,minimum width=3pt] (ball){};
			\node at (6,0)[circle,shading=ball,minimum width=3pt] (ball){};
			\node at (7,-1)[circle,shading=ball,minimum width=3pt] (ball){};
			\node at (8,-2)[circle,shading=ball,minimum width=3pt] (ball){};
			\node at (10,-1)[circle,shading=ball,minimum width=3pt] (ball){};
			\node at (11,0)[circle,shading=ball,minimum width=3pt] (ball){};
			\node at (12,1)[circle,shading=ball,minimum width=3pt] (ball){};
			\node at (14,0)[circle,shading=ball,minimum width=3pt] (ball){};
			\node at (15,-1)[circle,shading=ball,minimum width=3pt] (ball){};
			\node at (16,-2)[circle,shading=ball,minimum width=3pt] (ball){};
			\node at (18,-1)[circle,shading=ball,minimum width=3pt] (ball){};
			\node at (19,0)[circle,shading=ball,minimum width=3pt] (ball){};
			\node at (20,-1)[circle,shading=ball,minimum width=3pt] (ball){};
			\node at (22,-1)[circle,shading=ball,minimum width=3pt] (ball){};
			\node at (23,-2)[circle,shading=ball,minimum width=3pt] (ball){};
			\node at (24,-1)[circle,shading=ball,minimum width=3pt] (ball){};
			\node at (26.5,0){$\dots$};
			\draw[thin] plot [smooth] coordinates{(0,-2)(25.5,-2)};
			\draw[thin] plot [smooth] coordinates{(0,-1)(25.5,-1)};
			\draw[thin] plot [smooth] coordinates{(0,0)(25.5,0)};
			\draw[thin] plot [smooth] coordinates{(0,1)(25.5,1)};
			\draw[thin] plot [smooth] coordinates{(0,2)(25.5,2)};
		%	\node at (0,-3){$0$};
		%	\node at (1,-3){$1$};
		%	\node at (2,-3){$2$};
		%	\node at (3,-3){$3$};
		%	\node at (4,-3){$4$};
		%	\node at (5,-3){$5$};
		%	\node at (6,-3){$6$};
		%	\node at (7,-3){$7$};
		%	\node at (8,-3){$8$};
		%	\node at (9,-3){$9$};
			%\node at (-1,-2){$3$};
		%	\node at (-1,-1){$2$};
		%	\node at (-1,0){$1$};
		%	\node at (-1,1){$0$};
		%	\node at (-1,2){$4$};
			
			\end{tikzpicture}
		\end{center}
		\end{eg}

We let $w_e(\la)$ denote the total number of vacant spots which have a bead to their right and refer to this as the $e$-weight.   If $w(\rho)=0$ then we say that $\rho$ is an $e$-core.
Given a partition $\lambda$, we define the $e$-core of $\lambda$ 
to be the partition obtained by moving all beads on $\cA_e(\lambda)$
as far left as possible.   
%Then  $\la$ and $\mu$ have the same $e$-core if and only if $\lambda \in \widehat{{S}}_h \cdot \mu$.  
  We let 
  $$\Lambda(\rho, w) := \{ \mu \vdash |\rho| + we \mid e\text{-core}(\mu) = \rho\} $$   
  for $\rho$ an $e$-core.  
   	\begin{eg}  
	The $4$-abaci with 3 beads of  $ (3^4,1)$, $(3^3,2,1^2)$, $(3,2^5)$  
 and $(3^3,1^4) \in \mptn 1 {13}(3)$ are as follows 

		
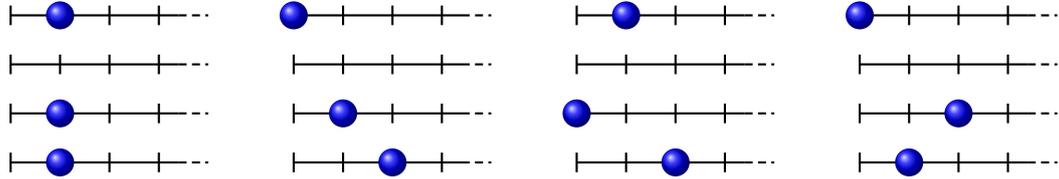
\begin{center}
			\begin{tikzpicture}[scale=1.3,every text node part/.style={align=center}]
			
\foreach\x in{0,0.5,1,...,1.5}{\draw[thick] (\x,-0.1)--(\x,0.1);
\draw[thick] (\x,-0.1+1.5)--(\x,0.1+1.5);
 \draw[thick] (\x,-0.1+0.5)--(\x,0.1+0.5);  \draw[thick] (\x,-0.1+1)--(\x,0.1+1); 
}
%\draw(0,0)--(0,1.5);
\draw[thick](0,0.5)--(1.7,0.5);
 \draw[thick](0,1.5)--(1.7,1.5);
  \draw[thick](0,1)--(1.7,1);
  \draw[thick](0,0)--(1.7,0);
\draw[thick,densely dashed](1.7,1)--(2,1);
\draw[thick,densely dashed](1.7,0.5)--(2,0.5);
 \draw[thick,densely dashed](1.7,1.5)--(2,1.5);
  \draw[thick,densely dashed](1.7,0)--(2,0);

       \node[circle,shading=ball,minimum width=3pt] (ball) at (0.5,0) {};
      \node[circle,shading=ball,minimum width=3pt] (ball) at (0.5,0.5) {};
      \node[circle,shading=ball,minimum width=3pt] (ball) at (0.5,1.5) {};

%      \draw(-0.25,0) node  {1};
%            \draw(-0.25,0.5) node  {2};
%                  \draw(-0.25,1.5) node  {3};
 
			\end{tikzpicture}
\qquad
	\begin{tikzpicture}[scale=1.3,every text node part/.style={align=center}]
			
\foreach\x in{0,0.5,1,...,1.5}{\draw[thick] (\x,-0.1)--(\x,0.1);
\draw[thick] (\x,-0.1+1.5)--(\x,0.1+1.5);
 \draw[thick] (\x,-0.1+0.5)--(\x,0.1+0.5);  \draw[thick] (\x,-0.1+1)--(\x,0.1+1); 
}
%\draw(0,0)--(0,1.5);
\draw[thick](0,0.5)--(1.7,0.5);
 \draw[thick](0,1.5)--(1.7,1.5);
  \draw[thick](0,1)--(1.7,1);
  \draw[thick](0,0)--(1.7,0);
\draw[thick,densely dashed](1.7,1)--(2,1);
\draw[thick,densely dashed](1.7,0.5)--(2,0.5);
 \draw[thick,densely dashed](1.7,1.5)--(2,1.5);
  \draw[thick,densely dashed](1.7,0)--(2,0);

       \node[circle,shading=ball,minimum width=3pt] (ball) at (1,0) {};
      \node[circle,shading=ball,minimum width=3pt] (ball) at (0.5,0.5) {};
      \node[circle,shading=ball,minimum width=3pt] (ball) at (0,1.5) {};

%      \draw(-0.25,0) node  {1};
%            \draw(-0.25,0.5) node  {2};
%                  \draw(-0.25,1.5) node  {3};
 
			\end{tikzpicture}	\qquad		
	\begin{tikzpicture}[scale=1.3,every text node part/.style={align=center}]
			
\foreach\x in{0,0.5,1,...,1.5}{\draw[thick] (\x,-0.1)--(\x,0.1);
\draw[thick] (\x,-0.1+1.5)--(\x,0.1+1.5);
 \draw[thick] (\x,-0.1+0.5)--(\x,0.1+0.5);  \draw[thick] (\x,-0.1+1)--(\x,0.1+1); 
}
%\draw(0,0)--(0,1.5);
\draw[thick](0,0.5)--(1.7,0.5);
 \draw[thick](0,1.5)--(1.7,1.5);
  \draw[thick](0,1)--(1.7,1);
  \draw[thick](0,0)--(1.7,0);
\draw[thick,densely dashed](1.7,1)--(2,1);
\draw[thick,densely dashed](1.7,0.5)--(2,0.5);
 \draw[thick,densely dashed](1.7,1.5)--(2,1.5);
  \draw[thick,densely dashed](1.7,0)--(2,0);

       \node[circle,shading=ball,minimum width=3pt] (ball) at (1,0) {};
      \node[circle,shading=ball,minimum width=3pt] (ball) at (0,0.5) {};
      \node[circle,shading=ball,minimum width=3pt] (ball) at (0.5,1.5) {};

%      \draw(-0.25,0) node  {1};
%            \draw(-0.25,0.5) node  {2};
%                  \draw(-0.25,1.5) node  {3};
 
			\end{tikzpicture}		
\qquad
	\begin{tikzpicture}[scale=1.3,every text node part/.style={align=center}]
			
\foreach\x in{0,0.5,1,...,1.5}{\draw[thick] (\x,-0.1)--(\x,0.1);
\draw[thick] (\x,-0.1+1.5)--(\x,0.1+1.5);
 \draw[thick] (\x,-0.1+0.5)--(\x,0.1+0.5);  \draw[thick] (\x,-0.1+1)--(\x,0.1+1); 
}
%\draw(0,0)--(0,1.5);
\draw[thick](0,0.5)--(1.7,0.5);
 \draw[thick](0,1.5)--(1.7,1.5);
  \draw[thick](0,1)--(1.7,1);
  \draw[thick](0,0)--(1.7,0);
\draw[thick,densely dashed](1.7,1)--(2,1);
\draw[thick,densely dashed](1.7,0.5)--(2,0.5);
 \draw[thick,densely dashed](1.7,1.5)--(2,1.5);
  \draw[thick,densely dashed](1.7,0)--(2,0);

       \node[circle,shading=ball,minimum width=3pt] (ball) at (0,1.5) {};
      \node[circle,shading=ball,minimum width=3pt] (ball) at (1,0.5) {};
      \node[circle,shading=ball,minimum width=3pt] (ball) at (0.5,0) {};

%      \draw(-0.25,0) node  {1};
%            \draw(-0.25,0.5) node  {2};
%                  \draw(-0.25,1.5) node  {3};
% 
			\end{tikzpicture}			\end{center}
			We have that $w_4(\la)=3$ and $4\text{-core}(\la)=(1)$ for each of these examples.  	\end{eg}

\begin{rmk}\label{conv:runnerlabels}	
%Let   $\mu\trianglelefteq \la$.  
Note that for $\mu \in \mptn 1 n (h)$, its removable box of highest content has content $h-k$, where $k$ is the position of the first bead in the $\Z$-abacus $\cA^{h}_{\Z}(\mu)$. In particular, this bead sits in runner $k$ mod $e$ in the $e$-abacus $\cA^{h}_{e}(\mu)$. Thus, in order to make the labels of the runners of the $e$-abaci in $\mu$
%in  $\mu\in\Po_e( \lambda)$
  correspond in a nice way to the contents of addable and removable $i$-boxes of the partitions, %in $\Po_e( \lambda)$, 
  one should label the runners from bottom to top by $h-1,h-2,\dots,1,0,e-1,e-2,\dots,h+1,h$. With this convention, removing a box of content $i$ mod $e$ corresponds to moving a bead on runner $i-1$ down to runner $i$; and adding a box of content $i$ mod $e$ corresponds to moving a bead on runner $i$ up to runner $i-1$.
%Here $h=\#\{\text{nonzero columns of }\lambda\}$.
\end{rmk}

\subsection{$e$-unitary partitions and    posets} We   recall the definition of   $e$-unitary partitions from \cite{conjecture} and show that these   are precisely the partitions in $\mathcal{F}^1_n=\cup_{h\geq 1}\mathcal{F}^1_n(h)$ studied in this paper.  
 
\begin{defn}\label{unitarypartitions}\cite{MR2534594,conjecture}  Fix $e\geq 2$. Suppose $\lambda$ has exactly $h$ columns and form  $\cA^h_\Z(\lambda)$ the abacus on $h$ beads.  
 We call $\lambda$ an $e$-unitary partition if all the beads on $\cA^h_\Z(\lambda)$ lie in an interval of width $e$.  In particular,  $\cA_e(\lambda)$ has at most one bead on each runner. Given an $e$-unitary partition $\lambda$, we let $\Po_e(\lambda)$ denote the set of all the $e$-abaci obtained from $\lambda$ by successively moving a bead on some runner one step to the right so long as we also move a bead on a different runner one step to the left.
\end{defn}
 
 \begin{prop}\label{unitary versus F}
The set $\mathcal{F}^1_n=\cup_{h\geq 1}\mathcal{F}^1_n(h)$ is precisely equal to the set of $e$-unitary partitions. 
\end{prop}
\begin{proof}
Suppose $\lambda\in\mathcal{P}^1_n$ has exactly $h$ columns and let $\gamma_h, \gamma_1$ denote the positions of the leftmost and rightmost beads on  $\cA^h_\Z(\lambda)$.  
% $h$'th bead from the right and the first bead from the right, respectively, in some $\ZZ$-abacus of $\lambda$ with at least $h$ beads.
  Now simply note that 
 $\gamma_1-\gamma_h\leq e-1$ if and only if $\langle \la+\rho ,  \varepsilon_1 - \varepsilon_h \rangle < e -1$ if and only if $\la \in\mathcal{F}^1_n(h)$.  \end{proof}

\begin{eg} When $e=4$, $(3^4,1)$ is a $4$-unitary partition, and $(3^3,2,1^2),(3,2^5),(3^3,1^4)\in\Po_4(\lambda)$.
\end{eg}
\begin{rmk}\label{rem:e-restricted} 
If $e=h$ then $\lambda$ is $e$-unitary if and only if $\lambda=(e^k)$ for some $k\in\N$. If $\lambda$ is an $e$-unitary partition, then any $\mu\in\Po_e(\lambda)$ is always $e$-restricted unless $\lambda=(e^k)$ and $\mu=\lambda$.
\end{rmk}

If an $e$-abacus $\cA_e(\mu)$ has at most one bead on each runner, let $b_i$ be the unique bead on the runner labeled $i$ if such a bead exists, and let $\beta_i\in\Z_{\geq 0}$ be the horizontal position of $b_i$. Sometimes by abuse of notation we might just refer to $\beta_i$ as a bead.

 \subsection{The affine and extended affine symmetric group actions}
%	The affine symmetric group $\widetilde{S}_h$ is generated by simple reflections $s_i$, $i\in\Z/h\Z$, subject to the relations $s_i^2=1$, $s_is_j=s_js_i$ if $|i-j|>1$, and $s_is_{i+1}s_i=s_{i+1}s_is_{i+1}$ (where all subscripts are taken mod $h$). 
%	
	
%	\begin{figure}[ht!]
	
%	\caption{Insert figure here explaining how to go between two notions of action (partitions versus abaci) }
%	\end{figure}

	There is a natural action of the affine symmetric group $\widetilde{\mathfrak{S}}_h$ on $\Po_e(\lambda)$ when we take the presentation of $\widetilde{\mathfrak{S}}_h$ given by generators $s_i$, $i\in\Z/h\Z$, subject to the relations $s_i^2=1$, $s_is_j=s_js_i$ if $|i-j|>1$, and $s_is_{i+1}s_i=s_{i+1}s_is_{i+1}$ (where all subscripts are taken mod $h$). %Suppose $\lambda\in\Po_e^u$ has $h$ columns.
	  $\mathfrak{S}_h=\langle s_1,\dots,s_{h-1}\rangle$ acts by permutation of the $h$ runners containing beads, while $s_0$ switches the top and bottom beads in the abacus, then moves the bottom bead one step to the right and the top bead one step to the left. 
	From the description of $\Po_e(\lambda)$ in Definition \ref{unitarypartitions}, $\widetilde{\mathfrak{S}}_h$ acts transitively on $\Po_e(\lambda)$.

	 \begin{eg}\label{exl:girlwiththeflaxenhair}
		Illustration of the action of $s_0$:
		
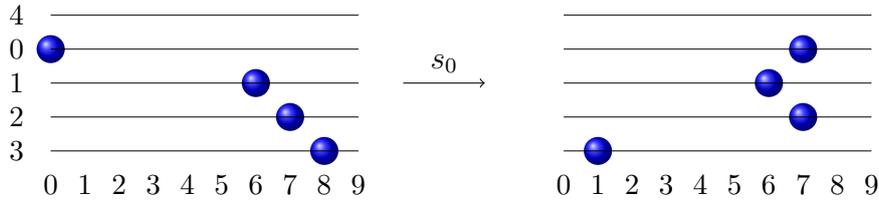
\begin{center}
			\begin{tikzpicture}[scale=.45,every text node part/.style={align=center}]
			\node at (0,1)[circle,shading=ball,minimum width=3pt] (ball){};
			\node at (6,0)[circle,shading=ball,minimum width=3pt] (ball){};
			\node at (7,-1)[circle,shading=ball,minimum width=3pt] (ball){};
			\node at (8,-2)[circle,shading=ball,minimum width=3pt] (ball){};
			\draw[thin] plot [smooth] coordinates{(0,-2)(9,-2)};
			\draw[thin] plot [smooth] coordinates{(0,-1)(9,-1)};
			\draw[thin] plot [smooth] coordinates{(0,0)(9,0)};
			\draw[thin] plot [smooth] coordinates{(0,1)(9,1)};
			\draw[thin] plot [smooth] coordinates{(0,2)(9,2)};
			\node at (0,-3){$0$};
			\node at (1,-3){$1$};
			\node at (2,-3){$2$};
			\node at (3,-3){$3$};
			\node at (4,-3){$4$};
			\node at (5,-3){$5$};
			\node at (6,-3){$6$};
			\node at (7,-3){$7$};
			\node at (8,-3){$8$};
			\node at (9,-3){$9$};
			\node at (-1,-2){$3$};
			\node at (-1,-1){$2$};
			\node at (-1,0){$1$};
			\node at (-1,1){$0$};
			\node at (-1,2){$4$};
			\node(A) at (10,0){};
			\node(B) at (13,0){};
			\draw[->] (A) -- node [font=\large, above] {$s_0$} (B) ;
			\node at (16,-2)[circle,shading=ball,minimum width=3pt] (ball){};
			\node at (21,0)[circle,shading=ball,minimum width=3pt] (ball){};
			\node at (22,-1)[circle,shading=ball,minimum width=3pt] (ball){};
			\node at (22,1)[circle,shading=ball,minimum width=3pt] (ball){};
			\draw[thin] plot [smooth] coordinates{(15,-2)(24,-2)};
			\draw[thin] plot [smooth] coordinates{(15,-1)(24,-1)};
			\draw[thin] plot [smooth] coordinates{(15,0)(24,0)};
			\draw[thin] plot [smooth] coordinates{(15,1)(24,1)};
			\draw[thin] plot [smooth] coordinates{(15,2)(24,2)};
			\node at (15,-3){$0$};
			\node at (16,-3){$1$};
			\node at (17,-3){$2$};
			\node at (18,-3){$3$};
			\node at (19,-3){$4$};
			\node at (20,-3){$5$};
			\node at (21,-3){$6$};
			\node at (22,-3){$7$};
			\node at (23,-3){$8$};
			\node at (24,-3){$9$};
			%\node at (14,-2){$3$};
			%\node at (14,-1){$2$};
			%\node at (14,0){$1$};
			%\node at (14,1){$0$};
			%\node at (14,2){$4$};
			\end{tikzpicture}
		\end{center}
	\end{eg}

	The extended affine symmetric group $\widehat{\mathfrak{S}}_h$ is the semidirect product $\Z^h\rtimes \mathfrak{S}_h$. There is a natural action of $\widehat{\mathfrak{S}}_h$ on the set of $e$-abaci 
	with exactly one bead on a fixed subset of $h$ runners, and no beads on the other runners: $\Z^h$ acts as the group of horizontal translations of the beads on their runners, and $\mathfrak{S}_h$ as permutations of the $h$ runners containing the beads. This action is locally nilpotent for the subgroup $\Z^h_{<0}$ consisting of left translations of the beads.
	In terms of partitions, the meaning is as follows: let $\rho$ be an $e$-core of some unitary partition; equivalently, $\cA_e(\rho)$ has its beads pushed all the way to the left and they are concentrated in the leftmost column of the $e$-abacus. Let $\cP_e(\rho)_h$ be the union of all $\Po_e(\lambda)$, $\lambda$ an $e$-unitary partition such that the $e$-core of $\lambda$ is $\rho$ and $\lambda$ has $h$ columns.  Let $\epsilon_i=(0,\dots,1,\dots,0)\in\Z^h$ with the $1$ in the $i$'th position. Then $\epsilon_i$ acts on $\mu\in\cP_e(\rho)_h$ by shifting the bead on the $i$'th runner containing a bead one unit to the right; on the Young diagram of $\mu$ it adds an $e$-rimhook whose arm-length is at most $h-1$. 
	 Observe that $\cP_e(\rho)_h$ is generated by $\cA_e(\rho)$ under the action of $\widehat{\mathfrak{S}}_h$:
	$$\widehat{\mathfrak{S}}_h\cdot \cA_e(\rho)=\cP_e(\rho)_h %=
	%\bigoplus_{\substack{\lambda\in\cP \\ e\text{-core}(\lambda)=\rho \\ L(\lambda)\in\cO_{1/e}\text{ is unitary} \\ \#\{\text{columns}(\lambda)\}=h}}[\cO_{1/e}(S_{|\lambda|})_{\leq \lambda}]
%Cherednik algebra not yet defined
	$$
	$\cP_e(\rho)_h$ is naturally identified with the monoid $\Z^{h}_{\geq 0}$ as a left $\widehat{\mathfrak{S}}_h$-module by identifying an abacus $\cA\in\cP_e(\rho)_h$ with the $h$-tuple of its beads' positions $(\beta_1,\dots,\beta_h)\in\Z^h_{\geq 0}$.
	
	\subsection{The homological degree statistic}
	
	Let $\lambda$ be an $e$-unitary partition. We recall the homological degree statistic on $\Po_e(\lambda)$ introduced by Berkesch-Griffeth-Sam. 
	
	\begin{defn}
		Suppose $\cA$ is an $e$-abacus with at most one bead on each runner. A \textit{disorder} of $\cA$ is an unordered pair $\{i,j\}$ such that runners $i$ and $j$ both contain a bead, satisfying $\beta_i>\beta_j$ and $b_j$ is above $b_i$. In other words, a pair of beads of $\cA$ yields a disorder if one bead is above and strictly to the left of the other bead.
	\end{defn}
	
%	\begin{eg}
%		In Example \ref{exl:girlwiththeflaxenhair}, every pair of beads yields a disorder; the disorders are $\{3,2\}$, $\{3,1\}$, $\{3,0\}$, $\{2,1\}$, $\{2,0\}$, $\{1,0\}$. 
%	\end{eg}
	
	\begin{defn}\label{def:hd} \cite[Definition 4.3]{conjecture}
		Let $\mu\in\Po_e(\lambda)$. The \textit{homological degree} of $\mu$, written $\hd(\mu)$, is the sum of the differences of all horizontal positions of beads in $\cA_e(\mu)$ minus the number of disorders of $\cA_e(\mu)$:
		$$\hd(\mu)=\sum_{\substack{i,j\in\Z/eZ \\ b_i,b_j\neq\emptyset \\ b_j \text{ is below }b_i}}|\beta_i-\beta_j|-\#\{\text{disorders of }\cA_e(\mu)\}$$
	\end{defn}
	
\begin{eg}
In Example \ref{exl:girlwiththeflaxenhair}, let $\nu$ denote the partition whose abacus is on the left, and let $\mu=s_0(\nu)$ as in the picture. Then $\cA_5(\nu)$ has $6$ disorders and $\hd(\nu)=1+2+8+1+7+6-6=19$; $\cA_5(\mu)$ has $1$ disorder and $\hd(\mu)=6+5+6+1+1-1=18$. Observe that $s_0$ changed the homological degree by $1$.
\end{eg}

	\subsection{Homological degree produced recursively by elements of $\widehat{\mathfrak{S}}_e$}
	Notice that empty runners of $\cA_e(\mu)$ play no role in $\hd(\mu)$; 
	if the empty runners are removed from $\cA_e(\mu)$, the homological degree remains the same. For simplicity of the formulas and exposition, we therefore work in the case that there are no empty runners, that is, $h=e$ columns and $\lambda=(e^k)$ for some $k\in\N$. The empty runners can be put back in at the end. %So in what follows, we are working with $\cP_e(\emptyset)_e$, the set of all $e$-abaci with exactly one bead on each runner, on which there is an action of $\widehat{S}_e$. Each abacus in $\cP_e(\emptyset)_e$ belongs to $\Po_e(e^k)$ for a unique $k\in\Z_{\geq 0}$, where $k$ is the total number of spaces to the left of beads in the abacus, that is, $k=\sum_{i=0}^{e-1}\beta_i$ where $\beta_i$ is the horizontal position of the bead on runner $i$.
	
	Our first characterization of the homological degree produces this statistic recursively starting from the empty partition by applying sequences of special elements $\tau_i\in\widehat{\mathfrak{S}}_e$, $i=e-1,e-2,\dots,1,0$, in a non-increasing order with respect to $i$. 
	\begin{defn}
		Let $\tau_i\in\widehat{\mathfrak{S}}_e$ be defined as follows: $\tau_i$ fixes the bottom $i$ runners; on the top $e-i$ runners, it first cyclically rotates the beads in the topwards direction, then shifts one space to the right the bead on the $(e-i)$'th runner from the top.
	\end{defn}
	Each $\tau_i$ is the ``affine generator" of the subgroup $\widehat{\mathfrak{S}}_{e-i}$ of $\widehat{\mathfrak{S}}_e$ which fixes the bottom $i$ runners: $\tau_i$ together with $\mathfrak{S}_{e-i}$ generates $\widehat{\mathfrak{S}}_{e-i}$ \cite[Section 2.1]{LeclercThibon}. 
	We are interested in applying $\tau_i$ to abaci whose bottom $i$ runners have their beads pushed all the way to the left.
	\begin{eg} Consider the $5$-abacus of $(3^{11},2^3,1^{11})$.
		Then $\tau_2$ acts as follows:
		
		
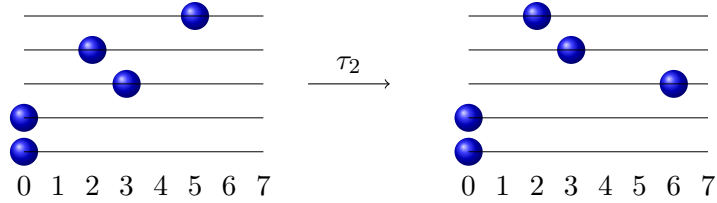
\begin{center}
			\begin{tikzpicture}[scale=.45,every text node part/.style={align=center}]
			\node at (5,2)[circle,shading=ball,minimum width=3pt] (ball){};
			\node at (2,1)[circle,shading=ball,minimum width=3pt] (ball){};
			\node at (3,0)[circle,shading=ball,minimum width=3pt] (ball){};
			\node at (0,-1)[circle,shading=ball,minimum width=3pt] (ball){};
			\node at (0,-2)[circle,shading=ball,minimum width=3pt] (ball){};
			\draw[thin] plot [smooth] coordinates{(0,-2)(7,-2)};
			\draw[thin] plot [smooth] coordinates{(0,-1)(7,-1)};
			\draw[thin] plot [smooth] coordinates{(0,0)(7,0)};
			\draw[thin] plot [smooth] coordinates{(0,1)(7,1)};
			\draw[thin] plot [smooth] coordinates{(0,2)(7,2)};
			\node at (0,-3){$0$};
			\node at (1,-3){$1$};
			\node at (2,-3){$2$};
			\node at (3,-3){$3$};
			\node at (4,-3){$4$};
			\node at (5,-3){$5$};
			\node at (6,-3){$6$};
			\node at (7,-3){$7$};
		%	\node at (-1,-2){$4$};
		%	\node at (-1,-1){$3$};
		%	\node at (-1,0){$2$};
		%	\node at (-1,1){$1$};
		%	\node at (-1,2){$0$};
			\node(A) at (8,0){};
			\node(B) at (11,0){};
			\draw[->] (A) -- node [above] {$\tau_2$} (B) ;
			\node at (19,0)[circle,shading=ball,minimum width=3pt] (ball){};
			\node at (15,2)[circle,shading=ball,minimum width=3pt] (ball){};
			\node at (16,1)[circle,shading=ball,minimum width=3pt] (ball){};
			\node at (13,-1)[circle,shading=ball,minimum width=3pt] (ball){};
			\node at (13,-2)[circle,shading=ball,minimum width=3pt] (ball){};
			\draw[thin] plot [smooth] coordinates{(13,-2)(20,-2)};
			\draw[thin] plot [smooth] coordinates{(13,-1)(20,-1)};
			\draw[thin] plot [smooth] coordinates{(13,0)(20,0)};
			\draw[thin] plot [smooth] coordinates{(13,1)(20,1)};
			\draw[thin] plot [smooth] coordinates{(13,2)(20,2)};
			\node at (13,-3){$0$};
			\node at (14,-3){$1$};
			\node at (15,-3){$2$};
			\node at (16,-3){$3$};
			\node at (17,-3){$4$};
			\node at (18,-3){$5$};
			\node at (19,-3){$6$};
			\node at (20,-3){$7$};
		%	\node at (12,-2){$4$};
		%	\node at (12,-1){$3$};
		%	\node at (12,0){$2$};
		%	\node at (12,1){$1$};
		%	\node at (12,2){$0$};
			\end{tikzpicture}
		\end{center}
Observe that $\tau_2$ increased the homological degree of the abacus by $2$.	
	\end{eg}
	
	Suppose $\tau$ is a partition all of whose parts are of size at most $e-1$, and which may contain parts of size $0$, so $\tau=((e-1)^{a_{e-1}},(e-2)^{a_{e-2}},\dots,1^{a_1},0^{a_0})$. Thus $\tau$ fits inside an $e-1$ by $k$ box, where $k$ is the total number of parts of $\tau$. Now identify $\tau$ with the element of $\widehat{\mathfrak{S}}_e$ given by the composition of operators $\tau_0^{a_0}\tau_1^{a_1}\dots\tau_{e-1}^{a_{e-1}}$. By abuse of notation  we will also call this element $\tau$. The proof of the following lemma is straightforward:
	\begin{lem}\label{thm:grassmannstyle}
		Let $\tau=((e-1)^{a_{e-1}},(e-2)^{a_{e-2}},\dots,1^{a_1},0^{a_0})$ with  $\sum_{i=0}^{e-1}a_i=k$, $k\in\Z_{\geq 0}$. Then $\tau(\cA_e(\emptyset))=\cA_e(\mu)$ with $\mu\in\Po_e(e^k)$. Any $\mu\in\Po_e(e^k)$ is produced in this way from a unique such $\tau$, and we have: $$\hd(\mu)=\sum_{i=0}^{e-1}ia_i=|\tau|$$
	\end{lem}
	
	Let $\lambda$ be an arbitrary $e$-unitary partition. By removing the empty runners from the $e$-abaci in $\Po_e(\lambda)$, there is likewise a natural bijection between the partitions $\mu$ in $\Po_e(\lambda)$ and partitions $\tau$ which fit inside an $(h-1)$ by $k$ box,

\!\!\!\!\!
$$\{\mu\in\Po_e(\lambda)\}\overset{\substack{\Phi\\\simeq }}{\longleftarrow}\{\tau\subset(h-1)^k\}, $$	
	given by $\Phi(\tau)=\tau(\emptyset)$ (where $\tau$ on the right-hand-side is the corresponding element of $\widehat{\mathfrak{S}}_h$ as described above). This bijection identifies $\hd(\mu)$ with $|\tau|$.
	
	\begin{rmk}
		Such a bijection turns up elsewhere in representation theory: notably, partitions $\tau$ which fit inside an $(h-1)$ by $k$ box also parametrize (1) the simple and standard modules of a regular block $\cB^\mathfrak{p}$ of parabolic category $\cO^\mathfrak{p}$ for $\mathfrak{gl}(h-1+k)$ with respect to the maximal parabolic $\mathfrak{gl}(h-1)\times\mathfrak{gl}(k)$ \cite{Stroppel}; (2) the Schubert cells in the Grassmannian $\mathsf{Gr}(k,h-1+k)=\mathsf{Gr}(h-1,h-1+k)$. The category $\cB^\mathfrak{p}$ is equivalent to perverse sheaves on the Grassmannian 
		%\cite{Braden},\cite[Theorem 5.8.1]{Stroppel},
		\cite{Braden,Stroppel},
		 explaining the coincidence of (1) and (2). Let $L_{(h-1)^k}$ denote the simple module in $\cB^\mathfrak{p}$ labeled by $\tau=(h-1)^k$, the unique maximal element of the poset (the poset structure is given by inclusion of Young diagrams). The bijection following Lemma \ref{thm:grassmannstyle} identifies the characters of unitary $L(\lambda)\in[\cO_{1/e}(\fS_n)]$ and $L_{(h-1)^k} \in[\cO^\mathfrak{p}]$. Moreover, $L_{(h-1)^k}$ has a BGG resolution \cite{BoeHunziker}
		and via the bijection $\Phi$ we obtain a natural bijection between the Verma modules appearing in degree $i$ of the respective resolutions in $\cO_{1/e}(S_n)$ and $\mathcal{O}^\mathfrak{p}$. However, the categories $\cO_{1/e}(\fS_n)_{\leq\lambda}$ and $\cB^\mathfrak{p}$ are not equivalent if $k>2$, and as a poset $\Po_e(\lambda)$ has ``extra edges" coming from the $\widetilde{\mathfrak{S}}_h$-action if $k>2$. The difference between the resolutions of $L(\lambda)$ and $L_{(h-1)^k}$ thus manifests itself in the maps in the complex.
	\end{rmk}
	
%	We also observe that Lemma \ref{thm:grassmannstyle} yields a symmetry of the poset $\Po_e(\lambda)$ that is not immediately apparent from the definition of the homological degree:
%	\begin{lem} Let $\lambda\in\cP_e^{u}$ have $h$ columns and $e$-weight $k$. Then:
%		$$\#\{\mu\in\Po_e(\lambda)\mid\hd(\mu)=i\}=\#\{\mu\in\Po_e(\lambda)\mid \hd(\mu)=(h-1)k-i\}$$
%	\end{lem}
%	\begin{proof}
%		The number of partitions $\tau$ of size $i$ that fit in an $h-1$ by $k$ box is obviously equal to the number of partitions $\tau$ of size $(h-1)k-i$ that fit in an $h-1$ by $k$ box. 
%	\end{proof}
	
	\subsection{Homological degree via rimhooks of minimal leg-length}\label{subsect:hd rimhooks}
	Consider again the case there is exactly one bead on every runner of the abacus. By the definition of $\tau_i$, it follows that the effect of applying $\tau=((e-1)^{a_{e-1}}, (e-2)^{a_{e-2}},\dots, 1^{a_1},0^{a_0})$ to the empty partition is to build a Young diagram $\lambda$
	by successively dropping $e$-rimhooks which meet the leftmost column (with leg-lengths $e-1$ ($a_{e-1}$ times), $e-2$ ($a_{e-2}$ times) and so on) Tetris-style on top of the partition constructed so far, then letting the boxes slide down the columns so that the result is a partition. This can change the shape of the previous rimhooks that were added, but not the set of their leg-lengths. 
	Thus we obtain a second combinatorial explanation of the homological degree:
	\textit{if $e=h$ then $\hd(\lambda)$ is the sum of the leg-lengths of the $e$-rimhooks of minimal leg-length composing $\lambda$}. \textit{If $e>h$ then $\hd(\lambda)$ is the sum of the leg-lengths of the $e$-rimhooks of minimal leg-length composing $\lambda$ minus $(e-h)k$, where $k=e\text{-weight}(\lambda)$.} This can be restated in a uniform way by considering the arm-lengths instead of the leg-lengths of the rimhooks: $\hd(\mu)$ is equal to $(h-1)k$ minus the sum of armlengths of the (minimal leg-length) rimhooks removed.
	
	\begin{eg}Let $e=h=5$ and $\tau=(3,3,1,0)$. Then $\tau(\emptyset)=(5,4,2^5,1)=:\lambda$ and $\hd(\lambda)=7$. We show the process of applying $\tau$ on abaci and partitions and the four $5$-rimhooks of minimal leg-lengths $3,3,1,0$ which compose $\lambda$:
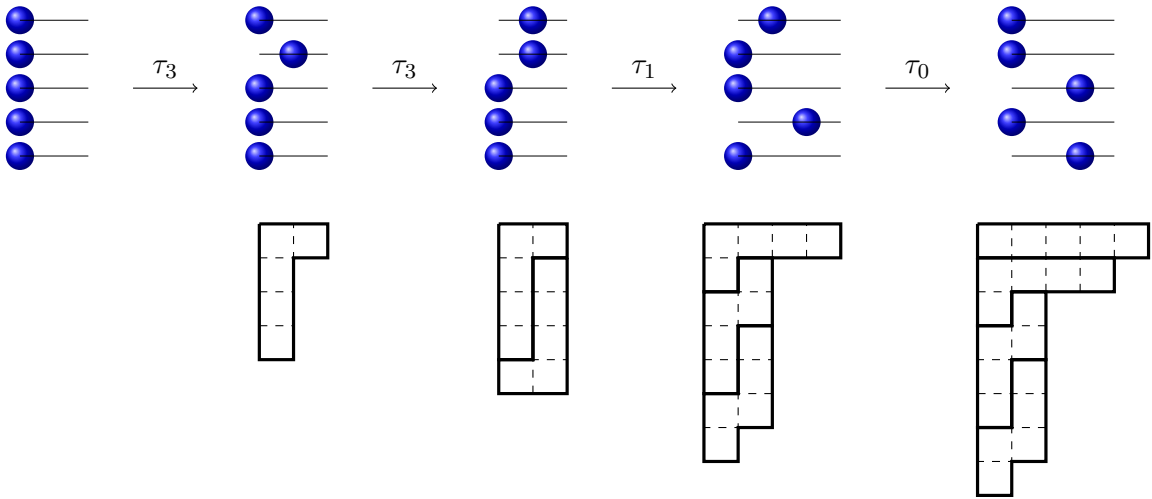
		\begin{center}
			\begin{tikzpicture}[scale=.45,every text node part/.style={align=center}]
			\node at (-20,2)[circle,shading=ball,minimum width=3pt] (ball){};
			\node at (-20,1)[circle,shading=ball,minimum width=3pt] (ball){};
			\node at (-20,0)[circle,shading=ball,minimum width=3pt] (ball){};
			\node at (-20,-1)[circle,shading=ball,minimum width=3pt] (ball){};
			\node at (-20,-2)[circle,shading=ball,minimum width=3pt] (ball){};
			\draw[thin] plot [smooth] coordinates{(-20,-2)(-18,-2)};
			\draw[thin] plot [smooth] coordinates{(-20,-1)(-18,-1)};
			\draw[thin] plot [smooth] coordinates{(-20,0)(-18,0)};
			\draw[thin] plot [smooth] coordinates{(-20,1)(-18,1)};
			\draw[thin] plot [smooth] coordinates{(-20,2)(-18,2)};
			\node(A1) at (-17,0){};
			\node(B1) at (-14.5,0){};
			\draw[->] (A1) -- node [above] {$\tau_3$} (B1) ;
			\node at (-13,2)[circle,shading=ball,minimum width=3pt] (ball){};
			\node at (-12,1)[circle,shading=ball,minimum width=3pt] (ball){};
			\node at (-13,0)[circle,shading=ball,minimum width=3pt] (ball){};
			\node at (-13,-1)[circle,shading=ball,minimum width=3pt] (ball){};
			\node at (-13,-2)[circle,shading=ball,minimum width=3pt] (ball){};
			\draw[thin] plot [smooth] coordinates{(-13,-2)(-11,-2)};
			\draw[thin] plot [smooth] coordinates{(-13,-1)(-11,-1)};
			\draw[thin] plot [smooth] coordinates{(-13,0)(-11,0)};
			\draw[thin] plot [smooth] coordinates{(-13,1)(-11,1)};
			\draw[thin] plot [smooth] coordinates{(-13,2)(-11,2)};
			\node(A2) at (-10,0){};
			\node(B2) at (-7.5,0){};
			\draw[->] (A2) -- node [above] {$\tau_3$} (B2) ;
			\node at (-5,2)[circle,shading=ball,minimum width=3pt] (ball){};
			\node at (-5,1)[circle,shading=ball,minimum width=3pt] (ball){};
			\node at (-6,0)[circle,shading=ball,minimum width=3pt] (ball){};
			\node at (-6,-1)[circle,shading=ball,minimum width=3pt] (ball){};
			\node at (-6,-2)[circle,shading=ball,minimum width=3pt] (ball){};
			\draw[thin] plot [smooth] coordinates{(-6,-2)(-4,-2)};
			\draw[thin] plot [smooth] coordinates{(-6,-1)(-4,-1)};
			\draw[thin] plot [smooth] coordinates{(-6,0)(-4,0)};
			\draw[thin] plot [smooth] coordinates{(-6,1)(-4,1)};
			\draw[thin] plot [smooth] coordinates{(-6,2)(-4,2)};
			\node(A3) at (-3,0){};
			\node(B3) at (-.5,0){};
			\draw[->] (A3) -- node [above] {$\tau_1$} (B3) ;
			\node at (2,2)[circle,shading=ball,minimum width=3pt] (ball){};
			\node at (1,1)[circle,shading=ball,minimum width=3pt] (ball){};
			\node at (1,0)[circle,shading=ball,minimum width=3pt] (ball){};
			\node at (3,-1)[circle,shading=ball,minimum width=3pt] (ball){};
			\node at (1,-2)[circle,shading=ball,minimum width=3pt] (ball){};
			\draw[thin] plot [smooth] coordinates{(1,-2)(4,-2)};
			\draw[thin] plot [smooth] coordinates{(1,-1)(4,-1)};
			\draw[thin] plot [smooth] coordinates{(1,0)(4,0)};
			\draw[thin] plot [smooth] coordinates{(1,1)(4,1)};
			\draw[thin] plot [smooth] coordinates{(1,2)(4,2)};
			\node(A4) at (5,0){};
			\node(B4) at (7.5,0){};
			\draw[->] (A4) -- node [above] {$\tau_0$} (B4) ;
			\node at (9,2)[circle,shading=ball,minimum width=3pt] (ball){};
			\node at (9,1)[circle,shading=ball,minimum width=3pt] (ball){};
			\node at (11,0)[circle,shading=ball,minimum width=3pt] (ball){};
			\node at (9,-1)[circle,shading=ball,minimum width=3pt] (ball){};
			\node at (11,-2)[circle,shading=ball,minimum width=3pt] (ball){};
			\draw[thin] plot [smooth] coordinates{(9,-2)(12,-2)};
			\draw[thin] plot [smooth] coordinates{(9,-1)(12,-1)};
			\draw[thin] plot [smooth] coordinates{(9,0)(12,0)};
			\draw[thin] plot [smooth] coordinates{(9,1)(12,1)};
			\draw[thin] plot [smooth] coordinates{(9,2)(12,2)};
			\draw[very thick] plot coordinates{(-13,-4)(-11,-4)(-11,-5)(-12,-5)(-12,-8)(-13,-8)(-13,-4)};
			\draw[dashed] plot coordinates{(-13,-5)(-12,-5)(-12,-4)};
			\draw[dashed] plot coordinates{(-13,-6)(-12,-6)};
			\draw[dashed] plot coordinates{(-13,-7)(-12,-7)};
			\draw[very thick] plot coordinates{(-6,-4)(-4,-4)(-4,-5)(-5,-5)(-5,-8)(-6,-8)(-6,-4)};
			\draw[very thick] plot coordinates{(-4,-5)(-4,-9)(-6,-9)(-6,-8)(-5,-8)(-5,-5)(-4,-5)};
			\draw[dashed] plot coordinates{(-6,-5)(-5,-5)(-5,-4)};
			\draw[dashed] plot coordinates{(-6,-6)(-4,-6)};
			\draw[dashed] plot coordinates{(-6,-7)(-4,-7)};
			\draw[dashed] plot coordinates{(-5,-9)(-5,-8)(-4,-8)};
			\draw[very thick] plot coordinates{(0,-4)(4,-4)(4,-5)(1,-5)(1,-6)(0,-6)(0,-4)};
			\draw[very thick] plot coordinates{(1,-5)(2,-5)(2,-7)(1,-7)(1,-9)(0,-9)(0,-6)(1,-6)(1,-5)};
			\draw[very thick] plot coordinates{(1,-7)(2,-7)(2,-10)(1,-10)(1,-11)(0,-11)(0,-9)(1,-9)(1,-7)};
			\draw[dashed] plot coordinates{(0,-5)(1,-5)(1,-4)};
			\draw[dashed] plot coordinates{(2,-5)(2,-4)};
			\draw[dashed] plot coordinates{(3,-5)(3,-4)};
			\draw[dashed] plot coordinates{(0,-7)(1,-7)(1,-6)(2,-6)};
			\draw[dashed] plot coordinates{(0,-8)(2,-8)};
			\draw[dashed] plot coordinates{(2,-9)(1,-9)(1,-10)(0,-10)};
			\draw[very thick] plot coordinates{(8,-4)(13,-4)(13,-5)(8,-5)(8,-4)};
			\draw[very thick] plot coordinates{(8,-5)(12,-5)(12,-6)(9,-6)(9,-7)(8,-7)(8,-5)};
			\draw[very thick] plot coordinates{(9,-6)(10,-6)(10,-8)(9,-8)(9,-10)(8,-10)(8,-7)(9,-7)(9,-6)};
			\draw[very thick] plot coordinates{(9,-8)(10,-8)(10,-11)(9,-11)(9,-12)(8,-12)(8,-10)(9,-10)(9,-8)};
			\draw[dashed] plot coordinates{(8,-6)(9,-6)(9,-4)};
			\draw[dashed] plot coordinates{(10,-6)(10,-4)};
			\draw[dashed] plot coordinates{(11,-6)(11,-4)};
			\draw[dashed] plot coordinates{(12,-5)(12,-4)};
			\draw[dashed] plot coordinates{(10,-7)(9,-7)(9,-8)(8,-8)};
			\draw[dashed] plot coordinates{(8,-9)(10,-9)};
			\draw[dashed] plot coordinates{(10,-10)(9,-10)(9,-11)(8,-11)};
			\end{tikzpicture}
		\end{center}
	\end{eg}

	\subsection{Homological degree is the length function.} % on $\widetilde{\mathfrak{S}}_h$}
	We now give a third combinatorial description of the homological degree by identifying it with the length function on $\widetilde{\mathfrak{S}}_h$. This unifies the   combinatorics   of abaci with that of alcove geometries and allows us to describe the BGG complex in type $A$ in terms of abaci. Let $\lambda$ be a unitary partition, and suppose $\lambda$ has $h<e$ columns.
	\begin{lem}\label{lem:hd=ell}
		The following are equivalent for $\nu,\mu\in\Po_e(\lambda)$, $\mu \rhd \nu$:
		\begin{itemize}[leftmargin=*]
		\item $\ell(\nu)=\ell(\mu)+1$ and $\nu$ is obtained from $\mu$ by moving a column of boxes as in Theorem \ref{onecol};
			\item $t\mu=\nu$ for some transposition $t\in\widetilde{\mathfrak{S}}_h$ acting on abaci as above, subject to the following conditions on the beads $\beta_i$ of $\cA_e(\mu)$: \begin{enumerate}[leftmargin=*]\item if $t\in \mathfrak{S}_h$ and swaps runners $i$ and $j$, then for each runner $k$ between runners $i$ and $j$, $\beta_k\notin[\beta_i,\beta_j]$;
				\item if $t$ is conjugate to $s_0$ and acts nontrivially on runners $i$ and $j$, runner $i$ below runner $j$, then: for each runner $k$ below runner $i$, $\beta_k\notin[\beta_i,\beta_j+1]$ and for each runner $\ell$ above runner $j$, $\beta_\ell\notin[\beta_{i-1},\beta_j]$.
			\end{enumerate}
		\end{itemize}
		Therefore, the homological degree statistic on $\Po_e(\lambda)$ coincides with the length function on $\Po_e(\lambda)$ coming from the $\widetilde{\mathfrak{S}}_h$ alcove geometry.
	\end{lem}
	
	\begin{proof}
	This is a  translation of 1-column moves from the language of Young diagrams into the language of abaci. It follows from a direct computation using Definition \ref{def:hd} that the conditions for a transposition $t$ to increase the homological degree by $1$ are exactly those given by (1) and (2).
	\end{proof}

 \!\!\!\!\!\!\!\!
 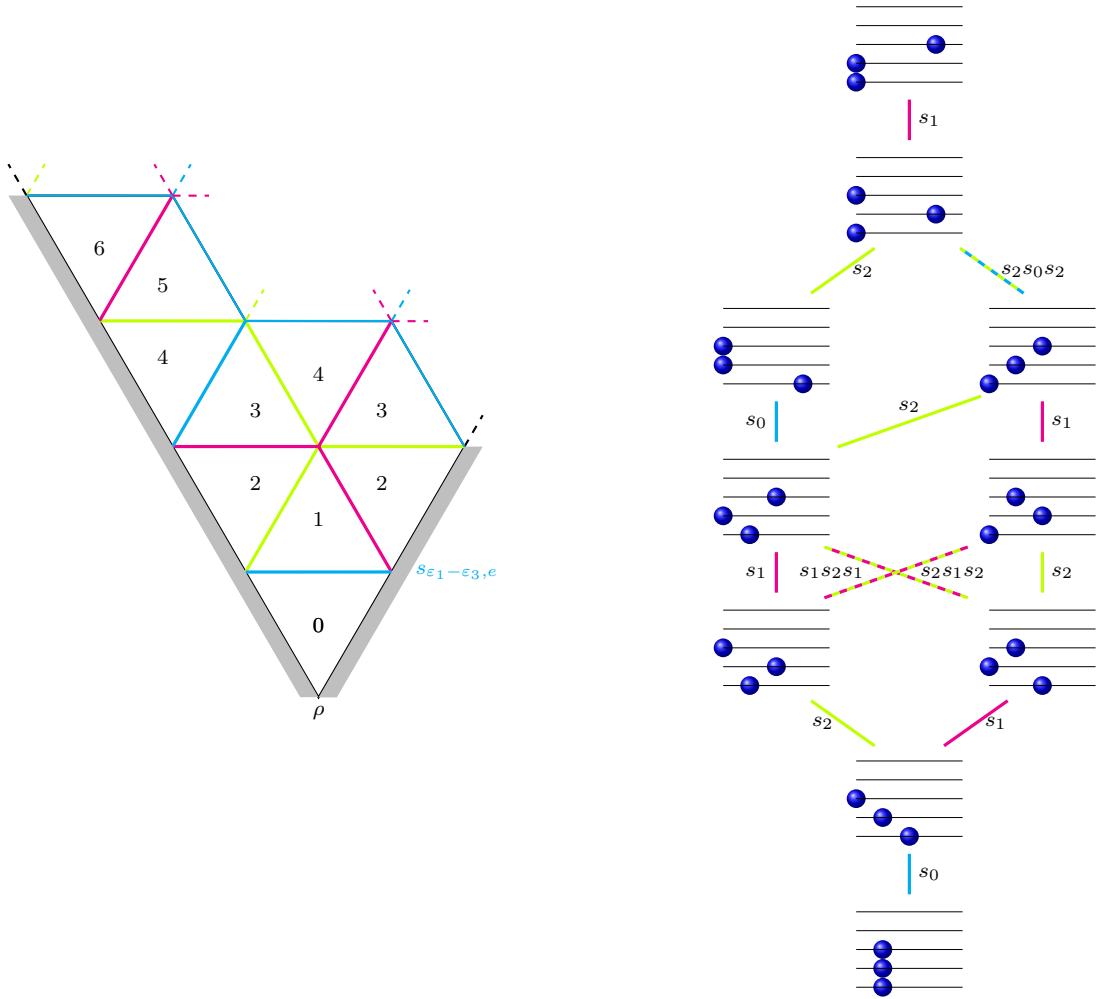
\begin{figure}[ht!]\label{tricolor}
 \scalefont{0.75}  
 $$\begin{minipage}{10cm}
 \begin{tikzpicture}    [scale=1.2]%[scale=1.4]  
     \clip (-3.2,-0.3) rectangle ++(7.4,7.7);
      \path (120:2.4cm)++(60:2.4cm) coordinate (BB);                    
       \path (0,0) coordinate (origin);

   \path[dotted] (origin)++(60:6.4cm) coordinate (aa);
      \path[dotted] (origin)++(60:4.8cm) coordinate (ab);
      
            \path[dotted] (origin)++(60:4.8cm)++(120:1.6) coordinate (ba);
   \path[dotted] (origin)++(120:6.4cm) coordinate (bb);
   \path[dotted] (origin)++(120:4.8cm) coordinate (bb1);
   \path[dotted] (bb1)++(60:3.2cm) coordinate (bb2);
      \path[dotted] (bb2)++(180:1.6cm) coordinate (bb25);
            \path[dotted] (ab)++(0:1.6cm) coordinate (bb325);
      \path[dotted] (bb25)++(-60:1.6cm) coordinate (bb252);
      
      \path[dotted] (bb3)++(180:1.6cm) coordinate (bb4);

      \path[dotted] (bb2)++(-60:1.6cm) coordinate (bb3);

	\path[dotted] (ab)++(180:1.6cm) coordinate (ee);

            \path[dotted] (origin)++(120:1.6cm)++(60:1.6cm) coordinate (dd);
				\path[dotted] (dd)++(-60:1.6cm) coordinate (dd2);
				\path[dotted] (dd)++(120:1.6cm) coordinate (dd0);
                        \path[dotted] (dd)++(0:1.6cm)  coordinate (dd3);

   \draw[thick]  (ab)--(ba)--(bb3)--(bb2)--(bb25)--(bb252)--(dd2)--(dd3)
  --(ab) ;
  \draw[very thick, lime] (dd)--(ee)--(bb3)--(bb4);
  \draw[very thick, lime] (ab)--(ee);
  \draw[very thick, magenta] (bb2)--(bb4);
  \draw[very thick, magenta] (ba)--(ee)--(dd3);
  \draw[very thick, magenta] (dd0)--(ee);
  \draw[very thick, cyan] (dd)--(dd3) node[right] {$\ \  s_{\varepsilon_1-\varepsilon_3,e} $};
  \draw[very thick, cyan] (dd0)--(bb3);
     
                      %%%CHRIS 

                        \draw[ thick,dashed](ab)--++(60:0.4);
                                                \draw[cyan,  thick,dashed](ba)--++(60:0.4);
   \draw[magenta, thick,dashed](ba)--++(0:0.4);
      \draw[magenta, thick,dashed](ba)--++(120:0.4);

         \draw[lime, thick,dashed](bb3)--++(60:0.4);         \draw[cyan, thick](bb3)--++(120:1.6);
              \draw[cyan, thick](bb3)--++(0:1.6);

   \draw[cyan, thick,dashed](bb2)--++(60:0.4);   \draw[magenta, thick,dashed](bb2)--++(120:0.4);   \draw[magenta, thick,dashed](bb2)--++(0:0.4);

   \draw[lime, thick,dashed](bb25)--++(60:0.4);   \draw[ thick,dashed](bb25)--++(120:0.4);

   \fill[gray!50](bb25)--++(180:0.2)--++(-60:1.6*4) --++(0:0.2);

 \fill[gray!50](ab)  --++(0:0.2) --++(-120:3.2)--++(-180:0.2) ; % --(ab)     --++(180:0.2)--++(-60:6.4);                   

\draw(dd2) node [below] {$\rho$};

 \path(dd2)--++(120:3.2)--++(0:1.6)  coordinate (stein);
  \path(dd2)--++(90:0.8)  coordinate (zero);
  
    \path(stein)--++(-90:0.8)  coordinate (st1);
        \path(stein)--++(-30:0.8)  coordinate (st2);
                \path(stein)--++(210:0.8)  coordinate (st25);
                \path(stein)--++(30:0.8)  coordinate (st3);                \path(stein)--++(150:0.8)  coordinate (st35);
                  \path(stein)--++(90:0.8)  coordinate (st4); 
  \draw(zero) node {$0$};
  
  \draw(st1) node {$1$};  \draw(st2) node {$2$}; 
   \draw(st25) node {$2$};\draw(st3) node {$3$}; 
   \draw(st35) node {$3$};   \draw(st4) node {$4$};

 \path(stein)--++(120:3.2)--++(-120:1.6)  coordinate (stein);

 \path(stein)--++(-30:0.8)  coordinate (st2);
                \path(stein)--++(210:0.8)  coordinate (st25);
                \path(stein)--++(30:0.8)  coordinate (st3);                \path(stein)--++(150:0.8)  coordinate (st35);
                  \path(stein)--++(90:0.8)  coordinate (st4); 
  \draw(zero) node {$0$};
  
   \draw(st2) node {$4$}; 
     \draw(st3) node {$5$}; 
    \draw(st4) node {$6$};
%%%

%re-colouring
    \draw[cyan,  thick ](ba)--++(-60:1.6);

   \draw[cyan, thick ](bb2)--++(180:1.6);
%%%END CHRIS colours

%%partitions

%
%  
%  
%  
 \end{tikzpicture}
 \end{minipage}
% \caption{The alcoves corresponding to partitions in $\Po_e(\lambda)$ when $h=k=3$ and $e\geq 3$; $\lambda$ can be any $e$-unitary partition with $3$ columns and $e$-weight $3$ and belongs to the bottom alcove (the fundamental alcove). Crossing a wall of color $i$ corresponds to applying $s_i$ to the partition in that alcove. The colors corresponding to $i$: ${\color{cyan}0}$, ${\color{magenta}1}$, ${\color{lime}2}$.}
%\end{SCfigure*}
%
%\begin{SCfigure*}
%
\begin{minipage}{7cm}	\begin{tikzpicture}[xscale=.35,yscale=.25,every text node part/.style={align=center}]
			\node at (3,0)[circle,shading=ball,inner sep=2.5pt](ball) {};
			\node at (0,-1)[circle,shading=ball,inner sep=2.5pt] (ball){};
			\node at (0,-2)[circle,shading=ball,inner sep=2.5pt] (ball){};
			\draw[thin] plot [smooth] coordinates{(0,-2)(4,-2)};
			\draw[thin] plot [smooth] coordinates{(0,-1)(4,-1)};
			\draw[thin] plot [smooth] coordinates{(0,0)(4,0)};
			\draw[thin] plot [smooth] coordinates{(0,1)(4,1)};
			\draw[thin] plot [smooth] coordinates{(0,2)(4,2)};
			
			\node(A) at (2,-2.5){};
			\node(B) at (2,-5.5){};
			\draw[very thick, color=magenta] (B) -- node [right] {${\color{black}s_1}$} (A) ;
			
				\node at (0,-8)[circle,shading=ball,inner sep=2.5pt] (ball){};
			\node at (3,-9)[circle,shading=ball,inner sep=2.5pt] (ball){};
			\node at (0,-10)[circle,shading=ball,inner sep=2.5pt] (ball){};
			\draw[thin] plot [smooth] coordinates{(0,-6)(4,-6)};
			\draw[thin] plot [smooth] coordinates{(0,-7)(4,-7)};
			\draw[thin] plot [smooth] coordinates{(0,-8)(4,-8)};
			\draw[thin] plot [smooth] coordinates{(0,-9)(4,-9)};
			\draw[thin] plot [smooth] coordinates{(0,-10)(4,-10)};
			
						\node(A) at (1,-10.5){};
			\node(B) at (-2,-13.5){};
			\draw[very thick, color=lime] (B) -- node [right] {${\color{black}s_2}$} (A) ;
					\node at (-5,-16)[circle,shading=ball,inner sep=2.5pt] (ball){};
			\node at (-5,-17)[circle,shading=ball,inner sep=2.5pt] (ball){};
			\node at (-2,-18)[circle,shading=ball,inner sep=2.5pt] (ball){};
			\draw[thin] plot [smooth] coordinates{(-1,-14)(-5,-14)};
			\draw[thin] plot [smooth] coordinates{(-1,-15)(-5,-15)};
			\draw[thin] plot [smooth] coordinates{(-1,-16)(-5,-16)};
			\draw[thin] plot [smooth] coordinates{(-1,-17)(-5,-17)};
			\draw[thin] plot [smooth] coordinates{(-1,-18)(-5,-18)};

			\node(A) at (3,-10.5){};
			\node(B) at (6,-13.5){};
			\node(C) at (3.3,-10.5){};
			\node(D) at (6.3,-13.5){};
			\node(E) at (3.6,-10.5){};
			\node(F) at (6.6,-13.5){};
%			\draw[thick, color=lime] (B) -- node [right] {} (A) ;
%			\draw[thick, color=cyan] (D) -- node [right] {} (C) ;
%			\draw[thick, color=lime] (F) -- node [right] {${\color{black}s_2s_0s_2}$} (E) ;
			\draw[very thick,cyan,dash pattern= on 3pt off 3pt] (F)   to  (E);
 \draw[very thick, lime,dash pattern= on 3pt off 3pt,dash phase=3pt] (F)  -- node [right] {${\color{black}s_2s_0s_2}$}  (E);

					\node at (7,-16)[circle,shading=ball,inner sep=2.5pt] (ball){};
			\node at (6,-17)[circle,shading=ball,inner sep=2.5pt] (ball){};
			\node at (5,-18)[circle,shading=ball,inner sep=2.5pt] (ball){};
			\draw[thin] plot [smooth] coordinates{(9,-14)(5,-14)};
			\draw[thin] plot [smooth] coordinates{(9,-15)(5,-15)};
			\draw[thin] plot [smooth] coordinates{(9,-16)(5,-16)};
			\draw[thin] plot [smooth] coordinates{(9,-17)(5,-17)};
			\draw[thin] plot [smooth] coordinates{(9,-18)(5,-18)};
						\node(A) at (-3,-18.5){};
			\node(B) at (-3,-21.5){};
			\draw[very thick, color=cyan] (B) -- node [left] {${\color{black}s_0}$} (A) ;

				\node at (-3,-24)[circle,shading=ball,inner sep=2.5pt] (ball){};
			\node at (-5,-25)[circle,shading=ball,inner sep=2.5pt] (ball){};
			\node at (-4,-26)[circle,shading=ball,inner sep=2.5pt] (ball){};
			\draw[thin] plot [smooth] coordinates{(-1,-22)(-5,-22)};
			\draw[thin] plot [smooth] coordinates{(-1,-23)(-5,-23)};
			\draw[thin] plot [smooth] coordinates{(-1,-24)(-5,-24)};
			\draw[thin] plot [smooth] coordinates{(-1,-25)(-5,-25)};
			\draw[thin] plot [smooth] coordinates{(-1,-26)(-5,-26)};
			
				\node(A) at (7,-18.5){};
			\node(B) at (7,-21.5){};
			\draw[very thick, color=magenta] (B) -- node [right] {${\color{black}s_1}$} (A) ;
			
						\node(A) at (5,-18.5){};
			\node(B) at (-1,-21.5){};
			\draw[very thick, color=lime] (B) -- node[left, above] {${\color{black}s_2}$} (A) ;
			
			\node at (6,-24)[circle,shading=ball,inner sep=2.5pt] (ball){};
			\node at (7,-25)[circle,shading=ball,inner sep=2.5pt] (ball){};
			\node at (5,-26)[circle,shading=ball,inner sep=2.5pt] (ball){};
			\draw[thin] plot [smooth] coordinates{(9,-22)(5,-22)};
			\draw[thin] plot [smooth] coordinates{(9,-23)(5,-23)};
			\draw[thin] plot [smooth] coordinates{(9,-24)(5,-24)};
			\draw[thin] plot [smooth] coordinates{(9,-25)(5,-25)};
			\draw[thin] plot [smooth] coordinates{(9,-26)(5,-26)};
			
						\node(A) at (7,-26.5){};
			\node(B) at (7,-29.5){};
			\draw[very thick, color=lime] (B) -- node [right] {${\color{black}s_2}$} (A) ;

		\node(Y) at (-1.5,-26.5){};
				\node(Z) at (4.5,-29.5){};			
				\node(A) at (-1,-26.5){};
			\node(B) at (5,-29.5){};
			\node(C) at (-0.5,-26.5){};
			\node(D) at (5.5,-29.5){};
%			\draw[thick, color=magenta] (Z) -- node [left] {{{\small ${\color{black}s_1s_2s_1}$}}} (Y);
%			\draw[thick, color=lime] (B) -- node {} (A) ;
%			\draw[thick, color=magenta] (D) -- node {} (C);	

		\draw[very thick,lime,dash pattern= on 3pt off 3pt] (Z)   to  (Y);
 \draw[very thick, magenta,dash pattern= on 3pt off 3pt,dash phase=3pt] (Z)  -- node [left] {${\color{black}s_1s_2s_1}$ \ \ }  (Y); 
				
					\node(Y) at (4.5,-26.5){};
				\node(Z) at (-1.5,-29.5){};			
				\node(A) at (5,-26.5){};
			\node(B) at (-1,-29.5){};
			\node(C) at (5.5,-26.5){};
			\node(D) at (-0.5,-29.5){};
%			\draw[thick, color=lime] (Z) -- node [above] {} (Y);
%			\draw[thick, color=magenta] (B) -- node {} (A) ;
%			\draw[thick, color=lime] (D) -- node [right]{{\small ${\color{black}s_2s_1s_2}$}} (C);	

			\draw[very thick,magenta,dash pattern= on 3pt off 3pt] (Z)   to  (Y);
 \draw[very thick, lime,dash pattern= on 3pt off 3pt,dash phase=3pt] (Z)  -- node [right] {\ \ ${\color{black}s_2s_1s_2}$}  (Y);

			\node at (6,-32)[circle,shading=ball,inner sep=2.5pt] (ball){};
			\node at (5,-33)[circle,shading=ball,inner sep=2.5pt] (ball){};
			\node at (7,-34)[circle,shading=ball,inner sep=2.5pt] (ball){};
			\draw[thin] plot [smooth] coordinates{(9,-32)(5,-32)};
			\draw[thin] plot [smooth] coordinates{(9,-33)(5,-33)};
			\draw[thin] plot [smooth] coordinates{(9,-34)(5,-34)};
			\draw[thin] plot [smooth] coordinates{(9,-30)(5,-30)};
			\draw[thin] plot [smooth] coordinates{(9,-31)(5,-31)};
			
				\node(A) at (-3,-26.5){};
			\node(B) at (-3,-29.5){};
			\draw[very thick, color=magenta] (B) -- node [left] {${\color{black}s_1}$} (A) ;
			
			\node at (-5,-32)[circle,shading=ball,inner sep=2.5pt] (ball){};
			\node at (-3,-33)[circle,shading=ball,inner sep=2.5pt] (ball){};
			\node at (-4,-34)[circle,shading=ball,inner sep=2.5pt] (ball){};
			\draw[thin] plot [smooth] coordinates{(-1,-32)(-5,-32)};
			\draw[thin] plot [smooth] coordinates{(-1,-33)(-5,-33)};
			\draw[thin] plot [smooth] coordinates{(-1,-34)(-5,-34)};
			\draw[thin] plot [smooth] coordinates{(-1,-30)(-5,-30)};
			\draw[thin] plot [smooth] coordinates{(-1,-31)(-5,-31)};
			
				\node(A) at (-2,-34.5){};
			\node(B) at (1,-37.5){};
			\draw[very thick, color=lime] (B) -- node [left] {${\color{black}s_2}$} (A) ;
				
			\node(A) at (6,-34.5){};
			\node(B) at (3,-37.5){};
			\draw[very thick, color=magenta] (B) -- node [right] {${\color{black}s_1}$} (A) ;
			
			\node at (0,-40)[circle,shading=ball,inner sep=2.5pt] (ball){};
			\node at (1,-41)[circle,shading=ball,inner sep=2.5pt] (ball){};
			\node at (2,-42)[circle,shading=ball,inner sep=2.5pt] (ball){};
			\draw[thin] plot [smooth] coordinates{(0,-42)(4,-42)};
			\draw[thin] plot [smooth] coordinates{(0,-39)(4,-39)};
			\draw[thin] plot [smooth] coordinates{(0,-38)(4,-38)};
			\draw[thin] plot [smooth] coordinates{(0,-40)(4,-40)};
			\draw[thin] plot [smooth] coordinates{(0,-41)(4,-41)};
		
				\node(A) at (2,-42.5){};
			\node(B) at (2,-45.5){};
			\draw[very thick, color=cyan] (B) -- node [right] {${\color{black}s_0}$} (A) ;
			
			\node at (1,-48)[circle,shading=ball,inner sep=2.5pt] (ball){};
			\node at (1,-49)[circle,shading=ball,inner sep=2.5pt] (ball){};
			\node at (1,-50)[circle,shading=ball,inner sep=2.5pt] (ball){};
			\draw[thin] plot [smooth] coordinates{(0,-48)(4,-48)};
			\draw[thin] plot [smooth] coordinates{(0,-49)(4,-49)};
			\draw[thin] plot [smooth] coordinates{(0,-47)(4,-47)};
			\draw[thin] plot [smooth] coordinates{(0,-50)(4,-50)};
			\draw[thin] plot [smooth] coordinates{(0,-46)(4,-46)};
			\end{tikzpicture}
			\end{minipage}$$
			\caption{
On the left we have the alcoves corresponding to partitions in $\Po_e(\lambda)$ when $h=3$ and $e =5$.  The fundamental alcove is at the bottom and contains  $(3^5)\in \mathcal{F}_{15}(3)$. 
Each alcove contains a number   indicateing the length/homological degree   for a point in that alcove.  
 The grey region denotes the non-dominant region.  The dotted lines indicate that we tile one sixth of $\mathbb{R}^2$ when we let $n\rightarrow \infty$.  
 Crossing a wall of color $i$ corresponds to applying $s_i$ to the partition in that alcove, with $i$: ${\color{cyan}0}$, ${\color{magenta}1}$, ${\color{lime}2}$.	 
On the right-hand side, we have extracted  the poset $\Po_e(\lambda)$.  The homological degree increases from the bottom (where it is zero) to the top (where it is 6). The edges of the poset are coloured and decorated so as to facilitate comparison between the two pictures.  }
\end{figure}

The superficial difference between the two $\widetilde{\mathfrak{S}}_h$ actions on $\Po_e(\lambda)$, the one coming from abaci, the other from alcove geometry, is simply the difference between generating sets for $\widetilde{\mathfrak{S}}_h$: the generators $s_i$ in the former case may be identified with the reflections across the hyperplane walls bordering the fundamental alcove %or unitary chamber
 in the latter case. Since all $s_i$ play a symmetrical role, we can cyclically relabel the simple reflections $s_i$ so that $s_0$ is the reflection across the unique wall that must be crossed to get out of the fundamental alcove while staying in $\Po_e(\lambda)$.

	\begin{eg} The conditions (1) and (2) on abaci in Lemma \ref{lem:hd=ell} in a picture: applying $t$ will increase the homological degree by $1$ if and only if no bead lies in the red regions of the runners.
$$
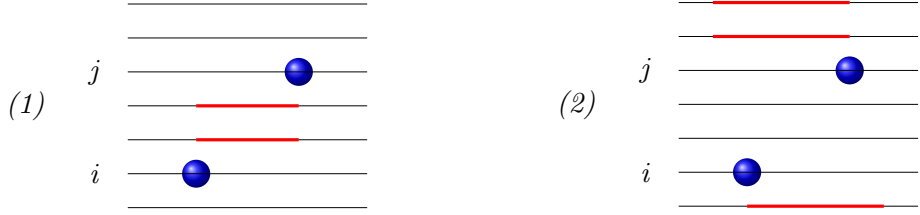
\begin{tikzpicture}[scale=.45,every text node part/.style={align=center}]
			\node at (-3,0){(1)};
			\node at (-1,1){j};
			\node at (-1,-2){i};
			\node at (5,1)[circle,shading=ball,minimum width=3pt] (ball){};
			\node at (2,-2)[circle,shading=ball,minimum width=3pt] (ball){};	
			\draw[very thick,red] plot [smooth] coordinates	{(2,-1)(5,-1)};		
			\draw[very thick,red] plot [smooth] coordinates	{(2,0)(5,0)};	
			\draw[thin] plot [smooth] coordinates{(0,0)(2,0)};
			\draw[thin] plot [smooth] coordinates{(5,0)(7,0)};
			\draw[thin] plot [smooth] coordinates{(0,-1)(2,-1)};
			\draw[thin] plot [smooth] coordinates{(5,-1)(7,-1)};
			\draw[thin] plot [smooth] coordinates{(0,1)(7,1)};
			\draw[thin] plot [smooth] coordinates{(0,2)(7,2)};
			\draw[thin] plot [smooth] coordinates{(0,3)(7,3)};
			\draw[thin] plot [smooth] coordinates{(0,-3)(7,-3)};
			\draw[thin] plot [smooth] coordinates{(0,-2)(7,-2)};
			\end{tikzpicture}		
		\qquad\qquad\qquad 	
			\begin{tikzpicture}[scale=.45,every text node part/.style={align=center}]
			\node at (-3,0){(2)};
			\node at (-1,1){j};
			\node at (-1,-2){i};
			\node at (5,1)[circle,shading=ball,minimum width=3pt] (ball){};
			\node at (2,-2)[circle,shading=ball,minimum width=3pt] (ball){};
			\draw[very thick,red] plot [smooth] coordinates	{(1,2)(5,2)};
			\draw[very thick,red] plot [smooth] coordinates	{(1,3)(5,3)};		
			\draw[thin] plot [smooth] coordinates{(0,-2)(7,-2)};			
			\draw[thin] plot [smooth] coordinates{(0,-1)(7,-1)};
			\draw[thin] plot [smooth] coordinates{(0,0)(7,0)};
			\draw[thin] plot [smooth] coordinates{(0,1)(7,1)};
			\draw[thin] plot [smooth] coordinates{(0,2)(1,2)};
			\draw[thin] plot [smooth] coordinates{(5,2)(7,2)};
			\draw[thin] plot [smooth] coordinates{(0,3)(1,3)};
			\draw[thin] plot [smooth] coordinates{(5,3)(7,3)};
			\draw[very thick,red] plot [smooth] coordinates	{(2,-3)(6,-3)};	
			\draw[thin] plot [smooth] coordinates{(0,-3)(2,-3)};
			\draw[thin] plot [smooth] coordinates{(6,-3)(7,-3)};
			\end{tikzpicture}$$
		\end{eg}

\section{The Mullineux map on unitary simple modules }\label{mulling}

% The purpose of this section is to examine the effect of the Mullineux map on the  simples modules $D^\Bbbk_n(\la)$
% for $\la \in \mathcal{F}_n^1$.  We show that these   simples are preserved under the Mullineux map.  Moreover, we construct an explicit Mullineux isomorphism in terms of the bases and representing matrices of these simples given in  \cref{maintheomre,action}.  
 
 We first recall the  Mullineux involution on the quiver Hecke algebra of the symmetric group:  
 Let $\Mull$ denote the $R_n$-automorphism determined by 
\begin{equation}\label{mulldef}
\Mull: e(i_1,i_2,\dots,i_n) \mapsto e(-i_1,-i_2,\dots,-i_n)
\quad
\Mull: \psi_r \mapsto \psi_r
\quad
\Mull: y_k \mapsto y_k
\end{equation}
for $0\leq k \leq n$ and $0\leq r <n$ and $\underline{i}=(i_1,\dots,i_n)\in I^n$.    
Given a   simple module  $D^\Bbbk_n(\la)$, we let   $D^\Bbbk_n(\la)^\Mull$ denote the module with the same underlying vector space but with the multiplication defined by twisting the action with the involution $\Mull$.   The relationship between these two simples was the subject of a conjecture of Mullineux \cite{MR545202}.  The combinatorics of this relationship is fiendishly complicated in general and is only understood on the level of the {\em labels} of simple modules.
  The purpose of this section is to examine the effect of the Mullineux map on the  simple modules $D^\Bbbk_n(\la)$
 for $\la \in \mathcal{F}_n^1$.  We show that the set of these simples is preserved under the Mullineux involution.  Moreover, we construct an explicit Mullineux isomorphism in terms of the bases and representing matrices of these simples given in  \cref{maintheomre,action} --- we remark that this is the time %a % first class of simples for which such 
 the Mullineux isomorphism has been explicitly constructed (outside of the trivial semisimple case).  Furthermore we shall see that the Mullineux combinatorics drastically simplifies on unitary $e$-regular partitions $\lambda$ and that we can easily compute $\Mull(\lambda)$ on the $e$-abacus of $\lambda$.  
 We define the {\sf unitary branching   graph}, $\mathcal{Y}$, to have  vertices on level $k$ given  by 
 $$\mathcal{Y}_k = \{  \lambda \mid    \lambda \text{ is $e$-restricted and $\lambda \in  \mathcal{F}^1_k$   } \}$$ and edges connecting levels $k$ and $k+1$ given by 
$$\mathcal{E}_{k,k+1}=\{ \lambda \to \mu  \mid 
 \lambda\in \mathcal{Y}_k, \mu\in \mathcal{Y}_{k+1} \text{ and }\la=\mu-\square \text{ for $\square$ a good node}\}.$$

We first discuss how the abaci of an $e$-core $\rho$ and its transpose $\rho^t$ are obtained from one another when $\rho$ has at most $e-1$ columns. Recall the basics of abaci from Section \ref{subsect:abacus}. First, note that if $\rho$ has at most $h<e$ columns then $\rho^t$ has at most $e-h$ columns. 
Now, let $\cA_e^{h}(\rho)$ denote the $e$-abacus of $\rho$ written with $h$ beads, and perform the following procedure on it: (1) swap the empty spots and the beads in the first column (so that the resulting abacus has $e-h$ beads), then (2) flip this abacus upside down. The resulting abacus, $\cA_{e}^{e-h}(\rho^t)$,  is   the $e$-abacus of $\rho^t$ written with $e-h$ beads. %, which we will denote by $\cA_{e}^{e-h}(\rho^t)$ .

\begin{defn}\label{def:mullineux} 
%Let $\lambda \in \mathcal{F}_n(h)$, $h<e$, and let $\cA_e(\lambda)$ be the $e$-abacus of $\lambda$ with $h$ beads, let $\rho$ be the $e$-core of $\lambda$ and $\cA_e(\rho)$ its abacus with $h$ beads. 
Let $\lambda \in \mathcal{F}_n(h)$ for some $1\leq h<e$ and  let $\rho$ be the $e$-core of $\lambda$. Write $w(\lambda) = (e-h)q + r$ for some $q\geq 0$, $0\leq r <e-h$. Define $ \lambda_\Mull$ to be the partition with abacus obtained from  $\cA^{e-h}_e(\rho^t)$ by 
 moving the bottom $r$ beads $(q+1)$-units to the right, and the top $e - h - r$ beads $q$ units to the right.  \end{defn}

\begin{prop}\label{easypeasy} If $\lambda\in\mathcal{Y}_n$, then $\lambda_\Mull\in\mathcal{Y}_n$. Specifically: in the case $\lambda=\rho$, we have $\rho_\Mull=\rho^t$. Otherwise, we have $\lambda_\Mull\in\cF_n(e-h)$.
\end{prop} 
\begin{proof}
If $\lambda=\rho$ is an $e$-core, then $w(\lambda)=0$ and algorithm just stops after the step where we take the transpose of $\rho$. The abacus $\cA_e^{e-h}(\rho^t)$ clearly satisfies the criterion for unitarity (Definition \ref{unitarypartitions}) since all of its beads are concentrated in the first column. If $w(\lambda)>0$, so $\lambda$ is not an $e$-core, we must move the bottom-most bead of $\cA_e^{e-h}(\rho^t)$ at least one unit to the right to obtain $\cA_e^{e-h}(\lambda_\Mull)$. This guarantees that $\cA_e^{e-h}(\lambda_\Mull)$ does not start with a bead, and since $\cA_e^{e-h}(\lambda_\Mull)$ has $e-h$ beads, we conclude that $\lambda_{\Mull}$ has precisely $e - h$ columns. Finally, by construction, $\lambda_\Mull$ satisfies the conditions of Definition \ref{unitarypartitions}.
\end{proof}

\begin{eg}
Take $e=5$, $h=2$, $\lambda=(2^{28},1^3)$, so $w(\lambda)=11$. We obtain $\lambda_\Mull=(3^{19},1^2)$ as follows:

		
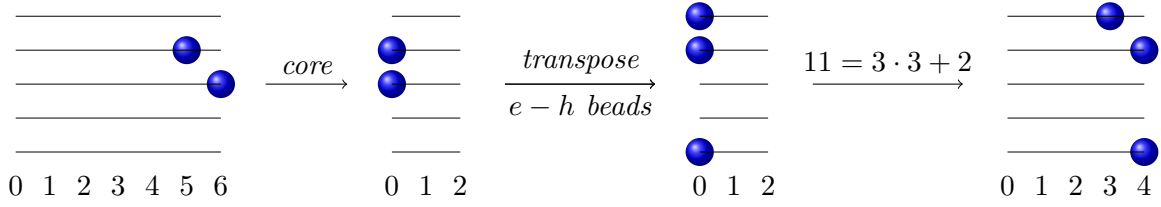
\begin{center}
			\begin{tikzpicture}[scale=.45,every text node part/.style={align=center}]
			\node at (-3,1)[circle,shading=ball,minimum width=3pt] (ball){};
			\node at (-2,0)[circle,shading=ball,minimum width=3pt] (ball){};
			\draw[thin] plot [smooth] coordinates{(-8,-2)(-2,-2)};
			\draw[thin] plot [smooth] coordinates{(-8,-1)(-2,-1)};
			\draw[thin] plot [smooth] coordinates{(-8,0)(-2,0)};
			\draw[thin] plot [smooth] coordinates{(-8,1)(-2,1)};
			\draw[thin] plot [smooth] coordinates{(-8,2)(-2,2)};
			\node at (-8,-3){$0$};
			\node at (-7,-3){$1$};
			\node at (-6,-3){$2$};
			\node at (-5,-3){$3$};
			\node at (-4,-3){$4$};
			\node at (-3,-3){$5$};
			\node at (-2,-3){$6$};
			\node(A) at (-1,0){};
			\node(B) at (2,0){};
\draw[->] (A) -- node[above]{core}(B) ;	
			\node at (3,0)[circle,shading=ball,minimum width=3pt] (ball){};
			\node at (3,1)[circle,shading=ball,minimum width=3pt] (ball){};
			\draw[thin] plot [smooth] coordinates{(3,-2)(5,-2)};
			\draw[thin] plot [smooth] coordinates{(3,-1)(5,-1)};
			\draw[thin] plot [smooth] coordinates{(3,0)(5,0)};
			\draw[thin] plot [smooth] coordinates{(3,1)(5,1)};
			\draw[thin] plot [smooth] coordinates{(3,2)(5,2)};
			\node at (3,-3){$0$};
			\node at (4,-3){$1$};
			\node at (5,-3){$2$};
			\node(C) at (6,0){};
			\node(D) at (11,0){};
\draw[->] (C) -- node[above]{transpose}(D);
\draw[->] (C) -- node[below]{$e-h$ beads}(D);
\node at (12,2)[circle,shading=ball,minimum width=3pt] (ball){};
			\node at (12,1)[circle,shading=ball,minimum width=3pt] (ball){};
			\node at (12,-2)[circle,shading=ball,minimum width=3pt] (ball){};
			\draw[thin] plot [smooth] coordinates{(12,-2)(14,-2)};
			\draw[thin] plot [smooth] coordinates{(12,-1)(14,-1)};
			\draw[thin] plot [smooth] coordinates{(12,0)(14,0)};
			\draw[thin] plot [smooth] coordinates{(12,1)(14,1)};
			\draw[thin] plot [smooth] coordinates{(12,2)(14,2)};
			\node at (12,-3){$0$};
			\node at (13,-3){$1$};
			\node at (14,-3){$2$};
			\node(E) at (15,0){};
			\node(F) at (20,0){};
\draw[->] (E) -- node[above]{$11=3\cdot 3+2$}(F);
\node at (24,2)[circle,shading=ball,minimum width=3pt] (ball){};
			\node at (25,1)[circle,shading=ball,minimum width=3pt] (ball){};
			\node at (25,-2)[circle,shading=ball,minimum width=3pt] (ball){};
			\draw[thin] plot [smooth] coordinates{(21,-2)(25,-2)};
			\draw[thin] plot [smooth] coordinates{(21,-1)(25,-1)};
			\draw[thin] plot [smooth] coordinates{(21,0)(25,0)};
			\draw[thin] plot [smooth] coordinates{(21,1)(25,1)};
			\draw[thin] plot [smooth] coordinates{(21,2)(25,2)};
			\node at (21,-3){$0$};
			\node at (22,-3){$1$};
			\node at (23,-3){$2$};
			\node at (24,-3){$3$};
			\node at (25,-3){$4$};
			\end{tikzpicture}
		\end{center}
\end{eg}

\begin{thm}\label{mull}
The map $\Mull: \mathcal{Y}_k \to  \mathcal{Y}_k$ for $k\geq 0$  is  a well-defined graph involution.  Given 
 $$\sts=(\lambda^{(0)}\xrightarrow{ r_1 } \lambda^{(1)}\xrightarrow{  r_2 } \dots \xrightarrow{  r_n } \lambda^{(n)})
$$ we let $\sts_\Mull$ denote the path 
$$ \sts_\Mull=(\lambda^{(0)}_\Mull\xrightarrow{ -r_1 } \lambda^{(1)}_\Mull\xrightarrow{  -r_2 } \dots \xrightarrow{  -r_n } \lambda^{(n)}_\Mull).
 $$
We have that 
 $
D^\Bbbk_n(\la_\Mull) \cong (D^\Bbbk_n(\la))^\Mull
 $
and that the isomorphism is determined by $: c_\sts \mapsto c_{\sts_M}$.  
\end{thm}
\begin{proof}
The Mullineux involution $\Mull$ is characterized as the unique involution on $e$-regular partitions mapping $\emptyset$ to $\emptyset$ and such that $\Mull(\tilde{f}_i(\lambda))=\tilde{f}_{-i}(\Mull(\lambda))$ \cite{MR1395065,MR1477629,MR1616083}. We want to identify $\lambda_\Mull$ with $\Mull(\lambda)$ for all vertices $\lambda$ of $\mathcal{Y}$. By construction we have that $\lambda_\Mull$ is also vertex of $\mathcal{Y}_k$ whenever $\lambda$ is, and that $(\lambda_\Mull)_\Mull=\lambda$. It is clear that $\emptyset_\Mull=\emptyset$. Thus if $i \in \{0, \dots, e-1\}$ is such that $\tilde{f}_{i}(\lambda)\in \mathcal{F}_n^1$, we need to show that 
	\[
\left(\tilde{f}_{i}(\lambda)\right)_\Mull = \tilde{f}_{-i}(\lambda_\Mull).
	\]

We remark that, if $\lambda \in \mathcal{F}^{1}_{n}$, then $\tilde{f}_{i}(\lambda)$ adds the leftmost addable box of content residue $i$, if any. In order to keep track of the action of $\tilde{f}_i$ on abaci, we follow the conventions of Remark \ref{conv:runnerlabels}, so we label the runners of an $e$-abacus with $h$ beads, at most one bead per runner, from bottom-to-top by $h-1, h-2, \dots, 1, 0, e-1, \dots, h$. This is done so that the labels of the runners correspond nicely to the contents of addable/removable boxes. Note that the labeling of runners changes in the process of constructing $\lambda_\Mull$, when $\cA_e(\rho)$ with $h$ beads is replaced by $\cA_e(\rho^t)$ with $e-h$ beads. The abacus $\cA_e(\lambda)$ has a bead (resp. empty space) on runner $i$ if and only if $\cA_e(\lambda_\Mull)$ has an empty space (resp. bead) on runner $-i - 1$. Finally, observe that if the top runner is labeled $m$ in these conventions, that $\tilde{f}_m$ increases the weight $w$ of a partition by at most $1$ but all $\tilde{f}_i$, $i\neq m$, do not increase the weight.

%Let us first consider the case where the box that $\tilde{f}_{i}$ adds to $\lambda$ is the addable box of highest content so that, in particular, $i = h$ and $\tilde{f}_{i}(\lambda) \in \mptn ^1 _n(h+1)$. Note that the assumption that $\tilde{f}_{i}(\lambda)$ is unitary implies that $\lambda$ is, in fact, a core. If $\tilde{f}_{i}(\lambda)$ is also a core

Set $\hat{\rho}$ to be the core of  $\tilde{f}_{i}(\lambda)$.
We consider two cases.
	
	{\it Case 1. $w(\tilde{f}_{i}(\lambda)) = w(\lambda)$.} So either $i \neq h$ or $i = h$ and $\lambda$ is a core. In the latter case, $\tilde{f}_{i}(\lambda)$ is also a core, and both $(\tilde{f}_{i}(\lambda))_{\Mull}$ and $\tilde{f}_{-i}(\lambda_{\Mull})$ coincide with the transpose of $\tilde{f}_{i}(\lambda)$. In the former case, the abaci $\cA_{e}(\rho)$ and $\cA_{e}(\hat{\rho})$ coincide on all runners except those labeled by $i$ and $i - 1$. Thus, $\cA_{e}(\rho^t)$ and $\cA_{e}(\hat{\rho}^t)$ only differ on runners $-i - 1$ and $-i$: $\cA_{e}(\hat{\rho}^t)$ has a bead on runner $-i-1$ and an empty space on runner $-i$, while the opposite is true for $\cA_{e}(\rho^t)$. Thus, $(\tilde{f}_{i}(\lambda))_\Mull$ is obtained from $\lambda_\Mull$ by sliding the bead on runner $-i$ up runner $-i-1$. But this is exactly how we obtain $\tilde{f}_{-i}(\lambda_\Mull)$ from $\lambda_\Mull$. We are done in this case. 
	
	{\it Case 2. $w(\tilde{f}_{i}(\lambda)) = w(\lambda) + 1$.} So $i = h$, and the abacus of $\tilde{f}_{i}(\lambda)$ is obtained from that of $\lambda$ by moving the bead on the top runner (labeled $h$) down to the bottom runner (labeled $h - 1$) and then one unit right. Just as in the first case, the abaci $\cA_{e}(\rho^t)$ and $\cA_{e}(\hat{\rho}^t)$ only differ on runners $-h$ and $-h-1$. Note that these are the top and bottom runners of the abacus, respectively. Write division with remainder $w = w(\lambda) = (e-h)q + r$, so that $\lambda_\Mull$ is obtained from $\cA_{e}(\rho^t)$ by moving the bottom $r$ beads $q + 1$ units to the right, and the remainder $e-h-r$ beads $q$ units to the right. We have a subdivision into two further cases. 
	
	{\it Case 2.1.} $r < e-h-1$. So $w + 1 = (e-h)q + (r+1)$ is division with remainder, and $(\tilde{f}_{i}(\lambda))_\Mull$ is obtained from  $\cA_{e}(\hat{\rho}^t)$ by moving the bottom $r + 1$ beads $q + 1$ units to the right, and the remaining beads $q$ units to the right. Note that the beads $2, \dots, r+1$ of $\cA_{e}(\hat{\rho}^t)$ coincide with the beads $1, \dots, r$ of $\cA_{e}(\rho^t)$. Thus, $(\tilde{f}_{i}(\lambda))_\Mull$ is obtained from $\lambda_\Mull$ by taking the bead in the top runner, moving it down to the bottom runner and sliding one unit to the right. This is precisely $\tilde{f}_{-i}(\lambda_\Mull)$. 
	
	{\it Case 2.2.} $r = e - h -1$. So $w + 1 = (e-h)(q+1)$. Here, $(\tilde{f}_{i}(\lambda))_\Mull$ is obtained from $\cA_{e}(\hat{\rho}^t)$ by moving \emph{all} beads $q+1$ units to the right, while $\lambda_\Mull$ is obtained from $\cA_{e}(\rho^t)$ by moving all beads $q+1$ units to the right, except the one in the top runner, that we only move $q$ units to the right. So we see that, again, $(f_{i}(\lambda))_\Mull$ is obtained from $\lambda_\Mull$ by taking the bead in the top runner, moving it down to the bottom runner and sliding one unit to the right. So $(\tilde{f}_{i}(\lambda))_\Mull = \tilde{f}_{-i}(\lambda_\Mull)$. 
	
This proves that the involution $\mathcal{Y} \to \mathcal{Y}$ given by $\lambda\mapsto\lambda_\Mull$ coincides with the Mullineux involution restricted to $\mathcal{Y}$. 
   Now, the bases of $D^\Bbbk_n(\la)$ for $\la \in \mathcal{F}^1_n$  are given by the paths in the unitary branching graph terminating at said vertices.   By \cref{action} we can match up these bases through the action of the idempotents under the twisting by the Mullineux map (see \cref{mulldef}).  The result follows.    \end{proof}

\begin{eg}
Let $e=7$.  
We have that  $\Mull(3^{10},2^4)=(4^8,1^3)$.  
We depict these partitions,  and the manner in which they can constructed via adding rim $7$-hooks in 
\cref{Mullfig}.  
Furthermore, we provide an example of   $\stt \in \Std_7(3^{10},2^4)$ and 
$\stt_\Mull \in \Std_7(4^8,1^3)$.  
Note that the map on the level of tableau preserves the rim hooks drawn in the two diagrams!!

\!\!\!\!\!\!\!\!\!
 \begin{figure}[ht!]
$$
\begin{tikzpicture}[scale=0.5]
\clip(-0.5,0.5)--(11,0.5)--(11,-7.5-6)--(-0.5,-7.5-6)--(-0.5,0.5);
\draw[thick](0,0)--(3,0)--(3,-1)--(1,-1)--(1,-5)--(0,-5)--(0,0);
\draw[thick,fill=gray!30]   (0,-5)--(0,-6)--(2,-6)--(2,-2)--(3,-2)--(3,-1)--(1,-1)--(1,-5)--(0,-5);
\draw[thick](0,-6)--(0,-7)--(3,-7)--(3,-2);
\draw[thick,fill=gray!30](0,0-7)--(3,0-7)--(3,-1-7)--(1,-1-7)--(1,-5-7)--(0,-5-7)--(0,0-7);
\draw[thick]  (0,0-7)--(0,-5-7)--(0,-6-7)--(2,-6-7)--(2,-2-7)--(3,-2-7)--(3,-1-7)--(3,0-7)--(0,0-7);
 \path(0.5,-0.5) coordinate (origin); 
 \foreach \i in {0,...,19}
  {
    \path (origin)++(-90:1*\i cm)  coordinate (a\i);
  \path (a\i)++(0:1cm)  coordinate (b\i);
        \path (b\i)++(0:1 cm)  coordinate (c\i);
            \path (c\i)++(0:1  cm)  coordinate (d\i);
      } 
\draw(a0) node {$1$};
\draw(a1) node {$4$};
\draw(a2) node {$5$}; 
\draw(a3) node {$6$}; 
\draw(a4) node {$7$}; 
\draw(b0) node {$2$};
\draw(c0) node {$3$};

 \draw(b1) node {$8$};
\draw(c1) node {$9$}; 
\draw(a5) node {$10$}; 
\draw(b2) node {$11$};
\draw(b3) node {$12$};
\draw(b4) node {$13$};
\draw(b5) node {$14$};

\draw(c2) node {$15$};
\draw(c3) node {$18$}; 
\draw(c4) node {$19$};
\draw(c5) node {$20$};
\draw(a6) node {$16$};
\draw(b6) node {$17$};
\draw(c6) node {$21$};

\draw(a7) node {$22$};
\draw(a8) node {$25$};
\draw(a9) node {$26$}; 
\draw(a10) node {$27$}; 
\draw(a11) node {$28$}; 
\draw(b7) node {$23$};
\draw(c7) node {$24$};

 \draw(b8) node {$29$};
\draw(c8) node {$30$}; 
\draw(a12) node {$31$}; 
\draw(b9) node {$32$};
\draw(b10) node {$33$};
\draw(b11) node {$34$};
\draw(b12) node {$35$};

\draw[thick](6,0)--(10,0)--(10,-1)--(7,-1)--(7,-4)--(6,-4)--(6,0); 
\draw[thick,fill=gray!30](6,-4)--(6,-5)--(8,-5)--(8,-2)--(10,-2)--(10,-1)--(7,-1)--(7,-4)--(6,-4); 
\draw[thick](6,-5)--(6,-6)--(9,-6)--(9,-3)--(10,-3)--(10,-2)--(8,-2)--(8,-5)--(6,-5);
\draw[thick,fill=gray!30] (6,-6)--(6,-7)--(10,-7)--(10,-3)--(9,-3)--(9,-6)--(6,-6);
\draw[thick](6,0-7)--(10,0-7)--(10,-1-7)--(7,-1-7)--(7,-4-7)--(6,-4-7)--(6,0-7);

 \path(6.5,-0.5) coordinate (origin); 
 \foreach \i in {0,...,19}
  {
    \path (origin)++(-90:1*\i cm)  coordinate (a\i);
    \path (a\i)++(0:1cm)  coordinate (b\i);
        \path (b\i)++(0:1 cm)  coordinate (c\i);
            \path (c\i)++(0:1  cm)  coordinate (d\i);
      } 
\draw(a0) node {$1$};
\draw(a1) node {$2$};
\draw(a2) node {$3$}; 
\draw(a3) node {$7$}; 
\draw(b0) node {$4$};
\draw(c0) node {$5$};
\draw(d0) node {$6$};

\draw(b1) node {$8$};
\draw(b2) node {$9$}; 
\draw(b3) node {$10$}; 
\draw(c1) node {$11$};
\draw(d1) node {$12$};
\draw(a4) node {$13$};
\draw(b4) node {$14$};

\draw(c2) node {$15$};
\draw(c3) node {$16$}; 
\draw(c4) node {$17$}; 
\draw(d2) node {$18$};
\draw(a5) node {$19$};
\draw(b5) node {$20$};
\draw(c5) node {$21$};

\draw(d3) node {$22$};
\draw(d4) node {$23$}; 
\draw(d5) node {$24$}; 
\draw(a6) node {$25$};
\draw(b6) node {$26$};
\draw(c6) node {$27$};
\draw(d6) node {$28$};

\draw(a7) node {$29$};
\draw(a8) node {$30$};
\draw(a9) node {$31$}; 
\draw(a10) node {$32$}; 
\draw(b7) node {$33$};
\draw(c7) node {$34$};
\draw(d7) node {$35$};

\end{tikzpicture}
$$

\!\!\!\!\!\!\! 
\caption{
 A pair of tableaux $\stt \in \Std_7(3^{10},2^4)$ and 
$\stt_\Mull \in \Std_7(4^8,1^3)$ indexing basis elements  swapped under the isomorphism $D^\Bbbk_{35}(3^{10},2^4)^\Mull  \cong D^\Bbbk_{35}(4^8,1^3)$.  
}
\label{Mullfig}
\end{figure}
 \end{eg}

 \!\!\!\!\!\!
\section{The rational Cherednik algebra of the symmetric group over $\mathbb{C}$}

For the remainder of the paper, we  restrict our attention to the field $\mathbb{C}$ and rational Cherednik algebras of type $G(1,1,n)$.  	 
%We wish to apply the results of the previous sections to both the representation theory of the Cherednik algebra itself (which we now define) 
%and 
%We now recall the original definition of the Cherednik algebra.  
%We shall use the results proven in the diagrammatic setting to deduce results concerning 
%
 Let   $\mathfrak{S}_n$ be the symmetric group on $n$ elements. The group $\mathfrak{S}_n  $ acts on the algebra of polynomials in $2n$ non-commuting variables $\C\langle x_{1}, \dots, x_{n}, y_{1}, \dots y_{n}\rangle$. 
 Fix a number $c \in \C$. The \emph{rational Cherednik algebra} $H_{c}(\mathfrak{S}_{n})$ is the quotient of the semidirect product algebra $\C\langle x_{1}, \dots, x_{n}, y_{1}, \dots, y_{n}\rangle\rtimes \mathfrak{S}_{n}$ by the relations 
	\[
	{[x_{i}, x_{j}] = 0, \hspace{1cm} [y_{i}, y_{j}] = 0, \hspace{1cm} [y_{i}, x_{j}] = c(ij)\; (i \neq j), \hspace{1cm} [y_{i}, x_{i}] = 1 - c \sum_{j \neq i}(ij)}
	\]	
	\noindent where $(ij)$ denotes the transposition in $\mathfrak{S}_{n}$ that switches $i$ and $j$, see \cite{EtingofGinzburg}.
 % The \emph{rational Cherednik algebra} $H_{1/e}(\mathfrak{S}_n)$ is the quotient of the semidirect product algebra $$\C\langle x_{1}, \dots, x_{n}, y_{1}, \dots, y_{n}\rangle\rtimes \mathfrak{S}_n$$ by commutation relations in the $x$'s and $y$'s that are similar to those of the Weyl algebra but involve an error term in $\C \mathfrak{S}_n$ determined by  $e\in \mathbb{N}$ \cite{EtingofGinzburg}.   In particular, these relations tell us that the $x$'s commute with each other and so do the $y$'s. 
$H_c$ has three distinguished subalgebras: $\C[\underline{y}] := \C[y_{1}, \dots, y_{n}]$, $\C[\underline{x}] := \C[x_{1}, \dots, x_{n}]$, and the group algebra $\C\mathfrak{S}_n$. The \emph{PBW theorem}  \cite[Theorem 1.3]{EtingofGinzburg} asserts that multiplication gives a vector space isomorphism
	\[
	\C[\underline{x}] \otimes \C\mathfrak{S}_n \otimes \C[\underline{y}] \buildrel \cong \over \rightarrow H_c
	\]
	\noindent called the \emph{triangular decomposition} of $H_c$, by analogy with the triangular decomposition of the universal enveloping algebra of a semisimple Lie algebra.

	We define the category $\cO_c(\mathfrak{S}_n)$ to be the full subcategory consisting of all finitely generated $H_c$-modules on which $y_{1}, \dots, y_{n}$ act locally nilpotently. Category $\cO_{c}$ is not always very interesting. By \cite{DunklDeJeuOpdam}, see also \cite[Section 3.9]{BezrukavnikovEtingof}, $\cO_{c}$ is semisimple (and equivalent to the category of representations of $\mathfrak{S}_{n}$) unless $c = r/e$, with $\gcd(r;e) = 1$ and $1 < e \leq n$. Equivalences of categories reduce the study of $\cO_{r/e}(\mathfrak{S}_n)$ to  $\cO_{1/e}(\mathfrak{S}_{n})$, for $1 < e \leq n$ \cite{rouquier}. For the rest of the paper we work with $\cO_{1/e}(\mathfrak{S}_n)$. It will be convenient to set
	\[
	\cO_{1/e} := \bigoplus_{n \geq 0} \cO_{1/e}(\mathfrak{S}_{n}).
	\]
 	The category  $\cO_{1/e}(\mathfrak{S}_n)$ is Morita equivalent (over $\mathbb{C}$) to $A(n,\kappa)$ for any  value of $\kappa \in I^\ell$ \cite[Theorem A]{Webster}.  
	Thus,  $\cO_{1/e}(\mathfrak{S}_n)$ is a highest weight category with respect to the dominance 
	 ordering  $\rhd$ on $\mptn 1 n$.   The standard modules are constructed as follows.  
 Extend the action of $\mathfrak{S}_n$ on $S_n(\lambda)$ to an action of $\C[\underline{y}]\rtimes \mathfrak{S}_n$ by letting $y_{1}, \dots, y_{n}$ act by $0$. The algebra $\mathbb{C}[\underline{y}] \rtimes \mathfrak{S}_n$ is a subalgebra of $H_{1/e}$ and we define
	\[
	\Delta(\lambda) := {\rm Ind}_{\C[\underline{y}]\rtimes \mathfrak{S}_n}^{H_{1/e}} S_n( \lambda) := 
	H_{1/e} \otimes_{\C[\underline{y}] \rtimes \mathfrak{S}_n} S_{n}(\lambda) = \C[\underline{x}] \otimes S_n(\lambda)
	\]
  where the last equality is \emph{only as $\C[\underline{x}]$-modules} and follows from the triangular decomposition.  
	 We let $L(\lambda)$ denote the unique irreducible quotient of $\Delta(\lambda)$.

\newcommand{\supp}{{\rm supp}}
\newcommand{\Spec}{{\rm Spec}}

%\subsubsection{Supports}
 Any module $M \in \cO_{1/e}(\mathfrak{S}_n)$ is finitely generated over the algebra $\C[\underline{x}]$ and, as such, it has a well-defined support $\supp(M) \subseteq \C^{n} = \Spec(\C[\underline{x}])$. We now explain a way to compute the supports of simple modules in $\cO_{1/e}(\mathfrak{S}_n)$ that was obtained in \cite{MR3658722}. To do this, for any $i = 0, \dots, \lfloor n/e\rfloor$, denote by $X_{i}$ the variety
	\[
	X_{i} := \mathfrak{S}_n\{(z_{1}, \dots, z_{n}) \in \C^{n} : z_{1} = z_{2} = \cdots z_{e}, z_{e+1} = \cdots = z_{2e}, \cdots, z_{(i-1)e + 1} = \cdots = z_{ie}\}
	\]
	By its definition, $X_{i}$ is a $\mathfrak{S}_n$-stable subvariety of $\C^{n}$. Note that $X_{0} = \C^{n}$, and these subvarieties form a chain $X_{0} \supsetneq X_{1} \supsetneq \cdots \supsetneq X_{\lfloor n/e\rfloor}$. 
		Now recall that a partition $\lambda$ is said to be \emph{$e$-restricted} if $\lambda_{i} - \lambda_{i+1} < e$ for every $i \geq 0$, that is, if no two consecutive parts of $\lambda$ differ by more than $e-1$ parts. By the division algorithm, for any partition $\lambda$ there exist unique partitions $\mu, \nu$ such that $\lambda = e\mu + \nu$ and $\nu$ is $e$-restricted. Then, according to \cite[Theorem 1.6]{MR3658722},
	\[
	\supp(L_{1/e}(\lambda)) = X_{|\mu|}
	\]
	\noindent So, for example, $L_{1/e}(\lambda)$ has full support if and only if $\lambda$ is $e$-restricted. On the other hand, if $e$ divides $n$, then $L_{1/e}(\lambda)$ has minimal support if and only if $\lambda = e\mu$, where $\mu$ is a partition of $n/e$.

	The categories $\cO_{1/e}(\mathfrak{S}_{n})$ come equipped with induction and restriction functors $$\Res^{n}_{n-1}: \cO_{1/e}(\mathfrak{S}_{n}) \rightleftarrows \cO_{1/e}(\mathfrak{S}_{n-1}): \Ind^{n}_{n-1}$$ that were constructed by Bezrukavnikov and Etingof in \cite{BezrukavnikovEtingof}. Their definition is quite technical and will not be needed. 
 In fact, 
		Bezrukavnikov and Etingof constructed restriction functors for any parabolic subgroup of $\mathfrak{S}_{n}$, \cite{BezrukavnikovEtingof}.  It follows from their construction that $M$ has full support if and only if it is not killed by restriction to any parabolic subgroup. We will use this property below without further mention.

 \newcommand{\opp}{{\rm opp}}
 \newcommand{\ue}{h}
\newcommand{\uq}{\underline{q}}
\newcommand{\um}{\underline{m}}
\renewcommand{\RR}{\mathsf{R}}
\newcommand{\DD}{\mathsf{D}}

\subsection{Unitary modules}\label{unitarysec}

For   $\lambda \in \mptn 1n$, fix a positive-definite, $\mathfrak{S}_{n}$-invariant Hermitian form on the irreducible representation $S_{n}(\lambda)$. A standard argument shows that this form can be extended to a Hermitian form $(\cdot, \cdot)$ on the standard module $\Delta_{1/e}(\lambda)$, which is $H_{c}$-invariant in   that $(y_{i}v, v') = (v, x_{i}v')$ for every $v, v' \in \Delta_{1/e}(\lambda)$ and $i = 1,\dots, n$. 
Moreover, the simple module $L_{1/e}(\lambda)$ is the quotient of $\Delta_{1/e}(\lambda)$ by the radical of this form. In particular, $L_{1/e}(\lambda)$ is equipped with a $H_{c}$-invariant, non-degenerate Hermitian form. We say that $L_{1/e}(\lambda)$ is \emph{unitary} if this form is positive-definite. 
%Recall that $\mathcal{F}^1_{n}$ denotes the set of partitions $\lambda \in \mptn1n$ such that $\ell(\lambda) = 0$.  

	  	\begin{thm} \cite{MR2534594}
%Suppose that $\lambda \in \mptn 1n$.
 The Hermitian form on $L(\lambda)$ is positive-definite if and only if 
  $\la $ is an $e$-unitary partition.  
Thus $L(\lambda)$ is unitary   if and only if   	$\lambda \in \mathcal{F}^1_n$ \end{thm} 

Applying the {\sf KZ} functor to these simples, we obtain the complete set of   simple unitary modules for the Hecke algebra.   
We emphasise that the  simples  labelled by $\la=(e^k)$ for some $k\geq 0$ do not survive under the   {\sf KZ} functor and so there are fewer unitary simples for the Hecke algebras.  

\begin{thm}[{\cite[Corollary 4.5]{unitaryhecke}}]
  The simple $R_n$-module   $D_n^\CC(\lambda)$ is    unitary   if and only if   $\lambda$ is $e$-restricted and	$\lambda \in \mathcal{F}^1_n$
\end{thm}

\begin{rmk}
We would like to say some words on higher levels. Associated to the group $G(\ell, 1, n)$ there is   a rational Cherednik algebra $H_{c}(G (\ell, 1, n))$, where $c = (c_{0}, c_{1}, \dots, c_{\ell - 1})$ is now a collection of $\ell$ complex numbers. 
The definition of a unitary module goes through unchanged.  
%Using results from \cite{griffeth}, we have checked that 
We let 
 $\ell = 2$ and taking the charge $c_0=1/e $ and $c_{1} = 0$ so that we are, essentially, working with rational Cherednik algebras associated to the Weyl group of type $D$.
 We have checked (using Griffeth's classification of unitary modules \cite{griffeth}) that if $\lambda \in \mathcal{F}^2_n$ then $L(\lambda)$ is indeed a unitary module.

\end{rmk}

	\subsection{Changing quantum characteristics }\label{Cahracteristicchange} 
	Having constructed a BGG resolution for any unitary module with $h<e$ columns, we proceed to relate these complexes to each other for various $e$, and to construct the complex in the special case $h=e$ for the unitary module $L(e^k)$. 
 % Having constructed a BGG resolution for any unitary module,
% we proceed to relate these complexes to each other for various quantum characteristics over $\mathbb{C}$.   
% In some sense, we shall see that all these complexes 
%  are, in some sense, the same and arise from the $e=h$ case.   
%
	As observed in   \cref{abacus}, the $e$-abacus of any unitary module which is not of the form $L(e^{k})$ will contain empty runners; removing the empty runners produces the $h$-abacus of a partition of the form $(h^{k})$, with $h < e$ and $k$ equal to the weight of the block containing $\lambda$. So we may try using the runner removal Morita equivalences of Chuang-Miyachi which upgrade the combinatorial operation ``removing runners" to an equivalence of highest weight categories \cite{ChuangMiyachi}.
  Given an $e$-core partition $\rho$ and $k \in \mathbb{N}$, let $n := |\rho| + ek$ and set 
  $$\Lambda({\rho,k}):= \{ \lambda \mid \la \in \mptn  1 n
, e{\rm-core}(\lambda)=\rho, w(\lambda)=k\}\subseteq \mptn 1n,
  $$
  $$
  \Lambda_{h}^{+}(\rho, k) := \{\lambda \in \Lambda(\rho, k) \mid \la \in \mptn 1 n (h)\} \subseteq \Lambda(\rho, k)\}
  $$
  $$
  \Lambda_{h}^{-}(\rho, k) := \{\lambda \mid \lambda^{T} \in \mptn 1 n (h), e{\rm-core}(\la) = \rho^{T}, w(\lambda) = k\} \subseteq \Lambda(\rho^{T}, k).
  $$
  Notice that the transpose map gives a bijection between the sets $\Lambda^{-}_{h}(\rho, k)$ and $\Lambda^{+}_{h}(\rho, k)$; 
under this map the partial ordering on the sets is reversed.  
 Let 
$\mathcal{O}_{1/e}({\rho,k})$
% and   $\mathcal{O}^-_h({\rho,k})	$ 
denote the block of category $\mathcal{O}_{1/e}$ corresponding to $\Lambda(\rho, k)$. Note that the set $\Lambda^{-}_{h}(\rho, k)$ is co-saturated in $\Lambda(\rho^{T}, k)$ so we can consider the quotient category of $\cO_{1/e}(\rho^{T}, k)$ by the Serre subcategory spanned by simples whose label does \emph{not} belong to $\Lambda^{-}_{h}(\rho, k)$. We denote this quotient by $\cO_{1/e, h}^{-}(\rho, k)$ This is a highest weight category, with standard objects $\Delta^{-}_{h}(\nu) := \pi(\Delta_{1/e}(\nu))$, where $\nu \in \Lambda^{-}_{h}(\rho, k)$ and $\pi: \cO_{1/e}(\rho^{T}, k) \rightarrow \cO_{1/e,h}^{-}(\rho, k)$ is the quotient functor. We remark that $\pi$ admits a left adjoint $\pi^{!}: \cO_{1/e, h}^{-}(\rho, k) \rightarrow \cO_{1/e}(\rho^{T}, k)$, and $\pi^{!}(\Delta^{-}(\nu)) = \Delta(\nu)$ for $\nu \in \Lambda^{-}_{h}(\rho, k)$. %sub/quotient category corresponding  to $\Lambda^+_h({\rho,k})$ and $\Lambda^-_h({\rho,k})$ respectively.  
%We let $ \Delta_h^+(\nu)$, $\nu\in \Lambda^+_h({\rho,k})$  and $ \Delta_h^-(\nu)$, $\nu \in \Lambda^-_h({\rho,k})$ denote the standard objects in these categories.  

% 	For $\nu \in \Lambda^-_h({\rho,k})$, consider the 
%$h$-abacus of $\nu$ -- but now read off its sequence of spaces and beads according to the dual convention which produces from this abacus the partition $\nu^t$, see Remark \ref{rmk:cmabaci}. We denote this abacus by
%  $\cA_{h}^{CM}(\nu^t)$. We re-emphasize: as diagrams, $\cA_{h}^{CM}(\nu^t)$ and $\cA_{h}(\nu)$ are identical. It is the partitions and objects they label, $\nu^t$ versus $\nu$, objects of a quotient category versus objects of the dual subcategory, which differ.
%	 Now, let ${\sf r}=(r_{0}, \dots, r_{\ue -1}, r_{\ue}) \in \Z_{\geq 0}^{\ue +1}$, and construct a partition $\nu^{\sf R}$ as follows. 
%	 In the abacus $\cA^t _{\ue}(\nu)$, insert $r_{i}$ empty runners between runners $i-1$ and $i$ (so $r_{0}$ and $r_{\ue}$ are the number of empty runners inserted at the top and bottom of the abacus). 
%	 This creates a new $e$-abacus, $\cA$, with   $e:= \ue + r_{0} + \cdots + r_{\ue}$ runners.
%	 We denote by $\nu^{\sf R}$ the unique partition such that $\cA = \cA^{t}_{e}(\nu^{\sf R})$. 
%We let  $\rho = \varnothing^{\sf R}$. We have a bijection 
%	 $$
%	{\sf R}:  \Lambda^-_{\varnothing, k} \to 	 \Lambda^-_{\rho, k}
%	 $$
%  given by ${\sf R}:\nu \mapsto \nu^{\sf R}$ and we let ${\sf R}^{-1}$ denote the inverse.  
%  

Given $\nu\in \Lambda^-_h({\rho,k})$ we set $ \cA^{CM}_{h}(\nu):=\cA_{h}(\nu^{T})$.  
%Let $\nu \in \Lambda^-_h({\rho,k})$ 
 	 Let ${\sf r}=(r_{0}, \dots, r_{\ue -1}, r_{\ue}) \in \Z_{\geq 0}^{\ue +1}$, and construct a partition $\nu^{+}$ as follows. 
	 In the abacus $\cA _{h}^{CM}(\nu)$, insert $r_{i}$ empty runners between runners $i-1$ and $i$ (so $r_{0}$ and $r_{\ue}$ are the number of empty runners inserted at the top and bottom of the abacus, respectively). 
	 This creates a new $e$-abacus, $\cA$, with   $e:= \ue + r_{0} + \cdots + r_{\ue}$ runners.
	 We denote by $\nu^{+}$ the unique partition such that $\cA = \cA^{CM}_{e}(\nu^{+})$. 
We let  $\rho = \varnothing^{+}$. We have a bijection 
	 $$
	{\sf R}:  \Lambda^-(\varnothing, k) \to 	 \Lambda^{-}(\rho, k)
	 $$
  given by ${\sf R}:\nu \mapsto \nu^{+}$ and we let ${\sf R}^{-1}: \nu \mapsto \nu^{-}$ denote the inverse.  
 We are now able to recall the main result of Chuang--Miyachi.

 	\begin{thm}[\cite{ChuangMiyachi}]\label{thm:rrme}
		The categories $\mathcal{O}^-_{1/h,h}({\varnothing,k})$ and $\mathcal{O}^-_{1/e,h}({\rho,k})	$ are equivalent as highest weight categories. Moreover, the equivalence $$\underline{\RR}: 
		\mathcal{O}^-_{1/h,h}({\varnothing,k})
		\rightarrow \mathcal{O}^-_{1/e,h}({\rho,k})$$ sends the standard module $ \Delta_h^-(\nu)$ to the standard module $\Delta_h^-(\nu^{\sf R})$.
	\end{thm}

  Note, however, that we cannot apply the above theorem directly since we are interested in the subcategories  $\cO_{{1/h},h}^+(\varnothing, k)$ and $\cO_{{1/e},h}^+(\rho^T,k)$ rather than the quotient categories $\mathcal{O}^-_{1/h,h}({\varnothing,k})$ and $\mathcal{O}^-_{1/e,h}({\rho,k})$, where $\cO_{1/h, h}^{+}(\varnothing, k)$ denotes the Serre subcategory spanned by the simples whose label belongs to $\Lambda^{+}(\varnothing, k)$, and similarly for $\cO_{1/e, h}^{+}({\rho, k})$. %the abacus $\cA_{h}^{CM}(\nu)$ is different from $\cA_{h}(\nu)$. %--no not really. the problem is with categories not abaci. the combinatorics matches. its meaning does not, we have subs they have quotients.
  Let us fix this. %We now recall how the above theorem can be of use to use via  Ringel duality. 
Following \cite[Section 4]{GGOR}, we note that the rational Cherednik  algebra $H_{1/e} := H_{1/e}(\fS_{n})$ has finite global dimension and is isomorphic to its opposite algebra; an explicit isomorphism is given by $w \mapsto w^{-1},\; x \mapsto x,\; y \mapsto -y$. In particular, the functor $\text{RHom}_{H_{1/e}}(\bullet, H_{1/e})$ gives an equivalence $D^{b}(H_{1/e}\text{-mod}) \rightarrow D^{b}(H_{1/e}\text{-mod}^{\text{opp}})$. Let us denote by $\DD$ the functor $\text{RHom}_{H_{1/e}}(\bullet, H_{1/e})[n]$. The following theorem summarizes various results of \cite[Section 4.3.2]{GGOR}. We denote by $D^{b}(\cO_{1/e}(\fS_{n}))$ the subcategory of $D^{b}(H_{1/e}\text{-mod})$ consisting of complexes with homology in $\cO_{1/e}$, and by $\cO_{1/e}^{\Delta}$ the category of objects in $\cO_{1/e}$ that admit a $\Delta$-filtration.
	
	\begin{thm}\label{thm:GGOR}
		The functor $\DD$ induces a derived equivalence $\DD: D^{b}(\cO_{1/e}(\fS_{n})) \rightarrow D^{b}(\cO_{1/e}(\fS_{n})^{\opp})$ as well as an equivalence of exact categories $\DD: \cO_{1/e}(\fS_{n})^{\Delta} \rightarrow (\cO_{1/e}(\fS_{n})^{\Delta})^{\opp}$. For a partition $\lambda \vdash n$, $\DD(\Delta(\lambda)) = \Delta(\lambda^{T})$ (where both sides of the equation are interpreted as complexes concentrated in degree $0$).
	\end{thm}
	
	By abuse of notation, we will write $\DD: D^{b}(\cO_{1/e}) \to D^{b}(\cO_{1/e}^{\opp})$ for $\bigoplus_{n \geq 0}\text{RHom}_{H_{1/e}(\fS_{n})}(\bullet, H_{1/e}(\fS_{n}))[n]$.
 	Let us   mention a property of $\DD$ that will be important later. The following is an immediate consequence of \cite[Lemma 2.5]{Losev2017}  and the definition of a perverse equivalence  \cite[Section 1.4]{Losev2017}. 
	
	\begin{lem}
		For every $n \geq 0$, the functor $\DD$ induces a (contravariant!) abelian autoequivalence of the category of \emph{minimally supported} modules in category $\cO_{1/e}(\fS_{n})$.
	\end{lem} 
	
%	
%In what follows, we shall abuse notation by writing 	
%  $$\DD: \mathcal{O}^+_h({\rho,k}) \to   
% \mathcal{O}^-_h({\rho,k})
  
Let $\la \in \mathcal{F}^1_n\subseteq  \mptn 1 n(h)$ be such that $(h^{k})^{+} = \lambda$, where $k$ is the $e$-weight of $\lambda$ and $h$ the number of nonempty runners in $\cA_{e}(\la)$. 
 	  Define the functor $\tilde{\RR}^{-}$ via the following composition
 	   	  $$\tilde{\RR}^-:=\DD\pi^!\underline{\RR}^{-1}\pi\DD:D^b(\cO_{1/e}(\rho, k))\rightarrow D^b(\cO_{1/h}(\varnothing, k))).$$
 	  
 	   Each functor in the composition defining $\tilde{\RR}^-$ takes Vermas to Vermas, and is either an equivalence of $\Delta$-filtered categories or exact on $\Delta$-filtered categories while being an isomorphism on spaces of homomorphisms between Vermas. It follows that for $\mu\in\Po_e(\lambda)$, $\tilde{\RR}^-\Delta(\mu)=\Delta(\mu^-)$, and that $\tilde{\RR}^-$ takes a complex to a complex and sends nonzero maps to nonzero maps (however, we cannot conclude from this that $\tilde{\RR}^-$ takes a resolution to a resolution). Define $C_\bullet(h^k)=\tilde{\RR}^-(C_{\bullet}(\lambda))$. By construction, this is a complex whose $\ell$-th term is given by
	\[
	C_\ell(h^k)=\bigoplus_{
 \begin{subarray}c \mu \in \Po_e(\lambda)\\
\hd(\mu)= \ell
\end{subarray}
}
	\Delta(\mu^-)=
	\bigoplus_{
%	\substack{\tau\vdash hk \\ h\text{-core}(\tau)=\emptyset \\ \tau\leq (h^k) \\ \hd\tau=i}}
\begin{subarray}c \tau \in \Po_h(h^k)\\
\hd(\tau)= \ell
\end{subarray}}
	\Delta(\tau)
	\]
	and which has a map $\Delta(\tau)\rightarrow\Delta(\tau')$ whenever $\hd(\tau)=\ell$, $\hd(\tau')=\ell-1$, and $\cA_h(\tau)=t\cA_h(\tau')$ for some transposition $t\in\widetilde{\mathfrak{S}}_h$. $C_\bullet(h^k)$ is a complex that looks identical to $C_\bullet(\lambda)$ but with the partitions $\mu$ relabeled by $\mu^-$, and in particular $L(h^k)$ is the head of $C_0(h^k)=\Delta(h^k)$.

	%%%%%%%%%%%%%%%%%%%%%%%%%%%%%%%%%%%%%
	%%%%%%%%%%%%%%%%%%%%%%%%%%%%%%%%%%%%%
	%%%%%%                         %%%%%%
	%%%%%%     %%%      %%%        %%%%%%
	%%%%%%      %%  %   %%         %%%%%%
	%%%%%%          %              %%%%%%
	%%%%%%         %%              %%%%%%   
	%%%%%%       %      %          %%%%%%
	%%%%%%        %%%%%%           %%%%%%
	%%%%%%                         %%%%%%
	%%%%%%%%%%%%%%%%%%%%%%%%%%%%%%%%%%%%%

	The following theorem answers \cite[Conjecture 4.5]{conjecture} in the affirmative.
	
	\begin{thm}\label{thm:exact} If $L(\lambda)\in\cO_{1/e}(\mathfrak{S}_n)$ is unitary then $L(\lambda)$ has a BGG resolution $C_\bullet(\lambda)$ whose $\ell$th term is given by $$C_\ell(\lambda)=\bigoplus_{\substack{\mu\in\Po_e(\lambda) \\ \hd(\mu)=\ell}}\Delta(\mu)$$
	\end{thm}
	\begin{proof}
		We have already shown the conjecture holds if $h<e$. 		
		Let $n=ke$ for some $k>0$ and take $(e^k)$, the unique unitary partition of $n$ with $e$ columns. Choose any $e'>e$ and any unitary partition $\lambda\in\cP_{e'}^u$ with $e$ columns and $e'$-weight $k$. Let $C_\bullet(\lambda)$ be the BGG resolution of $L(\lambda)$ and apply $\tilde{\RR}^-$ to it.  By the remarks above, $\tilde{\RR}^-(C_\bullet(\lambda))=C_\bullet(e^k)$ is the desired complex and $L(e^k)$ is the head of $C_0(e^k)$. We need to show that $C_\bullet(e^k)$ is exact except in degree $0$, where $H_0(C_\bullet)=L(e^k)$.
		
		As in the proof of the $h<e$ case, if $\lambda\in\Po_e(e^k)\setminus\{(e^k)\}$, then $\lambda$ is $e$-restricted. Thus $E_i(L(\lambda))\neq 0$ for some $i\in\Z/e\Z$, so if $L(\lambda)$ is a composition factor of a homology group $H_j(C_\bullet)$ then $E_i(C_\bullet)$ will fail to be exact. Similarly, it holds (by basic properties of highest weight categories) that $L(e^k)$ occurs exactly once in the composition series of all the $C_j$, when $j=0$.
		
		Next, $E_i(L(e^k))=0$ since $e^k$ has a single removable box and it is never a good removable box. Thus, it suffices to check that  $E_i(C_\bullet)$ is exact for each $i\in\Z/e\Z$. This is identical to the argument used in the $h<e$ cases for those $i$ such that $E_i(L(\lambda))=0$.
	\end{proof}
	
	We also make the observation that resolutions of unitary modules are, in a manner of speaking, independent of $e$. Let $h$ be the number of columns of $\lambda$ and let $k$ be the $e$-weight of $\lambda$.
	
	\begin{cor} Let $\la\in \mathcal{F}^1_n $.
	 The shape of the BGG complex $C_\bullet(\lambda)$ depends only on $h,k\in\N$. 
	\end{cor}
	\begin{proof}
		$\tilde{\RR}^-$ identifies $C_\bullet(\lambda)$ with $C_\bullet(h^k)$, thus sends a resolution of $L(\lambda)$ to a resolution of $L(h^k)$.
	\end{proof}
	
	%%%%%%%%%%%%%%%%%%%%%%%%%%%%%%%%%%%%%
	%%%%%%%%%%%%%%%%%%%%%%%%%%%%%%%%%%%%%
	%%%%%%                         %%%%%%
	%%%%%%     %%%      %%%        %%%%%%
	%%%%%%      %%  %   %%         %%%%%%
	%%%%%%          %              %%%%%%
	%%%%%%         %%              %%%%%%   
	%%%%%%       %      %          %%%%%%
	%%%%%%        %%%%%%           %%%%%%
	%%%%%%                         %%%%%%
	%%%%%%%%%%%%%%%%%%%%%%%%%%%%%%%%%%%%% 
	
	\subsection{Ringel duality and   more BGG resolutions} We can also construct some new BGG resolutions as corollaries of Theorem \ref{thm:exact} via Ringel duality. 
These resolutions will also be used in the  study subspace arrangements in \cref{keideal}.    	
	The character of $L(ek)=L(\triv)\in\cO_{1/e}(\mathfrak{S}_{ek})$ is dual to the character of $L(e^k)$ in the sense that its character is obtained from that of $L(e^k)$ by taking the transpose of each partition labelling a Verma module \cite[Remark 5.1]{EGL}:
	$$L(ek)=\sum_{\substack{\mu\in\Po_e(e^k) \\ \hd(\mu)=\ell}}(-1)^\ell\Delta(\mu^T).$$
	This is every bit as much an alternating sum character formula as that of $L(e^k)$, so we may naturally ask whether its character formula also comes from a BGG resolution.
	
	Let $C_\bullet$ be the BGG resolution of $L(e^k)$. 
	We apply Ringel duality to  construct a complex, $\DD(C_\bullet)$, in the principal block $\mathcal{O}(\varnothing, k) \subset\cO_{1/e}(\mathfrak{S}_{ek})$.  The complex $\DD(C_\bullet)$ is obtained from  $C_\bullet$ by replacing $\Delta(\mu)$ with $\Delta(\mu^T)$ for all $\mu\in\Po_e(e^k)$ and reversing the direction of all the arrows (since $\DD$ is a contravariant functor which takes Vermas to Vermas). By \cite{EGL}, the alternating sum of the terms of $\DD (C_\bullet)$ in the Grothendieck group $[\mathcal{O}(\varnothing, k)]$ coincides with the character of $L(\mathrm{triv})=L(ek)$.
	
	\begin{cor}\label{cor:BGGtriv}
		$\DD(C_\bullet)$ is a BGG resolution of $L(\triv)=L(ek)$.
	\end{cor}
	
	\begin{proof}
		A resolution is quasi-isomorphic to the module it resolves, so in $D^b(\cO_{1/e}(\mathfrak{S}_{ek}))$, $L(e^k)$ is isomorphic to its resolution $C_\bullet$. Since the Ringel duality $\DD$ is a derived self-equivalence of $D^b(\cO_{1/e}(\mathfrak{S}_{ek}))$ \cite{GGOR}, this implies $\DD(C_\bullet)\simeq\DD(L(e^k))$ in $D^b(\cO_{1/e}(\mathfrak{S}_{ek}))$. We know that at the end of the complex we have: $\Delta(ek-1,1)\rightarrow\Delta(ek)\rightarrow 0$, and so $L(ek)= \mathrm{Head}(\Delta(ek))$ must occur in the homology of $\DD(C_\bullet)$. Therefore $L(ek)$ is a composition factor of $\DD(L(e^k))$.
		
		We claim that $\DD(L(e^k))=L(ek)$. This follows from \cite[Lemma 2.5]{Losev2017} which states that $\DD$ is a perverse equivalence with respect to the filtration by dimensions of support: in particular, $\DD$ is a self-equivalence of the semi-simple subcategory spanned by the minimal support modules $L(e\sigma)$.
		Since $\DD^2=Id$, it follows that $\DD$ must permute the minimal support simple modules $L(e\sigma)$, $\sigma\vdash k$. We have already seen that $\DD(L(e^k)=\DD(L(e(1^k)))$ contains $L(ek)=L(e(k))$ as a composition factor; it follows that $\DD(L(e^k))=L(ek)$.
		
		To conclude, $\DD(C_\bullet)$ is equivalent to $L(ek)$ in $D^b(\cO_{1/e}(\mathfrak{S}_{ek}))$, where $L(ek)$ is considered as a complex concentrated in degree $0$.  Hence $H_i(\DD(C_\bullet))=\delta_{i0}L(ek)$, as required.  
	\end{proof}
	
	Let $\pi$ denote the quotient functor which kills the subcategory generated by $\{L(\nu)\mid \nu\hbox{ has more than }e\hbox{ rows}\}$.  
	
	\begin{cor} $\pi \DD(C_\bullet)$ is a BGG resolution of $\underline{L}(\triv)=\underline{L}(ek)$ in the quotient category  $\pi(\cO_{1/e}(\mathfrak{S}_{ke}))  $. By adding an arbitrary configuration of $a\in \ZZ_{>0}$ empty runners to the abacus, $\RR\pi\DD(C_\bullet)$ is a BGG resolution of $\RR\underline{L}(\triv)$ in $\RR\pi(\cO_{1/e}(\mathfrak{S}_{ke})) $. % ($r$ being determined by the runners that were added).
	\end{cor}
	
	\begin{proof}
		The quotient functor $\pi$ is %given by an idempotent truncation and is therefore exact \cite[Appendix]{Donkin}; $\pi$
		exact, sends $\Delta(\mu)$ to the standard module $\underline{\Delta}(\mu)$, and sends $L(\mu)$ to the simple module $\underline{L}(\mu)$. The first claim then follows from Corollary \ref{cor:BGGtriv}, implying the second claim by Theorem \ref{thm:rrme}.
	\end{proof}

	\subsection{Computation of Lie algebra and Dirac cohomology}	
	BGG resolutions for classical and affine Lie algebras over $\C$ are closely related to the computation of Lie algebra cohomology \cite{Kostant,MR0414645,BoeHunziker}. 
	Recently, a version of Lie algebra cohomology (and homology) for  rational Cherednik  algebras over $\C$ was constructed in \cite{HuangWong}; $\hstar := \bigoplus \mathbb{C}x_{i}$ plays the role of the nilradical $\mathfrak{n}\subset\mathfrak{b}\subset\mathfrak{g}$, and the complex reflection group $W$ plays the role of the Cartan subalgebra.

	%	 The analogue of the classical result holds: if $L\in\cO_c(W)$ has a BGG resolution $C_\bullet$ in which $\Delta(\mu)$ appears at most once for a given $\mu\in{\rm Irr} W$, then the $i$'th homology group is given by $$H_i(\hstar, L)=\bigoplus\limits_{\Delta(\mu)\subset C_i}S_n(\mu)$$ \cite[Proposition 6.1]{HuangWong}.  It then follows that the BGG resolutions we have constructed in this paper for unitary modules $L(\lambda)\in\cO_{1/e}(S_n)$ yield that
	\begin{thm} Let $\lambda \in \mathcal{F}^1_n\subseteq \mptn 1n(h)$.
		We have that 
		$$H_i(\hstar, L(\lambda))=\bigoplus\limits_{\substack{\mu\in\Po_e(\lambda)\\\hd(\mu)=i}}S_n(\mu).$$  
	\end{thm}
	\begin{proof}
		This follows immediately from  our main theorem and \cite[Proposition 6.1]{HuangWong}.  
	\end{proof}
	
	Likewise, if $L(\lambda)\in\cO_c(G(\ell,1,n))$ where $c$ corresponds to the rank $e$ and charge $\underline{s}=(\kappa_1,\kappa_2,\dots,\kappa_\ell)\in\ZZ^\ell$ for the Fock space, and  $\lambda \in \mathcal{F}^\ell_n$, then $$H_i(\hstar,L(\lambda))=\bigoplus\limits_{\substack{\mu\trianglelefteq  \lambda\\ \ell(\mu)=i}}S_n(\mu).$$ 
	This also computes the Lie algebra cohomology $H^{i}(\hstar, L(\lambda))$. Indeed, by Poincar\'e duality (cf. \cite[Proposition 2.7]{HuangWong}), we get $$H^{i}(\hstar, L(\lambda)) = H_{n-i}(\hstar, L(\lambda)) \otimes \wedge^{n}\mathfrak{h},$$ where $n := \dim\mathfrak{h}$.
	A consequence of the computation of Lie algebra cohomology for unitary modules admitting a BGG resolution is that this immediately gives the computation of the Dirac cohomology $H_{D}(L(\lambda))$. This is defined as the usual Dirac cohomology, where the Dirac operator $D \in H_{1/e}(\mathfrak{S}_{n}) \otimes \fc$ has been constructed in \cite{Ciubotaru}. Here $\fc$ is the Clifford algebra associated to $\bigoplus \mathbb{C}x_{i} \oplus \bigoplus \mathbb{C}y_{j}$ with its natural nondegenerate bilinear form $(x_{i}, y_{j}) = \delta_{ij}$. For a module $M \in \cO_{1/e}(\mathfrak{S}_{n})$, the algebra $H_{1/e}(\mathfrak{S}_{n}) \otimes \fc$ acts on the space $M \otimes \wedge^{\bullet}\fh$, and the Dirac cohomology is defined to be, as usual, $\ker(D)/\ker(D)\cap\text{im}(D)$. This is a representation of $\widetilde{W}$, a certain double-cover of the group $W$. Then, by \cite[Theorem 5.1]{HuangWong}, $H_{D}(L(\lambda)) = \bigoplus_{\mu \leq \lambda}S_n(\mu) \otimes \chi$, where $\chi$ is a 1-dimensional character of the double cover $\widetilde{W}$. We refer to \cite{HuangWong} for details.

	\section{Graded  free resolutions of algebraic varieties, \\ Betti numbers,  and Castelnuovo--Mumford regularity    }

We now consider the consequences of our results for computing minimal resolutions of linear subspace arrangements.  Easy examples of ideals whose resolutions we compute include the braid arrangements of type $A$ and type $D$.  
Such minimal resolutions are  difficult to   compute geometrically \cite{MR520233}. 
 As a consequence, we prove a combinatorial formula for the 
Betti numbers of the ideal of the $m$-equals arrangement  predicted in \cite{conjecture}.  
We also calculate the   Castelnuovo--Mumford regularity 
 for the coordinate ring of these arrangements, a notoriously difficult problem in general (see 
 \cite{MR1942401, MR3299723}).

It is pointed out in \cite{MR2037715} that BGG resolutions via parabolic Verma modules for Lie algebras can be used to provide commutative algebra  resolutions of determinantal ideals by viewing the coordinate ring  as a unitarizable  highest weight module.    
  We employ our  Cherednik algebra resolutions in an analogous fashion.   
  The first of these commutative algebra resolutions, given in \cref{eideal}, was predicted in \cite{conjecture} and concerns the smallest ideal, $I_{e,1,n}$, of  the polynomial representation (this is the vanishing ideal of the subspace arrangement, $X_{e,1,n}$, consisting of $e$ equal coordinates for $e\leq  n$).  We then provide a cyclotomic generalisation   of this resolution in \cref{eideal2}.  
  The third resolution, given in \cref{keideal}, concerns the smallest quotient, $\CC[X_{e,k,n}]$, of the polynomial representation (this is the coordinate ring of the subspace arrangement, $X_{e,k,n}$,  consisting of $k$ clusters of $e$ equal coordinates for $ke=n$); 
 the ideal  vanishing  on this space was studied in \cite{conjecture}, 
 however since neither this ideal nor its quotient is unitary (in general) 
 the authors did not predict any resolution arising via Cherednik algebras.  
 
% We remark that the resolutions of $\CC[X_{e,k,n}]$  and $I_{e,1,n}$ 
%  are swapped (on the level of Cherednik algebras) by   Ringel duality (modulo applying the  Morita equivalence of \cref{Cahracteristicchange}). 
 
%
%It is pointed out in \cite{MR2037715} that BGG resolutions via parabolic Verma modules for Lie algebras can be used to provide commutative algebras resolutions of determinantal ideals by viewing the coordinate ring  as a unitarizable  highest weight module.    
%We employ our  Cherednik algebra resolutions in an analogous fashion.   
%The first of these commutative algebra resolutions, given in \cref{eideal}, was predicted in \cite{conjecture} and concerns the smallest ideal of  the polynomial representation (this is the coordinate ring of the subspace arrangement consisting of $e$ equal coordinates for $e< n$).  We then provide a cyclotomic generalisation   of this resolution in \cref{eideal2}.  
%The third resolution, given in \cref{keideal}, concerns the smallest quotient of the polynomial representation (this is the coordinate ring of the subspace arrangement consisting of $k$ clusters of $e$ equal coordinates for $ke=n$); 
%the ideal  vanishing  on this space was studied in \cite{conjecture}, 
%however since neither this ideal nor its quotient is unitary (in general) 
%the authors did not predict any resolution arising via Cherednik algebras.  

	\subsection{Commutative algebra} 
	
	Let us discuss the consequences that the existence of the BGG resolution has for the study of graded modules over $\C[x_{1}, \dots, x_{n}] = \C[\underline{x}]$. First of all, for every $\mu \vdash n$, the standard module $\Delta_{1/e}(\mu)$ is free as a $\C[\underline{x}]$-module. So the resolution $C_{\bullet}(\lambda)$ is, in fact, a free resolution of $L_{1/e}(\lambda)$ when we view all involved modules as $\C[\underline{x}]$-modules.
	
	An observation now is that every module in category $\cO_{1/e}(\mathfrak{S}_{n})$ automatically acquires a grading compatible with the usual grading on $\C[\underline{x}]$, as follows. Consider the deformed Euler element\footnote{We remark that our Euler element \lq$\mathsf{eu}$\rq \; differs from the one used in \cite{conjecture} by the constant $n(e - n + 1)/2e$} $\mathsf{eu} := \frac{1}{2}\sum_{i = 1}^{n}x_{i}y_{i} + y_{i}x_{i} \in H_{1/e}$. This is a grading element of $H_{1/e}$ in the sense that $[\mathsf{eu}, x_{i}] = x_{i}, [\mathsf{eu}, y_{i}] = -y_{i}$, and $[\mathsf{eu}, w]= 0$ for $w \in \mathfrak{S}_{n}$. Any module in category $\cO_{1/e}(\mathfrak{S}_{n})$ is now graded by generalized eigenspaces for $\mathsf{eu}$:
	\[
	M = \bigoplus_{a \in \C}M_{a}, M_{a} := \{m \in M : (\mathsf{eu}-a)^{k}m = 0 \; \text{for} \; k \gg 0\}.
	\] 
	Note that, since the grading on $M$ was defined using an element of $H_{1/e}$, every morphism in category $\cO_{1/e}(\mathfrak{S}_{n})$ has degree $0$. In particular, this grading is \emph{different} from the grading of objects in $\cO_{1/e}(\mathfrak{S}_{n})$ that has been used so far in this paper. The grading by generalized eigenspaces of $\mathsf{eu}$, however, is better-suited for the purposes of commutative algebra. 
	
	A priori, $M \in \cO_{1/e}(\mathfrak{S}_{n})$ is only $\C$-graded, but in our case we can do better. Since $[\mathsf{eu}, w] = 0$ for $w \in \mathfrak{S}_{n}$, $\mathsf{eu}$ may be seen as an endomorphism of the $\mathfrak{S}_{n}$-module $S_{n}(\tau) \cong 1 \otimes S_{n}(\tau) \subseteq \C[\underline{x}] \otimes \tau = \Delta_{1/e}(\tau)$. Thus, $\mathsf{eu}$ acts by a scalar $c_{\tau}$ on $S_{n}(\tau)$, and by the definition of $\Delta_{1/e}(\tau)$ we get that $\Delta_{1/e}(\tau)_{a} \neq 0$ if and only if $a = c_{\tau} + k$ for some $k \in \Z_{\geq 0}$. Moreover,
	\[
	\Delta_{1/e}(\tau)_{c_{\tau} + k} = \C[\underline{x}]_{k} \otimes S_n(\tau)
	\]
	\noindent where $\C[\underline{x}]_{k}$ denotes the subspace of homogeneous polynomials of degree $k$ in the variables $x_{1}, \dots, x_{n}$. We will write $\C[\underline{x}] \otimes S_n(\lambda)$ to refer to the $\C[\underline{x}]$-module $\Delta_{1/e}(\lambda)[c_{\lambda}]$, where the brackets denote the usual grading shift. Thus, $\C[\underline{x}]\otimes S_n(\lambda)$ is $\Z_{\geq 0}$-graded, and $(\C[\underline{x}] \otimes S_n(\lambda))_{k} = \C[\underline{x}]_{k} \otimes S_n(\lambda)$. 
	
	Now consider the resolution of the graded $\C[\underline{x}]$-module $L_{1/e}(\lambda)[c_{\lambda}]$, where the $i$-th term of the complex is given by
	\[
	\bigoplus_{\substack{%\mu \vdash n \\ 
			\mu \in \Po_e(\lambda) \\ \hd(\mu) = i}}(\C[\underline{x}] \otimes S_n(\mu))[c_{\lambda} - c_{\mu}]
	\] 
	We remark that, since $\lambda$ and $\mu$ belong to the same block of category $\cO_{1/e}(\mathfrak{S}_{n})$, $c_{\lambda} - c_{\mu}$ is actually an integer. Of course, this is the same as the BGG resolution $C_{\bullet}(\lambda)$, but we write it in this way to emphasize that we are only interested in the $\C[\underline{x}]$-module structure. By abuse of notation, we also denote this complex by $C_{\bullet}(\lambda)$. Note that $(\C[\underline{x}] \otimes S_n(\mu))[c_{\lambda} - c_{\mu}] = \Delta_{1/e}(\mu)[c_{\lambda}]$, from where it follows that all the maps in the complex have degree $0$ as maps of graded $\C[\underline{x}]$-modules. In particular, $C_{\bullet}(\lambda)$ is a graded-free resolution of $L_{1/e}(\lambda)[c_{\lambda}]$. 
	
	The value of $c_{\lambda}$ can be expressed in terms of the content of the boxes of $\lambda$, namely
	$$
	c_{\lambda} = \frac{n}{2} -\frac{1}{e}\sum_{\square \in \lambda} \text{column}(\square) - \text{row}(\square)
	$$
	\noindent It follows from Section \ref{subsect:hd rimhooks} or from Lemma \ref{lem:hd=ell} 
	that if $\hd(\mu) < \hd(\nu)$ then $c_{\mu} < c_{\nu}$. In particular, when viewing the differential in the resolution $C_{\bullet}(\lambda)$ as matrices with coefficients in $\C[\underline{x}]$, no nonzero entry of the differential is a degree $0$ element of $\C[\underline{x}]$. It follows immediately that:
	
	\begin{lem}\label{lem:minimalres}
		The complex $C_{\bullet}(\lambda)$ is a minimal graded free resolution of $L_{1/e}(\lambda)[c_{\lambda}]$.
	\end{lem}
	
	Lemma \ref{lem:minimalres} implies a combinatorial formula for computing  many interesting invariants of the module $L_{1/e}(\lambda)[c_{\lambda}]$.  	
	In the rest of this section, if $L_{1/e}(\lambda)$ is unitary we write:
	$$n:=|\lambda|,\quad k:=e\text{-weight}(\lambda),\quad h:=\#\text{columns}(\lambda)$$
	Recalling the basics of abaci in Section \ref{subsect:abacus} this means that the abacus $\cA_e(\lambda)$ has $h$ nonempty runners and there are $k$ vacant spaces in $\cA_e(\lambda)$ with some bead to their right.
	\begin{prop}
		Suppose $L_{1/e}(\lambda)$ is unitary. Then,
				\begin{enumerate}	
			\item[(1)] $\displaystyle{\beta_{i,j} = \sum_{\substack{\mu \in \Po_e(\lambda) \\ c_{\lambda} - c_{\mu} = -j \\ \hd(\mu) = i}} \dim(S_n(\mu))}$, 
			\item[(2)] $\text{pdim}(L_{1/e}(\lambda)) = (h-1)k$,  \item[(3)] $\text{depth}(L_{1/e}(\lambda)) = n - (h-1)k$
		\end{enumerate}
		 where $\beta_{i,j}$ denotes the $(i,j)$-graded Betti number of $L_{1/e}(\lambda)[c_{\lambda}]$, and $\text{pdim}$ stands for the projective dimension \emph{as a graded $\C[\underline{x}]$-module}. 
	\end{prop}
	\begin{proof}
		Statement (1) is clear from the form of the resolution $C_{\bullet}(\lambda)$. The maximal homological degree of a partition in $\Po_e(\lambda)$ is acquired by sliding all the beads to the left and then sliding the highest bead $k$ spaces to the right. (2) follows from here. Finally, by the Auslander-Buchsbaum formula, (3) is equivalent to (2). 
	\end{proof}

	Another consequence of Lemma \ref{lem:minimalres} and the fact that the function $c_{\lambda}$ is strictly increasing on homological degrees, is the computation of the Castelnuovo-Mumford regularity of the module $L_{1/e}(\lambda)[c_{\lambda}]$. Recall that, by definition, the regularity of a module $M$ is 
	\[
	\reg(M) := \max\{j : \text{there exists} \; i \; \text{such that} \; \beta_{i, i+j}(M) \neq 0\}
	\]
	In other words, for a minimal graded-free resolution $C_{\bullet}$ of $M$, for each $i = 0, \dots, \text{pdim}(M)$, let $n_{i}$ be the maximum degree of a generator of $C_{i}$, and $m_{i} := n_{i} - i$. Then, $\reg(M) = \max_{i}\{m_{i}\}$. The Castelnuovo-Mumford regularity is a measure of the computational complexity of the module $M$ and it is, in general,   incredibly difficult to compute, cf. \cite{MR1942401, MR3299723}.

	\begin{prop}\label{prop:regularity}
		Suppose $L_{1/e}(\lambda)$ is unitary. Let $\mu_0 \in \Po_e(\lambda)$ be   obtained by, first, sliding all beads of $\cA_{e}(\lambda)$ to the left, and then, sliding the upmost beat $k$ spaces to the right. % This is the partition of maximal homological degree in $\Po_e(\lambda)$. %to the right in such a way that we still get a partition of $n$. 
		Then,
		\[
		\reg(L_{1/e}(\lambda)[c_{\lambda}]) = (c_{\mu_0} - c_{\lambda}) - (h-1)k
		\]
	\end{prop}
	\begin{proof}
		As in the paragraph above the statement of the proposition, let us denote by $n_{i}$ the maximum degree of a generator of $C(\lambda)_{i}$, and $m_{i} := n_{i} - i$. Note that $n_{i} := \max\{c_{\mu} - c_{\lambda} : \mu \in \Po_e(\lambda), \hd(\mu) = i\}$. Since the $c$-function is increasing in homological degree, the sequence $(n_{i})$ is increasing and therefore the sequence $(m_{i})$ is nondecreasing. So the regularity of $L_{1/e}(\lambda)[c_{\lambda}]$ is $m_{\text{pdim}(L_{1/e}(\lambda))}$. Since $\text{pdim}(L_{1/e}(\lambda)) = (h-1)k$, the result follows.
	\end{proof}
	
	%Let us assume that $L_{1/e}(\lambda)$ indeed has full support, the other case will be treated in Subsection ???. Since $\dim(\supp(L_{1/e}(\lambda)) = |\lambda|$, we get that $L_{1/e}(\lambda)$ is Cohen-Macaulay if and only if $\text{pdim}(L_{1/e}(\lambda)) = 0$, that is, if and only if $L_{1/e}(\lambda) = \Delta_{1/e}(\lambda)$, and we get the following.
	
	%\begin{prop}
	%	Let $\lambda\vdash n$ be such that $L_{1/e}(\lambda)$ is unitary and $\lambda \neq (e^{k})$ for $k \in \Z_{>0}$. Then, the following are equivalent.
	%	\begin{enumerate}
	%		\item $L_{1/e}(\lambda)$ is Cohen-Macaulay as a $\C[\underline{x}]$-module.
	%		\item $L_{1/e}(\lambda)$ is a free $\C[\underline{x}]$-module.
	%		\item $L_{1/e}(\lambda) = \Delta_{1/e}(\lambda)$.
	%	\end{enumerate}
	%\end{prop}
	
	\begin{eg}
		Consider $e = 5$, $n = 15$ and $\lambda = (3^4, 2, 1)$. Then, $\text{pdim}L_{1/e}(\lambda) = 4$, so $L_{1/e}(\lambda)$ is not Cohen-Macaulay and a minimal graded-free resolution of $L_{1/e}(\lambda)[c_{\lambda}]$ is 
		\begin{align*}
		0 &\to (3, 2, 1^{10})[-9] \rightarrow (3, 2^2, 1^7)[-5] &\to (3, 2^6)[-3] \oplus (3^3, 1^6)[-3]  &\to  (3^3, 2^2, 1^2)[-1] \rightarrow (3^4, 2, 1) \\ &\to L_{1/e}(\lambda)[c_{\lambda}] \rightarrow 0
		\end{align*}	
		\noindent where for brevity, we write $\mu[d]$ in place of $(\C[\underline{x}] \otimes S_n(\mu))[d]$. From the resolution, we see that $\reg(L_{1/e}(\lambda)[c_{\lambda}]) = 5$. 
	\end{eg}

	\subsubsection{The $e$-equals ideal} \label{eideal}
We examine these results in the situation where the modules $L_{1/e}(\lambda)$  
 have a clear geometric meaning.  
The representation theoretic import of $X_{e,1,n}$  was first noticed and explained 
in \cite{conjecture}.  Resolutions of the ideals vanishing on these subspace arrangements are given by BGG resolutions of the corresponding unitary module 
for $H_{1/e}(\mathfrak{S}_n)$.

	Let $n = (e-1)p + q$, with $0 \leq q < e-1$. Consider the partition $\lambda = ((e-1)^{p}, q)$ of $n$. Note that the $e$-abacus of $\lambda$ has exactly one empty runner, and the module $L_{1/e}(\lambda)$ is unitary. In fact, it follows from \cite{MR3658722} that $L_{1/e}(\lambda)$ is isomorphic to the socle of the polynomial representation 
	$$\Delta_{1/e}(\triv)\cong \mathbb{C}[x_1,x_2,\dots, x_n]$$ which by \cite[Theorem 5.10]{MR2534594} coincides with the {\sf $e$-equals ideal}  
	 $I_{e, 1, n}$ of functions vanishing on the set
	\[
	X_{e, 1, n} := \mathfrak{S}_{n}\{(z_{1}, \dots, z_{n}) \in \C^{n} : z_{1} = \dots = z_{e}\}.
	\]
	\noindent Note that $X_{e, 1, n}$ is an arrangement of $n \choose e$ linear subspaces of $\C^{n}$, each of dimension $n - e + 1$. When $e = 2$, $X_{2, 1, n}$ is nothing but the braid arrangement in $\C^{n}$, which consists of the reflection hyperplanes for the action of $\mathfrak{S}_{n}$ on $\C^{n}$. Let us give a set of generators for the ideal $I_{e, 1, n}$, following \cite{MR2534594, MR2959035}. Consider the partition $\lambda^{T} = ((p+1)^{q}, p^{e-1-q})$, which has exactly $e-1$ parts. Now consider the polynomial
	\[
	p_{\lambda^{T}}(x_{1}, \dots, x_{n}) = V(x_{1}, \dots, x_{\lambda^{T}_{1}})V(x_{\lambda^{T}_{1} + 1}, \dots, x_{\lambda^{T}_{1} + \lambda^{T}_{2}})\cdots V(x_{\lambda^{T}_{1} + \cdots + \lambda^{T}_{e-2} + 1}, \dots, x_{n})
	\]
	\noindent where $V(x_{1}, \dots, x_{k})$ is the Vandermonde determinant $\prod_{i < j}(x_{i} - x_{j})$. Then, the ideal $I_{e, 1, n}$ is generated by the $\mathfrak{S}_{n}$ images of the polynomial $p_{\lambda^{T}}$. 
	
	Since $L_{1/e}(\lambda)$ and $L_{1/e}(\triv)$ lie in the same block of category $\cO_{1/e}$, the weight of the partition $\lambda$ is $k = \lfloor n/e\rfloor$. Thus, as  was observed in \cite{conjecture}, the projective dimension of the algebra of functions $\C[X_{e, 1, n}] = \C[x_{1}, \dots, x_{n}]/I_{e, 1, n}$ is $\text{pdim}(\C[X_{e, 1, n}]) = \text{pdim}(L_{1/e}(\lambda)) + 1$  $= (e-2)\lfloor n/e\rfloor + 1$.  Since $\dim(X_{e, 1, n}) = n-e+1$, it follows that $\C[X_{e, 1, n}]$ is Cohen-Macaulay if and only if $e = 2$ or $\lfloor n/e\rfloor = 1$. This way, we recover part of \cite[Proposition 3.11]{EGL}.
	
	\begin{eg}
		Consider $e = 4$, $n = 10$. The minimal submodule in $\Delta_{1/e}(\triv)$ is $I_{4, 1, 10}$, and it is isomorphic to $L_{1/e}(3^{3}, 1)$. Note that $c_{\lambda} = 23/4$. The resolution of $L_{1/4}(3^3, 1)[-23/4]$ is   given by
		\begin{align*}
		0 &\to  (2, 1^{8})[-8] \rightarrow (2^{2}, 1^{6})[-6] \rightarrow (2^{5})[-3] \oplus (3^2, 1^{4})[-2] \rightarrow (3^{2}, 2, 1^{2})[-1] \rightarrow   (3^{3}, 1) \\ &\to L_{1/4}(3^{3}, 1)[-23/4] \rightarrow 0
		\end{align*}
		A resolution of the coordinate ring $\C[x_{1}, \dots, x_{10}]/I_{4, 1, 10}$ looks similar, but each term is further shifted by $-12$ (because $c_{\triv} - c_{\lambda} = -12$), and the end of the sequence is $(3^{3}, 1)[-12] \rightarrow (10) \rightarrow \C[\underline{x}]/I_{4, 1, 10} \rightarrow 0$. Note that the regularity of $\C[x_{1}, \dots, x_{10}]/I_{4, 1, 10}$ is $15$. 
	\end{eg}
	
	Let us now compute the regularity of the subspace arrangement $X_{e, 1, n}$. 
	
	\begin{prop}\label{prop:reg e-equals}
		The regularity of the $\C[\underline{x}]$-module $\C[X_{e, 1, n}]$ is given by
		\[
		\reg(\C[X_{e, 1, n}]) = \begin{cases} \lfloor n/e\rfloor(n-e + 1) - 1, & \text{if} \; n/e \in \Z  \\ \lfloor n/e\rfloor(n-e+2) - 1, & \text{else.} \end{cases}
		\]

	\end{prop}
	\begin{proof}
		Let us write $n = (e-1)p + q = ep_{1} + q_{1}$, with $0 \leq q < e-1$ and $0 \leq q_{1} < e$.  As above, let $\lambda = ((e-1)^p, q)$ be the partition such that $L_{1/e}(\lambda)$ is isomorphic to the socle of $\Delta_{1/e}(\triv)$. Note that the $e$-core of any partition in the block of $\triv=(n)$ is $(q_1)$ and the $e$-weight is $p_1$. It then follows from the rimhook description of homological degree in Section \ref{subsect:hd rimhooks}
		that the partition $\mu_0$ with highest homological degree in $\Po_e(\lambda)$ is given by adding $p_1$ vertical strips of length $e$ in the leftmost column to the $e$-core of $\lambda$: thus $\mu_0=(a,1^{n-a})$ where $a=q_1$ if $q_1>0$ and $a=1$ if $q_1=0$, and $\hd(\mu_0)=(e-2)p_1$.
		Now it follows by a direct computation that  		\[
		\reg(\C[X_{e, 1, n}]) = c_{\mu_0} - c_{\triv} - (e-2)p_{1} - 1 = \frac{(n-a)n}{e} - (e-2)p_{1} - 1
		\]
		\noindent which coincides with the formula in the statement of the proposition.
	\end{proof}
	
	\subsubsection{More BGG resolutions and  a generalisation of the $e$-equals ideal.}  \label{eideal2}
 	We take $\ell$th powers and obtain a generalisation of the $e$-equals ideal.   These subspaces arrangements admit   commutative algebra resolutions which can be constructed via   BGG-resolutions for the Cherednik algebra of $G(\ell,1,n)$ (which we also construct in this section).   
 	Consider the ideal $I_{e, 1, n}(\ell)$ of polynomials vanishing on the set		\[
	X_{e, 1, n}(\ell) := \mathfrak{S}_{n}\{(z_{1}, \dots, z_{n}) \in \C^{n} : z_{1}^{\ell} = z_{2}^{\ell} = \cdots = z_{e}^{\ell}\}.
	\]
		Note that $X_{e, 1, n}(\ell)$ is an arrangement of $\ell^{e}{n \choose e}$ linear subspaces of $\C^{n}$, each of dimension $n - e + 1$. When $e = \ell = 2$, $X_{2, 1, n}(2)$ is the braid arrangement of type $D_{n}$, consisting of reflection hyperplanes for the reflection representation of the Weyl group of type $D_{n}$ on $\C^{n}$. To give a set of generators for the ideal $I_{e, 1, n}(\ell)$, recall from the previous subsection the partition $\lambda = ((e-1)^{p}, q)$ and the polynomial $p_{\lambda^{T}} \in \C[x_{1}, \dots, x_{n}]$. According to \cite[Proposition 2.5]{MR2959035}, a set of generators of the ideal $I_{e, 1, n}(\ell)$ is given by the $\mathfrak{S}_{n}$-images of $p_{\lambda^{T}}(x_{1}^{\ell}, \dots, x_{n}^{\ell})$. 
	
	Our next goal is to construct a graded-free resolution of the algebra of functions $\C[x_{1}, \dots, x_{n}]/I_{e, 1, n}(\ell)$. In order to do this, we will use the following well-known commutative algebra result. 

	\begin{lem}\label{cor:powers}
		Let $F_{1}, F_{2}, F_{3}$ be free $\C[x_{1}, \dots, x_{n}]$-modules of finite rank, with bases $\{v^{1}_{1}, \dots, v^{1}_{i_{1}}\}$, $\{v^{2}_{1}, \dots, v^{2}_{i_{2}}\}$ and $\{v^{3}_{1}, \dots, v^{3}_{i_{3}}\}$, respectively. Let $A: F_{1} \to F_{2}$, $B: F_{2} \to F_{3}$ be morphisms defined in the given bases by matrices $(f_{ij}(x_{1}, \dots, x_{n})), (g_{jk}(x_{1}, \dots, x_{n}))$, respectively, and define new morphisms $\widetilde{A}, \widetilde{B}$ by the matrices $(f_{ij}(x_{1}^{\ell}, \dots, x_{n}^{\ell}))$, $(g_{jk}(x_{1}^{\ell}, \dots, x_{n}^{\ell}))$, respectively. If $\image(A) = \ker(B)$, then $\image(\widetilde{A}) = \ker(\widetilde{B})$. 
	\end{lem}
	
	Note that, for   $\mu\in \mptn 1 n$, the module $\C[x_{1}, \dots, x_{n}] \otimes S_n(\mu)$ has a distinguished basis indexed by   $\Std(\mu)$. Thus, if $\lambda$ is a unitary partition, we can apply Lemma \ref{cor:powers} to the complex $C_{\bullet}(\lambda)$ (viewed as a complex of free $\C[x_{1}, \dots, x_{n}]$-modules) to obtain a complex $\widetilde{C}_{\bullet}(\lambda)$, which is exact outside of degree $0$. By construction, thanks to \cite[Proposition 2.5]{MR2959035}, when $\lambda = ((e-1)^{p}, q)$, the zeroth homology of $\widetilde{C}_{\bullet}(\lambda)$ coincides with the ideal $I_{e, 1, \lambda}(\ell)$. Moreover, by multiplying the grading shifts of $C_{\bullet}(\lambda)$ by $\ell$, this obtains a minimal graded-free resolution of $I_{e, 1, \lambda}$, and extending by $\C[x_{1}, \dots, x_{n}]$, of the algebra of functions $\C[X_{e, 1, n}(\ell)]$. We then obtain the following result. 
	
	\begin{prop}
		The projective dimension of $\C[X_{e, 1, n}(\ell)]$ coincides with that of $\C[X_{e, 1, n}]$, which is $(e-2)\lfloor n/e \rfloor + 1$ so that, regardless of $\ell$, $\C[X_{e, 1, n}(\ell)]$ is Cohen-Macaulay if and only if $e = 2$ or $\lfloor n/e\rfloor = 1$. The regularity of $\C[X_{e, 1, n}(\ell)]$ is given by
 		$$
		\reg(\C[X_{e, 1, n}(\ell)]) = \begin{cases}\lfloor n/e\rfloor(\ell(n-1) - e + 2) -1 & \text{if} \; n/e \in \ZZ \\
		\lfloor n/e\rfloor(\ell n - e + 2) - 1 & \text{else} \end{cases}
		$$
	\end{prop}
	
	We have obtained the complex $\widetilde{C}_{\bullet}(\lambda)$ by means of pure commutative algebra. As it turns out, $\widetilde{C}_{\bullet}(\lambda)$ is a complex of standard modules for the rational Cherednik algebra of the group $G(\ell, 1, n) := \mathfrak{S}_{n} \ltimes (\ZZ/\ell\ZZ)^{n}$ under a special class of parameters. The group $G(\ell, 1, n)$ is a complex reflection group, acting naturally on $\C^{n}$, and the rational Cherednik algebra depends on a function $\widetilde{c}: S \to \C$, where $S \subseteq G(\ell, 1, n)$ is the set of reflections and $\widetilde{c}(s) = \widetilde{c}(wsw^{-1})$ for every $s \in S$, $w \in G(\ell, 1, n)$. Here, for a complex number $c \in \C$, we will take \emph{any} function $\widetilde{c}$ such that $\widetilde{c}(s) = c$, if $s \in G(\ell, 1, n)$ is conjugate to a reflection in $\mathfrak{S}_{n}$. Any other reflection in $G(\ell, 1, n)$ is conjugate to a nonzero element of, say, the first copy of $\Z/\ell\Z$, so we have $\ell - 1$ more parameters for $H_{\widetilde{c}}(G(\ell, 1, n))$, let us call them $c_{1}, \dots, c_{\ell -1}$.  
	
	The rational Cherednik algebra $H_{\widetilde{c}}(G(\ell, 1, n))$ admits a presentation very similar to that of the rational Cherednik algebra $H_{c}(\mathfrak{S}_{n})$ of the symmetric group. We will not give this presentation. Instead, we remark that $H_{\widetilde{c}}(G(\ell, 1, n))$ is the subalgebra of $\text{End}_{\C}(\C[x_{1}, \dots, x_{n}])$ generated by the functions $x_{i}$ of multiplication by $x_{i}$ ($i = 1, \dots, n$), the elements of $G(\ell, 1, n)$ (naturally viewed as automorphisms of $\C[x_{1}, \dots, x_{n}])$ and the \emph{Dunkl--Opdam operators}:
 	\[
	\widetilde{D}_{i} := \partial_{i} - c \sum_{j \neq i}\sum_{t = 0}^{\ell - 1}\frac{1}{x_{i} - \xi^{t}x_{j}}(1 - (ij)^{t}) - \sum_{k = 1}^{\ell - 1}\frac{2c_{k}}{(1 - \xi^{-k})x_{i}}(1 - \xi_{i}^{k})
	\] 
	 where $\xi := \exp(2\pi\sqrt{-1}/\ell)$, $\xi_{i} \in G(\ell, 1, n)$ is the element that acts by multiplication by $\xi$ on the $i$-th coordinate in $\C^{n}$, and $(ij)^{t} \in G(\ell, 1, n)$ is $(ij)^{t} = \xi_{i}^{t}\xi_{j}^{-t}(ij)$. Let us remark that a similar presentation exists for the algebra $H_{c}(\fS_{n})$, the Dunkl  operators are now given by
 	\[
	D_{i} = \partial_{i} - c\sum_{j \neq i}\frac{1}{x_{i} - x_{j}}(1 - (ij)).
	\] 
	
	We will need the following result, that relates the operators $D_{i}$ and $\widetilde{D}_{i}$.
	
	\begin{lem}\label{lemma:dunkl}
		For $g \in \C[x_{1}, \dots, x_{n}]$, denote by $\widetilde{g} := g(x_{1}^{\ell}, \dots, x_{n}^{\ell})$. Then, for any $i = 1, \dots, n$:
		\[
		\widetilde{D}_{i}(\widetilde{g}) = \ell x_{i}^{\ell - 1}\widetilde{D_{i}(g)}
		\]	
	\end{lem}
\begin{proof}
	First of all, note that $\widetilde{g}$ is invariant under the action of $(\Z/\ell\Z)^{n}$ on $\C[x_{1}, \dots, x_{n}]$, and so it follows that 
	\[
	\widetilde{D_{i}(\widetilde{g})}  = \partial_{i}(\widetilde{g}) - c\sum_{j \neq i}\sum_{t = 0}^{\ell - 1}\frac{\widetilde{g} - (ij)\widetilde{g}}{x_{i} - \xi^{t}x_{j}}.
	\]
 Now let $h(x_{1}, \dots, x_{n}) \in \C[x_{1}, \dots, x_{n}]$ be such that $g - (ij)g = (x_{i} - x_{j})h$. Note that it follows that $\widetilde{g} - (ij)\widetilde{g} = (x_{i}^{\ell} - x_{j}^{\ell})\widetilde{h}$, so
	\[
	\sum_{t = 0}^{\ell - 1}\frac{\widetilde{g} - (ij)\widetilde{g}}{x_{i} - \xi^{t}x_{j}} = \sum_{t = 0}^{\ell} \prod_{\substack{k = 0 \\ k \neq t}}^{\ell - 1}(x_{i} - \xi^{k}x_{j})\widetilde{h} = \ell x_{i}^{\ell - 1}\widetilde{h}
	\]	
 and the result follows.
\end{proof}

The algebra $H_{\widetilde{c}}(G(\ell, 1, n))$ still admits a triangular decomposition $H_{\widetilde{c}}(G(\ell, 1, n)) = \C[x_{1}, \dots, x_{n}] \otimes \C G(\ell, 1, n) \otimes \C[y_{1}, \dots, y_{n}]$, where $y_{i}$ is the Dunkl--Opdam operator $\widetilde{D}_{i}$. In particular, one can still define standard modules. For an irreducible representation $E$ of $G(\ell, 1, n)$, we have the standard module $\Delta_{\widetilde{c}}(E)$. As a $\C[x_{1}, \dots, x_{n}]$-module, $\Delta_{\widetilde{c}}(E) = \C[x_{1}, \dots, x_{n}] \otimes E$. 
	
	The irreducible representations, $S_n(\lambda)$, of $G(\ell, 1, n)$ are indexed by the set $\mptn \ell n$, and each $S_n(\lambda)$  has a natural basis indexed by  the set $\Std(\lambda)$. In particular, if $\lambda \in \mptn 1 n$, we can consider the $\ell$-partition $\widetilde{\lambda} \in \mptn \ell n$ given by $\widetilde{\lambda} = (\lambda, \emptyset, \dots, \emptyset)$. 
	The sets $\Std(\widetilde\lambda)$ and $\Std({\lambda})$ are obviously identified. Moreover, $G(\ell, 1, n)$ admits a natural surjection to $\mathfrak{S}_{n}$, and the irreducible representation $S_{n}(\widetilde{\lambda})$ of $G(\ell, 1, n)$ is simply given by the $\fS_{n}$-irreducible $S_{n}(\lambda)$ under this surjection.
	
	\begin{prop}\label{prop:powers}
		Let $c \in \C$. Then, for any   $\lambda, \mu \in \mptn 1 n$ and any parameter $\widetilde{c}$ as above, there is a natural identification	
		\[
		\begin{array}{rcl}
\sim \ :	\Hom_{H_{c}(\fS_{n})}(\Delta_{c}(\lambda), \Delta_{c}(\mu)) & \buildrel \cong \over \longrightarrow & \Hom_{H_{\widetilde{c}}(G(\ell, 1, n))}(\Delta_{\widetilde{c}}(\widetilde{\lambda}), \Delta_{\widetilde{c}}(\widetilde{\mu})) 
%\\
%		f & \mapsto & \widetilde{f}
		\end{array}
		\]
		 given as follows. For a standard Young tableau $\stt \in \Std(\lambda)$, if $f \in \Hom_{H_{c}(\fS_{n})}(\Delta_{c}(\lambda), \Delta_{c}(\mu))$ is given by $f(1 \otimes \stt) = \sum_{\sts \in \Std(\mu)}f_{\stt\sts}(x_{1}, \dots, x_{n}) \otimes \sts$, then $\widetilde{f}(1 \otimes \stt) = \sum_{\sts \in \Std(\mu)}f_{\stt\sts}(x_{1}^{\ell}, \dots, x_{n}^{\ell}) \otimes \sts$. 
	\end{prop}
	\begin{proof}
		We need to show, first, that $\widetilde{f}|_{1 \otimes S_{n}(\widetilde{\lambda})}$ is a map of $G(\ell, 1, n)$-representations. This follows from the fact that, for any polynomial $g \in \C[x_{1}, \dots, x_{n}]$, $g(x_{1}^{\ell}, \dots, x_{n}^{\ell})$ is invariant under the action of $(\Z/\ell\ZZ)^{n}$. Now we need to show that, for any standard Young tableau $\stt \in \Std(\lambda)$, $\widetilde{f}(1 \otimes \stt)$ is annihilated by all Dunkl operators $\widetilde{D}_{i}$. This is a direct consequence of Lemma \ref{lemma:dunkl}.
		
		This shows that $f \mapsto \widetilde{f}$ does define a morphism, which is clearly injective. To show that it is bijective, let $h: \Delta_{\widetilde{c}}(\widetilde{\lambda}) \to \Delta_{\widetilde{c}}(\widetilde{\mu})$ be a morphism. In particular, $h|_{1 \otimes S_{n}(\widetilde{\lambda})}$ is a map of $G(\ell, 1, n)$-modules. This implies that, if $h(1 \otimes \stt) = \sum_{\sts \in \Std(\mu)}h_{\stt\sts}(x_{1}, \dots, x_{n}) \otimes \sts$, then $h_{\stt\sts}(x_{1}, \dots, x_{n}) \in \C[x_{1}^{\ell}, \dots, x_{n}^{\ell}]$ for every $\sts \in \Std(\mu)$. Thanks to Lemma \ref{lemma:dunkl}, this implies that $h = \widetilde{f}$ for some $f: \Delta_{c}(\lambda) \to \Delta_{c}(\mu)$. % We are done. 
	\end{proof}
	
	\begin{rmk}
		If $c \not\in 1/2 + \Z$, then the existence of an isomorphism between $\Hom_{H_{c}(\fS_{n})}(\Delta_{c}(\lambda), \Delta_{c}(\mu))$ and $\Hom_{H_{\widetilde{c}}(G(\ell, 1, n))}(\Delta_{\widetilde{c}}(\widetilde{\lambda}), \Delta_{\widetilde{c}}(\widetilde{\mu}))$ follows from \cite[Proposition 5.9]{GGOR}.
	\end{rmk}
	
	By Proposition \ref{prop:powers} and Lemma \ref{cor:powers}, we have that if $\lambda$ is a unitary partition of $n$, then the complex $\widetilde{C}_{\bullet}(\lambda)$ is actually a complex of standard modules for $H_{\widetilde{c}}(G(\ell, 1, n))$, which is exact outside of degree zero, and thus it is a BGG resolution of its zeroth homology. 
	
	\begin{rmk}\label{rmk1}
		The zeroth homology of $\widetilde{C}_{\bullet}(\lambda)$ is \emph{not} necessarily   an irreducible $H_{\widetilde{c}}(G(\ell, 1, n))$-module. For example, if $\lambda = ((e-1)^{p}, q)$, we have seen that $H_{0}(\widetilde{C}_{\bullet}(\lambda))$ is the ideal $I_{e, 1, n}(\ell)$. When $\ell = 2$, $e < n$ is even and the parameter $\widetilde{c}$ is such that $\widetilde{c}(s) = 0$ if $s$ is not conjugate to a reflection in $\mathfrak{S}_{n}$, then this is an indecomposable, but not irreducible, $H_{\widetilde{c}}(G(\ell, 1, n))$-module.
	\end{rmk}
	
	\begin{rmk}
		Even if the zeroth homology of $\widetilde{C}_{\bullet}(\lambda)$ is irreducible (and thus it necessarily coincides with $L_{\widetilde{c}}(\widetilde{\lambda})$) the natural Hermitian form on $L_{\widetilde{c}}(\widetilde{\lambda})$ does not need to be positive-definite, even if that for $L_{c}(\lambda)$ is. An example of this is given by taking $\ell = 2$, odd $n$, $e = n$, $\lambda = (e-1, 1)$ and the parameter $\widetilde{c}$ as in Remark \ref{rmk1} . In this case, $L_{\widetilde{c}}(\widetilde{\lambda}) = I_{e, 1, n}(2)$, which does not admit an invariant positive-definite Hermitian form, cf. \cite[Proposition 7.1]{MR2959035} 
	\end{rmk}
	
%	\begin{rmk}
%		Finally, we remark that the one-column homomorphisms of this paper do not allow us to directly deal with \emph{all} parameters $\widetilde{c}$ for the rational Cherednik algebra of $G(\ell, 1, n)$. In order to deal with these parameters, we would have admit rational multicharges. Of course, these homomorphisms were essential in constructing the BGG resolution $C_{\bullet}(\lambda)$, and we could not have constructed the resolutions $\widetilde{C}_{\bullet}(\lambda)$ without them! 
%	\end{rmk}
	
	\subsubsection{The $(k, e)$-equals ideal} \label{keideal}
	We now consider the subspace arrangements of
	 $k$ distinct clusters of $e$ equal parameters for $n=ke$.  
	We show that the BGG resolution of $L(\triv)$ is a minimal resolution of the coordinate ring of this subspace arrangement and generalise this to type $G(\ell, 1,n)$ as before.  
	
	Let  $n = ke$,  as we have seen, in this case we can give a BGG resolution of $L_{1/e}(\triv)$. It follows from \cite[Theorem 5.10]{MR2534594} that   $\rad(\Delta_{1/e}(\triv))$ is the ideal $I_{e, k, n}$ of functions vanishing on
	\[
	X_{e, k, n} := \mathfrak{S}_{n}\{(z_{1}, \dots, z_{n}) \in \C^{n} : z_{1} = \cdots= z_{e}, z_{e+1} = \cdots = z_{2e}, \dots, z_{(k-1)e+1} = \cdots = z_{ke}\}.
	\]
	Recall that the resolution of $L_{1/e}(\triv)$ is obtained as the Ringel dual of the resolution of $L_{1/e}(e^{k})$. Thus, the projective dimension of the algebra of functions $\C[X_{e, k, n}] = \C[\underline{x}]/I_{e, k, n} \cong L_{1/e}(\triv)[c_{\triv}]$ is $(e - 1)k$. By the Auslander--Buchsbaum formula, the depth of $\C[X_{e, k, n}]$ is $n - (e-1)k = k$. So $\C[X_{e, k, n}]$ is always Cohen-Macaulay, and we recover a special case of \cite[Proposition 3.11]{EGL}.
	
	Let us now analyze the regularity of $L_{1/e}(\triv)[c_{\triv}]$. By an argument similar to the proof of Proposition \ref{prop:reg e-equals}, this is given by $c_{(k^{e})} - c_{\triv} - (e-1)k$. By a direct computation, this is
	\[
	\reg(L_{1/e}(\triv)[c_{\triv}]) = \frac{k(n-e-k+1)}{2}
	\]

	\begin{eg}
		Assume $e = 3$, $n = 6$. Then we have that a resolution of $L(\triv)[c_{\triv}] = \C[x_{1}, \cdots, x_{6}]/I_{3, 2, 6}$ is given by
		\[
		0 \rightarrow (2^3)[-6] \rightarrow (3, 2, 1)[-5] \rightarrow (3^2)[-4] \oplus (4, 1^2)[-4] \rightarrow (5,1)[-2] \rightarrow (6) \rightarrow \C[\underline{x}]/I_{3, 2, 6} \rightarrow 0
		\]		
		\noindent and $\reg(\C[\underline{x}]/I_{3,2,6}) = 2$. 
	\end{eg}

Of course, for $\ell \geq 1$ we also have the subspace arrangement
\[
X_{e, k, n}(\ell) :=  \mathfrak{S}_{n}\{(z_{1}, \dots, z_{n}) \in \C^{n} : z_{1}^{\ell} = \cdots= z_{e}^{\ell}, z_{e+1}^{\ell} = \cdots = z_{2e}^{\ell}, \dots, z_{(k-1)e+1}^{\ell} = \cdots = z_{ke}^{\ell}\}
\]
And its defining ideal $I_{e, k, n}(\ell)$. Since $I_{e, k, n}$ is the unique maximal submodule in $\Delta_{1/e}(\triv)$ and the submodules of this standard module are linearly ordered, the ideal $I_{e, k, n}$ is generated in a single degree. Thus, the exact same argument as that in the proof of \cite[Proposition 2.5]{MR2959035}, if $q_{1}(x_{1}, \dots, x_{n}), \dots, q_{t}(x_{1}, \dots, x_{n})$ are generators of $I_{e, k, n}$ of minimal degree, then $q_{1}(x_{1}^{\ell}, \dots, x_{n}^{\ell}), \dots,$ $ q_{t}(x_{1}^{\ell}, \dots, x_{n}^{\ell})$ are generators of $I_{e, k, n}(\ell)$. It follows that the complex $\widetilde{C}_{\bullet}(\triv)$ is a minimal graded-free resolution of the algebra of functions $\C[X_{e, k, n}(\ell)]$, and the variety $X_{e, k, n}(\ell)$ is always Cohen-Macaulay. Moreover, the regularity of $\C[X_{e, k, n}(\ell)]$ is given by $\ell(c_{(k^{e})} - c_{\triv}) - (e-1)k$, or more explicitly,
\[
\reg(\C[X_{e, k, n}(\ell)]) = \frac{k[\ell(n+e-k-1) -2(e-1)]}{2}.
\]
We remark that in general as $H_{\widetilde{c}}(G(\ell, 1, n))$-modules, $\C[X_{e, k, n}(n)]$ does \emph{not} coincide with $L_{\widetilde{c}}(\triv)$. For example, if $\ell = 2$, $e = n$ is even and $\widetilde{c}(s) = 0$ for a reflection $s$ not conjugate to an element of $\mathfrak{S}_{n}$, then $L_{\widetilde{1/e}}(\triv)$ is finite-dimensional, while $\C[X_{e, 1, n}(2)]$ is not.

	\begin{rmk}
		Changing the parameter of the rational Cherednik algebra to $c = a/e > 0$ with $\gcd(a;e) = 1$ does not change the shape of the resolution $C_{\bullet}(\lambda)$, so the projective dimension and depth of $L_{a/e}(\lambda)$ are independent of $a\in \ZZ_{>0}$ when $\lambda$ is $e$-unitary. However, the value of $c_{\lambda}$ is not independent of $a\in \ZZ_{>0}$, and we get
		\[
		\beta_{i,j}(L_{a/e}(\lambda) = \beta_{i, j/a}(L_{1/e}(\lambda))
		\] 
		\noindent where we implicitly agree that $\beta_{i, j/a} = 0$ if $j/a \not\in \mathbb{Z}$. For any such $a\in \ZZ_{>0}$ the module $L_{a/e}((m-1)^p, q)$ can be identified with an ideal of $\C[\underline{x}]$ whose vanishing set coincides with $X_{e, 1, n}$. This ideal is radical if and only if $a = 1$, cf. \cite[Theorem 5.10]{MR2534594}. Similar considerations apply to $L_{a/e}(\triv)$. 
	\end{rmk}

\subsection*{Acknowledgements} E. Norton gratefully acknowledges the support and perfect working conditions of the Max Planck Institute for Mathematics, and J. Simental thanks the Hausdorff Institute for Mathematics where he was a participant at the Junior Trimester Program on Symplectic geometry and Representation theory; these institutes made our collaboration possible. C. Bowman and E. Norton thank the organizers of the conference ``Representation theory of symmetric groups and related algebras" held in December 2017 at the  National University of Singapore, for providing an excellent venue for the fruitful exchange of mathematical ideas. We thank J.~Chuang, S.~Fishel, S.~Griffeth, H.~Ko, Z.~Lin, P.~Shan, C.~Stroppel, and J.~Torres for stimulating conversations  
 and to further thank S. Griffeth for his careful reading of an earlier version of this manuscript. 
S.~Fishel and S.~Griffeth have an alternative approach to \cite[Conjecture 4.5]{conjecture}; our work was independent of each other. E. Norton dedicates this paper to ``the Crew:" Danny, Jesse, and Joe.

\bibliographystyle{amsalpha}   
\bibliography{master} 

\end{document}